\documentclass[11pt,reqno]{amsproc}

\usepackage{amsmath,amsfonts,amssymb,amsthm}
\usepackage[abbrev,lite,nobysame]{amsrefs}
\usepackage{mathrsfs}
\usepackage{graphics,graphicx}
\usepackage[usenames,dvipsnames]{color}
\usepackage{times}
\usepackage{bbm}
\usepackage[margin=1in]{geometry}
\usepackage[colorlinks=true, pdfstartview=FitV, linkcolor=blue, citecolor=blue, urlcolor=blue]{hyperref}
\usepackage{subcaption}

\usepackage{enumerate}
\usepackage{enumitem}

\makeatletter

\renewcommand\subsection{\@startsection{subsection}{2}%
\normalparindent{.5\linespacing\@plus.7\linespacing}{-.5em}
{\normalfont\bfseries}}

\renewcommand\subsubsection{\@startsection{subsubsection}{3}%
\normalparindent{.5\linespacing\@plus.7\linespacing}{-.5em}
{\normalfont\bfseries}}

\newcommand{\bullpar}[1]{\vspace{.5em}\noindent $\bullet$ \normalfont {\bfseries #1.}}

\newcommand{\diampar}[1]{\vspace{.5em}\noindent $\diamond$ \normalfont {\itshape #1.}}

\def\@tocline#1#2#3#4#5#6#7{\relax
  \ifnum #1>\c@tocdepth 
  \else
    \par \addpenalty\@secpenalty\addvspace{#2}%
    \begingroup \hyphenpenalty\@M
    \@ifempty{#4}{%
      \@tempdima\csname r@tocindent\number#1\endcsname\relax
    }{%
      \@tempdima#4\relax
    }%
    \parindent\z@ \leftskip#3\relax \advance\leftskip\@tempdima\relax
    \rightskip\@pnumwidth plus4em \parfillskip-\@pnumwidth
    #5\leavevmode\hskip-\@tempdima
      \ifcase #1
       \or\or \hskip 1em \or \hskip 2em \else \hskip 3em \fi%
      #6\nobreak\relax
    \dotfill\hbox to\@pnumwidth{\@tocpagenum{#7}}\par
    \nobreak
    \endgroup
  \fi}
\makeatother

\newtheorem{theorem}{Theorem}
\newtheorem{proposition}{Proposition}[section]
\newtheorem{lemma}[proposition]{Lemma}
\newtheorem{corollary}[proposition]{Corollary}

\theoremstyle{definition}

\newtheorem{remark}[proposition]{Remark}

\numberwithin{equation}{section}

\newcommand{\green}[1]{\textcolor{green}{#1}}
\newcommand{\blue}[1]{\textcolor{blue}{#1}}
\newcommand{\yellow}[1]{\textcolor{yellow}{#1}}

\newcommand{\B}{{\mathcal{B}}}
\newcommand{\C}{{\mathbb C}}
\newcommand{\E}{{\mathcal E}}
\newcommand{\cI}{{\mathcal I}}
\renewcommand{\L}{{\mathcal{L}}}
\newcommand{\G}{{\mathcal{G}}}
\newcommand{\N}{{\mathbb N}}
\newcommand{\R}{{\mathbb R}}
\newcommand{\Z}{{\mathbb Z}}

\newcommand{\D}{{\Delta}}
\newcommand{\T}{{\mathbb T}}
\newcommand{\cT}{{\mathcal T}}
\newcommand{\cS}{{\mathcal S}}

\newcommand{\W}{{\mathcal{W}}}

\renewcommand{\O}{{\Omega}}
\renewcommand{\o}{{\omega}}

\renewcommand{\d}{\mathrm{d}}
\newcommand{\ep}{\varepsilon}
\newcommand\e{{\rm e}}

\renewcommand{\v}{\boldsymbol{v}}

\renewcommand{\b}{\beta}

\newcommand{\p}{\partial}

\renewcommand{\l}{\left}

\renewcommand{\r}{\right}

\def\Re{{\rm Re}}
\def\Im{{\rm Im}}

\begin{document}

\title[Limiting absorption principles and inviscid damping in Euler-Boussinesq]{Limiting absorption principles and linear inviscid damping in the Euler-Boussinesq system in the periodic channel}

\author[M. Coti Zelati]{Michele Coti Zelati}
\address{Department of Mathematics, Imperial College London, London, SW7 2AZ, UK}
\email{m.coti-zelati@imperial.ac.uk}

\author[M. Nualart]{Marc Nualart}
\address{Department of Mathematics, Imperial College London, London, SW7 2AZ, UK}
\email{m.nualart-batalla20@imperial.ac.uk}



\begin{abstract}
We consider the long-time behavior of solutions to the two dimensional non-homogeneous Euler equations under the Boussinesq approximation posed on a periodic channel. We study the linearized system near a linearly stratified  Couette flow and  prove inviscid damping of the perturbed density and velocity field for any positive Richardson number, with optimal rates. Our methods are based on time-decay properties of oscillatory integrals obtained using a limiting absorption principle, and require a careful understanding of the asymptotic expansion of the generalized eigenfunction near the critical layer. 

As a by-product of our analysis, we provide a precise description of the spectrum of the linearized operator, which, for sufficiently large Richardson number, consists of an essential spectrum (as expected according to classical hydrodynamic problems) as well as discrete neutral eigenvalues (giving rise to oscillatory modes) accumulating towards the endpoints of the essential spectrum. 
\end{abstract}

\maketitle

\setcounter{tocdepth}{2}
\tableofcontents

\section{Introduction}\label{inviscid}
Under the Boussinesq approximation, the motion of an incompressible, non-homogeneous, inviscid fluid is described  by the Euler equations
\begin{equation}\label{eq:EBintro}
\begin{aligned}
(\p_t+\tilde\v\cdot\nabla)\tilde\o  &= -\mathfrak{g}\p_x\tilde\rho,\\
(\p_t+\tilde\v\cdot\nabla)\tilde\rho&=0, 
\end{aligned}
\end{equation}
where  $\tilde\v=\nabla^\perp\Delta^{-1}\tilde\o$ denotes the velocity field of the fluid with vorticity $\tilde\o=\nabla^\perp \cdot \tilde\v$ and density $\tilde\rho$, and $\mathfrak{g}$ is the gravity constant. 

In the  periodic channel $\T\times[0,1]$, we are interested  in the linear asymptotic  stability of the special equilibrium solution 
\begin{equation}\label{eq:StratCouette}
\bar{\v}=(y,0), \qquad \bar{\rho}(y)=1-\vartheta y, \qquad \p_y p=-\mathfrak{g}\bar{\rho}(y),
\end{equation}
which describes a Couette flow that is linearly stratified by a density with slope $\vartheta>0$. We introduce the perturbed velocity 
$\tilde{\v}=\bar{\v}+\v$ and density profile $\tilde{\rho}=\bar{\rho}+\vartheta\rho$,  and define the corresponding vorticity perturbation $\o=\nabla^\perp\cdot \v$. After neglecting the nonlinear terms,
the linearized Euler-Boussinesq system \eqref{eq:EBintro}  near \eqref{eq:StratCouette} can be written as
\begin{equation}\label{eq:linEulerBouss}
\begin{cases}
\p_t\o + y\p_x\o=-\b^2\p_x\rho \\
\p_t\rho + y\p_x\rho =\p_x\psi,\\
\D\psi=\o,
\end{cases}
\end{equation}
with $\psi$ being the streamfunction and $\beta=\sqrt{\vartheta \mathfrak{g}} >0$. The understanding of the long-time dynamics of solutions to \eqref{eq:linEulerBouss} is very much related to the spectral properties of the associated linear operator
\begin{equation}\label{eq:LinSSC}
\L=\begin{pmatrix}
y\p_x & \b^2\p_x \\
-\D^{-1}\p_x & y\p_x
\end{pmatrix}.
\end{equation}
In the setting of the periodic channel, $\L$ can have quite interesting features: it has both continuous and point spectrum, with a sequence of eigenvalues accumulating to the endpoint of the spectrum. As a consequence, any
asymptotic stability result requires well-prepared initial data, whose projection onto the point spectrum vanishes. 

We summarize the main result of this article in the following theorem. There a few key assumptions on the initial data that we informally state in the theorem and comment on right after. 

\begin{theorem}\label{thm:mainchan}
Let $\b>0$ and assume that the initial data $(\o^0,\rho^0)$ vanish on the physical boundaries, is orthogonal to the subspace generated by the eigenfunctions of $\L$, satisfy an orthogonality condition at the endpoint of the essential spectrum and 
\begin{equation}\label{eq:zeroxave}
\int_\T \o^0(x,y)\d x = \int_\T\rho^0(x,y)\d x =0.
\end{equation}
Let $\v=(v^x,v^y)=\nabla^\perp\psi=(-\p_y\psi,\p_x\psi)$ be the corresponding velocity field. We have the following estimates. 
\begin{itemize}
\item If $\b^2\neq1/4$, let $\mu=\Re\sqrt{1/4-\b^2}$ and $\nu=\Im\sqrt{1/4-\b^2}$. Then,
 \begin{align}
\Vert v^x(t) \Vert_{L^2}&\lesssim \frac{1}{t^{\frac12-\mu}}\l( \Vert \rho^0 \Vert_{L^2_xH^3_y} + \Vert \o^0 \Vert_{L^2_xH^3_y}\r), \label{eq:decayvx} \\
\Vert v^y(t) \Vert_{L^2}&\lesssim \frac{1}{t^{\frac32-\mu}}\l( \Vert \rho^0 \Vert_{L^2_xH^4_y} + \Vert \o^0 \Vert_{L^2_xH^4_y}\r), \label{eq:decayvy}\\
\Vert \rho(t) \Vert_{L^2}&\lesssim \frac{1}{t^{\frac12-\mu}}\l( \Vert \rho^0 \Vert_{H^1_xH^3_y} + \Vert \o^0 \Vert_{H^1_xH^3_y}\r), \label{eq:decayrho} 
\end{align}
for all $t\geq 1$.
\item If $\b^2=1/4$, then 
\begin{align}
\Vert v^x(t) \Vert_{L^2}&\lesssim \frac{1+\log(t)}{t^\frac12}\l( \Vert \rho^0 \Vert_{L^2_xH^3_y} + \Vert \o^0 \Vert_{L^2_xH^3_y}\r), \label{eq:decayvxlog}\\
\Vert v^y(t) \Vert_{L^2}&\lesssim \frac{1+\log(t)}{t^\frac32}\l( \Vert \rho^0 \Vert_{L^2_xH^4_y} + \Vert \o^0 \Vert_{L^2_xH^4_y}\r), \label{eq:decayvylog}\\
\Vert \rho(t) \Vert_{L^2}&\lesssim\frac{1+\log(t)}{t^\frac12}\l( \Vert \rho^0 \Vert_{H^1_xH^3_y} + \Vert \o^0 \Vert_{H^1_xH^3_y}\r),   \label{eq:decayrholog} 
\end{align}
for all $t\geq 1$. 
\end{itemize}
\end{theorem}

\begin{remark}[Assumptions on data]\label{rmk:data}
The assumptions on the initial data are completely natural. The vanishing at the boundary points $y\in\{0,1\}$ is a typical requirement \cite{Jia20,IJ22}, while \eqref{eq:zeroxave} is inessential, as the $x$-average is a constant of motion for \eqref{eq:linEulerBouss}. The orthogonality to eigenfunctions of $\mathcal{L}$ is needed to avoid oscillatory, non-decaying modes (which are present for $\beta^2>1/4$, see Section \ref{sub:specpic}). Lastly, the precise meaning of the spectral assumption at the endpoints of the essential spectrum $\sigma_{ess}(\mathcal{L})=[0,1]$ is in condition \eqref{eq:specassume} in Section \ref{sub:LAPboundary} below. It requires orthogonality to certain generalized eigenfunctions that appear at $\p \sigma_{ess}(\mathcal{L})=\{0,1\}$.    
\end{remark}

The inviscid damping estimates \eqref{eq:decayvx}-\eqref{eq:decayrholog} encode the asymptotic stability of \eqref{eq:linEulerBouss} and precisely describe the long-time dynamics.
The decay is due to a combination of \emph{mixing} (due to the background Couette flow) and \emph{stratification} (due to the background density). The former has been extensively
studied in the homogeneous Euler equations both at the linear level \cite{BCZV19,CZZ19,GNRS20,Jia20,JiaGev20,WZZ18,WZZ19,WZZKolmo20,Zillinger16,Zillinger17,ZillingerCirc17,Zillinger21} and at the nonlinear level \cite{BM15,IJ22,IJ20,IJnon20,MZ20}.

In the presence of stratification, the spectral stability of the Euler-Boussinesq system has been address in the classical work of Miles \cite{Miles} and Howard \cite{Howard}. See  \cite{Yaglom}*{Section 3.2.3} for a survey on the literature regarding the spectral problem.  The first work in the
direction of asymptotic stability dates back to Hartman \cite{Hartman} in 1975, in which \eqref{eq:linEulerBouss} on $\T\times\R$ was solved explicitly on the Fourier side using hypergeometric functions. Moreover,
it was predicted the vorticity should be unstable in $L^2$, with a growth proportional to $\sqrt{t}$. This approach was used in \cite{YL18} to prove decay rates analogous to those in Theorem \ref{thm:mainchan} in $\T\times\R$. In this spatial setting, a different approach based on an energy method in Fourier space was used in \cite{BCZD22} to prove both inviscid damping and instability in the spectrally stable regime $\beta^2>1/4$,
confirming the predictions of \cite{Hartman}. The analysis has been extended in the full nonlinear setting in \cite{BBCZD21}. A third proof of linear inviscid damping on $\T\times\R$ can be found in our companion article \cite{CZN23strip}, in which the methods developed here can be used  to provide explicit solutions in physical variables to \eqref{eq:linEulerBouss}.

Our article constitutes the first result of (linear) asymptotic stability of a stably stratified shear flow for the Euler-Boussinesq equations in the periodic channel, as well as the first rigorous characterization of the spectrum of the linearized operator \eqref{eq:LinSSC} and in particular the existence of discrete neutral eigenvalues for $\b^2>1/4$. From a technical standpoint, the main difficulty lies in the stratification of the background density $\bar\rho$. This manifests itself in the equation that rules the underlying spectral problem (the Taylor-Goldstein equation, see \eqref{eq:TG} below), which becomes more singular than the usual Rayleigh equation for inviscid homogeneous fluids. 

This work also connects with the global well-posedness for the Euler-Boussinesq equations and, by extension, to the axisymmetric 3d Euler equations. Certain solutions to the Euler-Boussinesq and 3d Euler equations are known to blow up in finite time, see the ground-breaking work of Elgindi \cite{Elgindi21}, and related works \cites{ElgindiJeong19, ElgindiJeong20, ChenHou23}. On the other hand, there are examples where inviscid damping plays a key role in proving global well-posedness for the 3d Euler equations and for the inhomogeneous 2d Euler equations, see \cite{GPW} and \cite{CWZZ, Zhao23}, respectively. In the case of Euler-Boussinesq near stratified shear flows, a long-time existence result relying on inviscid damping estimates can be found in \cite{BBCZD21}.

\section{Main ideas and outline of the article}
In this section, we give a brief account of the strategy of proof of Theorem \ref{thm:mainchan}, recording the main steps that will be then expanded in the subsequent sections, and providing a quick reasoning behind the assumptions of Theorem \ref{thm:mainchan} on the initial data. We focus on the case $\b^2\neq 1/4$ for the sake of clarity. When $\b^2=1/4$, the strategy is the same, but the statements of the main results typically differ by a logarithmic correction, and we prefer to postpone them in the relevant Section \ref{sec: Bounds Greens Function special}.
We also set some of the notation and assumptions that will be used throughout the manuscript.

\subsection{Fourier decomposition and spectral representation}
The setting of the periodic channel $\T\times[0,1]$ considered in this article poses new challenges as it forbids the use of Fourier methods in the vertical direction $y$. However, we can decouple \eqref{eq:linEulerBouss}
 in Fourier modes in $x\in\T$, writing
\begin{equation*}
\o=\sum_{m\in\Z}\o_m(t,y)\e^{imx}, \qquad \rho=\sum_{m\in\Z}\rho_m(t,y)\e^{imx}, \qquad \psi=\sum_{m\in\Z}\psi_m(t,y)\e^{imx},
\end{equation*} 
so that 
\begin{equation*}
(\p_t+imy)\o_m=-im\b^2\rho_m, \qquad (\p_t+imy)\rho_m=im\psi_m,
\end{equation*}
for each $m\in\Z$, with 
\begin{equation*}
\begin{cases}
\D_m\psi_m =\o_m, \\
\psi_m|_{y=0,1}=0,
\end{cases}
\qquad  \D_m:= \p_y^2-m^2.
\end{equation*}  
The modes corresponding to the $x$-average, namely when $m=0$, are clearly conserved and therefore we will not consider them further (cf. \eqref{eq:zeroxave}). 
Moreover, since $\o$ and $\rho$ are real-valued, we necessarily have that $\overline{\o_{-m}}=\o_m$ and $\overline{\rho_{-m}}=\rho_m$. Without loss of generality, we take $m\geq 1$. 

For our purposes, it is more convenient to write \eqref{eq:linEulerBouss} in the compact  stream-function formulation 
\begin{equation*}
\p_t \begin{pmatrix} \psi_m \\ \rho_m\end{pmatrix}+imL_m\begin{pmatrix}
\psi_m \\ \rho_m
\end{pmatrix}=0,
\end{equation*}
and directly obtain its solution as
\begin{equation*}
\begin{pmatrix} \psi_m \\ \rho_m\end{pmatrix}=\e^{-imL_mt}\begin{pmatrix}
\psi_m^0 \\ \rho_m^0
\end{pmatrix}
\end{equation*}
where $L_m$ is the linear operator defined by
\begin{equation}\label{eq:linOP}
L_m =\begin{pmatrix}
\D_m^{-1}(y\D_m)  & \b^2\D_m^{-1} \\
-1 & y
\end{pmatrix}
\end{equation}
Using Dunford's formula \cite{Engel-Nagel, Taylor-11}, we have that
\begin{equation}\label{eq:Dunford}
\begin{pmatrix}
\psi_m(t,y) \\ \rho_m(t,y)
\end{pmatrix} = \frac{1}{2\pi i} \int_{\p\O}\e^{-imct} (c-L_m)^{-1}\begin{pmatrix}
\psi_m^0(y) \\ \rho_m^0(y)
\end{pmatrix} \,\d c, 
\end{equation}
where here $\O$ is any domain containing the spectrum $\sigma(L_m)$. Under suitable conditions on the initial data 
(see Proposition \ref{prop: contour limiting absorption principle} below),  we can reduce the contour of integration to  
\begin{equation}\label{eq:psirholim}
\begin{pmatrix}
\psi_m(t,y) \\ \rho_m(t,y)
\end{pmatrix} 
=\frac{1}{2\pi i }\lim_{\ep\rightarrow 0}\int_0^1 \e^{-imy_0t}\l[(-y_0-i\ep+L_m)^{-1}-(-y_0+i\ep+L_m)^{-1}\r]
\begin{pmatrix}
\psi_m^0 \\ \rho_m^0 
\end{pmatrix}\, \d y_0.
\end{equation}
In particular, the contour integral along the {essential spectrum of $L_m$, $\sigma_{ess}(L_m)=[0,1]$ is the only non-trivial contribution from $\sigma(L_m)$ to the Dunford's formula. For $\ep>0$, we denote
\begin{equation}\label{eq:geneigen}
\begin{pmatrix}
\psi^{\pm}_{m,\ep}(y,y_0) \\ \rho^\pm_{m,\ep}(y,y_0)
\end{pmatrix}:=\l( -y_0\pm i\ep+L_m\r)^{-1}\begin{pmatrix}
\psi_m^0(y) \\ \rho_m^0(y)
\end{pmatrix}
\end{equation}
and obtain the coupled system of equations
\begin{equation*}
\begin{aligned}
\o_m^0(y)&=(y-y_0\pm i\ep)\D_m\psi^\pm_{m,\ep}(y,y_0)+\b^2 \rho^\pm_{m,\ep}(y,y_0), \\
\rho_m^0(y)&=(y-y_0\pm i\ep)\rho^\pm_{m,\ep}(y,y_0) -\psi^\pm_{m,\ep}(y,y_0).
\end{aligned}
\end{equation*}
We first solve
\begin{equation}\label{eq rho m ep}
\rho^\pm_{m,\ep}(y,y_0)=\frac{1}{y-y_0\pm i\ep}\l( \rho_m^0(y)+\psi^\pm_{m,\ep}(y,y_0)\r)
\end{equation}
and from there we obtain the following inhomogeneous \emph{Taylor-Goldstein equation} for $\psi^\pm_{m,\ep}$,
\begin{equation}\tag{TG}\label{eq:TG}
\D_m\psi^{\pm}_{m,\ep}+\b^2\frac{\psi^\pm_{m,\ep}}{(y-y_0\pm i\ep)^2}=\frac{\o_m^0}{y-y_0\pm i\ep}-\b^2\frac{\rho_m^0}{(y-y_0\pm i\ep)^2},
\end{equation}
along with homogeneous Dirichlet boundary conditions at $y=0,1$.

\subsection{Notation and conventions}
Throughout the manuscript, we assume $\b>0$ and $m\geq 1$. We say that $A\lesssim B$ when there exists $C>0$ such that $A\leq CB$. 
Also, for $j\geq 0$ we define
\begin{equation*}
Q_{j,m}=\Vert \rho_m^0 \Vert_{H^{j+2}_y} + \Vert \o_m^0 \Vert_{H^{j+2}_y},
\end{equation*}
to quantify the regularity requirements on the initial data.

\subsection{Green's function for the Taylor-Goldstein equation}
Solutions to \eqref{eq:TG} are fundamental objects of study of this work. They can be constructed via the classical method of Green's functions, by first solving
the \emph{homogeneous} Taylor-Goldstein equation
\begin{equation}\label{eq:TGoperator}\tag{TGh}
\textsc{TG}_{m,\ep}^\pm\phi=0, \qquad \textsc{TG}_{m,\ep}^\pm := \D_m +\frac{\b^2}{(y-y_0\pm i\ep)^2},
\end{equation}
for $y\in (0,1)$. We refer to $\textsc{TG}_{m,\ep}^\pm$ as to the \emph{Taylor-Goldstein operator}. As in the statement of Theorem \ref{thm:mainchan}, we define throughout the
article the numbers
\begin{equation}\label{eq:munu}
\mu=\Re\l(\sqrt{1/4-\b^2}\r), \qquad  \nu=\Im\l(\sqrt{1/4-\b^2}\r),
\end{equation}
and we denote by $\G^\pm_{m,\ep}(y,y_0,z)$ the Green's function of the Taylor-Goldstein equation, which  satisfies
\begin{equation}\label{eq:TGgreens}
\textsc{TG}_{m,\ep}^\pm\G^\pm_{m,\ep}(y,y_0,z)=\delta(y-z).
\end{equation}
While $\G^\pm_{m,\ep}(y,y_0,z)$ has an explicit expression, reported in Proposition \ref{prop: def Green's Function}, we record its important properties as the key result.

\begin{theorem}\label{L2 bounds of G}
Let $\beta^2\neq 1/4$.
There exists $\ep_0>0$ such that for all $\ep\in(0, \ep_0)$ and for all $y,y_0\in[0,1]$ such that $m|y-y_0|\leq 3\b$, we have 
\begin{equation*}
|y-y_0+i\ep|^{-\frac12+\mu} \Vert\G^\pm_{m,\ep}(y,y_0,\cdot)\Vert_{L^2_z}+ |y-y_0+i\ep|^{\frac12+\mu} \Vert \p_y\G^\pm_{m,\ep}(y,y_0,\cdot)\Vert_{L^2_z}\lesssim \frac{1}{m^{1+\mu}}.
\end{equation*}
\end{theorem}
The theorem provides {sharp} bounds on the Green's function near the \emph{critical layer} $y=y_0$, where  \eqref{eq:TGoperator} is singular and \eqref{eq:TG}
has a regular singular point. The scale of the problem is crucially determined by $\beta$ and $m$.

The proof of Theorem \ref{L2 bounds of G} is carried out in Section \ref{sec: Bounds Greens Function}, while the analogous result for $\beta^2=1/4$ 
is stated in Theorem \ref{L2 bounds G special} and proven in Section \ref{sec: Bounds Greens Function special}. They are based on the asymptotic
properties of Whittaker functions  \cite{Whittaker03}, whose main properties can be found in Appendix \ref{app:Whittaker}. 

\subsection{Regularization of the generalized stream-functions}
The source term of \eqref{eq:TG} is, a priori, too singular for $\psi_{m,\ep}^\pm$ to be obtained as an application of the Green's function on \eqref{eq:TG}. However, the singularity of the source term is no worse than $\frac{\b^2}{(y-y_0\pm i\ep)^2}$, which is precisely the potential of the Taylor-Goldstein operator \eqref{eq:TGoperator}. Then, \eqref{eq:TG} may be written as
\begin{equation*}
\text{TG}_{m,\ep}^\pm\psi_m^\pm=\text{TG}_{m,\ep}^\pm \l( \frac{1}{\b^2}(y-y_0\pm i\ep)\o_m^0 -\rho_m^0\r) + \D_m\big(\rho_m^0(y)-\frac{1}{\b^2}(y-y_0\pm i\ep)\o_m^0(y)\big).
\end{equation*}
Hence, for $z,y_0\in[0,1]$ and $0\leq \ep\leq 1$, define
\begin{equation}\label{eq:Fdata}
F_{m,\ep}^\pm(z,y_0):=\D_m\rho_m^0(z)-\frac{1}{\b^2}\D_m\big((z-y_0\pm i\ep)\o_m^0(z)\big)
\end{equation}
and note that, since the pair of initial data vanish on the physical boundaries $y=0$ and $y=1$, the solution $\psi^\pm_{m,\ep}(y,y_0)$  to \eqref{eq:TG}
is given by
\begin{equation}\label{def psi}
\begin{aligned}
\psi^\pm_{m,\ep}(y,y_0)&= \frac{1}{\b^2}(y-y_0\pm i\ep)\o_m^0(y) -\rho_m^0(y) + \varphi_{m,\ep}^\pm(y,y_0),
\end{aligned}
\end{equation}
while 
\begin{equation}\label{def rho}
\rho_{m,\ep}^\pm(y,y_0)=\frac{1}{\b^2}\o_m^0(y) +\frac{1}{y-y_0\pm i\ep} \varphi_{m,\ep}^\pm(y,y_0).
\end{equation}
Here, $\varphi_{m,\ep}^\pm$ solves 
\begin{equation}\label{eq:varphieq}
\textsc{TG}_{m,\ep}^\pm\varphi_{m,\ep}^\pm  = F_{m,\ep}^\pm
\end{equation}
and is given by 
\begin{equation}\label{eq:defvarphi}
    \varphi_{m,\ep}^\pm(y,y_0)=\int_0^1 \G^\pm_{m,\ep}(y,y_0,z) F_{m,\ep}^\pm(z,y_0)\d z.
\end{equation}
The main reason to write $\psi_{m,\ep}^\pm$ and $\rho_{m,\ep}^\pm$ using \eqref{def psi} and \eqref{def rho} is that now $F_{m,\ep}^\pm\in L^2_z$ and we can use the bounds on the Green's function $\G_{m,\ep}^\pm$ from Theorem \ref{L2 bounds of G} in \eqref{eq:defvarphi} to estimate $\varphi_{m,\ep}^\pm$, and thus $\psi_{m,\ep}^\pm$ and $\rho_{m,\ep}^\pm$, near the critical layer.  The introduction of $F_{m,\ep}^\pm$ constitutes a first example of the underlying \emph{motif} of inviscid damping, namely that \emph{decay costs regularity}.

\subsection{Spectral picture}\label{sub:specpic}
The main assumption of Theorem \ref{thm:mainchan} consists in requiring that the initial data are orthogonal to the subspace generated by the eigenfunctions of $L_m$. 
Generically speaking, (embedded) eigenvalues may constitute an obstruction to damping phenomena, as they can give rise to oscillatory modes or even growing (hence unstable)
modes. The spectral picture here is quite intriguing and drastically different compared to the case of the periodic strip. The main result on the spectrum of $L_m$ is below.

\begin{theorem}\label{thm:spectralpic}
Let $\b>0$. Then the essential spectrum of $L_m$ is $\sigma_{ess}(L_m)=[0,1]$. Moreover,
\begin{itemize}
\item any eigenvalue $c\in\C$ such that $\l| \Re(c)-1/2\r| \geq 1/2$, must have $\Im(c)= 0$;
\item for $\b^2>1/4$,
\begin{itemize}
\item there are no eigenvalues $c\in\C$ such that $\Im(c)\neq0$ and $\Re(c)\in(0,1)$.
\item there are no real eigenvalues $c\in\R$ such that  $c<-\b/m$ or $c>1+\b/m$.
\item there is a countably infinite number of discrete eigenvalues $c\in\C$, with $\Im(c)=0$ and $\Re(c)\in \l(-\b/m,0\r)\cup \l(1,1+\b/m\r)$. Moreover, they accumulate towards $0$ and $1$.
\end{itemize}
\item for $\b^2\leq1/4$,
\begin{itemize}
\item there is no  eigenvalue $c\in\C$ such that $\Re(c)\leq0$ or $\Re(c)\geq1$.
\item there is no  eigenvalue $c\in\C$ such that $\l| \Im (c)\r|\geq \b/m$ or $\l| \Im (c)\r|\leq \ep_0$.
\end{itemize}
\end{itemize}
\end{theorem}

The three cases outlined above are depicted in Figure \ref{fig:spectralpic}. Unstable eigenmodes can be ruled out by the classical 
Miles-Howard stability criterion \cites{Miles,Howard}  when $\b^2\geq1/4$, so that any eigenvalue $c\in\C$ of $L_m$ must have $\Im(c)=0$.
However, spectral stability is typically not sufficient to deduce asymptotic stability. This is particularly clear when $\beta^2>1/4$, for which
infinitely many eigenvalues exist, corresponding to neutral (oscillatory) modes. This is a specific feature of the problem in the \emph{periodic channel}.
The same problem on the periodic strip does not have any of these modes, as the essential spectrum is the whole real line, and hence eigenvalues
are ``pushed away to infinity''. In the periodic channel, each of these discrete eigenvalues are found to be zeroes of the Wronskian of the Green's function and this is precisely how we characterize them in Proposition \ref{prop discrete eigenvalue and peak}.

\begin{figure}[h!]
  \centering
  \begin{subfigure}[b]{0.32\linewidth}
    \includegraphics[width=\linewidth]{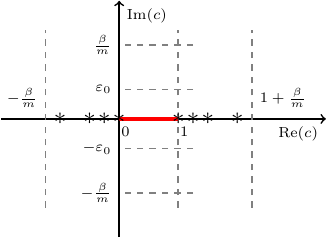}
     \caption{$\beta^2>1/4$}
  \end{subfigure}
  \begin{subfigure}[b]{0.32\linewidth}
    \includegraphics[width=\linewidth]{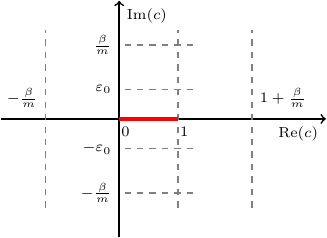}
     \caption{$\beta^2=1/4$}  
     \end{subfigure}
  \begin{subfigure}[b]{0.32\linewidth}
    \includegraphics[width=\linewidth]{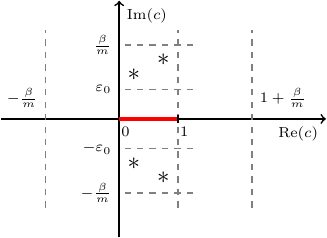}
     \caption{$\beta^2<1/4$}
       \end{subfigure}
  \caption{The essential spectrum   $\sigma_{ess}(L_m)=[0,1]$ is in red. Eigenvalues are denotes by  $*$. Theorem \ref{thm:spectralpic} shows their existence 
  for $\beta^2> 1/4$, while when $\beta^2<1/4$ we can only discern that they do not exist close to the essential spectrum.}
  \label{fig:spectralpic}
\end{figure}

When $\beta^2<1/4$, we are able to rule out the existence of eigenvalues in the proximity of the essential spectrum, which is a consequence of suitable lower bounds on the Wronskian. Nonetheless, isolated unstable eigenvalues in an intermediate region may exist in this case,
although their presence does not affect the conclusion of Theorem \ref{thm:mainchan} if the data are orthogonal to them. The proof of their existence is an interesting open
question.

The proof of Theorem \ref{thm:mainchan} is postponed to Section \ref{sec:spectral}. It requires an extensive analysis of the resolvent operator $(c-L_m)^{-1}$ and
of spectral integrals of the form \eqref{eq:Dunford}, where the domain of integration   containing the essential spectrum is carefully designed. 

\subsection{Solutions to the inhomogeneous Taylor-Goldstein equation} Once the Green's function is established and \eqref{eq:TG} is regularized due to the introduction of $F_{m,\ep}^\pm$ and $\varphi_{m,\ep}^\pm$, most of the analysis on $\psi_{m,\ep}^\pm$ will follow from the properties of generic solutions $\Phi_{m,\ep}^\pm$ to the general inhomogeneous Taylor-Goldstein equation
\begin{equation}\tag{TGf}\label{eq:inhom TG}
\textsc{TG}_{m,\ep}^\pm\Phi_{m,\ep}^\pm  = f,
\end{equation}
for some $f\in L^2$ and with boundary conditions $\Phi_{m,\ep}^\pm(0,y_0)=\Phi_{m,\ep}^\pm(1,y_0)=0$. To formally quantify the distance to the critical layer, for $y_0\in[0,1]$ and $n\geq 1$ we introduce the nested sets
\begin{equation*}
J_n=\lbrace y\in [0,1]: m|y-y_0|\leq n\b\rbrace
\end{equation*}
and $J_n^c=[0,1]\setminus J_n$. A direct consequence of Theorem \ref{L2 bounds of G} are the asymptotic expansions of $\Phi_{m,\ep}^\pm$ near the critical layer. That is, for all $y\in J_3$ we have
\begin{equation}\label{eq: locboundsPhi}
|y-y_0\pm i\ep|^{-\frac12+\mu} |\Phi_{m,\ep}^\pm(y,y_0)|+ |y-y_0\pm i\ep|^{\frac12+\mu} |\p_y \Phi_{m,\ep}^\pm(y,y_0)|\lesssim \frac{1}{m^{1+\mu}}\Vert f \Vert_{L^2_y}.
\end{equation}
Using the entanglement inequality
\begin{equation}\label{eq: entangleineq}
\Vert \p_y \Phi_{m,\ep}^\pm \Vert_{L^2_y(J_3^c)}^2 + m^2 \Vert  \Phi_{m,\ep}^\pm \Vert_{L^2_y(J_3^c)}^2 \lesssim m^2 \Vert \Phi_{m,\ep}^\pm \Vert_{L^2_y(J_2^c\cap J_3)}^2 + \frac{1}{m^2}\Vert f \Vert_{L^2_y(J_2^c)}^2,
\end{equation}
which is inspired from \cite{IIJ22} and proved in Lemma \ref{lemma:entangle ineq}, the localised asymptotic expansions \eqref{eq: locboundsPhi} provide integral estimates on $\Phi_{m,\ep}^\pm$ away from the critical layer,
\begin{equation}\label{eq:globboundsPhi}
\Vert\p_y \Phi_{m,\ep}^\pm(y,y_0)\Vert_{L^2_y(J_3^c)} + m\Vert \Phi_{m,\ep}^\pm(y,y_0)\Vert_{L^2_y(J_3^c)}\lesssim \frac{1}{m}\Vert f \Vert_{L^2_y}.
\end{equation}
The precise statements and proofs of \eqref{eq: locboundsPhi} and \eqref{eq:globboundsPhi}, as well as the corresponding versions for $\b^2=1/4$, can be found in Proposition \ref{L2 bounds inhom TG solution} in Section \ref{sec: bounds inhom TG}.

\subsection{Inviscid damping estimates through the limiting absorption principle}
The last step in the proof of Theorem \ref{thm:mainchan} is a stationary phase argument to deduce decay of $\psi_m$ and $\rho_m$ in \eqref{eq:psirholim}.
As customary, it involves an integration by parts in the spectral variable $y_0$ to gain time-decay from the oscillatory phase. The amount of decay that can be obtained is linked
to the regularity of the generalized streamfunctions $\psi^{\pm}_{m,\ep}$ in \eqref{eq:geneigen}, and even more crucially to their asymptotic expansion at the critical layer (matching
that of the Green's function in Theorem \ref{L2 bounds of G}, as can be seen from \eqref{def psi} and \eqref{eq:defvarphi}). Moreover, the integration leads to boundary terms at the endpoint of the spectrum that need to be treated \emph{ad hoc}.

To obtain the asymptotic expansions of $\psi_{m,\ep}^\pm$ near the critical layer, in Proposition \ref{stream derivatives formulae} we observe that $\p_y + \p_{y_0}$ commutes with the Taylor-Goldstein operator \eqref{eq:TGoperator} and we deduce formulas for $\p_{y_0}\psi_{m,\ep}^\pm$, and several other derivatives with respect to both $y$ and $y_0$. These formulas involve solutions $\Phi_{m,\ep}^\pm$ to \eqref{eq:inhom TG} for source terms $f$ given by derivatives of $F_{m,\ep}^\pm$. As is clear from \eqref{eq: locboundsPhi}, the asymptotic expansions of $\Phi_{m,\ep}^\pm$, and in turn of $\p_{y_0}\psi_{m,\ep}^\pm$ and related derivatives, are conditional to the $L^2$ boundedness of derivatives of $F_{m,\ep}^\pm$, constituting a further example of the fact that decay costs regularity.

Some formulas from Proposition \ref{stream derivatives formulae} involve as well terms related to $\p_y\varphi_{m,\ep}^\pm(z,y_0)$, and higher derivatives, evaluated at the physical boundaries $z=0$ and $z=1$. In general, these boundary terms arise when the Taylor-Goldstein operator \eqref{eq:TGoperator} acting on the $\p_y$ derivative of solutions to \eqref{eq:inhom TG} is inverted, and usually they do not vanish. See Proposition \ref{stream derivatives formulae} for more details. Near the critical layer, these boundary terms are studied in Section \ref{sec: boundary term estimates} and some require \eqref{eq:specassume} to be sufficiently regular, see Proposition \ref{prop: dyVarphi0 estimates} for more details.

Once the asymptotic expansions for $\psi_{m,\ep}^\pm$ near the critical layer are established via Proposition \ref{stream derivatives formulae} and Proposition \ref{L2 bounds inhom TG solution}, these are used through the entanglement inequality  \eqref{eq: entangleineq} to derive the regularity estimates of $\psi_{m,\ep}^\pm$ away from the critical layer. Additionally, asymptotic expansions and regularity estimates for $\rho_{m,\ep}^\pm$ are deduced accordingly thanks to \eqref{def rho}. The precise statements and proofs are found in Section \ref{sec: estimates generalized stream}. Both the asymptotic expansions and the regularity estimates are uniform in $\ep$ sufficiently small, so that the limiting functions in \eqref{eq:psirholim} retain the same properties.

\subsection{Limiting absorption principle for spectral boundary terms}\label{sub:LAPboundary}
The stationary phase argument employed in the proof of Theorem \ref{thm:mainchan} requires an integration by parts in the spectral variable $y_0$ in \eqref{eq:psirholim} regarding $\psi_m$ that involves spectral boundary terms evaluated at $y_0=0$ and $y_0=1$. These boundary terms are 
\begin{equation}\label{eq: psi solid term}
\begin{aligned}
-\frac{1}{2\pi i}\frac{1}{imt}\lim_{\ep\rightarrow 0}\Big[ \e^{-imy_0t}\l(\psi_{m,\ep}^-(y,y_0)-\psi_{m,\ep}^+(y,y_0)\r)\Big]_{y_0=0}^{y_0=1}.
\end{aligned}
\end{equation} 
For $y_0=0$, from \eqref{def psi} and \eqref{eq:defvarphi} we note that 
\begin{equation}\label{eq: dif psi boundary}
\begin{aligned}
\psi_{m,\ep}^-(y,0)-\psi_{m,\ep}^+(y,0)&= -\frac{2i\ep}{\b^2}\o_m^0 +\int_0^1 \l(\G_{m,\ep}^-(y,0,z)-\G_{m,\ep}^+(y,0,z)\r)F_{m}(z,0)\d z \\
&\quad +\frac{i\ep}{\b^2}\int_0^1 \l(\G_{m,\ep}^-(y,0,z)+\G_{m,\ep}^+(y,0,z)\r)\D_m\o_m^0\d z,
\end{aligned}
\end{equation}
where $F_m(z,0)=F_{m,0}^\pm(z,0)$. Moreover, for $\b^2>\frac14$, from Lemma \ref{Green's limit on boundary}, there exists $\ep_0>0$ and $C_\ep\geq C_0>0$ such that
\begin{equation}\label{eq: dif Greens boundary}
\left| \G_{m,\ep}^-(y,0,z)-\G_{m,\ep}^+(y,0,z)-C_\ep\phi_{u,m}(y)\phi_{u,m}(z)\right|\lesssim \ep^\frac12,
\end{equation}
for all $\ep\leq \ep_0$ and uniformly in $y,z\in[0,1]$. Here, $\phi_{u,m}$, given by \eqref{eq: phium}, denotes the generalized eigenfunction associated to the generalized eigenvalue $y_0=0$. Analogous expressions to \eqref{eq: dif psi boundary} and \eqref{eq: dif Greens boundary} can be deduced for the boundary term associated to $y_0=1$, now involving $\phi_{l,m}$, the generalized eigenfunction associated to the generalized eigenvalue $y_0=1$ and given by \eqref{eq: philm}.

In view of \eqref{eq: dif Greens boundary}, for \eqref{eq: psi solid term} to vanish we require the initial data $(\o_m^0, \rho_m^0)$ to be such that 
\begin{equation}\label{eq:specassume}\tag{H}
\int_0^1 \phi_{u,m}(z)F_m(z,0)\d z=\int_0^1 \phi_{l,m}(z)F_m(z,1)\d z=0.
\end{equation}
This is the key orthogonality assumption at the endpoint of the essential spectrum, which was discussed in Remark \ref{rmk:data}. Then, we are able to show
\begin{theorem}\label{thm: boundary limiting absorption principle}
We have that
\begin{equation*}
\lim_{\ep\rightarrow 0}\l\Vert\psi_{m,\ep}^-(\cdot,y_0)-\psi_{m,\ep}^+(\cdot,y_0)\r\Vert_{L^2_y}=0, \qquad  y_0\in\{ 0, 1\}.
\end{equation*}
\end{theorem}
The proof of Theorem \ref{thm: boundary limiting absorption principle} is carried out in Section \ref{sec:spectral}, where \eqref{eq: dif Greens boundary} is shown in Lemma \ref{Green's limit on boundary} for $\b^2>1/4$. For the case $\b^2\leq 1/4$, the difference of Green's functions at $y_0=0$ and $y_0=1$ vanish as $\ep\rightarrow 0$ and no orthogonality conditions are needed, see Lemma \ref{real Green's limit on boundary} and Lemma \ref{special Green's limit on boundary} for more details.

\section{Explicit solutions to the Taylor-Goldstein equation}\label{sec: explicit TG solutions}
The first step towards the proof of Theorem \ref{thm:mainchan} is to derive the expression of the Green's function associated to \eqref{eq:TG}. The building block consists of the so-called Whittaker functions \cite{Whittaker03}, a modified form of hypergeometric functions that solve equations of the form
\begin{equation}\label{eq:whit}
\p^2_{\zeta} M_{\kappa,\gamma}+\left(-{\frac  {1}{4}}+{\frac  {\kappa }{\zeta}}+{\frac  {1/4- \gamma^{2}}{\zeta^{2}}}\right)M_{\kappa,\gamma}=0, \qquad \zeta\in\C,
\end{equation}
for parameters $\kappa,\gamma\in\C$.
Their properties are reported in Appendix \ref{app:Whittaker}. 

\subsection{The case $\b^2\neq 1/4$}
We use Whittaker functions with $\gamma=\pm(\mu+i\nu)=\pm \sqrt{1/4-\b^2}$ and $b=0$, see \eqref{eq:munu}, and denote  by
$M_\pm(\zeta):=M_{0,\pm(\mu+i\nu)}(2m\zeta)$ the solution to the rescaled Whittaker equation
\begin{equation}\label{eq:Whittakereqn}
\p^2_{\zeta} M_\pm+\left(-{\frac  {1}{4}} +{\frac  {1/4- (1/4-\b^2)}{4m^2\zeta^{2}}}\right)M_\pm=0, \qquad \zeta\in\C.
\end{equation}
The construction of the Green's function is contained in the following result.
\begin{proposition}\label{prop: def Green's Function}
Let  $\ep\in(0,1)$ and $\b^2\neq1/4$. The Green's function $\G_{m,\ep}^\pm$ of  $\textsc{TG}_{m,\ep}^\pm$
is given by
\begin{equation}\label{Green Ray}
\G_{m,\ep}^\pm(y,y_0,z)=\frac{1}{\W_{m,\ep}^\pm(y_0)}
\begin{cases}
\phi_{u,m,\ep}^\pm(y,y_0)\phi_{l,m,\ep}^\pm(z,y_0), \quad &0\leq z\leq y\leq 1,\\
\phi_{u,m,\ep}^\pm(z,y_0)\phi_{l,m,\ep}^\pm(y,y_0), \quad &0\leq y\leq z\leq 1,
\end{cases}
\end{equation}
where  $\phi_{u,m,\ep}^\pm(\cdot,y_0)$ and $\phi_{l,m,\ep}^\pm(\cdot,y_0)$ are two homogeneous solutions to \eqref{eq:TGoperator} such that  
$\phi_{u,m,\ep}^\pm(1,y_0)=0$ and $\phi_{l,m,\ep}^\pm(0,y_0)=0$, respectively, for all $y_0\in [0,1]$. They are explicitly given by
\begin{equation}\label{eq:homoup}
\phi_{u,m,\ep}^\pm(y,y_0):=M_+(1-y_0\pm i \ep)M_-(y-y_0\pm i\ep)-M_-(1-y_0\pm i\ep)M_+(y-y_0\pm i\ep)
\end{equation}
and
\begin{equation}\label{eq:homolow}
\phi_{l,m,\ep}^\pm(y,y_0):=M_+(-y_0\pm i \ep)M_-(y-y_0\pm i\ep)-M_-(-y_0\pm i\ep)M_+(y-y_0\pm i\ep),
\end{equation}
with Wronskian 
\begin{equation}\label{def Wronskian}
\W_{m,\ep}^\pm(y_0):=4 (\mu+i\nu) m\Big(M_+(-y_0\pm i\ep)M_-(1-y_0\pm i\ep)-M_-(-y_0\pm i\ep)M_+(1-y_0\pm i\ep)\Big).
\end{equation} 
Furthermore, we have the relation
$\G_{m,\ep}^+(y,y_0,z)=\overline{\G_{m,\ep}^-(y,y_0,z)}$,
for all $y,y_0,z\in[0,1]$.
\end{proposition}
\begin{proof}
We introduce the variables $\tilde{y}_\pm=2m(y-y_0\pm i\ep)$ and $\tilde{z}_\pm=2m(z-y_0\pm i\ep)$ and we write $\G_{m,\ep}^\pm(y,y_0,z)=\G(\tilde{y}_\pm,\tilde{z}_\pm)$ and rewrite \eqref{eq:TGgreens}
as
\begin{equation}
\p^2_{\tilde{y}_\pm}\G + \left(-\frac{1}{4}+\frac{1/4-(1/4-\b^2)}{\tilde{y}_\pm^2}\right)\G=\frac{1}{4m^2}\delta\l(\frac{1}{2m}\l(\tilde{y}_\pm-\tilde{z}_\pm\r)\r).
\end{equation}
The left-hand side above has precisely the form of \eqref{eq:Whittakereqn}, and therefore  the general solution is given in terms of the homogeneous solutions 
in \eqref{eq:homoup}-\eqref{eq:homolow}
by
\begin{align*}
\G_{m,\ep}^\pm(y,y_0,z)=
\begin{cases}C_1(\widetilde{z}_\pm)\phi_{u,m,\ep}^\pm(y,y_0), &0\leq z\leq y\leq 1,\\
C_2(\widetilde{z}_\pm)\phi_{l,m,\ep}^\pm(y,y_0), &0\leq y\leq z\leq 1,
\end{cases}
\end{align*}
where $C_i$ are constants to be determined. Imposing the continuity and jump conditions of the Green's function, together with basic properties of the Whittaker functions \cite{NIST}, we obtain the desired result.
\end{proof}

We also record the following proposition regarding homogeneous solutions to \eqref{eq:TGoperator}. 
\begin{proposition}\label{prop: generalized eigenfunction}
The unique solutions to the homogeneous \eqref{eq:TGoperator} for $\ep=0$ and $y_0=0,1$ with homogeneous Dirichlet boundary conditions at $y=0,1$ are given by 
\begin{equation}\label{eq: phium}
\phi_{u,m}(y):=M_+(1)M_-(y)-M_-(1)M_+(y)
\end{equation}
and
\begin{equation}\label{eq: philm}
\phi_{l,m}(y):=M_+(1)M_-(1-y)-M_-(1)M_+(1-y).
\end{equation}
\end{proposition}

\subsection{The case $\b^2= 1/4$}
We next provide the Green's function to the Taylor-Goldstein equation in the case $\b^2=1/4$. In this case, the Whittaker equation \eqref{eq:whit} has to be taken for $a=b=0$, and 
 $M_0(\zeta):=M_{0,0}(2m\zeta)$ satisfies
\begin{equation}\label{eq:Whittakereqnbeq}
\p^2_{\zeta} M_0+\left(-{\frac  {1}{4}} +{\frac {1/4}{4m^2\zeta^{2}}}\right)M_0=0, \qquad \zeta\in\C.
\end{equation}
The second independent homogeneous solution from which we build the Green's function is given by $W_0(\zeta):=W_{0,0}(2m\zeta)$, defined to be the unique solution to \eqref{eq:Whittakereqnbeq} such that 
$$W_{0,0}(\zeta)=\sqrt{\tfrac{\zeta}{\pi}}\l( 2\log(2) + \varsigma -\log(\zeta)\r) + O\l(\zeta^\frac32\log(\zeta)\r),$$
as $\zeta\rightarrow 0$, where $\varsigma$ denotes the Euler constant. Apart from the introduction of $W_0$, the result here is similar to that in Proposition \ref{prop: def Green's Function}.
\begin{proposition}\label{Prop def special Green's function}
Let  $\ep\in(0,1)$. The Green's function $\G_{m,\ep}^\pm$ of  $\textsc{TG}_{m,\ep}^\pm$
is given by
\begin{equation}\label{Green Log TG}
\G_{m,\ep}^\pm(y,y_0,z)=\frac{1}{\W_{m,\ep}^\pm(y_0)}
\begin{cases}\phi_{u,m,\ep}^\pm(y,y_0)\phi_{l,m,\ep}^\pm(z,y_0)&0\leq z\leq y\leq 1,\\\phi_{u,m,\ep}^\pm(z,y_0)\phi_{l,m,\ep}^\pm(y,y_0)&0\leq y\leq z\leq 1,\end{cases}
\end{equation}
where  $\phi_{u,m,\ep}^\pm(\cdot,y_0)$ and $\phi_{l,m,\ep}^\pm(\cdot,y_0)$ are two homogeneous solutions to \eqref{eq:TGoperator} such that  
$\phi_{u,m,\ep}^\pm(1,y_0)=0$ and $\phi_{l,m,\ep}^\pm(0,y_0)=0$, respectively, for all $y_0\in [0,1]$. They are explicitly given by
\begin{equation}\label{eq:homoupspecial}
\phi_{u,m,\ep}^\pm(y,y_0):=W_0(1-y_0\pm i \ep)M_0(y-y_0\pm i\ep)-M_0(1-y_0\pm i\ep)W_0(y-y_0\pm i\ep)
\end{equation}
and
\begin{equation}\label{eq:homolowspecial}
\phi_{l,m,\ep}^\pm(y,y_0):=W_0(-y_0\pm i \ep)M_0(y-y_0\pm i\ep)-M_0(-y_0\pm i\ep)W_0(y-y_0\pm i\ep).
\end{equation}
with Wronskian 
\begin{equation}\label{def special Wronskian}
\W_{m,\ep}^\pm(y_0):=\frac{2m}{\sqrt{\pi}}\Big(W_0(-y_0\pm i\ep)M_0(1-y_0\pm i\ep)-M_0(-y_0\pm i\ep)W_0(1-y_0\pm i\ep)\Big).
\end{equation}
Furthermore, we have the relation
$\G_{m,\ep}^+(y,y_0,z)=\overline{\G_{m,\ep}^-(y,y_0,z)}$,
for all $y,y_0,z\in[0,1]$.
\end{proposition}

Similarly, we state the following proposition regarding homogeneous solutions to \eqref{eq:TGoperator} when $\b^2=1/4$.
\begin{proposition}\label{prop: special generalized eigenfunction}
The unique solutions to the homogeneous \eqref{eq:TGoperator} for $\ep=0$ and $y_0=0,1$ with homogeneous Dirichlet boundary conditions at $y=0,1$ are given by 
\begin{equation}\label{eq: phium special}
\phi_{u,m}(y):=W_0(1)M_0(y)-M_0(1)W_0(y)
\end{equation}
and
\begin{equation}\label{eq: philm special}
\phi_{l,m}(y):=W_0(1)M_0(1-y)-M_0(1)W_0(1-y).
\end{equation}
\end{proposition}

\subsection{Derivative formulae for solutions to the Taylor-Goldstein equation}
We finish this section by exhibiting the following useful expressions for various derivatives of $\psi_{m,\ep}^\pm$ and $\rho_{m,\ep}^\pm$.

\begin{proposition}\label{stream derivatives formulae}
Let $\ep\in(0,1)$. Then,
\begin{equation}\label{def p y0 psi}
\begin{aligned}
\p_{y_0}\psi_{m,\ep}^\pm(y,y_0)
&=-\frac{1}{\b^2}\o_m^0(y) + \Big. \B_{m,\ep}^\pm(y,y_0,z)\Big]_{z=0}^{z=1} -\int_0^1\p_y\G_{m,\ep}^\pm(y,y_0,z)F_{m,\ep}^\pm(z,y_0)\d z \\
&\quad+ \int_0^1\G_{m,\ep}^\pm(y,y_0,z)\l(\p_zF_{m,\ep}^\pm(z,y_0)+\p_{y_0}F_{m,\ep}^\pm(z,y_0)\r)\d z,
\end{aligned}
\end{equation}
where $\B_{m,\ep}^\pm(y,y_0):=(\p_y+\p_{y_0})^2\varphi_{m,\ep}^\pm(y,y_0)$. Moreover, 
\begin{equation}\label{def p y0 y0 psi}
\begin{aligned}
\p^2_{y_0}\psi_{m,\ep}(y,y_0)
&= F_{m,\ep}^\pm(y,y_0) -2\Big. \p_y \B_{m,\ep}^\pm(y,y_0,z)\Big]_{z=0}^{z=1} + \Big. \widetilde{\B_{m,\ep}^\pm}(y,y_0,z)\Big]_{z=0}^{z=1}\\
&\quad+ \l(m^2-\frac{\b^2}{(y-y_0\pm i\ep)^2}\r)\int_0^1\G_{m,\ep}^\pm(y,y_0,z)F_{m,\ep}^\pm(z,y_0)\d z \\
&\quad- 2\int_0^1\p_y\G_{m,\ep}^\pm(y,y_0,z)\l(\p_z+\p_{y_0}\r)F_{m,\ep}^\pm(z,y_0) \d z \\
&\quad+\int_0^1\G_{m,\ep}^\pm(y,y_0,z)\l( \p_z+\p_{y_0}\r)^2F_{m,\ep}^\pm(z,y_0) \d z,
\end{aligned}
\end{equation}
where $\widetilde{\B_{m,\ep}^\pm}(y,y_0,z):=\p_z\G_{m,\ep}^\pm(y,y_0,z)\l(\p_z + \p_{y_0}\r)^2\varphi_{m,\ep}^\pm(z,y_0)$.  Additionally,
\begin{equation}\label{def p y psi}
\begin{aligned}
\p_{y}\psi_{m,\ep}^\pm(y,y_0)&=\frac{1}{\b^2}\l(\o_m^0(y)+(y-y_0\pm i\ep)\p_y\o_m^0(y)\r)-\p_y\rho_m^0(y)+ \p_y\varphi_{m,\ep}^\pm(y,y_0)
\end{aligned}
\end{equation}
and 
\begin{equation}\label{def p y0 y psi}
\begin{aligned}
\p_{y_0,y}^2\psi_{m,\ep}^\pm(y,y_0)&=-\frac{1}{\b^2}\p_y\o_m^0(y) - F_{m,\ep}^\pm(y,y_0) +\p_y \B_{m,\ep}^\pm(y,y_0,z)\Big]_{z=0}^{z=1} \\
&\quad+ \int_0^1\p_y\G_{m,\ep}^\pm(y,y_0,z)\l(\p_zF_{m,\ep}^\pm(z,y_0)+\p_{y_0}F_{m,\ep}^\pm(z,y_0)\r)\d z \\
&\quad-\l(m^2-\frac{\b^2}{(y-y_0\pm i\ep)^2}\r)\int_0^1\G_{m,\ep}^\pm(y,y_0,z)F_{m,\ep}^\pm(z,y_0)\d z.
\end{aligned}
\end{equation}
\end{proposition}

\begin{proof}
The formula for $\p_y\psi_{m,\ep}^\pm$ follows from taking a $\p_y$ derivative in \eqref{def psi}. Similarly, once $\p_{y_0}\psi_{m,\ep}^\pm$ is established, the expression for $\p_{y_0,y}\psi_{m,\ep}^\pm$ follows from taking a $\p_y$ derivative in \eqref{def p y0 psi} and noting that $\G_{m,\ep}^\pm$ is the Green's function of the Taylor-Goldstein operator. As for $\p_{y_0}\psi_{m,\ep}^\pm$ and $\p_{y_0}^2\psi_{m,\ep}$ we show these expressions using the Taylor-Goldstein equation and taking $y_0$ and $y$ derivatives there. More precisely, note that $\p_y+\p_{y_0}$ commutes with the Taylor-Goldstein operator \eqref{eq:TGoperator}. As such, $\textsc{TG}_{m,\ep}^\pm \l(\p_y + \p_{y_0}\r)\varphi_{m,\ep}^\pm=\l(\p_y + \p_{y_0}\r)F_{m,\ep}^\pm $ and the first part of the lemma follows, upon noting that 
\begin{equation*}
\p_{y_0}\psi_{m,\ep}^\pm=-\frac{1}{\b^2}\o_m^0+\p_{y_0}\varphi_{m,\ep}^\pm
\end{equation*}
and that
\begin{equation*}
\int_0^1 \G_{m,\ep}^\pm(y,y_0,z)\textsc{TG}_{m,\ep}^\pm\p_z\varphi_{m,\ep}^\pm(z,y_0)\d z = {-\l. \p_z\G_{m,\ep}^\pm(y,y_0,z)\p_z\varphi_{m,\ep}^\pm(z,y_0)\r|_{z=0}^{z=1}} + \p_y\varphi_{m,\ep}^\pm(y,y_0).
\end{equation*}
As for the second part of the lemma, $\textsc{TG}_{m,\ep}^\pm\l( \p_y + \p_{y_0}\r)^2 \varphi_{m,\ep}^\pm=\l(\p_y + \p_{y_0}\r)^2F_{m,\ep}^\pm$,
from which we deduce that
\begin{align*}
\p_{y_0}\l( \p_y + \p_{y_0}\r)\varphi_{m,\ep}^\pm
&=-\p_{y}\l( \p_y + \p_{y_0}\r) \varphi_{m,\ep}^\pm  +\Big[ \p_z \G_{m,\ep}^\pm(y,y_0,z)(\p_{y_0}+\p_z)^2\varphi_{m,\ep}^\pm (z,y_0)\Big]_{z=0}^{z=1} \\
&\quad+ \int_0^1\G_{m,\ep}^\pm(y,y_0,z)\l(\p_z+\p_{y_0}\r)^2F_{m,\ep}^\pm(z,y_0)\d z.
\end{align*}
Now, since $(\p_{y_0}+\p_y)^2\varphi_{m,\ep}^\pm=\p_{y_0}^2\varphi_{m,\ep}^\pm+2\p_{y_0}\p_y\varphi_{m,\ep}^\pm + \p_y^2\varphi_{m,\ep}^\pm$, we observe that for $y=0$ and $y=1$,
\begin{equation*}
\begin{aligned}
(\p_{y_0}+\p_y)^2\varphi_{m,\ep}^\pm(y,y_0)&=-\frac{2}{\b^2}\p_y\o_m^0(y) -F_{m,\ep}^\pm(y,y_0) +2\p_y\Big[ \p_z\G_{m,\ep}^\pm(y,y_0,z)\p_z\varphi_{m,\ep}^\pm(z,y_0)\Big]_{z=0}^{z=1} \\
&\quad+2 \int_0^1\p_y\G_{m,\ep}^\pm(y,y_0,z)\l(\p_z+\p_{y_0}\r)F_{m,\ep}^\pm(z,y_0)\d z \\
\end{aligned}
\end{equation*}
Moreover, from $\textsc{TG}_{m,\ep}^\pm\l(\p_y + \p_{y_0}\r)\varphi_{m,\ep}^\pm=\l(\p_y+\p_{y_0}\r)F_{m,\ep}^\pm$ we can also obtain
\begin{equation*}
\l(\p_y + \p_{y_0}\r)\varphi_{m,\ep}^\pm=\Big[ \p_z\G_{m,\ep}^\pm(y,y_0,z)\p_y\varphi_{m,\ep}^\pm(y,y_0)\Big]_{z=0}^{z=1} + \int	_0^1 \G_{m,\ep}^\pm(y,y_0,z)(\p_{z}+\p_{y_0})F_{m,\ep}^\pm(z,y_0)\d z,
\end{equation*}
so that
\begin{align*}
\p_y \l(\p_y + \p_{y_0}\r)\varphi_{m,\ep}^\pm &=\p_y\Big[ \p_z\G_{m,\ep}^\pm(y,y_0,z)\p_y\varphi_{m,\ep}^\pm(y,y_0)\Big]_{z=0}^{z=1} \\
&\quad+ \int	_0^1 \p_y\G_{m,\ep}^\pm(y,y_0,z)(\p_{z}+\p_{y_0})F_{m,\ep}^\pm(z,y_0)\d z.
\end{align*}
We finish with the observation that $(\p_{y_0}-\p_y)\l(\p_y + \p_{y_0}\r)\varphi_{m,\ep}^\pm=(\p_{y_0}^2-\p_y^2)\varphi_{m,\ep}^\pm$,
that is, 
\begin{align*}
\p_{y_0}^2\varphi_{m,\ep}^\pm &=\p_y^2\varphi_{m,\ep}^\pm+ (\p_{y_0}-\p_y)\l(\p_y + \p_{y_0}\r)\varphi_{m,\ep}^\pm \\
&=\l(m^2 -\b^2\frac{1}{(y-y_0+i\ep)^2}\r)\varphi_{m,\ep}^\pm + F_{m,\ep}^\pm+ \p_{y_0}\l(\p_y + \p_{y_0}\r)\varphi_{m,\ep}^\pm - \p_y \l(\p_y + \p_{y_0}\r)\varphi_{m,\ep}^\pm.
\end{align*}
Gathering the previously obtained terms, the second part of the lemma follows, since $\p^2_{y_0}\psi_{m,\ep}^\pm=\p^2_{y_0}\varphi_{m,\ep}^\pm$.
\end{proof}

With the same ideas as above, we can also find useful expressions for $\p_{y_0}\rho_{m,\ep}^\pm$ and $\p_y\rho_{m,\ep}^\pm$, thanks again to \eqref{eq rho m ep}.
\begin{corollary}\label{density derivative formulae}
Let $\ep\in(0,1)$. Then,
\begin{equation}\label{def p y0 rho}
\begin{aligned}
\p_{y_0}\rho_{m,\ep}^\pm(y,y_0) &= \frac{1}{(y-y_0\pm i\ep)^2}\int_0^1 \G^\pm_{m,\ep}(y,y_0,z) F_{m,\ep}^\pm(z,y_0)\d z + \frac{1}{y-y_0\pm i\ep}\Big. \B_{m,\ep}^\pm(y,y_0,z)\Big]_{z=0}^{z=1}\\
&\quad+ \frac{1}{y-y_0\pm i\ep}\int_0^1 \G^\pm_{m,\ep}(y,y_0,z) \l( \p_zF_{m,\ep}^\pm(z,y_0) + \p_{y_0}F_{m,\ep}^\pm(z,y_0)\r)\d z \\
&\quad- \frac{1}{y-y_0\pm i\ep}\int_0^1 \p_y\G^\pm_{m,\ep}(y,y_0,z) F_{m,\ep}^\pm(z,y_0)\d z,
\end{aligned}
\end{equation} 
and
\begin{equation}\label{def p y rho}
\begin{aligned}
\p_{y}\rho_{m,\ep}^\pm(y,y_0) &= \frac{1}{\b^2}\p_y\o_m^0(y) -\frac{1}{(y-y_0\pm i\ep)^2}\int_0^1 \G^\pm_{m,\ep}(y,y_0,z) F_{m,\ep}^\pm(z,y_0)\d z \\
&\quad+ \frac{1}{y-y_0\pm i\ep}\int_0^1 \p_y\G^\pm_{m,\ep}(y,y_0,z) F_{m,\ep}^\pm(z,y_0)\d z.
\end{aligned}
\end{equation}
\end{corollary}

\section{Bounds on the Green's function  for $\b^2\neq 1/4$} \label{sec: Bounds Greens Function}
This section is devoted to the proof of Theorem \ref{L2 bounds of G}, which provides $L^2$ bounds on the Green's function $\G_{m,\ep}^\pm$. We separate the estimates into
bounds near the critical layer (Section \ref{sub:nearCL}) and away from the critical layer (Section \ref{sub:awayCL}). We wrap up the proof in Section \ref{sub:proofthm2}.

\subsection{Pointwise bounds near the critical layer}\label{sub:nearCL}
The aim is to provide pointwise bounds for the Green's function and its $\p_y$ derivative when both $y$ and $z$ variables are close to the spectral variable $y_0$. 
\begin{proposition}\label{Bounds Small y Small z}
Let  $y,y_0,z\in[0,1]$ such that $m|y-y_0+i\ep|\leq 10\b$ and $m|z-y_0+i\ep|\leq 10\b$.
There exists $\ep_0>0$ such that
\begin{equation*}
|\G_{m,\ep}^+(y,y_0,z)|\lesssim m^{-2\mu} |y-y_0+i\ep|^{\frac12-\mu}|z-y_0+i\ep|^{\frac12-\mu}
\end{equation*}  
and
\begin{equation*}
|\p_y\G_{m,\ep}^+(y,y_0,z)|\lesssim m^{-2\mu} |y-y_0+i\ep|^{-\frac12-\mu}|z-y_0+i\ep|^{\frac12-\mu},
\end{equation*}  
for all $\ep\leq\ep_0$.
\end{proposition}
The proofs depend heavily on the Wronskian associated to the Green's function $\G_{m,\ep}^\pm(y,y_0,z)$ and whether $\b^2>1/4$ or not. We begin with the case in which $\b^2>1/4$, for which $\mu=0$ and $\nu >0$.
\begin{proposition}
Let $\b^2>1/4$. Within the assumptions of Proposition \ref{Bounds Small y Small z}, there exists $\ep_0>0$ such that
\begin{equation*}
|\G_{m,\ep}^+(y,y_0,z)|\lesssim |y-y_0+i\ep|^{\frac12}|z-y_0+i\ep|^{\frac12}
\end{equation*}  
and
\begin{equation*}
|\p_y\G_{m,\ep}^+(y,y_0,z)|\lesssim |y-y_0+i\ep|^{-\frac12}|z-y_0+i\ep|^{\frac12},
\end{equation*}  
for all $\ep\leq\ep_0$.
\end{proposition}
\begin{proof}
Let us assume that $y\leq z$. Then \eqref{Green Ray} tells us that 
\begin{equation}\label{eq:greenylz}
\G_{m,\ep}^+(y,y_0,z)=\frac{1}{W_{m,\ep}^\pm(y_0)}\phi_{u,m,\ep}^+(z,y_0)\phi_{l,m,\ep}^+(y,y_0)
\end{equation}
and we have from Lemma \ref{asymptotic expansion M} that
\begin{equation*}
|\phi_{u,m,\ep}^+(z,y_0)|\lesssim m^\frac12|z-y_0+i\ep|^\frac12\big( |M_-(1-y_0+i\ep)| + |M_+(1-y_0+i\ep)|\big),
\end{equation*}
while
\begin{equation*}
|\phi_{l,m,\ep}^+(y,y_0)|\lesssim m^\frac12|y-y_0+i\ep|^\frac12\big( |M_-(y_0-i\ep)| + |M_+(y_0-i\ep)|\big).
\end{equation*}
The proof follows once we show that
\begin{equation}\label{Wronskian lower bound}
|\W_{m,\ep}^+(y_0)|\geq C_\nu m|M_-(y_0-i\ep)||M_+(1-y_0+i\ep)|
\end{equation}
and 
\begin{equation}\label{Quotient bounds}
\frac{|M_+(y_0-i\ep)|}{|M_-(y_0-i\ep)|}+ \frac{|M_-(1-y_0+i\ep)|}{|M_+(1-y_0+i\ep)|}\lesssim 1.
\end{equation}
To prove the lower bound on the Wronskian, we begin by writing out a suitable expression for $\W^+_{m,\ep}(y_0)$, where we have used the analytic continuation properties of the Whittaker functions $M_\pm$:
\begin{equation*}
\begin{aligned}
\W^+_{m,\ep}(y_0)&=4i\nu m\Big( \e^{\nu\pi}M_-(y_0-i\ep)M_+(1-y_0+i\ep) -\e^{-\nu\pi}M_+(y_0-i\ep)M_-(1-y_0+i\ep)\Big) \\
&=4i\nu mM_-(y_0-i\ep)M_+(1-y_0+i\ep) \l( \e^{\nu\pi} - \e^{-\nu\pi}\frac{M_+(y_0-i\ep)}{M_-(y_0-i\ep)}\frac{M_-(1-y_0+i\ep)}{M_+(1-y_0+i\ep)}\r).
\end{aligned}
\end{equation*}
The proof depends on the location of $y_0\in[0,1]$ as well as on the smallness of $m$. In this direction, let $N_\nu>0$ given in Lemma \ref{growth bounds M large argument}.

\bullpar{Case 1: $m< N_\nu$} Assume that $y_0\leq1/2$ (otherwise we would have $1-y_0\leq 1/2$ and the proof would carry over unaltered). Therefore, it follows that $2m y_0 < N_\nu$ and $m\leq 2m(1-y_0) < 2N_\nu$. Hence, there exists $\ep_0>0$ such that from Lemma \ref{Comparison bounds M small argument} 
\begin{equation*}
\l| \frac{M_+(y_0-i\ep)}{M_+(y_0+i\ep)}\r| \leq \e^{\frac54\nu\pi},
\end{equation*}
and from Lemma \ref{Comparison bounds M order one argument}
\begin{equation*}
\frac{|M_-(1-y_0+i\ep)|}{|M_+(1-y_0+i\ep)|}\leq \e^{\frac14\nu\pi}, 
\end{equation*} 
for all $\ep\leq\ep_0$. Moreover,
\begin{equation*}
|\W_{m,\ep}^+(y_0)|\geq C_\nu m|M_-(y_0-i\ep)||M_+(1-y_0+i\ep)|
\end{equation*}
for $C_\nu=4\nu(\e^{\nu\pi}-\e^{\nu\pi/2})$.

\bullpar{Case 2: $m\geq N_\nu$}
Assume now that $2m y_0\leq N_\nu$. Then, since $m\geq N_\nu$ we have that $2m(1-y_0)\geq N_\nu$. The other case is completely analogous and $m\geq N_\nu$ ensures that $2my_0<N_\nu$ and $2m(1-y_0)<N_\nu$ cannot hold simultaneously for any $y_0\in[0,1]$. Therefore, it follows from Lemma \ref{Comparison bounds M small argument} that 
\begin{equation*}
\l| \frac{M_+(y_0-i\ep)}{M_+(y_0+i\ep)}\r| \leq \e^{\frac54\nu\pi},
\end{equation*}
while from Lemma \ref{growth bounds M large argument} we obtain
\begin{equation*}
\l|\frac{M_-(1-y_0+i\ep)}{M_+(1-y_0+i\ep)}\r|\leq \e^{\frac14\nu\pi},
\end{equation*} 
for all $\ep\leq \ep_0$, for some $\ep_0>0$. The lower bound on $|\W^+_{m,\ep}(y_0)|$ holds for the same $C_\nu$ as above.
\end{proof}
We next consider the case $\b^2<1/4$, for which $\nu=0$ and $0<\mu<1/2$.
\begin{proposition}\label{bounds real G small y small z}
Let $\b^2<1/4$. Within the assumptions of Proposition \ref{Bounds Small y Small z}, there exists $\ep_0>0$ such that
\begin{equation*}
|\G_{m,\ep}^+(y,y_0,z)|\lesssim m^{-2\mu} |y-y_0+i\ep|^{\frac12-\mu}|z-y_0+i\ep|^{\frac12-\mu}
\end{equation*}  
and
\begin{equation*}
|\p_y\G_{m,\ep}^+(y,y_0,z)|\lesssim m^{-2\mu} |y-y_0+i\ep|^{-\frac12-\mu}|z-y_0+i\ep|^{\frac12-\mu},
\end{equation*}  
for all $\ep\leq\ep_0$.
\end{proposition}
\begin{proof}
Let us assume that $y\leq z$ and deal with the expression in \eqref{eq:greenylz}. From Lemma \ref{asymptotic expansion M} we have
\begin{equation*}
|\phi_{u,m,\ep}^+(z,y_0)|\lesssim m^{\frac12-\mu}|z-y_0+i\ep|^{\frac12-\mu}\big( |M_-(1-y_0+i\ep)| + |M_+(1-y_0+i\ep)|\big),
\end{equation*}
while
\begin{equation*}
|\phi_{l,m,\ep}^+(y,y_0)|\lesssim m^{\frac12-\mu}|y-y_0+i\ep|^{\frac12-\mu}\big( |M_-(y_0-i\ep)| + |M_+(y_0-i\ep)|\big),
\end{equation*}
both following from the observation that $(2m|y-y_0+i\ep|)^{2\mu}\leq 10$. Using the analytic continuation properties of the Whittaker functions $M_\pm$ we obtain 
\begin{equation*}
\begin{aligned}
\W_{m,\ep}^+(y_0)&=4\mu m\l( \e^{i\mu\pi}M_+(y_0-i\ep)M_-(1-y_0+i\ep)-\e^{-i\mu\pi}M_-(y_0-i\ep)M_+(1-y_0+i\ep)\r) \\
&=-4\mu m\l( \e^{-i\mu\pi}M_-(y_0-i\ep)M_+(1-y_0+i\ep) - \e^{i\mu\pi}M_+(y_0-i\ep)M_-(1-y_0+i\ep)\r)
\end{aligned}
\end{equation*}
One needs to obtain suitable estimates on several quotients. This is again done considering the location of $y_0\in[0,1]$ and the smallness of $m$. Thus, let $N_{\mu,0}>0$ given in Lemma \ref{growth bounds real M large argument}.

\bullpar{Case 1: $m\leq N_{\mu,0}$} 
Assume initially that $y_0\leq \frac12$. Then, $2my_0\leq N_{\mu,0}$ and $m\leq 2m(1-y_0)\leq 2N_{\mu,0}$.
Assume further that $2my_0\leq \delta_{\mu,1}$ as given in Lemma \ref{Comparison bounds real M small argument}. From Lemma \ref{Comparison bounds real M order one argument} choosing $N_{\mu,1}:=m$ and Lemma \ref{lower bounds M} we have that 
\begin{equation*}
\frac34 \leq \l| \frac{M_+(1-y_0+i\ep)}{M_+(1-y_0)} \r| \leq \frac54,
\end{equation*}
and from Lemma \ref{Comparison bounds real M small argument}
\begin{equation*}
\l| \frac{M_-(1-y_0+i\ep)}{M_+(1-y_0)} \r| = \l| \frac{M_-(1-y_0+i\ep)}{M_-(1-y_0)}\r| \l|  \frac{M_-(1-y_0)}{M_+(1-y_0)} \r| \leq \frac54 M\l(\frac12-\mu,1-2\mu,2N_{\mu,0}\r),
\end{equation*} 
for $\ep\leq \ep_0$ small enough. Additionally, since $2my_0\leq \delta_{\mu,1}$, we have from Lemma \ref{Comparison bounds real M small argument} that 
\begin{equation*}
\l| \frac{M_+(y_0-i\ep)}{M_-(y_0-i\ep)} \r| \leq \frac{1}{5M\l(\frac12-\mu,1-2\mu,2N_1\r)}.
\end{equation*}
With the above comparison estimates at hand, we note that
\begin{equation*}
\begin{aligned}
\e^{-i\mu\pi}&M_-(y_0-i\ep)M_+(1-y_0+i\ep) - \e^{i\mu\pi}M_+(y_0-i\ep)M_-(1-y_0+i\ep) \\
&=\e^{-i\mu\pi}M_-(y_0-i\ep)M_+(1-y_0)\l( \frac{M_+(1-y_0+i\ep)}{M_+(1-y_0)} - \e^{2i\mu\pi}\frac{M_+(y_0-i\ep)}{M_-(y_0-i\ep)}\frac{M_-(1-y_0+i\ep)}{M_+(1-y_0)} \r) \\
\end{aligned}
\end{equation*}
and therefore we can lower bound
\begin{equation*}
|\W_{m,\ep}^+(y_0)|\geq 2\mu m M_+(1-y_0)|M_-(y_0-i\ep)|.
\end{equation*}
The bounds on the Green's functions follow from the lower bound on the Wronskian and the comparison estimates stated above.

Assume now that $2my_0 > \delta_{\mu,1}$. Then, due to Lemma \ref{Comparison bounds real M order one argument} we have both
\begin{equation*}
\l| M_\pm(y_0-i\ep)-M_\pm(y_0)\r|\leq \frac{\sin\mu\pi}{4}\l|M_\pm(y_0)\r|
\end{equation*}
and
\begin{equation*}
\l| M_\pm(1-y_0+i\ep)-M_\pm(1-y_0)\r|\leq \frac{\sin\mu\pi}{4}\l|M_\pm(1-y_0)\r|,
\end{equation*}
for all $\ep\leq \ep_0$. With the observation that
\begin{equation*}
\begin{aligned}
&\l| \e^{-i\mu\pi}M_-(y_0)M_+(1-y_0) - \e^{i\mu\pi}M_+(y_0)M_-(1-y_0) \r| \\
 &\qquad\geq \sin\mu\pi \big( M_-(y_0)M_+(1-y_0) + M_+(y_0)M_-(1-y_0) \big),
\end{aligned}
\end{equation*}
and the expansion
\begin{equation*}
\begin{aligned}
M_-(y_0-i\ep)M_+(1-y_0+i\ep) &= M_-(y_0)M_+(1-y_0)\\
&\quad+\big( M_-(y_0-i\ep) - M_-(y_0) \big)M_+(1-y_0) \\
&\quad+ \big( M_-(y_0-i\ep) - M_-(y_0) \big) \big( M_+(1-y_0+i\ep)-M_+(1-y_0)\big) \\
&\quad + M_-(y_0) \big( M_+(1-y_0+i\ep)-M_+(1-y_0)\big),
\end{aligned}
\end{equation*}
one can lower bound 
\begin{equation*}
|\W_{m,\ep}^+(y_0)|\geq \mu m\sin\mu\pi \big( M_-(y_0)M_+(1-y_0) + M_+(y_0)M_-(1-y_0) \big).
\end{equation*}
As before, we obtain the bound on the Green's function combing the lower bound on the Wronskian with the above comparison estimates.

\bullpar{Case 2: $m\geq N_{\mu,0}$} Assume $2my_0<N_{\mu,0}$. Since $m\geq N_{\mu,0}$ then $2m(1-y_0)\geq N_{\mu,0}$. 
Assume further that $2my_0\leq \delta_{\mu,1}$ as given in Lemma \ref{Comparison bounds real M small argument} and let $C_\mu:=2^{-4\mu}\frac{\Gamma(1-\mu)}{\Gamma(1+\mu)}$. Then, from Lemma \ref{Comparison bounds M small argument} we have that 
\begin{equation*}
\l| \frac{M_+(y_0-i\ep)}{M_-(y_0-i\ep)} \r| \leq \frac13 C_\mu^{-1},
\end{equation*}
while from Lemma \ref{growth bounds real M large argument},
\begin{equation*}
\l| \frac{M_-(1-y_0+i\ep)}{M_+(1-y_0+i\ep)}\r| \leq \frac32 C_\mu.
\end{equation*}
Since we can write
\begin{equation*}
\W_{m,\ep}^+(y_0)=-4\mu m\e^{-i\mu\pi}M_-(y_0-i\ep)M_+(1-y_0+i\ep)\l( 1 - \e^{2i\mu\pi}\frac{M_+(y_0-i\ep)}{M_-(y_0-i\ep)}\frac{M_-(1-y_0+i\ep)}{M_+(1-y_0+i\ep)} \r)
\end{equation*} 
we are able to lower bound
\begin{equation*}
|\W_{m,\ep}^+(y_0)|\geq 2\mu m |M_-(y_0-i\ep)||M_+(1-y_0+i\ep)|,
\end{equation*}
and the estimates on the Green's function follow directly.

On the other hand, if $2my_0\geq \delta_{\mu,1}$, we shall write
\begin{equation*}
\begin{aligned}
\e^{-i\mu\pi}&M_-(y_0-i\ep)M_+(1-y_0+i\ep) - \e^{i\mu\pi}M_+(y_0-i\ep)M_-(1-y_0+i\ep) \\
&=\e^{-i\mu\pi}M_-(y_0)M_+(1-y_0+i\ep) - \e^{i\mu\pi}M_+(y_0)M_-(1-y_0+i\ep) \\
&\quad+ \e^{-i\mu\pi}\big(M_-(y_0-i\ep) -M_-(y_0)\big)M_+(1-y_0+i\ep) \\
&\quad - \e^{i\mu\pi}\big(M_+(y_0-i\ep) -M_+(y_0)\big)M_-(1-y_0+i\ep) \\
&= T_1 + T_2 + T_3.
\end{aligned}
\end{equation*}
and we note that
\begin{equation*}
T_1=M_+(1-y_0+i\ep)\l( \e^{-i\mu\pi}M_-(y_0) - \e^{i\mu\pi}M_+(y_0)\frac{M_-(1-y_0+i\ep)}{M_+(1-y_0+i\ep)} \r),
\end{equation*}
with 
\begin{equation*}
\begin{aligned}
&\Im \l( \e^{-i\mu\pi}M_-(y_0) - \e^{i\mu\pi}M_+(y_0)\frac{M_-(1-y_0+i\ep)}{M_+(1-y_0+i\ep)} \r) \\
&= -\sin\mu\pi \l( M_-(y_0) + M_+(y_0)\Re\l( \frac{M_-(1-y_0+i\ep)}{M_+(1-y_0+i\ep)} \r) + \frac{1}{\tan \mu\pi}M_+(y_0)\Im\l( \frac{M_-(1-y_0+i\ep)}{M_+(1-y_0+i\ep)}\r) \r).
\end{aligned}
\end{equation*}
Once again, from Lemma \ref{growth bounds real M large argument}, we have that 
\begin{equation*}
\l| \Re\l( \frac{M_-(1-y_0+i\ep)}{M_+(1-y_0+i\ep)} \r) - C_\mu  \r| +
\l| \frac{1}{\tan \mu\pi}\Im\l( \frac{M_-(1-y_0+i\ep)}{M_+(1-y_0+i\ep)} \r)  \r| \leq \frac{C_\mu}{4},
\end{equation*} 
so that we can lower bound 
\begin{equation*}
|T_1| \geq \sin\mu\pi |M_+(1-y_0+i\ep)|\l( M_-(y_0) + \frac{C_\mu}{2} M_+(y_0) \r).
\end{equation*}
Next, we shall see that the terms $T_2$ and $T_3$ are sufficiently small so that they can be absorbed by $T_1$. To this end, from Lemma \ref{Comparison bounds real M order one argument} we have that
\begin{equation*}
|T_2|\leq \frac{\sin\mu\pi}{2}|M_+(1-y_0+i\ep)|M_-(y_0),
\end{equation*}
and, combined with Lemma \ref{growth bounds real M large argument}, we also have that
\begin{equation*}
|T_3|\leq \sin\mu\pi \frac{C_\mu}{4}M_+(y_0) |M_+(1-y_0+i\ep)|,
\end{equation*}
for all $\ep\leq \ep_0$ small enough. Hence, we conclude that 
\begin{equation*}
|T_1+T_2+T_3|\geq \sin\mu\pi \l( \frac{1}{2}M_-(y_0) +\frac{C_\mu}{4}M_+(y_0)\r) |M_+(1-y_0+i\ep)|
\end{equation*}
and we lower bound
\begin{equation*}
|\W_{m,\ep}^+(y_0)|\geq \mu m\sin\mu\pi\l( 2M_-(y_0) + C_\mu M_+(y_0)\r) |M_+(1-y_0+i\ep)|.
\end{equation*}
The bounds on the Green's function are a straightforward consequence of the above lower bound $\W_{m,\ep}^+(y_0)$ and the comparison estimates.
\end{proof}

\subsection{Estimates for $\G_{m,\ep}$ away from the critical layer}\label{sub:awayCL}
Throughout this section, let $\ep_0$ be given by Proposition \ref{Bounds Small y Small z} and assume that $m>8\b$. Hence, both $y_0<\tfrac{4\b}{m}$ and $y_0>1-\tfrac{4\b}{m}$ cannot hold simultaneously and through the section we assume without loss of generality that $y0<1-\tfrac{4\b}{m}$. 

The proof of the following results combines an entanglement inequality inspired by \cite{IIJ22} and the estimates from Proposition \ref{Bounds Small y Small z}. Firstly we obtain estimates when $z$ is far from the critical layer, but $y$ is still near the spectral variable $y_0$.
\begin{lemma}\label{Bounds G Small z Large y}
Let $y_0\in[0,1]$ and $0<\ep\leq \ep_0$. For all $z\in[0,1]$ such that $m|z-y_0|\leq 9\b$ we have the following.
\begin{equation*}
\Vert \p_y\G_{m,\ep}^\pm(\cdot,y_0,z)\Vert_{L^2_y(J_3^c)}^2 + m^2\Vert \G_{m,\ep}^\pm(\cdot,y_0,z)\Vert_{L^2_y(J_3^c)}^2 \lesssim m^{-2\mu}|z-y_0\pm i\ep|^{1-2\mu}.
\end{equation*}
\end{lemma}

\begin{proof}
Assume without loss of generality that $y_0<1-\frac{3\b}{m}$. Let $y_2=y_0+\frac{2\b}{m}$ and take $\eta\in C_p^1([y_2,1])$, the space of piecewise continuously differentiable functions. To ease notation, we denote $h(y):=\G_{m,\ep}^+(y,y_0,z)$. Hence $h(y)$ solves
\begin{equation*}
\l( \D_m +\b^2\frac{1}{(y-y_0+i\ep)^2}\r) h=\delta(y-z).
\end{equation*}
Multiplying the equation by $\overline{h}\eta^2$ and integrating from $y_2$ to 1, we find that
\begin{equation*}
-\overline{h}(z)\eta^2(z)\mathcal{H}(z-y_2)=\int_{y_2}^1 |\p_y h|^2\eta^2 + 2\p_yh \overline{h}\p_y\eta\eta +m^2|h|^2\eta^2 -\b^2\frac{|h|^2\eta^2}{(y-y_0+i\ep)^2} \, \d y 
\end{equation*}
and thus
\begin{equation*}
|\overline{h}(z)\eta^2(z)\mathcal{H}(z-y_2)|
\geq \int_{y_2}^1 \frac12|\p_yh|^2\eta^2 +\l(\frac{m^2}{2}\eta^2-2(\p_y\eta)^2\r)|h|^2\, \d y,
\end{equation*}
where we have used Young's inequality and $m|y-y_0+i\ep|\geq2\b$, for all $y\geq y_2$. Here, $\mathcal{H}$ represents the Heavyside function. Now, we shall choose $\eta$ as follows:
\begin{equation*}
\eta(y)= \begin{cases} \frac{m}{\b}(y-y_2),&y\in(y_2,y_2+\frac{\b}{m}),\\1,&y\in(y_2+\frac{\b}{m}, 1)\end{cases}.
\end{equation*}
Note that $\eta$ is a piecewise $C^1$ function such that it is a linear function in $(y_2,y_2+\frac{\b}{m})$ and it is constant in $(y_2+\frac{\b}{m}, 1)$. Hence, 
\begin{equation*}
|h(z)|+\frac{m^2}{\b^2}\int_{y_2}^{y_2+\frac{\b}{m}}|h|^2\d y \geq \frac12 \int_{y_2+\frac{\b}{m}}^1 \l(|\p_y h|^2 + m^2 |h|^2\r) \d y.
\end{equation*}
Using Proposition \ref{Bounds Small y Small z}, we can estimate
\begin{equation*}
\begin{aligned}
|h(z)| + \frac{m^2}{\b^2}\int_{y_2}^{y_2+\frac{\b}{m}}|h(y,y_0,z)|^2\d y &\lesssim m^{-4\mu}|z-y_0+i\ep|^{1-2\mu}\l( 1+ \frac{m^2}{\b^2}\int_{y_2}^{y_2+\frac{\b}{m}}|y-y_0+i\ep|^{1-2\mu}\d y\r) \\
&\lesssim m^{-2\mu}|z-y_0+i\ep|^{1-2\mu}.
\end{aligned}
\end{equation*}
Therefore, since $y_2=y_0+\frac{2\b}{m}$ we have the bound
\begin{equation*}
\int_{y_0+\frac{3\b}{m}}^1 \l(|\p_y h|^2 + m^2 |h|^2\r) \d y \lesssim m^{-2\mu}|z-y_0+i\ep|^{1-2\mu}.
\end{equation*}
and the Lemma follows.
\end{proof}
 
We shall now deduce estimates for $\p_y \G_{m,\ep}^\pm(y,y_0,z)$ when $y$ is still near $y_0$ but $z$ is away from the critical layer. To this end, we shall use the symmetry of the Green's function and the following result.
\begin{lemma}\label{Bounds dzG Large y Small z}
Let $y_0\in[0,1]$ and $0<\ep\leq \ep_0$. For all $z\in[0,1]$ such that $m|z-y_0|\leq 3\b$ we have the following.
\begin{equation*}
\Vert \p_y\p_z\G_{m,\ep}^\pm(\cdot,y_0,z)\Vert_{L^2_y(J_4^c)}^2 + \Vert \p_z\G_{m,\ep}^\pm(\cdot,y_0,z)\Vert_{L^2_y(J_4^c)}^2 \lesssim m^{-2\mu}|z-y_0\pm i\ep|^{-1-2\mu}.
\end{equation*}
\end{lemma}

\begin{proof}
We assume without loss of generality that $y_0\leq 1-\frac{4\b}{m}$. For any $y>z$, we have that $g(y):=\p_z\G_{m,\ep}^+(y,y_0,z)$ solves
\begin{equation*}
\l( \D_m +\b^2\frac{1}{(y-y_0+i\ep)^2}\r) g=0,
\end{equation*}
with $g(1)=0$. Multiplying the equation by $\overline{g}\eta^2$ and integrating from $y_2=y_0+\frac{7\b}{2m}>z$ to 1, we find that
\begin{equation*}
\begin{aligned}
0&=\int_{y_2}^1 |\p_y g|^2\eta^2 + 2\p_yh \overline{g}\p_y\eta\eta +m^2|g|^2\eta^2 -\b^2\frac{|g|^2\eta^2}{(y-y_0+i\ep)^2} \, \d y \\
&\geq \int_{y_2}^1 \frac12|\p_yg|^2\eta^2 +\l(\frac{m^2}{2}\eta^2-2(\p_y\eta)^2\r)|g|^2\, \d y,
\end{aligned}
\end{equation*}
where we have used Young's inequality and $m|y-y_0|\geq 2\b$, for all $y\geq y_2$. For
\begin{equation*}
\eta(y)= \begin{cases}\frac{2m}{\b}(y-y_2),&y\in(y_2,y_2+\frac{\b}{2m}),\\1,&y\in(y_2+\frac{\b}{2m}, 1),\end{cases}.
\end{equation*}
we get
\begin{equation*}
\frac{m^2}{\b^2}\int_{y_2}^{y_2+\frac{\b}{2m}}|g|^2\d y \geq \frac12 \int_{y_2+\frac{\b}{2m}}^1 \l(|\p_y g|^2 + m^2 |g|^2\r) \d y.
\end{equation*}
Using Proposition \ref{Bounds Small y Small z}, we can estimate
\begin{equation*}
\begin{aligned}
\frac{m^2}{\b^2}\int_{y_2}^{y_2+\frac{\b}{2m}}|g(y,y_0,z)|^2\d y &\lesssim m^{-4\mu}|z-y_0+i\ep|^{-1-2\mu}\frac{m^2}{\b^2}\int_{y_2}^{y_2+\frac{\b}{2m}}|y-y_0+i\ep|^{1-2\mu}\d y\\
&\lesssim m^{-2\mu}|z-y_0+i\ep|^{-1-2\mu}.
\end{aligned}
\end{equation*}
Now, $y_2 +\frac{\b}{2m}=y_0+\frac{4\b}{m}$ so that
\begin{equation*}
\int_{y_0+\frac{4\b}{m}}^1 \l(|\p_y g|^2 + m^2 |g|^2\r) \d y \lesssim m^{-2\mu}|z-y_0+i\ep|^{-1-2\mu}.
\end{equation*}
The proof is finished.
\end{proof}
The next corollary is a direct consequence of the above Lemma together with the observation that once the estimate for $\p_z\G_{m,\ep}^\pm(y,y_0,z)$ is established, the estimate of $\p_y\G_{m,\ep}^\pm(y,y_0,z)$ follows from the fact that, since $\G_{m,\ep}^\pm(y,y_0,z)=\G_{m,\ep}^\pm(z,y_0,y)$, then $(\p_y\G_{m,\ep}^\pm)(y,y_0,z)=(\p_z\G_{m,\ep}^\pm)(z,y_0,y)$.

\begin{corollary}\label{Bounds dyG Small y Large z}
Let $y_0\in[0,1]$ and $0<\ep\leq \ep_0$. For all $y\in[0,1]$ such that $m|y-y_0|\leq 3\b$ we have that.
\begin{equation*}
\Vert \p_y\G_{m,\ep}^\pm(y,y_0,\cdot)\Vert_{L^2_z(J_4^c)}\lesssim \frac{1}{m^{1+\mu}}|y-y_0\pm i\ep|^{-\frac12-\mu}.
\end{equation*}
\end{corollary}

\subsection{Proof of Theorem \ref{L2 bounds of G}}\label{sub:proofthm2}
Let $0<\ep\leq{\ep_0}\leq\frac{\b}{m}$ and assume that $m|y-y_0|\leq 3\b$. For $m\leq 8\b$, the Theorem follows directly from Proposition \ref{Bounds Small y Small z}. Hence, we consider for $m>8\b$ and we note that 
\begin{equation*}
\Vert \G_{m,\ep}^\pm\Vert_{L^2_z}\leq \Vert \G_{m,\ep}^\pm\Vert_{L^2_z(J_3)} + \Vert \G_{m,\ep}^\pm\Vert_{L^2_z(J_3^c)}, \quad \Vert \p_y\G_{m,\ep}^\pm\Vert_{L^2_z}\leq \Vert \p_y\G_{m,\ep}^\pm\Vert_{L^2_z(J_4)} + \Vert \p_y\G_{m,\ep}^\pm\Vert_{L^2_z(J_4^c)}.
\end{equation*}
Now, the bounds for $\Vert \G_{m,\ep}^\pm\Vert_{L^2_z(J_3)}$ and $\Vert \p_y\G_{m,\ep}^\pm\Vert_{L^2_z(J_4)}$ follow from Proposition \ref{Bounds Small y Small z}, while the estimate for $\Vert \G_{m,\ep}^\pm\Vert_{L^2_z(J_3^c)}$ is given in Lemma \ref{Bounds G Small z Large y} due to the $y,z$ symmetry of the Green's function and the estimate for $\Vert \p_y\G_{m,\ep}^\pm\Vert_{L^2_z(J_4^c)}$ is given by Corollary \ref{Bounds dyG Small y Large z}. The theorem follows.

\section{Bounds on the Green's function for $\b^2= 1/4$}\label{sec: Bounds Greens Function special} 
This section studies and obtains $L^2$ bounds on the Green's function for the case $\b^2=1/4$. Most of the results and proof are analogous to the ones presented in Section \ref{sec: Bounds Greens Function} above, so we limit ourselves to present  the statements we use and comment on the main ingredients of the proof.
\begin{theorem}\label{L2 bounds G special}
There exists $\ep_0>0$ such that for all $\ep\in(0, \ep_0)$ and for all $y,y_0\in[0,1]$ such that $m|y-y_0|\leq 3\b$, we have 
\begin{equation*}
|y-y_0 \pm i\ep|^{-\frac12} \Vert\G^\pm_{m,\ep}(y,y_0,\cdot)\Vert_{L^2_z}+ |y-y_0 \pm i\ep|^\frac12 \Vert \p_y\G^\pm_{m,\ep}(y,y_0,\cdot)\Vert_{L^2_z}\lesssim \frac{1}{m}\l( 1 + \big| \log m\l|y-y_0 \pm i\ep\r|\big|\r)
\end{equation*}
\end{theorem}
In comparison with Theorem \ref{L2 bounds of G}, we have a logarithmic correction to the behavior near the critical layer.
The proof is omitted as it is analogous to that of Theorem \ref{L2 bounds of G}, once all the intermediate steps are established. 
The rest of this section is devoted to the proof of such steps, to be compared with the analogous one of Section \ref{sec: Bounds Greens Function}.

\subsection{Estimates near the critical layer}
Using the analytic continuation properties from Lemma \ref{analytic continuation}, we can write the Wronskian as
\begin{equation*}
\begin{aligned}
\W_{m,\ep}^+(y_0)&=\frac{2im}{\sqrt{\pi}}\big( M_0(1-y_0+i\ep)W_0(y_0-i\ep)-W_0(1-y_0+i\ep)M_0(y_0-i\ep)\big) \\
&\quad+ 2m M_0(1-y_0+i\ep)M_0(y_0-i\ep).
\end{aligned}
\end{equation*}
We then have the following.
\begin{proposition}\label{Bounds special G Small y Small z}
Let  $y,y_0,z\in[0,1]$ such that $m|y-y_0+i\ep|\leq 10\b$ and $m|z-y_0+i\ep|\leq 10\b$.
There exists $\ep_0>0$ such that
\begin{equation*}
|\G_{m,\ep}^+(y,y_0,z)|\lesssim |y-y_0+i\ep|^{\frac12}|z-y_0+i\ep|^{\frac12}\l( 1 +  \big|\log (m|y-y_0+i\ep)|\big|\r)\l( 1 +  \big|\log (m|z-y_0+i\ep)|\big|\r)
\end{equation*}  
and
\begin{equation*}
|\p_y\G_{m,\ep}^+(y,y_0,z)|\lesssim |y-y_0+i\ep|^{-\frac12}|z-y_0+i\ep|^{\frac12}\l( 1 + \big|\log (m|y-y_0+i\ep)|\big|\r)\l( 1 + \big|\log (m|z-y_0+i\ep)|\big|\r)
\end{equation*}  
for all $\ep\leq\ep_0$.
\end{proposition}
\begin{proof}
Assume without loss of generality that $y\leq z$. From the asymptotic expansions given by Lemma \ref{asymptotic expansion M}, we have that 
\begin{equation*}
|\phi_{u,m,\ep}^+(z,y_0)|\lesssim (2m|z-y_0+i\ep|)^\frac12\l( 1+ \log (2m|z-y_0+i\ep)|\r)  \big(|W_0(1-y_0+i\ep)|+|M_0(1-y_0+i\ep)|\big) , 
\end{equation*}
while
\begin{equation*}
|\phi_{l,m,\ep}^+(y,y_0)|\lesssim (2m|y-y_0+i\ep|)^\frac12\l( 1+ \log (2m|y-y_0+i\ep)|\r)  \big(|W_0(y_0-i\ep)|+|M_0(y_0-i\ep)|\big).
\end{equation*}
The proposition follows from the estimates on the Wronskian given in the lemma below.
\end{proof}

\begin{lemma}\label{lower bound special Wronskian}
Let  $y_0\in[0,1]$. There exists $0<\ep_0\leq\frac{\b}{m}$ and $C>0$ such that
\begin{equation*}
|\W_{m,\ep}^+(y_0)|\geq C_\nu m|M_0(y_0-i\ep)||M_0(1-y_0+i\ep)|,
\end{equation*}
 for all $\ep\leq \ep_0$.
\end{lemma}

\begin{proof}
The proof follows from treating the next two cases. Let $N_0>0$ be given as in Lemma \ref{growth bounds special M large argument}.

\bullpar{Case 1: $m<N_0$} Assume that $y_0\leq \frac12$. Then $2my_0< N_0$ and $m\leq 2m(1-y_0)< 2N_0$.
Assume further that $2my_0\leq \delta_1$ given by Lemma \ref{Comparison bounds special M small argument}. Then,
\begin{equation*}
\W_{m,\ep}^+(y_0) = \frac{2mi}{\sqrt{\pi}}M_0(1-y_0+i\ep)W_0(y_0-i\ep)\l(1-\frac{W_0(1-y_0+i\ep)}{M_0(1-y_0+i\ep)}\frac{M_0(y_0-i\ep)}{W_0(y_0-i\ep)} -i \sqrt{\pi}\frac{M_0(y_0-i\ep)}{W_0(y_0-i\ep)}\r)
\end{equation*}
and from Lemma \ref{Comparison bounds special M order one argument} and Lemma \ref{Comparison bounds special M small argument} we have  
\begin{equation*}
\l| \frac{W_0(1-y_0+i\ep)}{M_0(1-y_0+i\ep)} \r| \leq C_0, \quad \l| \frac{M_0(y_0-i\ep)}{W_0(y_0-i\ep)} \r| \leq \frac{1}{2(C_0+\sqrt\pi)},
\end{equation*}
for all $\ep\leq \ep_0$, from which the lower bounds on the Wronskian follows.

Similarly, assume now that $\delta_1 <2my_0< N_0$, in this case we write
\begin{equation*}
\begin{aligned}
\W_{m,\ep}^+(y_0)&=2m M_0(1-y_0+i\ep)M_0(y_0-i\ep)\l( 1 + \frac{i}{\sqrt\pi}\l(\frac{W_0(y_0-i\ep)}{M_0(y_0-i\ep)} - \frac{W_0(1-y_0+i\ep)}{M_0(1-y_0+i\ep}\r) \r)
\end{aligned}
\end{equation*}
and we further note that, for all $\ep\leq \ep_0$, 
\begin{equation*}
\begin{aligned}
&\l| 1 + \frac{i}{\sqrt\pi}\l(\frac{W_0(y_0-i\ep)}{M_0(y_0-i\ep)} - \frac{W_0(1-y_0+i\ep)}{M_0(1-y_0+i\ep}\r) \r| \\
&\qquad\geq  1 - \frac{1}{\sqrt\pi}\l( \l|\Im\l( \frac{W_0(y_0-i\ep)}{M_0(y_0-i\ep)} \r)\r| + \l|  \Im \l(\frac{W_0(1-y_0+i\ep)}{M_0(1-y_0+i\ep}\r)\r| \r) \\
&\qquad\geq \frac12,
\end{aligned}
\end{equation*}
due to the estimates from Lemma \ref{Comparison bounds special M order one argument}. The lower bound on the Wronskian follows as before.

\bullpar{Case 2: $m\geq N_0$}
Under the assumption that $2m(1-y_0)\geq N_0$ and that $2m(y_0-i\ep)\leq \delta_1$ we can write
\begin{equation*}
\W_{m,\ep}^+(y_0) = \frac{2mi}{\sqrt{\pi}}M_0(1-y_0+i\ep)W_0(y_0-i\ep)\l(1-\frac{W_0(1-y_0+i\ep)}{M_0(1-y_0+i\ep)}\frac{M_0(y_0-i\ep)}{W_0(y_0-i\ep)} -i \sqrt{\pi}\frac{M_0(y_0-i\ep)}{W_0(y_0-i\ep)}\r)
\end{equation*}
and we have that
\begin{equation*}
\l| \frac{W_0(1-y_0+i\ep)}{M_0(1-y_0+i\ep)} \r| \leq \sqrt\pi, \quad \l| \frac{M_0(y_0-i\ep)}{W_0(y_0-i\ep)} \r| \leq \frac{1}{4\sqrt\pi},
\end{equation*}
from which we obtain the lower bound
\begin{equation*}
|\W_{m,\ep}^+(y_0)|\geq \frac{4m}{\sqrt\pi}|M_0(1-y_0+i\ep)||W_0(y_0-i\ep)|.
\end{equation*}

Now, when $2my_0\geq \delta_1$, we write
\begin{equation*}
\begin{aligned}
\W_{m,\ep}^+(y_0)&=2m M_0(1-y_0+i\ep)M_0(y_0-i\ep)\l( 1 + \frac{i}{\sqrt\pi}\l(\frac{W_0(y_0-i\ep)}{M_0(y_0-i\ep)} - \frac{W_0(1-y_0+i\ep)}{M_0(1-y_0+i\ep}\r) \r)
\end{aligned}
\end{equation*}
and we further note that, for all $\ep\leq \ep_0$, 
\begin{equation*}
\begin{aligned}
&\l| 1 + \frac{i}{\sqrt\pi}\l(\frac{W_0(y_0-i\ep)}{M_0(y_0-i\ep)} - \frac{W_0(1-y_0+i\ep)}{M_0(1-y_0+i\ep}\r) \r| \\
&\qquad\qquad\geq  1 - \frac{1}{\sqrt\pi}\l( \l|\Im\l( \frac{W_0(y_0-i\ep)}{M_0(y_0-i\ep)} \r)\r| + \l|  \frac{W_0(1-y_0+i\ep)}{M_0(1-y_0+i\ep}\r| \r) \\
&\qquad\qquad\geq \frac12,
\end{aligned}
\end{equation*}
due to the estimates from Lemma \ref{Comparison bounds special M order one argument} and Lemma \ref{growth bounds special M large argument}.
\end{proof}
\subsection{Estimates for $\G_{m,\ep}$ away from the critical layer}
Throughout this section, let $\ep_0$ be given by Lemma \ref{lower bound special Wronskian} and let  $m>8\b$.
\begin{lemma}\label{Bounds special G Small z Large y}
Let $y_0\in[0,1]$ and $0<\ep\leq \ep_0$. For all $z\in[0,1]$ such that $m|z-y_0|\leq 9\b$ we have 
\begin{equation*}
\Vert \p_y\G_{m,\ep}^\pm(\cdot,y_0,z)\Vert_{L^2_y(J_3^c)}^2 +  m^2\Vert \G_{m,\ep}^\pm(\cdot,y_0,z)\Vert_{L^2_y(J_3^c)}^2 \lesssim |z-y_0\pm i\ep|\l( 1 + \big| \log (m|z-y_0 \pm i\ep|)\big|\r)^2.
\end{equation*}
\end{lemma}

\begin{proof}
We comment the case $y_0<1-\frac{3\b}{m}$. The proof goes on the same spirit as the one for Lemma \ref{Bounds G Small z Large y}. For $y_2=y_0+\frac{2\b}{m}$ and $h(y)=\G_{m,\ep}^+(y,y_0,z)$, introducing a suitable cut-off function we have that
\begin{equation*}
|h(z)|+\frac{m^2}{\b^2}\int_{y_2}^{y_2+\frac{\b}{m}}|h|^2\d y \geq \frac12 \int_{y_2+\frac{\b}{m}}^1 \l(|\p_y h|^2 + m^2 |h|^2\r) \d y.
\end{equation*}
Using Proposition \ref{Bounds special G Small y Small z}, we estimate
\begin{equation*}
\begin{aligned}
&|h(z)| + \frac{m^2}{\b^2}\int_{y_2}^{y_2+\frac{\b}{m}}|h(y,y_0,z)|^2\d y \\
&\lesssim \frac{m^2}{\b^2}|z-y_0+i\ep|\l( 1 + \big| \log (m|z-y_0+i\ep|)\big|\r)^2\int_{y_2}^{y_2+\frac{\b}{m}}|y-y_0+i\ep|\l( 1 + \big| \log (m|y-y_0+i\ep|)\big|\r)^2\d y \\
&\quad +|z-y_0+i\ep|\l( 1 + \big| \log (m|z-y_0+i\ep|)\big|\r)^2 \\
&\lesssim |z-y_0+i\ep|\l( 1 + \big| \log (m|z-y_0+i\ep|)\big|\r)^2,
\end{aligned}
\end{equation*}
since $1\leq m|y-y_0+i\ep|\leq 2$, for all $y\in\l[y_2, y_2+\frac{\b}{m}\r]$. Therefore, recalling $y_2=y_0+\frac{2\b}{m}$ we have the bound
\begin{equation*}
\int_{y_0+\frac{3\b}{m}}^1 \l(|\p_y h|^2 + m^2 |h|^2\r) \d y \lesssim |z-y_0+i\ep|\l( 1 + \big| \log (m|z-y_0+i\ep|)\big|\r)^2
\end{equation*}
and the proof follows.
\end{proof}

We next provide an intermediate result towards estimates for $\Vert \p_y\G_{m,\ep}^\pm(y,y_0,\cdot)\Vert_{L^2_z(J_4^c)}$.
\begin{lemma}\label{Bounds special dzG Large y Small z}
Let $y_0\in[0,1]$ and $0<\ep\leq \ep_0$. For all $z\in[0,1]$ such that $m|z-y_0|\leq 3\b$ we have 
\begin{equation*}
\Vert \p_y\p_z\G_{m,\ep}^\pm(\cdot,y_0,z)\Vert_{L^2_y(J_4^c)}^2 + \Vert \p_z\G_{m,\ep}^\pm(\cdot,y_0,z)\Vert_{L^2_y(J_4^c)}^2 \lesssim |z-y_0\pm i\ep|^{-1}\l( 1 + \big| \log (m|z-y_0+i\ep|)\big|\r)^2.
\end{equation*}
\end{lemma}

\begin{proof}
From the proof of Lemma \ref{Bounds dzG Large y Small z}, for $g(y):=\p_z\G_{m,\ep}^+(y,y_0,z)$ we have that
\begin{equation*}
\frac{m^2}{\b^2}\int_{y_2}^{y_2+\frac{\b}{2m}}|g|^2\d y \geq \frac12 \int_{y_2+\frac{\b}{2m}}^1 \l(|\p_y g|^2 + m^2 |g|^2\r) \d y.
\end{equation*}
Using Proposition \ref{Bounds special G Small y Small z} we estimate
\begin{equation*}
\begin{aligned}
\frac{m^2}{\b^2}\int_{y_2}^{y_2+\frac{\b}{2m}}|g(y,y_0,z)|^2\d y &\lesssim |z-y_0+i\ep|^{-1}\l( 1 + \big| \log (m|z-y_0+i\ep|)\big|\r)^2.
\end{aligned}
\end{equation*}
Therefore, since $y_2 +\frac{\b}{2m}=y_0+\frac{4\b}{m}$ we have the bound
\begin{equation*}
\int_{y_0+\frac{4\b}{m}}^1 \l(|\p_y g|^2 + m^2 |g|^2\r) \d y \lesssim |z-y_0+i\ep|^{-1}\l( 1 + \big| \log (m|z-y_0+i\ep)|\big|\r)^2,
\end{equation*}
and the lemma follows.
\end{proof}

We finish the section with the estimates for $\Vert \p_y\G_{m,\ep}^\pm(y,y_0,\cdot)\Vert_{L^2_z(J_4^c)}$, which are deduce using the symmetry properties of the Green's function as in Corollary \ref{Bounds dyG Small y Large z} and are given in the next result.

\begin{corollary}\label{Bounds special dyG Small y Large z}
Let $y_0\in[0,1]$ and $0<\ep\leq \ep_0$. For all $y\in[0,1]$ such that $m|y-y_0|\leq 3\b$ we have 
\begin{equation*}
\Vert \p_y\G_{m,\ep}^\pm(y,y_0,\cdot)\Vert_{L^2_z(J_4^c)}\lesssim \frac{1}{m}|y-y_0\pm i\ep|^{-\frac12}\l( 1 + \big| \log (m|y-y_0+i\ep|)\big|\r).
\end{equation*}
\end{corollary}

\section{Contour integral reduction}\label{sec:spectral}
In this section, we study the contour integration that is present in the Dunford's formula (see \eqref{eq:Dunford})
\begin{equation}\label{Dunford's Formula}
\begin{pmatrix}
\psi_m(t,y) \\ \rho_m(t,y)
\end{pmatrix} = \frac{1}{2\pi i} \int_{\p\O}\e^{-imct} \mathcal{R}(c,L_m)\begin{pmatrix}
\psi_m^0 \\ \rho_m^0
\end{pmatrix} \,\d c,
\qquad 
{\mathcal{R}(c,L_m):=(c-L_m)^{-1}},
\end{equation}
where $\O$ is any domain containing $\sigma(L_m)$, the spectrum of the linearized operator $L_m$ in \eqref{eq:linOP}. The main goal of this section is, under suitable conditions on the initial data, to reduce the above contour integration to a much simpler integration along the essential spectrum $\sigma_{ess}(L_m)=[0,1]$.

As the domain of integration  we take the rectangle $\Omega= [-\beta/m,1+\beta/m]\times [-\beta/m,\beta/m]$, and further split in the regions,
\begin{equation*}
\begin{aligned}
R_{0}&=\l\lbrace c=y_0 + is\in \C: -\frac{\b}{m}\leq y_0 \leq 0, \, 0\leq |s|\leq\frac{\b}{m} \r\rbrace, \\
R_{ess}&=\l\lbrace c=y_0 + is\in \C: 0\leq y_0 \leq 1, \, 0\leq |s|\leq\frac{\b}{m} \r\rbrace, \\
R_{1}&=\l\lbrace c=y_0 + is\in \C: 1\leq y_0 \leq 1+\frac{\b}{m}, \, 0\leq |s|\leq\frac{\b}{m} \r\rbrace, \\
\end{aligned}
\end{equation*}
so that \eqref{Dunford's Formula} becomes 
\begin{equation*}
\begin{pmatrix}
\psi_m(t,y) \\ \rho_m(t,y)
\end{pmatrix} = \frac{1}{2\pi i} \l(\int_{\p R_0} + \int_{\p R_{ess}} + \int_{\p R_1}\r)\e^{-imct} \mathcal{R}(c,L_m)\begin{pmatrix}
\psi_m^0 \\ \rho_m^0
\end{pmatrix} \,\d c.
\end{equation*}
The decomposition of $\Omega$ is depicted in  Figure \ref{fig:spectraldecomposition} below.
\begin{figure}[h!]
  \centering
    \includegraphics[width=0.5\linewidth]{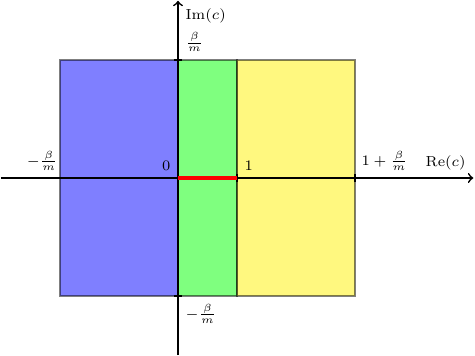}
  \caption{The domain of integration $\Omega$ and the decomposition into \blue{$R_0$}, \green{$R_{ess}$} and \yellow{$R_{1}$}.}
  \label{fig:spectraldecomposition}
\end{figure}

The goal of the next three sections is to show the following result, which amounts to reduce our contour integration as in \eqref{eq:psirholim}. 
\begin{proposition}\label{prop: contour limiting absorption principle}
Assume that the pair of initial data $(\o_m^0, \rho_m^0)$ is orthogonal to the subspace generated by the eigenfunctions of $L_m$. Then,
\begin{equation*}
\l\Vert \int_{\p R_0}\e^{-imct}\mathcal{R}(c,L_m)\d c \r\Vert_{L^2_y} + \l\Vert \int_{\p R_1}\e^{-imct}\mathcal{R}(c,L_m)\d c \r\Vert_{L^2_y}=0.
\end{equation*}
Moreover,
\begin{equation*}
\begin{aligned}
\int_{\p R_{ess}}\e^{-imct}\mathcal{R}(c,L_m)\d c &=\frac{1}{2\pi i }\lim_{\ep\rightarrow 0}\int_0^1 \e^{-imy_0t}\big((-y_0-i\ep+L_m)^{-1}\big. \\
&\qquad\qquad\qquad\qquad\qquad\big.-(-y_0+i\ep+L_m)^{-1}\big)\begin{pmatrix}
\psi_m^0 \\ \rho_m^0 
\end{pmatrix}\, \d y_0.
\end{aligned}
\end{equation*}
\end{proposition}
The description of the spectrum in Theorem \ref{thm:spectralpic} will then be clear from the following three sections. As a first step
towards proving Proposition \ref{prop: contour limiting absorption principle} we show that $\sigma(L_m)\subset \Omega$. 
\begin{lemma}
Let $c\in \C\setminus \Omega$. Then, $\l(c - L_m\r)^{-1}$ exists and
\begin{equation*}
\l\Vert \l(c - L_m\r)^{-1}\begin{pmatrix}
\psi_m^0 \\ \rho_m^0 
\end{pmatrix}\r\Vert_{L^2}\lesssim \Vert \o_m^0\Vert_{L^2_y} + \Vert \rho_m^0\Vert_{L^2_y}
\end{equation*}
\end{lemma}
\begin{proof}
For any $c\in \C\setminus \Omega$, we note that $\tfrac{1}{|y-c|}\leq \tfrac{m}{\b}$, for all $y\in[0,1]$ and we define $\psi_m(y,c)$ as the unique solution to
\begin{equation*}
\D_m\psi_{m}+\b^2\frac{\psi_{m}}{(y-c)^2}=\frac{\o_m^0}{y-c}-\b^2\frac{\rho_m^0}{(y-c)^2}.
\end{equation*}
with homogeneous boundary conditions $\psi_m(0,c)=\psi_m(1,c)=0$ given by standard ODE theory. We also define
\begin{equation*}
\rho_{m}(y,c)=\frac{1}{y-c}\l( \rho_m^0(y)+\psi_{m}(y,c)\r)
\end{equation*}
and it is straightforward to see that
\begin{equation*}
( -c+L_m)\begin{pmatrix}
\psi_{m}(y,c) \\ \rho_{m}(y,c)
\end{pmatrix}=\begin{pmatrix}
\psi_m^0(y) \\ \rho_m^0(y)
\end{pmatrix}.
\end{equation*}
Therefore, $(-c+L_m)$ is invertible and the desired resolvent estimates follow from usual energy estimates on the equation, that is, multiply the equation by $\overline{\psi_m(y,c)}$ integrate by parts and absorb the potential term.
\end{proof}
In order to show that the contributions from $\p R_0$ and $\p R_1$ vanish, we study the resolvent operator for the cases $\b^2>1/4$, $\b^2=1/4$ and $\b^2<1/4$ separately.

\subsection{Integral reduction for $\b^2>1/4$: discrete and embedded eigenvalues}
The classical Miles-Howard stability criterion \cite{Miles,Howard} rules out the existence of unstable modes when $\b^2\geq1/4$. That is, any eigenvalue $c\in\C$ of $L_m$ must have $\Im(c)=0$.
We next find and characterize the discrete set of real isolated eigenvalues that accumulate towards the essential spectrum, that is, towards 0 and 1. Our study involves a precise understanding of the Wronskian when $\ep=0$. For this, we denote
\begin{equation}\label{eq: Wc}
\W_m(c):=M_+(c)M_-(1+c)-M_-(c)M_+(1+c)
\end{equation}
for all $c>0$ and we note from \eqref{def Wronskian} that $\W_{m,0}^\pm(-c)=4i\nu m\W_m(c)$. We state the following.

\begin{proposition}\label{prop discrete eigenvalue and peak}
There exists sequences $\lbrace p_k\rbrace_{k\geq 1}$ and $\lbrace q_k\rbrace_{k\geq 1}$ of strictly positive real numbers such that $p_k,q_k\to 0$ as $k\to\infty$ and
\begin{equation*}
|\W_m(p_k)|=2|M_+(p_k)||M_+(1+p_k)|, \quad \W_m(q_k)=0,
\end{equation*}
for all $k\geq 1$. 
\end{proposition}
\begin{proof}
For any $c>0$, from \eqref{eq: Wc} we have that 
\begin{equation*}
\W_m(c)= -2i\Im\big(M_-(c)M_+(1+c)\big),
\end{equation*}
where further $M_-(c)=\overline{M_+(c)}=|M_+(c)|\e^{-i\text{Arg}(M_+(c))}$ and $M_+(1+c)=|M_+(1+c)|\e^{i\text{Arg}(M_+(1+c))}$. For $x>0$, we define $\Theta(x)=\text{Arg}(M_+(x))$ and we write
\begin{equation*}
\W_m(c)=-2i|M_+(c)||M_+(1+c)|\sin\l( \Theta(1+c)-\Theta(c)\r).
\end{equation*}
The proposition follows if we can find some integer $k_0\geq 0$ and two sequences $\lbrace p_k\rbrace_{k\geq 1}$ and $\lbrace q_k\rbrace_{k\geq 1}$ of strictly positive real numbers such that
\begin{equation*}
\Theta(1+p_k) - \Theta(p_k) = \l(k+k_0+\frac12\r)\pi, \qquad \Theta(1+q_k) - \Theta(q_k) = (k+k_0)\pi,
\end{equation*}
for all $k\geq 1$. To this end, given the Wronskian properties of the pair $M_+(x)$ and $M_-(x)$ from \cite{NIST}, we note that for all $x>0$
\begin{equation*}
\begin{aligned}
-2i\nu &= M_+(x)M_-'(x) - M_+'(x)M_-(x) \\
&= |M_+(x)|\e^{i\Theta(x)}\l( |M_+(x)|'\e^{-i\Theta(x)} -i\Theta'(x)|M_+(x)|\e^{-i\Theta(x)}\r) \\
&\quad - |M_+(x)|\e^{-i\Theta(x)}\l(|M_+(x)|'\e^{i\Theta(x)} + i\Theta'(x)|M_+(x)|\e^{i\Theta(x)}\r) \\ 
&=-2i\Theta'(x)|M_+(x)|^2,
\end{aligned}
\end{equation*}
and thus, $\Theta'(x)=\frac{\nu}{|M_+(x)|^2}>0$. Hence, for all $c>0$ we define
\begin{equation*}
r(c):=\Theta(1+c)-\Theta(c) = \nu \int_c^{1+c}\frac{1}{|M_+(x)|^2}\d x.
\end{equation*}
Note that $r(c)$ is continuous for all $c>0$ and strictly decreasing. This follows from $|M_+(x)|$ being strictly  increasing, see Lemma \ref{lower bounds M}. Moreover, also from Lemma \ref{lower bounds M}, we have that
\begin{equation*}
r(c)\gtrsim_\nu \int_c^{1+c} \frac{1}{x}\d x \gtrsim_\nu \ln \l(\frac{1+c}{c}\r),
\end{equation*}
which diverges as $c\rightarrow 0$, while from Lemma \ref{growth bounds M large argument} and since $|M_+(x)|$ is an increasing function of $x\geq 0$, we have  
\begin{equation*}
r(c)\leq \frac{\nu}{|M_+(c)|^2}\lesssim \e^{-c}.
\end{equation*}
Therefore, $r(c):(0,+\infty)\rightarrow(0,+\infty)$ is a bijection and we conclude the existence of two sequences of strictly positive real numbers $\lbrace p_k\rbrace_{k\geq 1}$ and $\lbrace q_k\rbrace_{k\geq 1}$ such that $q_{k+1}<p_k<q_{k}$, for all $k\geq 1$ and
\begin{equation*}
r(p_k)=\l(k+\frac12\r)\pi, \qquad r(q_k) = k\pi,
\end{equation*}
with the further property that $p_k,q_k\to 0$ as $k\to\infty$.
\end{proof}

\begin{corollary}
There are infinitely many eigenvalues $c_k:=-q_k<0$ of $L_m$ accumulating at $c=0$.
\end{corollary}

\begin{proof}
Any eigenvalue $c<0$ of $L_m$ is such that there exists a non-trivial solution $\psi=\psi_m(y,c)$ to 
\begin{equation*}
\D_m\psi + \b^2\frac{\psi}{(y-c)^2}=0
\end{equation*}
satisfying the boundary conditions $\psi(0)=\psi(1)=0$. We can write such solution as
\begin{equation*}
\psi(y)=AM_+(y-c) + BM_-(y-c)
\end{equation*}
and since $c<0$, it is smooth. Imposing the boundary conditions, we have non-trivial coefficients $A,B\in\C$ if and only if $\W_m(c)=M_+(c)M_-(1+c)-M_-(c)M_+(1+c)$ vanishes. This is the case for the sequence $\lbrace q_k\rbrace_{k\in\N}$. These $q_k$ are the discrete eigenvalues of $L_m$ and they accumulate towards 0.
\end{proof}

We shall next obtain suitable estimates on the contour integral of the resolvent. We focus on the integral along $\p R_0$, the results we show are the same for the integral along $\p R_1$. In this direction, we write
\begin{equation*}
\int_{\p R_0}\e^{-imct}\mathcal{R}(c,L_m)\d c = \int_{\p R_*}\e^{-imct}\mathcal{R}(c,L_m)\d c + \int_{\p (R_0\setminus R_*)}\e^{-imct}\mathcal{R}(c,L_m)\d c,
\end{equation*}
where $R_*=\l\lbrace c=y_0 + is\in \C: y_0\in \l[y_*,0\r], \, s\in[-\ep_*, \ep_*] \r\rbrace$, for some $y_*<0$ and $\ep_*>0$ that will be determined later. More precisely, we have that
\begin{equation*}
\begin{aligned}
\l\Vert \int_{\p R_*}\e^{-imct}\mathcal{R}(c,L_m)\d c \r\Vert_{L^2_y}&\leq \l\Vert \int_{y_*}^0 \e^{-imy_0t}\l(\e^{m\ep_*t}\psi_{m,\ep_*}^-(y,y_0)-\e^{-m\ep_*t}\psi_{m,\ep_*}^+(y,y_0) \r)\d y_0\r\Vert_{L^2_y}  \\
&\quad+ \l\Vert \int_0^{\ep_*} \e^{-imy_*t}\l(\e^{mst}\psi_{m,s}^-(y,y_*)+\e^{-mst}\psi_{m,s}^+(y,y_*) \r)\d s\r\Vert_{L^2_y} \\
&\quad+ \l\Vert \int_0^{\ep_*} \l(\e^{mst}\psi_{m,s}^-(y,0)+\e^{-mst}\psi_{m,s}^+(y,0) \r)\d s\r\Vert_{L^2_y} \\
&\quad+\l\Vert \int_{y_*}^0 \e^{-imy_0t}\l(\e^{m\ep_*t}\rho_{m,\ep_*}^-(y,y_0)-\e^{-m\ep_*t}\rho_{m,\ep_*}^+(y,y_0) \r)\d y_0\r\Vert_{L^2_y}  \\
&\quad+ \l\Vert \int_0^{\ep_*} \e^{-imy_*t}\l(\e^{mst}\rho_{m,s}^-(y,y_*)+\e^{-mst}\rho_{m,s}^+(y,y_*) \r)\d s\r\Vert_{L^2_y} \\
&\quad+ \l\Vert \int_0^{\ep_*} \l(\e^{mst}\rho_{m,s}^-(y,0)+\e^{-mst}\rho_{m,s}^+(y,0) \r)\d s\r\Vert_{L^2_y}.
\end{aligned}
\end{equation*}
In what follows, we obtain suitable estimates for each integral. We begin by obtaining bounds on the Green's functions $\G_{m,\ep}^\pm(y,y_*,z)$ when $y_*=-p_k$ for some $k\geq0$ small and for $\ep=\ep_*$ small.
 
\begin{lemma}\label{Green's limit on peak point}
Let $\ep>0$ and $p_k>0$ given by Proposition \ref{prop discrete eigenvalue and peak} for some $k\geq 1$. Then,
\begin{equation*}
|\G_{m,\ep}^\pm(y,-p_k,z)|\lesssim_m |y+p_k -i\ep|^\frac12,
\end{equation*}
uniformly for all $y,z\in[0,1]$, for $\frac{\ep}{p_k}$ sufficiently small.
\end{lemma}

\begin{proof}
We proceed similarly as in the proof of Lemma \ref{Green's limit on boundary}. That is, for $y\leq z$,
\begin{equation*}
\G_{m,\ep}^-(y,-p_k,z)=\frac{\phi_{l,m,\ep}^-(y,-p_k)\phi_{u,m,\ep}^-(z,-p_k)}{{\W_{m,\ep}^-(-p_k)}} 
\end{equation*} 
Due to the explicit solutions of the Taylor-Goldstein equation, we can find that
\begin{equation*}
\begin{aligned}
{\phi_{l,m,\ep}^-(y,-p_k)}&=  M_+(p_k- i \ep)M_-(y+p_k- i\ep)-M_-(p_k- i\ep)M_+(y+p_k- i\ep)\\
&= M_-(p_k)M_-(y+p_k-i\ep) -M_-(p_k)M_+(y+p_k-i\ep) + R_1(p_k,\ep)\\
\end{aligned}
\end{equation*}
where $R_1(p_k,\ep)\lesssim \frac{\ep}{p_k}|M_+(p_k)|\big( |M_+(y+p_k-i\ep)| + |M_-(y+p_k-i\ep)|\big)$. In particular,
\begin{equation*}
|\phi_{l,m,\ep}^-(y,-p_k)| \lesssim \l(1+\frac{\ep}{p_k}\r)|M_+(p_k)||y+p_k-i\ep|^\frac12.
\end{equation*}
On the other hand,
\begin{equation*}
\begin{aligned}
\phi_{u,m,\ep}^-(z,-p_k)&=M_+(1+p_k)M_-(z+p_k)-M_-(1+p_k)M_+(z+p_k) + R_2(p_k,\ep) \\
&=2i\Im\l(M_+(1+p_k)M_-(z+p_k)\r) + R_2(\ep)
\end{aligned}
\end{equation*}
where $|R_2(\ep)|\lesssim |M_+(1+p_k)|{\frac{\ep}{p_k}}$. In particular,
\begin{equation*}
|\phi_{u,\ep}^-(z,-p_k)|\lesssim \l(1+\frac{\ep}{p_k}\r)|M_+(1+p_k)|.
\end{equation*}
Now, let us now estimate the Wronskian. We trivially have that
\begin{equation*}
\begin{aligned}
\W_{m,\ep}^-(-p_k)&= 4i\nu m \Big( M_+(p_k-i\ep)M_-(1+p_k-i\ep) - M_+(p_k-i\ep)M_-(1+p_k-i\ep)\Big) \\
&=4i\nu m \Big( M_+(p_k)M_-(1+p_k) - M_-(p_k)M_+(1+p_k)\Big) + R_3(p_k,\ep) \\
&=8\nu m \Im\Big(M_-(p_k)M_+(1+p_k)\Big) + R_3(p_k,\ep) \\
&=8\nu m (-1)^k|M_+(p_k)||M_-(1+p_k)|+ R_3(p_k,\ep),
\end{aligned}
\end{equation*}
where $|R_3(p_k,\ep)|\lesssim \nu m|M_+(p_k)||M_+(1+p_k)|\frac{\ep}{p_k}$. In particular,
\begin{equation*}
\begin{aligned}
|\W_{m,\ep}(0)| \geq 8\nu m|M_+(p_k)||M_+(1+p_k)| - |R_3(p_k,\ep)|\gtrsim \nu m|M_+(p_k)||M_+(1+p_k)|,
\end{aligned}
\end{equation*}
for $\frac{\ep}{p_k}$ small enough. The bound on $\G_{m,\ep}^-(y,-p_k,z)$ follows directly.
\end{proof}
Once we have the pointwise bounds on the Green's function, we are able to prove the following.
\begin{proposition}\label{imag peak limiting absorption principle}
Let $y_*=-p_k$ be given by Proposition \ref{prop discrete eigenvalue and peak} for some $k\geq 0$ and let $\ep_*>0$ such that $\frac{\ep_*}{|y_*|}$ is small enough. Then,
\begin{equation*}
\l\Vert \int_0^{\ep_*} \e^{-imp_kt}\l(\e^{mst}\psi_{m,s}^-(y,y_*)+\e^{-mst}\psi_{m,s}^+(y,y_*) \r)\d s\r\Vert_{L^2_y}\lesssim \ep_*
\end{equation*}
and
\begin{equation*}
\l\Vert \int_0^{\ep_*} \e^{-imp_kt}\l(\e^{mst}\rho_{m,s}^-(y,y_*)+\e^{-mst}\rho_{m,s}^+(y,y_*) \r)\d s\r\Vert_{L^2_y}\lesssim \ep_*^\frac12
\end{equation*}
\end{proposition}

\begin{proof}
Firstly, Minkowski inequality provides
\begin{equation*}
\begin{aligned}
\l\Vert \int_0^{\ep_*} \e^{-imp_kt}\l(\e^{mst}\psi_{m,s}^-(y,y_*)+\e^{-mst}\psi_{m,s}^+(y,y_*) \r)\d s\r\Vert_{L^2_y} &\lesssim \sum_{\kappa\in \lbrace +, -\rbrace } \int_0^{\ep_*} \Vert\psi_{m,s}^\kappa(y,-p_k)\Vert_{L^2_y}\d s \\
\end{aligned}
\end{equation*}
and we have that
\begin{equation*}
\begin{aligned}
\psi^\pm_{m,s}(y,y_0)= \frac{1}{\b^2}(y-y_0\pm i\ep)\o_m^0(y) -\rho_m^0(y) + \int_0^1 \G^\pm_{m,\ep}(y,y_0,z) F_{m,\ep}^\pm(z,y_0) \d z.
\end{aligned}
\end{equation*}
Using Cauchy-Schwarz and the uniform estimates from Lemma \ref{Green's limit on peak point}, we bound
\begin{equation*}
\begin{aligned}
\sum_{\sigma\in\lbrace +,-\rbrace}\Vert \psi_{m,s}^\sigma(y,-p_k)\Vert_{L^2_y} &\lesssim \Vert \o_m^0\Vert_{L^2_y} + \Vert \rho_m^0\Vert_{L^2_y} +\sum_{\sigma\in\lbrace +,-\rbrace}\Vert \G_{m,s}^\sigma(y,-p_k,z)\Vert_{L^2_yL^2_z} \lesssim_m 1
\end{aligned}
\end{equation*}
and thus integrating in $s$ from 0 to $\ep_*$ we get the first part of the Proposition. For the perturbed density, we recall that
\begin{equation*}
\rho_{m,s}^\pm(y,y_0) = \frac{1}{\b^2}\o_m^0(y)+\frac{1}{y-y_0\pm is}\int_0^1\G_{m,s}^\pm(y,y_0,z)F_{m,s}^\pm(z,y_0)\d z.
\end{equation*}
In particular, from Lemma \ref{Green's limit on peak point} we have that
\begin{equation*}
\begin{aligned}
|\rho_{m,s}^\pm(y,-p_k)| \lesssim |\o_m(y)| + |y+p_k-is|^{-\frac12}\Vert F_{m,s}^\pm(z,y_0)\Vert_{L^2_z}
\end{aligned}
\end{equation*}
For $\Vert F_{m,s}^\pm(z,y_0)\Vert_{L^2_z}\lesssim 1$ uniformly in $s\in(0,\ep_*)$, we integrate in $s$ from 0 to $\ep_*$ to get the desired result.
\end{proof}

We next obtain bounds on the Green's function when the spectral parameter has non-zero imaginary part. These bounds are shown to depend both on the modulus and on the argument of the complex spectral parameter.

\begin{lemma}\label{Green's limit on angle}
Let $y_0<0$ and $\ep>0$. Denote $c=-y_0+i\ep=r\e^{i\theta}$, with $r>0$ and $\theta \in \l(0,\frac{\pi}{2}\r)$. Then,
\begin{equation*}
|\G_{m,\ep}^\pm(y,y_0,z)|\lesssim \frac{|y-y_0\pm i\ep|^\frac12}{\sinh^2(\nu\theta)}
\end{equation*}
and there exists $K_c>0$ such that
\begin{equation*}
|\G_{m,\ep}^-(y,y_0,z)-\G_{m,\ep}^+(y,y_0,z)-K_c\phi_{u,m}(y)\phi_{u,m}(z)|\lesssim \frac{r^\frac12}{\sinh^2(\nu\theta)},
\end{equation*}
uniformly for all $y,z\in[0,1]$.
\end{lemma}

\begin{proof}
For $y_0<0$ and $\ep>0$, we consider $c=-y_0+i\ep=r\e^{i\theta}$, with $r>0$ and $\theta\in \l(0,\frac{\pi}{2}\r)$. We next study $\G_{m,\ep}^+(y,y_0,z)$. For $y\leq z$, we write
\begin{equation*}
\G_{m,\ep}^+(y,y_0,z)=\frac{\phi_{l,m,\ep}^+(y,y_0)\phi_{u,m,\ep}^+(z,y_0)}{\W_{m,\ep}^+(y_0)} = \frac{\phi_{l,m,\ep}^+(y,y_0)\phi_{u,m,\ep}^+(y,y_0)\W_{m,\ep}^-(y_0)}{|\W_{m,\ep}^+(y_0)|^2}
\end{equation*}
The main difference with respect to the other estimates we have been carrying out is that now, we control $|\W_{m,\ep}^+(y_0)|^2$ as follows:
\begin{equation*}
\begin{aligned}
\W_{m,\ep}^+(y_0)&=4i\nu m\big( M_+(c)M_-(1+c)-M_-(c)M_+(1+c)\big) \\
&=4i\nu m \big(M_+(c)M_-(1)-M_-(c)M_+(1)\big) + R_1(c),
\end{aligned}
\end{equation*}
with $|R_1(c)|\lesssim r|M_+(c)||M_+(1)|$. For $c=r\e^{i\theta}$, a detailed asymptotic analysis of $M_+(c)$ and $M_-(c)$ shows that
\begin{equation*}
M_+(c)=r^\frac12\e^{-\nu\theta}\e^{i\frac{\theta}{2}}r^{i\nu} + R_2(c), \quad M_-(c)=r^\frac12\e^{\nu\theta}\e^{i\frac{\theta}{2}}r^{-i\nu} + R_3(c),
\end{equation*}
where $|R_2(c)|,\,|R_3(c)|\lesssim r^{\frac52}$. Hence, 
\begin{equation*}
\W_{m,\ep}^+(y_0) = 4i\nu m r^\frac12 \e^{i\frac{\theta}{2}}\big( \e^{-\nu\theta}r^{i\nu}M_-(1) - \e^{\nu\theta}r^{-i\nu}M_+(1)\big) + R_4(c),
\end{equation*}
with $|R_4(c)|\leq C_4r^{\frac32}|M_+(1)|$. In particular, for $r\leq\frac{4\nu m\sinh(\nu\theta)}{C_4}$ small enough, we estimate
\begin{equation*}
|\W_{m,\ep}^+(y_0)|\geq 4\nu m r^\frac12|M_+(1)|\sinh(\nu\theta).
\end{equation*}
As expected, the bound degenerates as $\theta\rightarrow 0^+$. With this lower bound we are able to prove the first part of the proposition, using the asymptotic expansions of $M_+(y-y_0+i\ep)$ and $M_-(y-y_0+i\ep)$, see Lemma \ref{asymptotic expansion M}. Nevertheless, to obtain the second part of the proposition, we continue by estimating
\begin{equation*}
\begin{aligned}
\phi_{u,m,\ep}^+(z,y_0)&=M_+(1)M_-(z)-M_-(1)M_+(z) +  R_5(c)=2i\Im\l(M_+(1)M_-(z)\r) +  R_5(c),
\end{aligned}
\end{equation*}
where $|R_5(c)|\lesssim r^\frac12|M_+(1)|$. Similarly, we have 
\begin{equation*}
\begin{aligned}
{\phi_{l,m,\ep}^+(y,0)}{{\W_{m,\ep}^-(y_0)}}&= 4i\nu m\Big( M_+(c)M_-(y)M_+(\overline{c})M_-(1) -M_+(c)M_-(y)M_-(\overline{c})M_+(1)\Big. \\
&\qquad\qquad- \Big. M_-(c)M_+(y)M_+(\overline{c})M_-(1) + M_-(c)M_+(y)M_-(\overline{c})M_+(1)\Big) + R_6(c),
\end{aligned}
\end{equation*}
with $|R_6(c)|\lesssim r^\frac12|M_+(c)|^2|M_+(1)|$. In fact, we can recognize
\begin{equation*}
\begin{aligned}
{\phi_{l,m,\ep}^+(y,0)}{{\W_{m,\ep}^-}(y_0)}&= 4i\nu m\Big( 2\Re\big(M_+(c)M_-(y)M_+(\overline{c})M_-(1)\big)\Big. \\
&\qquad\qquad-\Big.|M_+(c)|^2M_-(y)M_+(1) -|M_-(c)|^2M_+(y)M_-(1) \Big) + R_6(c).
\end{aligned}
\end{equation*}
Hence, we obtain
\begin{equation*}
\begin{aligned}
\phi_{l,m,\ep}^+(y,y_0)\phi_{u,m,\ep}^+(y,y_0)\W_{m,\ep}^-(y_0) &= -8\nu m \Im\l(M_+(1)M_-(z)\r)\Big( 2\Re\big(M_+(c)M_-(y)M_+(\overline{c})M_-(1)\big)\Big. \\
&\qquad-\Big.|M_+(c)|^2M_-(y)M_+(1) -|M_-(c)|^2M_+(y)M_-(1) \Big) + R_7(c),
\end{aligned}
\end{equation*}
where now $|R_7(c)|\lesssim r^\frac12|M_+(c)|^2|M_+(1)|^2$. In particular, 
\begin{align*}
\Im&\l( \phi_{l,m,\ep}^+(y,y_0)\phi_{u,m,\ep}^+(y,y_0)\W_{m,\ep}^-(y_0)\r)\\
&= -2\nu m \phi_{u,m}(z)\phi_{u,m}(y)\l( |M_+(c)|^2-|M_-(c)|^2\r) + \Im\l(R_7(c)\r).
\end{align*}
and
\begin{equation*}
\l|\phi_{l,m,\ep}^+(y,y_0)\phi_{u,m,\ep}^+(y,y_0)\W_{m,\ep}^-(y_0)\r|\lesssim |M_+(c)|^2|M_+(1)|^2.
\end{equation*}
Together with the lower bound on the Wronskian, we conclude the proof.
\end{proof}
With the above bounds, we are able to estimate the contribution of the integral along the horizontal boundary.

\begin{proposition}\label{imag angle limiting absorption principle}
For $y_*<0$ small enough, let $r_*\e^{i\theta_*}=y_*+i\ep_*$. We have that
\begin{equation*}
\l\Vert \int_0^{y_*} \e^{-imy_0t}\big(\e^{m\ep_*t}\psi_{m,\ep_*}^-(y,y_0)-\e^{-m\ep_*t}\psi_{m,\ep_*}^+(y,y_0) \big)\d y_0\r\Vert_{L^2_y}\lesssim \frac{r_*^\frac32}{\sinh^2(\nu\theta_*)}
\end{equation*}
and
\begin{equation*}
\l\Vert \int_0^{y_*} \e^{-imy_0t}\big(\e^{m\ep_*t}\rho_{m,\ep_*}^-(y,y_0)-\e^{-m\ep_*t}\rho_{m,\ep_*}^+(y,y_0) \big)\d y_0\r\Vert_{L^2_y}\lesssim \frac{r_*^\frac12}{\sinh^2(\nu\theta_*)}
\end{equation*}
\end{proposition}

\begin{proof}
Firstly, note that
\begin{equation*}
\begin{aligned}
\e^{m\ep_*t}\psi_{m,\ep_*}^-(y,y_0)-\e^{-m\ep_*t}\psi_{m,\ep_*}^+(y,y_0) &= \e^{m\ep_*t}\l(\psi_{m,\ep_*}^-(y,y_0)-\psi_{m,\ep_*}^+(y,y_0) \r) \\
&\quad+\l(\e^{m\ep_*t}-\e^{-m\ep_*t} \r)\psi_{m,\ep_*}^+(y,y_0)
\end{aligned}
\end{equation*} 
while
\begin{equation*}
\begin{aligned}
\psi_{m,\ep_*}^-(y,y_0)-\psi_{m,\ep_*}^+(y,y_0)&= -\frac{2i\ep_*}{\b^2}\o_m^0 +\int_0^1 \l(\G_{m,\ep_*}^-(y,y_0,z)-\G_{m,\ep_*}^+(y,y_0,z)\r)F_{m}(z,0)\d z \\
&\quad +\frac{i\ep_*}{\b^2}\int_0^1 \l(\G_{m,\ep_*}^-(y,y_0,z)+\G_{m,\ep_*}^+(y,y_0,z)\r)\D_m\o_m^0\d z.
\end{aligned}
\end{equation*}
Now, for $r\e^{i\theta}=-y_0+i\ep_*$, we use Lemma \ref{Green's limit on angle} to bound
\begin{equation*}
\ep_* \l\Vert  \int_0^1 \l(\G_{m,\ep_*}^-(y,y_0,z)+\G_{m,\ep_*}^+(y,y_0,z)\r)\D_m\o_m^0\d z \r\Vert_{L^2_y}\lesssim \frac{\ep}{\sinh^2(\nu\theta)}\lesssim \frac{\ep}{\sinh^2(\nu\theta_*)}
\end{equation*}
and, together with the orthogonality condition of the initial data,
\begin{equation*}
\l\Vert \int_0^1 \l(\G_{m,\ep_*}^-(y,y_0,z)-\G_{m,\ep_*}^+(y,y_0,z)\r)F_{m}(z,0)\d z \r\Vert_{L^2_y}\lesssim \frac{r^\frac12}{\sinh^2(\nu\theta)}\lesssim \frac{r_*^\frac{1}{2}}{\sinh^2(\nu\theta_*)},
\end{equation*}
where $r_*\e^{i\theta_*}=y_*+i\ep_*$. With this bound uniform in $y_0\in[y_*,0]$, we obtain
\begin{equation*}
\int_0^{y_*}\Vert \psi_{m,\ep}^-(y,y_0)-\psi_{m,\ep}^+(y,y_0)\Vert_{L^2_y}\d y_0 \lesssim \ep |y_*| + \frac{r_*^\frac12}{\sinh^2(\nu\theta_*)}|y_*| + \frac{\ep}{\sinh^2(\nu\theta_*)}|y_*|\lesssim r_*^\frac32\frac{\cos(\theta_*)}{\sinh^2(\nu\theta_*)}.
\end{equation*}
On the other hand,
\begin{equation*}
\l\Vert \int_0^{|y_*|} \l(\e^{m\ep_*t}-\e^{-m\ep_*t}\r)\psi_{m,\ep_*}^+(y,y_0)\d y_0 \r\Vert_{L^2_y}\lesssim \frac{\ep_*|y_*|}{\sinh^2(\nu\theta)}\lesssim \frac{r_*^2}{\sinh^2(\nu\theta_*)},
\end{equation*}
For the second part of the proposition, we recall that
\begin{equation*}
\begin{aligned}
\rho_{m,\ep_*}^-(y,y_0) &- \rho_{m,\ep_*}^+(y,y_0) \\
&= \frac{1}{y-y_0-i\ep_*}\int_0^1\big( \G_{m,\ep_*}^-(y,y_0,z)-\G_{m,\ep_*}^+(y,y_0,z)\big)F_m(z,y_0)\d z \\
&\quad+ \frac{2i\ep_*}{(y-y_0)^2+\ep_*^2}\int_0^1\G_{m,\ep_*}^+(y,y_0,z)F_m(z,y_0)\d z \\
&\quad+ \frac{i\ep_*}{\b}\int_0^1\l( \frac{1}{y-y_0-i\ep_*}\G_{m,\ep_*}^-(y,y_0,z)+\frac{1}{y-y_0+i\ep_*}\G_{m,\ep_*}^+(y,y_0,z)\r)\D_m\o_m^0\d z.
\end{aligned}
\end{equation*}
Using the bounds of Lemma \ref{Green's limit on angle} and the orthogonal condition on the initial data, we bound
\begin{equation*}
\begin{aligned}
|\rho_{m,\ep_*}^-(y,y_0) - \rho_{m,\ep_*}^+(y,y_0)|&\lesssim \frac{1}{|y-y_0-i\ep_*|}\frac{|y-y_0+i\ep_*|^\frac12}{\sinh^2(\nu\theta)} + \frac{\ep_*}{|y-y_0+i\ep_*|^2}\frac{|y-y_0+i\ep_*|^\frac12}{\sinh^2(\nu\theta)} \\
&\quad + \frac{|y-y_0+i\ep_*|^\frac12}{\sinh^2(\nu\theta)} \\
&\lesssim \frac{1}{\sinh^2(\nu\theta_*)}\frac{1}{|y-y_0+i\ep_*|^\frac12}.
\end{aligned}
\end{equation*}
Hence,
\begin{equation*}
\int_0^{y_*} |\rho_{m,\ep_*}^-(y,y_0) - \rho_{m,\ep_*}^+(y,y_0)| \d y_0 \lesssim \frac{1}{\sinh^2(\nu\theta_*)}\int_0^{|y_*|}\frac{1}{|y_0|^\frac12}\d y_0 \lesssim \frac{|y_*|^\frac12}{\sinh^2(\nu\theta_*)}
\end{equation*}
and similarly
\begin{equation*}
\l\Vert \int_0^{y_*} \l(\e^{m\ep_*t}-\e^{-m\ep_*t}\r)|\rho_{m,\ep_*}^+(y,y_0)| \d y_0 \r\Vert_{L^2_y}\lesssim \frac{\ep_*|y_*|^\frac12}{\sinh^2(\nu\theta_*)}\lesssim \frac{r_*^\frac32}{\sinh^2(\nu\theta_*)}.
\end{equation*}
The proof is concluded.
\end{proof}

We combine the estimates from Proposition \ref{imag peak limiting absorption principle} and Proposition \ref{imag angle limiting absorption principle} to obtain the following result.

\begin{proposition}\label{imag small box limiting absorption principle}
For all $\delta>0$, there exists $\theta_*\in\l(0,\frac{\pi}{2}\r)$ such that, for $r_*=\sinh^8(\nu\theta_*)$, $y_*=-r_*\cos(\theta_*)$ and $\ep_*=r_*\sin(\theta_*)$, we denote $R_*:=\l\lbrace c=y_0 + is\in \C: y_0\in \l[y_*,0\r], \, s\in[-\ep_*, \ep_*] \r\rbrace$. Then, there holds
\begin{equation*}
\l\Vert \int_{\p R_*}\e^{-imct}\mathcal{R}(c,L_m)\d c \r\Vert_{L^2_y}\leq \delta.
\end{equation*}
\end{proposition}

\begin{proof}
We choose $\theta_*>0$ such that $y_*=-r_*\cos(\theta_*)=-\sinh^8(\nu\theta_*)\cos(\theta_*)=-p_k$, for some $k>0$, where $p_k$ is given by Proposition \ref{prop discrete eigenvalue and peak}. This is possible because for $\theta_*$ small enough, $g(\theta_*):=\sinh^8(\nu\theta_*)\cos(\theta_*)$ is a continuous strictly monotone increasing function of $\theta_*$ such that $g(0)=0$. Moreover, since $p_k\rightarrow 0^+$ for $k\rightarrow \infty$, we may assume $\theta_*$ is sufficiently small. Hence, $\frac{\ep_*}{y_*}=\tan(\theta_*)$ is sufficiently small and we use Proposition \ref{imag peak limiting absorption principle} to bound
\begin{equation*}
\l\Vert \int_0^{\ep_*} \e^{-imy_*t}\l(\e^{mst}\psi_{m,s}^-(y,y_*)+\e^{-mst}\psi_{m,s}^+(y,y_*) \r)\d s\r\Vert_{L^2_y} \lesssim \ep_* \lesssim \sinh^8(\nu\theta_*),
\end{equation*}
and
\begin{equation*}
\l\Vert \int_0^{\ep_*} \e^{-imy_*t}\l(\e^{mst}\rho_{m,s}^-(y,y_*)+\e^{-mst}\rho_{m,s}^+(y,y_*) \r)\d s\r\Vert_{L^2_y} \lesssim \ep_*^\frac12\lesssim \sinh^4(\nu\theta_*).
\end{equation*}
Now, we use Proposition \ref{imag angle limiting absorption principle} to bound
\begin{equation*}
\begin{aligned}
\l\Vert \int_0^{y_*} \e^{-imy_0t}\l(\e^{m\ep_*t}\psi_{m,\ep_*}^-(y,y_0)-\e^{-m\ep_*t}\psi_{m,\ep_*}^+(y,y_0) \r)\d y_0 \r\Vert_{L^2_y} \lesssim \frac{r_*^\frac32}{\sinh^2(\nu\theta_*)} \lesssim \sinh^{10}(\nu\theta_*)
\end{aligned}
\end{equation*}
and
\begin{equation*}
\l\Vert \int_0^{y_*} \e^{-imy_0t}\big(\e^{m\ep_*t}\rho_{m,\ep_*}^-(y,y_0)-\e^{-m\ep_*t}\rho_{m,\ep_*}^+(y,y_0)\big)\d y_0\r\Vert_{L^2_y}\lesssim \frac{r_*^\frac12}{\sinh^2(\nu\theta_*)}\lesssim \sinh^2(\nu\theta_*)
\end{equation*}
Finally, we use \eqref{def psi}, \eqref{def rho} and the bounds from Proposition \ref{L2 bounds inhom TG solution} to estimate
\begin{equation*}
\l\Vert \int_0^{\ep_*} \l(\e^{mst}\psi_{m,s}^-(y,0)+\e^{-mst}\psi_{m,s}^+(y,0) \r)\d s\r\Vert_{L^2_y} \lesssim \ep_*\lesssim \sinh^{8}(\nu\theta_*)
\end{equation*}
and
\begin{equation*}
\l\Vert \int_0^{\ep_*} \l(\e^{mst}\rho_{m,s}^-(y,0)+\e^{-mst}\rho_{m,s}^+(y,0) \r)\d s\r\Vert_{L^2_y} \lesssim \ep_*^\frac12\lesssim \sinh^4(\nu\theta),
\end{equation*}
The proposition follows choosing $\theta_*$ small enough (that is, $\theta_*=g^{-1}(p_k)$ for $k>0$ sufficiently large). 
\end{proof}

We are finally in position to prove Proposition \ref{prop: contour limiting absorption principle} for the case $\b^2>1/4$.

\begin{proof}[Proof of Proposition \ref{prop: contour limiting absorption principle}]
We shall see that $\l\Vert \int_{\p R_0}\e^{-imct}\mathcal{R}(c,L_m)\d c \r\Vert_{L^2_y}\leq \delta$, for all $\delta>0$. Indeed, given $\delta>0$, from Proposition \ref{imag small box limiting absorption principle} there exists $\theta_*$ such that $y_*=-\sinh^8(\nu\theta_*)\cos(\theta_*)=-p_k$, for some $k>0$ large enough and such that
\begin{equation*}
\l\Vert \int_{\p R_*}\e^{-imct}\mathcal{R}(c,L_m)\d c \r\Vert_{L^2_y}\leq \delta.
\end{equation*}
Now, there are finitely many isolated eigenvalues in $R_0\setminus R_*$, they are real and between $-\frac{\b}{m}$ and $y_*<0$. Moreover, $\e^{-imct}\mathcal{R}(c,L_m)$ is an holomorphic function, for all $c\in R\setminus R_*$  such that $c\neq q_j$, for any $0\leq j\leq k$. Thus, 
\begin{equation*}
\frac{1}{2\pi i}\int_{\p R\setminus \p R_*}\e^{-imct}\mathcal{R}(c,L_m)\d c = \sum_{j=0}^k  \e^{-imq_jt} \mathbb{P}_{q_j} \begin{pmatrix}
\o_m^0 \\ \rho_m^0
\end{pmatrix}=0,
\end{equation*}
where $\mathbb{P}_{q_j} \begin{pmatrix}
\o_m^0 \\ \rho_m^0
\end{pmatrix}$ denotes the $L^2$-projection of $\begin{pmatrix}
\o_m^0 \\ \rho_m^0\end{pmatrix}$ to the generalized eigenspace $E_{q_j}$ associated to the eigenvalue $q_j$. With this, the proof is finished.
\end{proof}

The next proposition shows that the generalized eigenspace associated to any discrete eigenvalues is, in fact, simple.

\begin{proposition}
Let $c\in\R$ be a discrete eigenvalue of $L_m$. Then
$\ker\l(L_m -c\r)^2 = \ker \l(L_m-c\r)$.
In particular, $c$ is a semi-simple eigenvalue.
\end{proposition}

\begin{proof}
Note that the pair $(\o,\rho)\in\ker\l(L_m-c\r)$ if and only if
\begin{equation*}
\begin{aligned}
(y-c)\o + \b^2\rho=0, \qquad -\psi + (y-c)\rho =0,
\end{aligned}
\end{equation*}
where, as usual, we denote $\psi=\Delta_m^{-1}\o$. Hence, $\rho=\frac{\psi}{y-c}$ and the equation
\begin{equation*}
\Delta_m\psi + \b^2\frac{\psi}{(y-c)^2}=0
\end{equation*}
characterizes the eigenfunctions of $L_m$ of eigenvalue $c\in\R$. Now, the pair $(\o,\rho)\in\ker\l(L_m-c\r)^2$ if and only if
\begin{equation*}
\begin{aligned}
(y-c)^2\o - \b^2\psi + 2(y-c)\b^2\rho &= 0, \\
-\Delta_m^{-1}((y-c)\o) - (y-c)\psi + (y-c)^2\rho - \b^2\Delta_m^{-1}\rho &=0.
\end{aligned}
\end{equation*}
Obtaining $\rho$ in terms of $\o$ from the first equation and plugging it into the second one, we see that $\o$ solves
\begin{equation}
0=-\Delta_m^{-1}((y-c)\o) -(y-c)\psi - \frac{(y-c)^3}{2\b^2}\o + \frac{y-c}{2}\psi +\frac12\Delta_m^{-1}\l( \frac{1}{y-c}\l[ (y-c)^2\o -\b^2\psi\r]\r).
\end{equation}
Multiplying by $\overline{\o}$ and integrating by parts, we see that
\begin{equation*}
\begin{aligned}
0&=-\int_0^1\overline{\psi}(y-c)\o -\int_0^1 \overline{\o}(y-c)\psi -\int_0^1 \frac{(y-c)^3}{2\b^2}|\o|^2  + \int_0^1\overline{\o}\frac{y-c}{2}\psi + \int_0^1 \overline{\psi}\frac{y-c}{2}\o - \int_0^1 \frac{\b^2}{2}|\psi|^2 \\
&=-\frac{1}{2\b^2}\int_0^1 (y-c)^3 \l( |\o|^2 + \b^2\frac{\overline{\psi}\o + \overline{\o}\psi}{(y-c)^2} + \b^4 \frac{|\psi|^2}{(y-c)^4}\r) \\
&=-\frac{1}{2\b^2}\int_0^1 (y-c)^3\l| \D_m\psi + \b^2\frac{\psi}{(y-c)^2}\r|^2.
\end{aligned}
\end{equation*}
Hence, since either $c<0$ or $c>1$, we conclude that the pair $(\o,\rho)$ satisfies
\begin{equation*}
\Delta_m\psi + \b^2\frac{\psi}{(y-c)^2}=0,
\end{equation*}
with $\D_m\psi=\o$. That is, $(\o,\rho)\in \ker\l(L_m -c\r)$.
\end{proof}

We finish the subsection proving Theorem \ref{thm: boundary limiting absorption principle} for $\b^2> 1/4$. We need the following lemma, which shows that for $y_0\in\lbrace 0, 1\rbrace$, the difference $\G_{m,\ep}^-(y,y_0)-\G_{m,\ep}^+(y,y_0)$ approaches a varying multiple of a generalized eigenfunction of the linearized operator as $\ep\rightarrow 0$ associated to the ``embedded eigenvalue" $c=0$. We state the result for $y_0=0$, since the case $y_0=1$ is analogous.

\begin{lemma}\label{Green's limit on boundary}
Let $y_0=0$ and $0<\ep\ll 1$ sufficiently small. Then, there exists $C_\ep\in\R$ such that 
\begin{equation*}
\left| \G_{m,\ep}^-(y,0,z)-\G_{m,\ep}^+(y,0,z)-C_\ep\phi_{u,m}(y)\phi_{u,m}(z)\right|\lesssim \ep^\frac12,
\end{equation*}
where $\phi_{u,m}$ is given in \eqref{eq: phium}.
\end{lemma}

\begin{proof}
The result is trivially true for $y=0$ and $y=1$ because both $\phi_u$ and $\G^{\pm}$ vanish there. Therefore, in the sequel we consider $0<y<1$. Due to the complex conjugation property of the Green's function, 
\begin{equation*}
\G_{m,\ep}^-(y,0,z)- \G_{m,\ep}^+(y,0,z)=2i\Im \l( \G_{m,\ep}^-(y,0,z)\r)
\end{equation*}
Assuming initially that $y\leq z$, we have
\begin{align*}
2i\Im \l( \G_{m,\ep}^-(y,0,z)\r)&=2i\Im\l(\frac{\phi_{l,m,\ep}^-(y,0)\phi_{u,m,\ep}^-(z,0)}{{\W_{m,\ep}^-}(0)}\r)\\
&=\frac{2i}{|{\W_{m,\ep}^-}(0)|^2}\Im\l(\phi_{l,m,\ep}^-(y,0)\phi_{u,m,\ep}^-(z,0)\W_{m,\ep}^+(0)\r).
\end{align*} 
Due to the explicit solutions of the Taylor-Goldstein equation, we write
\begin{equation*}
\begin{aligned}
{\phi_{l,m,\ep}^-(y,0)}{\W_{m,\ep}^+(0)}&= 4i\nu m\Big( M_+(-i\ep)M_-(y)M_+(i\ep)M_-(1)  + M_-(-i\ep)M_+(y)M_-(i\ep)M_+(1) \Big. \\
&\quad- \Big. M_+(-i\ep)M_-(y)M_-(i\ep)M_+(1)  - M_-(-i\ep)M_+(y)M_+(i\ep)M_-(1) + R_1(\ep)\Big),
\end{aligned}
\end{equation*}
where $|R_1(\ep)|\lesssim_{m,\nu} \ep^\frac12 |M_+(i\ep)|^2$. Observe also that since $\overline{M_\pm(\zeta)}=M_\mp(\overline{\zeta})$ for all $\zeta\in\C$, we can write
\begin{equation*}
\begin{aligned}
{\phi_{l,m,\ep}^-(y,0)}{\W_{m,\ep}^+(0)}&= 4i\nu m \Big( 2\Re\big(M_+(-i\ep)M_-(y)M_+(i\ep)M_-(1)\big) \Big. \\
&\qquad\qquad -\Big. |M_+(-i\ep)|^2M_-(y)M_+(1) - |M_+(i\ep)|^2M_+(y)M_-(1)+R_1(\ep)\Big)
\end{aligned}
\end{equation*}
Since $M_+(-i\ep)=-i\e^{\nu\pi}M_+(i\ep)$, we obtain
\begin{equation*}
\begin{aligned}
{\phi_{l,m,\ep}^-(y,0)}{\W_{m,\ep}^+(0)}&= 4i\nu m \Big( 2\Re\big(M_+(-i\ep)M_-(y)M_+(i\ep)M_-(1)\big) \Big.\\
&\qquad\qquad-\Big. \e^{2\nu\pi}|M_+(i\ep)|^2M_-(y)M_+(1) - |M_+(i\ep)|^2M_+(y)M_-(1) + R_1(\ep) \Big).
\end{aligned}
\end{equation*}
Now,  $\overline{M_-(y)M_+(1)}=M_+(y)M_-(1)$ and we further observe that
\begin{equation*}
\Im \big( M_+(y)M_-(1)\big) = - \Im \big( M_-(y)M_+(1)\big) = -\frac{1}{2i}\big( M_-(y)M_+(1) - M_+(y)M_-(1)\big) =-\frac{1}{2i}\phi_{u,m}(y).
\end{equation*}
On the other hand,
\begin{equation*}
\begin{aligned}
\phi_{u,m,\ep}^-(z,0)&=M_+(1)M_-(z)-M_-(1)M_+(z) +  R_2(\ep)  =2i\Im\l(M_+(1)M_-(z)\r) +  R_2(\ep) 
\end{aligned}
\end{equation*}
with $|R_2(\ep)|\lesssim \ep^\frac12$. Thus,
\begin{equation*}
\begin{aligned}
\phi_{l,m,\ep}^-(y,0)&\phi_{u,m,\ep}^-(z,0)\W_{m,\ep}^+(0) \\
&=-8\nu m\Im\l(M_+(1)M_-(z)\r)\Big( 2\Re\big(M_+(-i\ep)M_-(y)M_+(i\ep)M_-(1)\big)   \Big.\\
&\qquad\qquad\qquad-\Big. \e^{2\nu\pi}|M_+(i\ep)|^2M_-(y)M_+(1) - |M_+(i\ep)|^2M_+(y)M_-(1) \Big) + R_3(\ep) \\
&=-8\nu m\Im\l(M_+(1)M_-(z)\r)\Big( 2\Re\big(M_+(-i\ep)M_-(y)M_+(i\ep)M_-(1)\big)\Big. \\
&\qquad\qquad\qquad -\Big. |M_+(i\ep)|^2\l( \e^{2\nu\pi}+1\r) \Re\l(M_-(y)M_+(1)\r)\Big) \\
&\quad-8i\nu m\Im\l(M_+(1)M_-(z)\r)|M_+(i\ep)|^2\l(\e^{2\nu\pi}-1\r)\Im\l(M_-(y)M_+(1)\r) + R_3(\ep),
\end{aligned}
\end{equation*}
where $|R_3(\ep)|\lesssim \ep^\frac12|M_+(i\ep)|^2$, uniformly in $y,z\in[0,1]$.
In particular,
\begin{equation*}
\Im\l(\phi_{l,m,\ep}^-(y,0)\phi_{u,m,\ep}^-(z,0)\W_{m,\ep}^+(0)\r)=2\nu m|M_+(i\ep)|^2\l(\e^{2\nu\pi}-1\r)\phi_{u,m}(y)\phi_{u,m}(z) + \Im\l(R_3(\ep)\r).
\end{equation*}
Moreover, due to the symmetry of the Green's function with respect to $y$ and $z$, we also have that
\begin{equation*}
\Im\l(\phi_{l,m,\ep}^-(z,0)\phi_{u,m,\ep}^-(y,0)\W_{m,\ep}^+(0)\r)=2\nu m|M_+(i\ep)|^2\l(\e^{2\nu\pi}-1\r)\phi_{u,m}(y)\phi_{u,m}(z) + \Im\l(\widetilde{R}_3(\ep)\r).
\end{equation*}
Let us now estimate the modulus squared of the Wronskian, that is, $|\W_{m,\ep}(0)|^2$. We trivially have that
\begin{equation*}
\begin{aligned}
|\W_{m,\ep}(0)|^2 &= 16\nu^2 m^2 \Big[ |M_+(-i\ep)|^2|M_-(1-i\ep)|^2 + |M_+(i\ep)|^2|M_-(1+i\ep)|^2\Big. \\
&\qquad\Big. -2\Re\Big(M_+(-i\ep)M_-(1-i\ep)M_+(i\ep)M_-(1+i\ep)\Big)\Big]  \\
&=16\nu^2 m^2 \Big[ \e^{2\nu\pi}|M_+(i\ep)|^2|M_-(1-i\ep)|^2+|M_+(i\ep)|^2|M_-(1+i\ep)|^2 \Big. \\
&\qquad \Big. -2\e^{\nu\pi}|M_+(i\ep)|^2|M_-(1-i\ep)||M_-(1+i\ep)|\cos(\theta_\ep)\Big]
\end{aligned}
\end{equation*}
where $\theta_\ep= \text{Arg}\big(M_+(-i\ep)M_-(1-i\ep)M_+(i\ep)M_-(1+i\ep)\big)$. As before, since $M_-(\zeta)$ is smooth at $\zeta=1$, we can further write
\begin{equation*}
\begin{aligned}
|\W_{m,\ep}(0)|^2 &=16\nu^2m^2|M_+(i\ep)|^2|M_-(1)|^2\l(\e^{2\nu\pi}-2\e^{\nu\pi}\cos(\theta_\ep)+1\r)+R_4(\ep)
\end{aligned}
\end{equation*}
where $|R_4(\ep)|\lesssim |\ep|^\frac12 |M_+(i\ep)|^2$. With this, we are able to write
\begin{equation*}
\begin{aligned}
\G_{m,\ep}^-(y,0,z)-\G_{m,\ep}^+(y,0,z)&=\frac{2i}{{|\W_{m,\ep}^-(0)|^2}}\Im\l(\phi_{l,m,\ep}^-(y,0)\phi_{u,m,\ep}^-(z,0)\W_{m,\ep}^+(0)\r) \\
&= -\frac{i}{4\nu m}\frac{|M_+(i\ep)|^2|\l( \e^{2\nu\pi}-1\r)\phi_u(y)\phi_u(z)}{|M_+(i\ep)|^2|M_-(1)|^2\l(\e^{2\nu\pi}-2\e^{\nu\pi}\cos(\theta_\ep)+1\r)+R_4(\ep)}\\
&\quad+\frac{i}{8\nu^2m^2}\frac{\Im \l( R_3(\ep)\r)}{|M_+(i\ep)|^2|M_-(1)|^2\l(\e^{2\nu\pi}-2\e^{\nu\pi}\cos(\theta_\ep)+1\r)+R_4(\ep)}
\end{aligned}
\end{equation*}
The lemma follows for 
$$C_\ep:=-\frac{i}{4\nu m}\frac{\e^{2\nu\pi}-1}{|M_-(1)|^2\l(\e^{2\nu\pi}-2\e^{\nu\pi}\cos(\theta_\ep)+1\r)}.$$
and recalling that $\l| \Im \l( R_3(\ep)\r) \r| \lesssim \ep^\frac12 |M_+(i\ep)|^2$,
\end{proof}

\begin{remark}
Note that $C_\ep$ is bounded but $\lim_{\ep\rightarrow0}C_\ep$ does not exist. Indeed, as can be seen from the asymptotic expansions of Lemma \ref{asymptotic expansion M}, we have that $\theta_\ep=2\nu\log(\ep) + 2\text{Arg}(M_-(1)) + O(\ep)$, as $\ep\rightarrow 0$. Thus, $\theta_\ep$ diverges to $-\infty$ and $\cos(\theta_\ep)$ does not converge. Hence, assumption \eqref{eq:specassume} becomes necessary in order to have a well defined pointwise limiting absorption principle for $y_0=0$.
\end{remark}

We are now in position to prove Theorem \ref{thm: boundary limiting absorption principle} for $\b^2>1/4$.
\begin{proof}[Proof of Theorem \ref{thm: boundary limiting absorption principle}]
For $y_0=0$ we have that
\begin{equation*}
\begin{aligned}
\psi_{m,\ep}^-(y,y_0)-\psi_{m,\ep}^+(y,y_0)&= -\frac{2i\ep}{\b^2}\o_m^0 +\int_0^1 \l(\G_{m,\ep}^-(y,0,z)-\G_{m,\ep}^+(y,0,z)\r)F_{m}(z,0)\d z \\
&\quad +\frac{i\ep}{\b^2}\int_0^1 \l(\G_{m,\ep}^-(y,0,z)+\G_{m,\ep}^+(y,0,z)\r)\D_m\o_m^0\d z.
\end{aligned}
\end{equation*}
Since $\o_m^0\in H_y^1$, the first term vanishes easily, while the third term also tends to zero when $\ep\rightarrow 0$, after a direct application of the Cauchy-Schwarz inequality and the facts that $\Vert \G_{m,\ep}^\pm\Vert_{L^2_z}$ is uniformly bounded in $\ep$ due to Theorem \ref{L2 bounds of G} and $\o_m^0\in H_y^2$. 

As for the second term, we invoke Lemma \ref{Green's limit on boundary} to show that 
\begin{equation*}
\begin{aligned}
\l|\int_0^1 (\G^-_{m,\ep} -\G_{m,\ep}^+)F_{m}(z) \d z\r| &\leq \l|C_\ep\phi_u(y,0)\int_0^1 \phi_u(z,0)F_m(z,0) \d z \r| \\
&+ \l|\int_0^1 (\G_{m,\ep}^- -\G_{m,\ep}^+ -C_\ep\phi_{u,m}(y)\phi_{u,m}(z))F_m(z) \d z \r| \\
&\lesssim \ep^\frac12\int_0^1 |F_m(z)|\d z \\
&\lesssim \ep^\frac12 \l( \Vert \rho_m^0 \Vert_{H^{2}_y} + \Vert \o_m^0 \Vert_{H^{2}_y} \r),
\end{aligned}
\end{equation*}
which vanishes as $\ep\rightarrow 0$. The proof for $y_0=1$ follows similarly from Lemma \ref{Green's limit on boundary}. We thus omit the details.
\end{proof}

\subsection{Integral reduction for $\b^2<1/4$: no discrete eigenvalues}
Thanks to the Hardy inequality \cite{HLP88}
\begin{equation}\label{Hardy Ineq}
\int_0^1 \l( \frac{1}{x}\int_0^x f(t)\d t\r) ^2  \d x\leq 4\int_0^1 |f(t)|^2 \d t, \qquad f\in L^2(0,1),
\end{equation}
we are able to prove $H^1$ bounds for the generalized stream functions $\psi_{m,\ep}^\pm(y,y_0)$ that are uniform in $\ep>0$. 
\begin{proposition}\label{H1 bounds boundary streamfunction}
Let   $0\leq \ep\leq 1$. Then,
\begin{equation}\label{eq:H1estbsmall}
\Vert \p_y \psi_{m,\ep}^\pm(y,y_0)\Vert_{L^2}^2 + m^2 \Vert \psi_{m,\ep}^\pm(y,y_0)\Vert_{L^2}^2 \lesssim \Vert \o_m^0 \Vert_{H^2}^2 + \Vert \rho_m^0 \Vert^2_{H^2}.
\end{equation}
Moreover,
\begin{equation*}
\begin{aligned}
\int_0^1\frac{\ep (y-y_0)}{((y-y_0)^2+\ep^2)^2}|\varphi_{m,\ep}^\pm|^2 \d y \lesssim \Vert \o_m^0 \Vert_{H^2}^2 + \Vert \rho_m^0 \Vert_{H^2}^2,
\end{aligned}
\end{equation*}
and
\begin{equation*}
\Vert \rho_{m,\ep}^\pm(y,y_0)\Vert_{L^2}^2 \lesssim \Vert \o_m^0 \Vert_{H^2}^2 + \Vert \rho_m^0 \Vert_{H^2}^2.
\end{equation*}
If we further assume that $|-y_0\pm i\ep|\geq c_0$, for some $c_0>0$, then
\begin{equation*}
\Vert \p_y \psi_{m,\ep}^\pm(y,y_0)\Vert_{L^2}^2 + m^2 \Vert \psi_{m,\ep}^\pm(y,y_0)\Vert_{L^2}^2 \lesssim \frac{1}{c_0^2}\Vert \o_m^0 \Vert_{L^2}^2 + \frac{1}{c_0^4}\Vert \rho_m^0 \Vert_{L^2}
\end{equation*}
and
\begin{equation*}
\Vert \rho_{m,\ep}^\pm(y,y_0)\Vert_{L^2}^2 \lesssim \frac{1}{c_0^2}\Vert \o_m^0 \Vert_{L^2}^2 + \frac{1}{c_0^4}\Vert \rho_m^0 \Vert_{L^2}^2.
\end{equation*}
In particular, $c=-y_0\pm i\ep$ belongs to the resolvent set of the operator $L_m$.
\end{proposition}
\begin{proof}
Multiplying \eqref{eq:varphieq} by $\overline{\varphi_{m,\ep}^\pm(y,y_0)}$ and integrating by parts, we obtain
\begin{equation*}
\int_0^1 \l(|\p_y\varphi_{m,\ep}^\pm(y,y_0)|^2 + m^2 |\p_y\varphi_{m,\ep}^\pm(y,y_0)|^2 -\b^2\frac{|\varphi_{m,\ep}^\pm(y,y_0)|^2}{(y-y_0\pm i\ep)^2}\r)\d y = -\int_0^1 F^\pm_{m,\ep}(y,y_0)\varphi_{m,\ep}^\pm(y,y_0)\d y
\end{equation*}
Assume now that $y_0\leq 0$ (the case $y_0\geq 1$ would be done similarly) and observe that, thanks to the Hardy inequality \eqref{Hardy Ineq} and $\varphi_{m,\ep}^\pm(0,y_0)=0$, 
\begin{equation*}
\begin{aligned}
\l| \b^2\int_0^1 \frac{|\varphi_{m,\ep}^\pm(y,y_0)|^2}{(y-y_0\pm i\ep)^2}\d y\r|\leq  \b^2 \int_0^1 \frac{|\varphi_{m,\ep}^\pm(y,y_0)|^2}{y^2}\d y &\leq \b^2 \int_0^1 \l(\frac{1}{y}\int_0^y |\p_y\varphi_{m,\ep}^\pm(y',y_0)|\d y'\r)^2\d y \\
&\leq 4\b^2\int_0^1 |\p_y\varphi_{m,\ep}^\pm(y,y_0)|^2 \d y.
\end{aligned}
\end{equation*}
Therefore, we conclude that
\begin{equation*}
(1-4\b^2)\int_0^1 |\p_y\varphi_{m,\ep}^\pm(y,y_0)|^2\d y + m^2 \int_0^1 |\varphi_{m,\ep}^\pm(y,y_0)|^2\d y \lesssim \frac{1}{m^2}\int_0^1 |F_{m,\ep}^\pm(y,y_0)|^2\d y.
\end{equation*}
Thus \eqref{eq:H1estbsmall} follows from \eqref{def psi}, the observation that $4\b^2<1$ and $\Vert F_{m,\ep}^\pm(y,y_0)\Vert_{L^2_y}^2\lesssim \Vert \o_m^0 \Vert_{H^2}^2 + \Vert \rho_m^0 \Vert_{H^2}$.
For the second statement, we take the real and imaginary part of \eqref{eq:varphieq}, for which we get
\begin{equation*}
\Delta_m \Re(\varphi_{m,\ep}^\pm)+\frac{1}{4((y-y_0)^2+\ep^2)^2}\l( ((y-y_0)^2-\ep^2)\Re(\varphi_{m,\ep}^\pm)\pm 2\ep (y-y_0)\Im(\varphi_{m,\ep}^\pm)\r)=\Re(F_{m,\ep}^\pm)
\end{equation*}
and
\begin{equation*}
\Delta_m\Im(\varphi_{m,\ep}^\pm)+\frac{1}{4((y-y_0)^2+\ep^2)^2}\l( ((y-y_0)^2-\ep^2)\Im(\varphi_{m,\ep}^\pm)\mp 2\ep (y-y_0)\Re(\varphi_{m,\ep}^\pm)\r)=\Im(F_{m,\ep}^\pm).
\end{equation*}
Cross multiplying the equations by $\Im(\varphi_{m,\ep}^\pm)$ and $\Re(\varphi_{m,\ep}^\pm)$, respectively, subtracting  them and integrating, we obtain
\begin{equation*}
\pm\int_0^1\frac{\ep (y-y_0)}{2((y-y_0)^2+\ep^2)^2}|\varphi_{m,\ep}^\pm|^2 \d y = \int_0^1 \Im(\varphi_{m,\ep}^\pm)\Re(F_{m,\ep}^\pm)-\Re(\varphi_{m,\ep}^\pm)\Im(F_{m,\ep}^\pm)\d y
\end{equation*}
so that
\begin{equation*}
\int_0^1\frac{\ep (y-y_0)}{((y-y_0)^2+\ep^2)^2}|\varphi_{m,\ep}^\pm|^2 \d y \lesssim \int_0^1 |\varphi_{m,\ep}^\pm|^2 + |F_{m,\ep}^\pm|^2 \d y \lesssim \Vert \o_m^0 \Vert_{H^2}^2 + \Vert \rho_m^0 \Vert^2_{H^2}.
\end{equation*}
The third statement of the proposition follows from the density formula \eqref{def rho}, the Hardy-type inequality and the uniform bounds from the first statement of the proposition. The proof is finished.
\end{proof}

From the arguments of the proof, one can directly obtain the following result.

\begin{corollary}
Let $y_0\leq 0$ or $y_0\geq 1$. Then, $y_0+ic$ is not an eigenvalue of $L_m$, for any $c\in\R$.
\end{corollary}

With the $\ep$-uniform $H^1_0$ for $\psi_{m,\ep}^\pm(y,y_0)$ at hand we are now able to prove Theorem \ref{thm: boundary limiting absorption principle}.

\begin{proposition}\label{real horizontal limiting absorption principle} 
We have that
\begin{equation*}
\lim_{\ep\rightarrow 0}\l\Vert \int_{-\frac{\b}{m}}^0 \e^{-imy_0t}\big(\e^{m\ep_*t}\psi_{m,\ep_*}^-(y,y_0)-\e^{-m\ep_*t}\psi_{m,\ep_*}^+(y,y_0) \big)\d y_0\r\Vert_{L^2_y}=0
\end{equation*}
and 
\begin{equation*}
\lim_{\ep\rightarrow 0}\l\Vert \int_{-\frac{\b}{m}}^0 \e^{-imy_0t}\big(\e^{m\ep_*t}\rho_{m,\ep_*}^-(y,y_0)-\e^{-m\ep_*t}\rho_{m,\ep_*}^+(y,y_0) \big)\d y_0\r\Vert_{L^2_y}=0.
\end{equation*}
\end{proposition}

\begin{proof}
Let us denote $\psi_{m,\ep}(y,y_0) =\psi_{m,\ep}^+(y,y_0) - \psi_{m,\ep}^-(y,y_0)$. Using \eqref{def psi}, we have
\begin{equation*}
\psi_{m,\ep}(y,y_0) = \frac{2i\ep}{\b^2}\o_m^0(y) +\varphi_{m,\ep}^+(y,y_0) - \varphi_{m,\ep}^-(y,y_0)
\end{equation*} 
and we further denote $\varphi_{m,\ep}(y,y_0)= \varphi_{m,\ep}^+(y,y_0) - \varphi_{m,\ep}^-(y,y_0)$, which solves
\begin{equation*}
\l( \p_y^2 -m^2 +\b^2\frac{1}{(y-y_0+i\ep)^2}\r) \varphi_{m,\ep}(y,y_0) = \b^2\frac{4i\ep (y-y_0)}{((y-y_0)^2 + \ep^2)^2}\varphi_{m,\ep}^-(y,y_0) -\frac{2i\ep}{\b^2}\D_m\o_m^0(y).
\end{equation*}
Multiplying by $\overline{\varphi_{m,\ep}(y,y_0)}$ integrating by parts and proceeding as before, we see
\begin{equation*}
(1-4\b^2)\int_0^1 |\p_y\varphi_{m,\ep}|^2 \d y + m^2 \int_0^1 |\varphi_{m,\ep}|^2 \d y \lesssim \frac{2\ep}{\b^2}\Vert \o_m^0\Vert_{H^2}^2 + \b^2\int_0^1 \frac{4\ep (y-y_0)}{((y-y_0)^2+\ep^2)^2}|\varphi_{m,\ep}^-| |\varphi_{m,\ep}|\d y.
\end{equation*}
Moreover, using Young's inequality we can bound
\begin{equation*}
\begin{aligned}
\b^2\int_0^1 \frac{4\ep (y-y_0)}{((y-y_0)^2+\ep^2)^2}|\varphi_{m,\ep}^-| |\varphi_{m,\ep}|\d y &\leq \b^2 \int_0^1 \frac{2\ep (y-y_0)}{((y-y_0)^2+\ep^2)^2}\l(c_0^2|\varphi_{m,\ep}|^2 +\frac{1}{c_0^2}|\varphi_{m,\ep}^-|^2\r)\d y \\
&\leq c_0^2\b^2 \int_0^1 \frac{1}{y^2}||\varphi_{m,\ep}|^2\d y +\frac{\b^2}{c_0^2}\int_0^1 \frac{2\ep (y-y_0)}{((y-y_0)^2+\ep^2)^2}|\varphi_{m,\ep}^-|^2\d y \\
&\leq 4c_0^2\b^2 \int_0^1 |\p_y\varphi_{m,\ep}|^2 \d y +\frac{\b^2}{c_0^2}\int_0^1 \frac{2\ep (y-y_0)}{((y-y_0)^2+\ep^2)^2}|\varphi_{m,\ep}^-|^2\d y.
\end{aligned}
\end{equation*}
Therefore, absorbing the derivative term in the left hand side for some $c_0$ small enough, we obtain 
\begin{equation}\label{H1 bounds dif stream}
\int_0^1 |\p_y\varphi_{m,\ep}|^2 \d y + m^2\int_0^1 |\varphi_{m,\ep}|^2 \d y\lesssim \ep\Vert \o_m^0\Vert_{H^2}^2 + \int_0^1 \frac{2\ep (y-y_0)}{((y-y_0)^2+\ep^2)^2}|\varphi_{m,\ep}^-|^2\d y.
\end{equation}
Given the uniform bounds in $\ep>0$ from Proposition \ref{H1 bounds boundary streamfunction}, we have that
\begin{equation*}
\begin{aligned}
\lim_{\ep\rightarrow 0}&\l\Vert \int_{-\frac{\b}{m}}^0 \e^{-imy_0t}\big(\e^{m\ep t}\psi_{m,\ep}^-(y,y_0)-\e^{-m\ep t}\psi_{m,\ep}^+(y,y_0) \big)\d y_0\r\Vert_{L^2_y}\\
\quad&=\lim_{\ep\rightarrow 0}\l\Vert \int_{-\frac{\b}{m}}^0 \e^{-imy_0t}\e^{m\ep t}\big(\psi_{m,\ep}^-(y,y_0)-\psi_{m,\ep}^+(y,y_0) \big)\d y_0\r\Vert_{L^2_y} \\
\quad &\leq \lim_{\ep\rightarrow 0} \int_{-\frac{\b}{m}}^0 \l(\frac{2\ep}{\b}\Vert \o_m^0\Vert_{L^2_y} + \l\Vert\varphi_{m,\ep }(y,y_0)\r\Vert_{L^2_y}\r)\d y_0 \\
\quad &=  \lim_{\ep\rightarrow 0} \int_{-\frac{\b}{m}}^0  \l\Vert\varphi_{m,\ep }(y,y_0)\r\Vert_{L^2_y}\d y_0.
\end{aligned}
\end{equation*}
Now, note that with \eqref{H1 bounds dif stream} we can estimate
\begin{equation*}
\begin{aligned}
\lim_{\ep\rightarrow 0} \int_{-\frac{\b}{m}}^0  \l\Vert\varphi_{m,\ep }(y,y_0)\r\Vert_{L^2_y}\d y_0 &\lesssim \lim_{\ep\rightarrow 0} \int_{-\frac{\b}{m}}^0  \l( \int_0^1 \frac{2\ep (y-y_0)}{((y-y_0)^2+\ep^2)^2}|\varphi_{m,\ep}^-|^2\d y\r)^\frac12 \d y_0 \\
&\lesssim \lim_{\ep\rightarrow 0} \ep^\frac12\int_{-\frac{\b}{m}}^0  \l( \int_0^1 \frac{1}{(y-y_0)^2}\d y\r)^\frac12 \d y_0, \\
\end{aligned}
\end{equation*}
where we have used the pointwise bound $|\varphi_{m,\ep}^-(y,y_0)|^2\lesssim y$, obtained from Proposition \ref{H1 bounds boundary streamfunction}. The conclusion follows.

For the density statement, we recall that
\begin{equation*}
\rho_{m,\ep}^-(y,y_0)-\rho_{m,\ep}^+(y,y_0) = \frac{\varphi_{m,\ep}(y,y_0)}{y-y_0-i\ep} + \frac{2i\ep}{(y-y_0)^2+\ep^2}\varphi_{m,\ep}^+(y,y_0),
\end{equation*}
from which, together with Proposition \ref{H1 bounds boundary streamfunction}, we deduce that 
\begin{equation*}
\begin{aligned}
\lim_{\ep\rightarrow 0}&\l\Vert \int_{-\frac{\b}{m}}^0 \e^{-imy_0t}\big(\e^{m\ep_*t}\rho_{m,\ep_*}^-(y,y_0)-\e^{-m\ep_*t}\rho_{m,\ep_*}^+(y,y_0) \big)\d y_0\r\Vert_{L^2_y}\\
\quad&\leq\lim_{\ep\rightarrow 0}\l\Vert \int_{-\frac{\b}{m}}^0 \e^{-imy_0t}\e^{m\ep_*t}\big(\rho_{m,\ep_*}^-(y,y_0)-\rho_{m,\ep_*}^+(y,y_0) \big)\d y_0\r\Vert_{L^2_y} \\
\quad &\lesssim \lim_{\ep\rightarrow 0} \int_{-\frac{\b}{m}}^0 \l(\l\Vert\frac{\varphi_{m,\ep}(y,y_0)}{y-y_0-i\ep}\r\Vert_{L^2_y} + \l\Vert\frac{2i\ep \varphi_{m,\ep}^+(y,y_0)}{(y-y_0)^2+\ep^2}\r\Vert_{L^2_y}\r)\d y_0 \\
\end{aligned}
\end{equation*}
Using the Hardy inequality \eqref{Hardy Ineq}, the estimates from \eqref{H1 bounds dif stream} and the above arguments, we have
\begin{equation*}
\lim_{\ep\rightarrow 0} \int_{-\frac{\b}{m}}^0 \l\Vert\frac{\varphi_{m,\ep}(y,y_0)}{y-y_0-i\ep}\r\Vert_{L^2_y}\d y_0 \lesssim \lim_{\ep\rightarrow 0} \int_{-\frac{\b}{m}}^0 \l\Vert \p_y\varphi_{m,\ep}(y,y_0)\r\Vert_{L^2_y}\d y_0 =0.
\end{equation*}
On the other hand, thanks to the bounds from Proposition \ref{H1 bounds boundary streamfunction}, we also have
\begin{equation*}
\lim_{\ep\rightarrow 0} \int_{-\frac{\b}{m}}^0 \l\Vert\frac{2i\ep \varphi_{m,\ep}^+(y,y_0)}{(y-y_0)^2+\ep^2}\r\Vert_{L^2_y}\d y_0\leq  \lim_{\ep\rightarrow 0} \ep^\frac12 \int_{-\frac{\b}{m}}^0 \l\Vert\frac{2\ep^\frac12(y-y_0)^\frac12 \varphi_{m,\ep}^+(y,y_0)}{(y-y_0)^2+\ep^2}\r\Vert_{L^2_y}\frac{1}{(-y_0)^\frac12}\d y_0=0.
\end{equation*}
With this, the proof is finished.
\end{proof}

We next show that the contribution from the vertical boundaries of the contour integral is also negligible.

\begin{proposition}\label{real vertical limiting absorption principle}
Let $y_0\in\l[-\frac{\b}{m}, 0\r]$. We have that
\begin{equation*}
\lim_{\ep\rightarrow 0}\l\Vert \int_0^\ep \e^{-imy_0t}\big(\e^{mst}\psi_{m,s}^-(y,y_0)+\e^{-ms_*t}\psi_{m,s}^+(y,y_0) \big)\d s\r\Vert_{L^2_y}=0
\end{equation*}
and 
\begin{equation*}
\lim_{\ep\rightarrow 0}\l\Vert \int_0^\ep \e^{-imy_0t}\big(\e^{mst}\rho_{m,\ep_*}^-(y,y_0)+\e^{-mst}\rho_{m,s}^+(y,y_0) \big)\d s\r\Vert_{L^2_y}=0.
\end{equation*}
\end{proposition}

\begin{proof}
The statement follows from Minkowski inequality and the fact that
\begin{equation*}
\lim_{\ep\rightarrow 0} \int_0^\ep  \l\Vert\psi_{m,s}^-(y,y_0)\r\Vert_{L^2_y}+ \l\Vert\psi_{m,s}^+(y,y_0)\r\Vert_{L^2_y} + \l\Vert\rho_{m,s}^-(y,y_0)\r\Vert_{L^2_y} + \l\Vert\rho_{m,s}^+(y,y_0)\r\Vert_{L^2_y}\d s=0,
\end{equation*}
due to the uniform bounds in $s\in[0,\ep]$ of these quantities from Proposition \ref{H1 bounds boundary streamfunction}. 
\end{proof}

We are now in position to carry out the proof of Proposition \ref{prop: contour limiting absorption principle} for the case $\b^2<1/4$.
\begin{proof}[Proof of Proposition \ref{prop: contour limiting absorption principle} ]
Since the resolvent operator $\mathcal{R}(c,L_m)$ is invertible for all $c\in \C$ with $\Re(c)\leq 0$ and $|c|\geq c_0$ for some $c_0>0$, confer Proposition \ref{H1 bounds boundary streamfunction}, we can reduce the contour integral to the boundary of the set $R_\ep:=\l\lbrace c=y_0 + is\in \C: y_0\in \l[-\b/m,0\r], \, s\in[-\ep, \ep] \r\rbrace$. Now, Proposition \ref{real horizontal limiting absorption principle} and Proposition \ref{real vertical limiting absorption principle} show that the integral along $\p R_\ep$ is negligible as $\ep\rightarrow 0$. The Proposition follows.
\end{proof}

We finish the subsection with the proof of Theorem \ref{thm: boundary limiting absorption principle} for $b^2<1/4$, which is a direct consequence of the following Lemma.

\begin{lemma}\label{real Green's limit on boundary}
Let $y,z\in[0,1]$. There exists $\ep_0>0$  such that
\begin{equation*}
\sup_{y,z\in[0,1]}\left| \G_{m,\ep}^-(y,y_0,z)-\G_{m,\ep}^+(y,y_0,z)\right|\lesssim \ep^{2\mu} + \ep^\frac12, \qquad y_0\in \{0,1\},
\end{equation*}
for all $\ep\leq \ep_0$.
\end{lemma}
\begin{proof}
We take $y_0=0$, the other case is analogous. The argument is similar to the one presented for Lemma \ref{Green's limit on boundary}. As before, we need to understand
\begin{equation*}
\G_{m,\ep}^-(y,0,z)- \G_{m,\ep}^+(y,0,z)=2i\Im \l( \G_{m,\ep}^-(y,0,z)\r)
\end{equation*}
For $y\leq z$, we have
\begin{align*}
2i\Im \l( \G_{m,\ep}^-(y,0,z)\r)=\frac{2i}{|{\W_{m,\ep}^-}(0)|^2}\Im\l(\phi_{l,m,\ep}^-(y,0)\phi_{u,m\ep}^-(z,0)\W_{m,\ep}^+(0)\r).
\end{align*} 
Due to the explicit solutions of the Taylor-Goldstein equation, we can find that
\begin{equation*}
\begin{aligned}
{\phi_{l,m,\ep}^-(y,0)}&{\W_{m,\ep}^+(0)} \\
&= 4\mu m\Big( M_+(-i\ep)M_-(y)M_+(i\ep)M_-(1)  + M_-(-i\ep)M_+(y)M_-(i\ep)M_+(1) \Big. \\
&\qquad\qquad - \Big. M_+(-i\ep)M_-(y)M_-(i\ep)M_+(1) - \Big.M_-(-i\ep)M_+(y)M_+(i\ep)M_-(1) +R_1(\ep)\Big),
\end{aligned}
\end{equation*}
where $|R_1(\ep)|\lesssim_{m,\mu} \ep^{\frac12-\mu} |M_-(i\ep)|^2$. Moreover, since now $\overline{M_\pm(\zeta)}=M_\pm(\overline{\zeta})$ for all $\zeta\in\C$, we can write
\begin{equation*}
\begin{aligned}
{\phi_{l,m,\ep}^-(y,0)}{\W_{m,\ep}^+(0)}&= 4\mu m \Big( |M_+(i\ep)|^2M_-(y)M_-(1) + |M_-(i\ep)|^2M_+(y)M_+(1) \Big.\\
&\quad- \Big. M_+(-i\ep)M_-(i\ep)M_-(y)M_+(1) - M_+(i\ep)M_-(-i\ep)M_+(y)M_-(1) + R_1(\ep)\Big).
\end{aligned}
\end{equation*}
On the other hand,
\begin{equation*}
\begin{aligned}
\phi_{u,m,\ep}^-(z,0)&=M_+(1)M_-(z)-M_-(1)M_+(z) +  R_2(\ep) = \phi_u(z) + R_2(\ep)\\
\end{aligned}
\end{equation*}
with $|R_2(\ep)|\lesssim \ep^{\frac12-\mu}$. Thus,
\begin{equation*}
\begin{aligned}
\phi_{l,m,\ep}^-(y,0)&\phi_{u,m,\ep}^-(z,0)\W_{m,\ep}^+(0) \\
&=-4\mu m \phi_u(z) \Big( M_+(-i\ep)M_-(i\ep)M_-(y)M_+(1) + M_+(i\ep)M_-(-i\ep)M_+(y)M_-(1) \Big) \\
&\quad+4\mu m \phi_u(z)\Big(|M_+(i\ep)|^2M_-(y)M_-(1) + |M_-(i\ep)|^2M_+(y)M_+(1)\Big) + R_3(\ep),
\end{aligned}
\end{equation*}
where $|R_3(\ep)|\lesssim \mu m \ep^{\frac12-\mu}|M_-(i\ep)|^2$, uniformly in $y,z\in[0,1]$.
In particular, since $M_\pm(y)\in\R$, for all $y\in[0,1]$, we have
\begin{equation*}
\Im\l(\phi_{l,m,\ep}^-(y,0)\phi_{u,m,\ep}^-(z,0)\W_{m,\ep}^+(0)\r)=-4\mu m \phi_{u,m}(z)\phi_{u,m}(y)\Im\l( M_+(-i\ep)M_-(i\ep)\r) + \Im\l(R_3(\ep)\r).
\end{equation*}
Moreover, due to the symmetry of the Green's function with respect to $y$ and $z$, we also have that
\begin{equation*}
\Im\l(\phi_{l,m,\ep}^-(z,0)\phi_{u,m,\ep}^-(y,0)\W_{m,\ep}^+(0)\r)=-4\mu m \phi_{u,m}(z)\phi_{u,m}(y)\Im\l( M_+(-i\ep)M_-(i\ep)\r) + \Im\l(\widetilde{R}_3(\ep)\r).
\end{equation*}
For the Wronskian, we have from \eqref{def Wronskian} that
\begin{equation*}
\begin{aligned}
\l|\W_{m,\ep}^+(0)\r|=4\mu m|M_-(i\ep)||M_+(1)|\l| 1- \frac{M_+(i\ep)}{M_-(i\ep)}\frac{M_-(1)}{M_+(1)} + R_4(\ep)\r|.
\end{aligned}
\end{equation*}
where $|R_4(\ep)|\lesssim \ep^{\frac12-\mu}$. In particular, for $\ep\leq \ep_0$ small enough we have from Lemma \ref{Comparison bounds real M small argument} that 
\begin{equation*}
\begin{aligned}
\l|\W_{m,\ep}^+(0)\r|\geq 2\mu m|M_-(i\ep)||M_+(1)|.
\end{aligned}
\end{equation*}
Therefore,
\begin{equation*}
\begin{aligned}
\l|\G_{m,\ep}^-(y,0,z)-\G_{m,\ep}^+(y,0,z)\r|&=\frac{2}{|\W_{m,\ep}^+(0)|^2}\l|\Im\l(\phi_{l,m,\ep}^-(y,0)\phi_{u,m,\ep}^-(z,0)\W_{m,\ep}^+(0)\r)\r| \\
&\leq \frac{2}{\mu m}\frac{\l|\phi_{u,m}(z)\phi_{u,m}(y)\Im\l( M_+(-i\ep)M_-(i\ep)\r)\r|}{|M_-(i\ep)|^2|M_+(1)|^2} + R_5(\ep)\\
&\lesssim \ep^{2\mu} + \ep^{\frac12-\mu},
\end{aligned}
\end{equation*}
and the lemma follows.
\end{proof}

\subsection{Integral reduction for $\b^2=1/4$}
The special case in which $\b^2=1/4$ is critical in the sense that the Hardy inequality \eqref{Hardy Ineq} may saturate and thus the derivative bounds in Proposition \ref{H1 bounds boundary streamfunction} are no longer uniform in $\ep>0$. Still, we are able to prove the following result.
\begin{proposition}\label{H1 bounds boundary special streamfunction}
Let $y_0\leq 0$ and $0<\ep\leq 1$. Then, 
\begin{equation*}
\frac{\ep^2}{1+\ep^2}\Vert \p_y \psi_{m,\ep}^\pm(\cdot,y_0)\Vert_{L^2}^2 + m^2 \Vert \psi_{m,\ep}^\pm(\cdot,y_0)\Vert_{L^2}^2 \lesssim \Vert \o_m^0 \Vert_{H^2}^2 + \Vert \rho_m^0 \Vert_{H^2}^2.
\end{equation*}
Moreover,
\begin{equation*}
\begin{aligned}
\int_0^1\frac{\ep (y-y_0)}{((y-y_0)^2+\ep^2)^2}|\varphi_{m,\ep}^\pm|^2 \d y \lesssim \Vert \o_m^0 \Vert_{H^2}^2 + \Vert \rho_m^0 \Vert_{H^2}^2.
\end{aligned}
\end{equation*}
If we further assume that $|-y_0\pm i\ep|\geq c_0$, for some $c_0>0$, then
\begin{equation*}
\frac{c_0^2}{1+c_0^2}\Vert \p_y \psi_{m,\ep}^\pm(y,y_0)\Vert_{L^2}^2 + m^2 \Vert \psi_{m,\ep}^\pm(y,y_0)\Vert_{L^2}^2 \lesssim \frac{1}{c_0^2}\Vert \o_m^0 \Vert_{L^2}^2 + \frac{1}{c_0^4}\Vert \rho_m^0 \Vert_{L^2}^2.
\end{equation*}
In particular, $c=-y_0\pm i\ep$ belongs to the resolvent set  of $L_m$.
\end{proposition}
\begin{proof}
The proof is similar to the one for Proposition \ref{H1 bounds boundary streamfunction}. Here, since $\b^2=1/4$ we estimate
\begin{equation*}
 \frac14\int_0^1\frac{|\psi_{m,\ep}^\pm|^2}{|y-y_0\pm i\ep|^2}\d y \leq \frac{1}{1+c_0^2}\int_0^1\frac{|\psi_{m,\ep}^\pm|^2}{4y^2}\d y \leq \frac{1}{1+c_0^2}\int_0^1 |\p_y\psi_{m,\ep}|^2\d y,
\end{equation*}
which can be absorbed by $\int_0^1 |\p_y\psi_{m,\ep}|^2\d y$, thus producing the desired $H^1$ estimates. 
\end{proof}
The estimate on the $L^2$ norm of the derivative degenerates as $\ep$ becomes small. We may lose pointwise bounds on the solution, and for this reason we investigate more thoroughly the Green's function $\G_{m,\ep}^\pm(y,y_0,z)$ when $-1\ll y_0 \leq 0$. In particular, we have that 
\begin{proposition}\label{Bounds special G small spectral boundary}
Let $y,z\in[0,1]$. There exists $\delta>0$ such that
\begin{equation*}
|\G_{m,\ep}^+(y,y_0,z)|\lesssim |y-y_0+i\ep|^{\frac12}|z-y_0+i\ep|^{\frac12}\l( 1 +  \big|\log m|y-y_0+i\ep|\big|\r)\l( 1 +  \big|\log m|z-y_0+i\ep|\big|\r)
\end{equation*}  
and
\begin{equation*}
|\p_y\G_{m,\ep}^+(y,y_0,z)|\lesssim |y-y_0+i\ep|^{-\frac12}|z-y_0+i\ep|^{\frac12}\l( 1 + \big|\log m|y-y_0+i\ep|\big|\r)\l( 1 + \big|\log m|z-y_0+i\ep|\big|\r)
\end{equation*}  
for all $y_0<0$ and $\ep>0$ with  $|-y_0\pm i\ep|\leq \delta$.
\end{proposition}
We remark that the hidden implicit constant may depend on $m$, but for our purposes this is unimportant.
\begin{proof}
The proof follows the same steps as the one for Proposition \ref{Bounds special G Small y Small z}. We shall obtain suitable estimates on the Wrosnkian. Now, since $y_0<0$, we recall
\begin{equation*}
\W_{m,\ep}^\pm(y_0):=\frac{2m}{\sqrt{\pi}}\Big(W_0(-y_0\pm i\ep)M_0(1-y_0\pm i\ep)-M_0(-y_0\pm i\ep)W_0(1-y_0\pm i\ep)\Big).
\end{equation*}
Using Lemma \ref{Comparison bounds special M small argument} and Lemma \ref{Comparison bounds special M order one argument}, there exists $C>0$ and $\delta>0$ such that 
\begin{equation*}
\l| \frac{W_0(1-y_0\pm i\ep)}{M_0(1-y_0\pm i\ep)}\r|\leq C, \quad \l| \frac{M_0(-y_0\pm i\ep)}{W_0(-y_0\pm i\ep)}\r|\leq \frac{1}{2C},
\end{equation*}
for all $|-y_0\pm i\ep|\leq \delta$. Hence, we can lower bound
\begin{equation*}
\l| \W_{m,\ep}^{\pm}(y_0)\r| \geq \frac{m}{\sqrt{\pi}}\l|W_0(-y_0\pm i\ep)\r|\l|M_0(1-y_0\pm i\ep)\r|
\end{equation*}
and the proposition follows from the asymptotic expansions of the homogeneous solutions that conform the Green's function.
\end{proof}
With the above asymptotics at hand, we are now able to prove the following result.
\begin{proposition}\label{special horizontal limiting absorption principle}
Let $\delta>0$ be given by Proposition \ref{Bounds special G small spectral boundary} and let  $y_0<0$ such that $|y_0|\leq \frac{\delta}{2}$. We have that
\begin{equation*}
\l\Vert \int_{-\frac{\delta}{2}}^0 \e^{-imy_0t}\big(\e^{m\ep t}\psi_{m,\ep}^-(y,y_0)-\e^{-m\ep t}\psi_{m,\ep}^+(y,y_0) \big)\d y_0\r\Vert_{L^2_y}\lesssim \ep^\frac12
\end{equation*}
and also
\begin{equation*}
\l\Vert \int_{-\frac{\delta	}{2}}^0 \e^{-imy_0t}\big(\e^{m\ep t}\rho_{m,\ep }^-(y,y_0)-\e^{-m\ep t}\rho_{m,\ep }^+(y,y_0) \big)\d y_0\r\Vert_{L^2_y}\lesssim \ep^\frac12,
\end{equation*}
for all $\ep>0$ such that $|-y_0 +i\ep|\leq \delta$.
\end{proposition}
\begin{proof}
Following the same strategy as in the proof of Proposition \ref{real horizontal limiting absorption principle}, we see that $\varphi_{m,\ep}(y,y_0)$ satisfies
\begin{equation*}
m^2 \Vert \varphi_{m,\ep}\Vert_{L^2}^2 \lesssim \ep \Vert \o_m^0\Vert_{H^2}^2 + \int_0^1 \frac{2\ep(y-y_0)}{((y-y_0)^2 + \ep^2)^2}\l(|\varphi_{m,\ep}^-|^2 + |\varphi_{m,\ep}^+|^2\r)\d y.
\end{equation*}
In particular, using the asymptotic bounds from Proposition \ref{Bounds special G small spectral boundary} we can estimate 
\begin{equation*}
\begin{aligned}
\int_0^1 \frac{2\ep(y-y_0)}{((y-y_0)^2 + \ep^2)^2}\l(|\varphi_{m,\ep}^-|^2 + |\varphi_{m,\ep}^+|^2\r)\d y &\lesssim \int_0^1 \frac{\ep(y-y_0)}{|y-y_0 + i\ep|^3}\l(1 + \l|\log |y-y_0+i\ep|\r|\r)^2 \d y \\
&\lesssim \int_0^1 \frac{\ep(y-y_0)}{|y-y_0 + i\ep|^{\frac72}} \d y \\
&\lesssim \ep \int_0^1 \frac{1}{(y-y_0)^\frac52} \d y \\
&\lesssim \ep\l(1 + (-y_0)^{-\frac32}\r).
\end{aligned}
\end{equation*}
We conclude the first part of the proof upon noting that
\begin{equation*}
 \int_{-\frac{\delta}{2}}^0 \l\Vert \varphi_{m,\ep}(y,y_0) \r\Vert_{L^2_y}\d y_0\lesssim \ep^\frac12\Vert \o_m^0\Vert_{H^2} + \ep^\frac12 \int_{-\frac{\delta}{2}}^0 \l(1 + (-y_0)^{-\frac32}\r)^\frac12\d y_0 \lesssim \ep^\frac12.
\end{equation*}
For the second part of the proposition, from \eqref{def rho} we have
\begin{equation*}
\rho_{m,\ep}^-(y,y_0)-\rho_{m,\ep}^+(y,y_0) = \frac{\varphi_{m,\ep}(y,y_0)}{y-y_0-i\ep} + \frac{2i\ep}{(y-y_0)^2+\ep^2}\varphi_{m,\ep}^+(y,y_0)
\end{equation*}
and we write 
\begin{equation*}
\varphi_{m,\ep}(y,y_0)=\int_0^1 \G_{m,\ep}^-(y,y_0,z) \l( \frac{i\ep(z-y_0)}{((z-y_0)^2 + \ep^2)^2}\varphi_{m,\ep}^-(z,y_0)-8i\ep\D_m\o_m^0(z)\r) \d z.
\end{equation*}
In particular, using Proposition  \ref{H1 bounds boundary special streamfunction} and Proposition \ref{Bounds special G small spectral boundary} we estimate 
\begin{equation*}
\begin{aligned}
&\l| \int_0^1 \G_{m,\ep}^-(y,y_0,z) \frac{i\ep(z-y_0)}{((z-y_0)^2 + \ep^2)^2}\varphi_{m,\ep}^-(z,y_0) \d z \r| \\
&\qquad\lesssim \ep^\frac12 \l(\int_0^1 \frac{(z-y_0)}{((z-y_0)^2 + \ep^2)^2}|\G_{m,\ep}^-(y,y_0,z)|^2\d z\r)^\frac12 \\
&\qquad \lesssim \ep^\frac12|y-y_0-i\ep|^\frac12 (1+ \l| \log	|y-y_0-i\ep|\r|) \l( \int_0^1 \frac{(z-y_0)}{|z-y_0-i\ep|^3}(1+ \l| \log	|z-y_0-i\ep|\r|)^2 \d z \r)^\frac12 \\
&\qquad \lesssim \ep^\frac12|y-y_0-i\ep|^\frac12 (1+ \l| \log	|y-y_0-i\ep|\r|)\l( \int_0^1 \frac{1}{(z-y_0)^{2+\frac14}} \d z \r)^\frac12 \\
&\qquad \lesssim \ep^\frac12|y-y_0-i\ep|^\frac12 (1+ \l| \log	|y-y_0-i\ep|\r|)\l( 1+(-y_0)^{-\frac12-\frac18} \r).
\end{aligned}
\end{equation*}
With this pointwise bound, we obtain
\begin{equation*}
\begin{aligned}
\l\Vert \frac{\varphi_{m,\ep}(y,y_0)}{y-y_0-i\ep} \r\Vert_{L^2_y} &\lesssim \ep^\frac12\l( 1+(-y_0)^{-\frac12-\frac18} \r) \l( \int_0^1 |y-y_0-i\ep|^{-1} (1+ \l| \log	|y-y_0-i\ep|\r|)^2 \d y\r)^\frac12 \\
&\lesssim \ep^\frac12\l( 1+(-y_0)^{-\frac12-\frac18} \r)\l( \int_0^1 |y-y_0-i\ep|^{-1-\frac14} \d y\r)^\frac12 \\
&\lesssim \ep^\frac12\l( 1+(-y_0)^{-\frac34} \r)
\end{aligned}
\end{equation*}
and thus
\begin{equation*}
\int_{-\frac{\delta}{2}}^0 \l\Vert \frac{\varphi_{m,\ep}(y,y_0)}{y-y_0-i\ep} \r\Vert_{L^2_y} \d y_0 \lesssim \ep^\frac12.
\end{equation*}
On the other hand, from the bounds obtained in Proposition \ref{H1 bounds boundary special streamfunction}, we have
\begin{equation*}
\int_{-\frac{\delta}{2}}^0 \l\Vert \frac{2i\ep}{(y-y_0)^2+\ep^2}\varphi_{m,\ep}^+(y,y_0) \r\Vert_{L^2_y} \d y_0 \lesssim \ep^\frac12 \int_{-\frac{\delta}{2}}^0 \frac{1}{(-y_0)^\frac12}\l\Vert \frac{\ep^\frac12 (y-y_0)^\frac12}{(y-y_0)^2+\ep^2}\varphi_{m,\ep}^+(y,y_0) \r\Vert_{L^2_y} \d y_0\lesssim \ep^\frac12,
\end{equation*}
and the proof is concluded.
\end{proof}
Similarly, the contribution from the resolvent integral along the vertical boundaries of the contour is also negligible. 
\begin{proposition}\label{special vertical limiting absorption principle}
Let $y_0\in\l[-\b/m, 0\r]$. We have that
\begin{equation*}
\l\Vert \int_0^\ep \e^{-imy_0t}\big(\e^{mst}\psi_{m,s}^-(y,y_0)+\e^{-ms_*t}\psi_{m,s}^+(y,y_0) \big)\d s\r\Vert_{L^2_y}\lesssim \ep
\end{equation*}
and 
\begin{equation*}
\l\Vert \int_0^\ep \e^{-imy_0t}\big(\e^{mst}\rho_{m,s}^-(y,y_0)+\e^{-mst}\rho_{m,s}^+(y,y_0) \big)\d s\r\Vert_{L^2_y}\lesssim \ep^\frac14.
\end{equation*}
\end{proposition}
\begin{proof}
The first part concerning the stream-functions $\psi_{m,\ep}^\pm(y,y_0)$ is a direct consequence of the uniform $L^2$ bounds of $\psi_{m,\ep}^\pm(y,y_0)$ obtained in Proposition \ref{H1 bounds boundary special streamfunction}. As for the density statement, we use \eqref{def rho}; thanks to the asymptotic bounds from Proposition \ref{Bounds special G small spectral boundary} we further observe that
\begin{equation*}
\begin{aligned}
\int_0^\ep \l\Vert \frac{\varphi_{m,\ep}^\pm(y,y_0)}{y-y_0\pm is}\r\Vert_{L^2_y}\d s &\lesssim \int_0^\ep \l( \int_0^1 |y-y_0\pm is|^{-1}(1+ \l| \log	|y-y_0-i\ep|\r|)^2 \d y\r)^\frac12 \d s \\
&\lesssim \int_0^\ep \l( \int_0^1 |y-y_0\pm is|^{-\frac32}\d y\r)^\frac12 \d s 
\lesssim \int_0^\ep |s|^{-\frac34}\d s.
\end{aligned}
\end{equation*}
With the above estimate, the bound follows swiftly.
\end{proof}
We are now in position to prove Proposition \ref{prop: contour limiting absorption principle} for the special case $\b^2=1/4$.

\begin{proof}[Proof of Proposition \ref{prop: contour limiting absorption principle}]
Let $\delta>0$ be given by Proposition \ref{Bounds special G small spectral boundary}. For all $\ep<\frac{\delta}{2}$, we introduce the rectangular region $R_\ep:=\l\lbrace c=y_0 + is\in \C: y_0\in \l[-\delta/2,0\r], \, s\in[-\ep, \ep] \r\rbrace$. From Proposition \ref{special horizontal limiting absorption principle} and Proposition \ref{special vertical limiting absorption principle} we conclude that 
\begin{equation*}
\l\Vert \int_{\p R_\ep}\e^{-imct}\mathcal{R}(c,L_m)\d c \r\Vert_{L^2_y}\lesssim \ep^\frac14,
\qquad
\l\Vert \int_{\p (R\setminus R_\ep)}\e^{-imct}\mathcal{R}(c,L_m)\d c \r\Vert_{L^2_y}=0,
\end{equation*}
because any $c\in R\setminus R_\ep$ belongs to the resolvent set of the operator $L_m$. Indeed, any $c\in R\setminus R_\ep$ is such that $\Re(c)\leq 0$, $|c|\geq \frac{\delta}{2}$, and we can see from Proposition \ref{H1 bounds boundary special streamfunction} that $\mathcal{R}(c,L_m)$ is invertible.
\end{proof}

Finally, in order to prove Theorem \ref{thm: boundary limiting absorption principle} for $\b^2=1/4$, we  state and prove the following key Lemma, from which the Theorem easily follows.
\begin{lemma}\label{special Green's limit on boundary}
Let $y_0=0$ and $y,z\in[0,1]$. Then, there exists $\ep_0>0$ such that 
\begin{equation*}
\sup_{y,z\in[0,1]}\left| \G_{m,\ep}^-(y,0,z)-\G_{m,\ep}^+(y,0,z)\right|\lesssim \frac{1}{\log\l(\frac{4}{\ep}\r)} + \ep^\frac14,
\end{equation*}
for all $\ep\leq \ep_0$.
\end{lemma}

\begin{proof}
We have $\G_{m,\ep}^-(y,0,z)- \G_{m,\ep}^+(y,0,z)=2i\Im \l( \G_{m,\ep}^-(y,0,z)\r)$ and for $y\leq z$, 
\begin{equation*}
2i\Im \l( \G_{m,\ep}^-(y,0,z)\r)=\frac{2i}{|{\W_{m,\ep}^-}(0)|^2}\Im\l(\phi_{l,m,\ep}^-(y,0)\phi_{u,m,\ep}^-(z,0)\W_{m,\ep}^+(0)\r).
\end{equation*} 
Now, using Proposition \ref{Prop def special Green's function}, Lemma \ref{asymptotic expansion W} and Lemma \ref{Comparison bounds special M small argument},
\begin{equation*}
\begin{aligned}
{\phi_{l,m,\ep}^-(y,0)}&{\W_{m,\ep}^+(0)} \\
&= \frac{2m}{\sqrt{\pi}}\Big( |W_0(i\ep)|^2M_0(y)M_0(1) + |M_0(i\ep)|^2W_0(y)W_0(1)\Big. \\
&\qquad\qquad - \Big.W_0(i\ep)M_0(-i\ep)W_0(y)M_0(1) -  W_0(-i\ep)M_0(i\ep)M_0(y)W_0(1) +R_1(\ep)\Big),
\end{aligned}
\end{equation*}
where $|R_1(\ep)|\lesssim_{m,\mu} \ep^{\frac14} |W_0(i\ep)|^2$. Similarly,
\begin{equation*}
\begin{aligned}
\phi_{u,m,\ep}^-(z,0)&=W_0(1)M_0(z)-M_0(1)W_0(z) +  R_2(\ep) =: \phi_{u,m}(z) + R_2(\ep)\\
\end{aligned}
\end{equation*}
with $|R_2(\ep)|\lesssim \ep^{\frac14}$. Thus,
\begin{equation*}
\begin{aligned}
\phi_{l,m,\ep}^-(y,0)&\phi_{u,m,\ep}^-(z,0)\W_{m,\ep}^+(0) \\
&=-\frac{2m}{\sqrt{\pi}} \phi_{u,m}(z) \Big( W_0(-i\ep)M_0(i\ep)M_0(y)W_0(1) + W_0(i\ep)M_0(-i\ep)W_0(y)M_0(1) \Big) \\
&\quad+\frac{2m}{\sqrt{\pi}} \phi_{u,m}(z)\Big( |W_0(i\ep)|^2M_0(y)M_0(1) + |M_0(i\ep)|^2W_0(y)W_0(1)\Big) + R_3(\ep),
\end{aligned}
\end{equation*}
where $|R_3(\ep)|\lesssim m \ep^{\frac14}|W_0(i\ep)|^2$, uniformly in $y,z\in[0,1]$. In particular, since $M_0(y)\in\R$ and $W_0(y)\in \R$, for all $y\in[0,1]$, we have
\begin{equation*}
\Im\l(\phi_{l,m,\ep}^-(y,0)\phi_{u,m,\ep}^-(z,0)\W_{m,\ep}^+(0)\r)=-\frac{2m}{\sqrt{\pi}} \phi_{u,m}(z)\phi_{u,m}(y)\Im\l( W_0(-i\ep)M_0(i\ep)\r) + \Im\l(R_3(\ep)\r).
\end{equation*}
Due to symmetry, we also have
\begin{equation*}
\Im\l(\phi_{l,m,\ep}^-(z,0)\phi_{u,m,\ep}^-(y,0)\W_{m,\ep}^+(0)\r)=-\frac{2m}{\sqrt{\pi}} \phi_{u,m}(z)\phi_{u,m}(y)\Im\l( W_0(-i\ep)M_0(i\ep)\r) + \Im\l(\widetilde{R}_3(\ep)\r).
\end{equation*}
For the Wronskian, we have from \eqref{def special Wronskian} that
\begin{equation*}
\begin{aligned}
\l|\W_{m,\ep}^+(0)\r|=\frac{2m}{\sqrt{\pi}}|W_0(i\ep)||M_0(1)|\l| 1- \frac{M_0(i\ep)}{W_0(i\ep)}\frac{W_0(1)}{M_0(1)} + R_4(\ep)\r|.
\end{aligned}
\end{equation*}
where $|R_4(\ep)|\lesssim \ep$. In particular, for $\ep\leq \ep_0$ small enough we have from Lemma \ref{Comparison bounds special M small argument} that 
\begin{equation*}
\begin{aligned}
\l|\W_{m,\ep}^+(0)\r|\geq \frac{m}{\sqrt{\pi}}|W_0(i\ep)||M_0(1)|.
\end{aligned}
\end{equation*}
Therefore,
\begin{equation*}
\begin{aligned}
\l|\G_{m,\ep}^-(y,0,z)-\G_{m,\ep}^+(y,0,z)\r|&=\frac{2}{|\W_{m,\ep}^+(0)|^2}\l|\Im\l(\phi_{l,m,\ep}^-(y,0)\phi_{u,m,\ep}^-(z,0)\W_{m,\ep}^+(0)\r)\r| \\
&\leq \frac{4\sqrt{\pi}}{m}\frac{\l|\phi_{u,m}(z)\phi_{u,m}(y)\Im\l( W_0(-i\ep)M_0(i\ep)\r)\r|}{|W_0(i\ep)|^2|M_0(1)|^2} + R_5(\ep)\\
&\lesssim \l|\frac{M_0(i\ep)}{W_0(i\ep)}\r| + \ep^\frac14,
\end{aligned}
\end{equation*}
and the conclusion follows from Lemma \ref{Comparison bounds special M small argument}.
\end{proof}

\section{Bounds on solutions to the inhomogeneous Taylor-Goldstein equation}\label{sec: bounds inhom TG}
This section provides bounds for solutions $\Phi_{m,\ep}$ to the inhomogeneous Taylor-Goldstein equation \eqref{eq:inhom TG} with boundary conditions $\Phi_{m,\ep}(0,y_0)=\Phi_{m,\ep}(1,y_0)=0$. The following lemma relates regions of the interval $(0,1)$ that are far away from a fixed $y_0\in[0,1]$ to nearby regions of $y_0$. 
\begin{lemma}\label{lemma:entangle ineq}
Let $y_0\in[0,1]$, $n\geq 1$ and $\Phi_{m,\ep}$ be the solution to \eqref{eq:inhom TG}.
Then, we have that
\begin{equation*}
\Vert \p_y \Phi_{m,\ep} \Vert_{L^2_y(J_3^c)}^2 + m^2 \Vert  \Phi_{m,\ep} \Vert_{L^2_y(J_3^c)}^2 \lesssim m^2 \Vert \Phi_{m,\ep} \Vert_{L^2_y(J_2^c\cap J_3)}^2 + \frac{1}{m^2}\Vert f \Vert_{L^2_y(J_2^c)}^2.
\end{equation*}
\end{lemma}

\begin{proof}
For $y_n=y_0+\frac{n\b}{m}$, the lemma follows from the energy inequality 
\begin{equation*}
\frac{1}{2}\int_{y_3}^1 \left[|\p_y \Phi_{m,\ep}|^2 + m^2|\Phi_{m,\ep}|^2 \right]\d y \leq \frac{m^2}{\b^2}\int_{y_2}^{y_3}|\Phi_{m,\ep}|^2 \d y + \int_{y_2}^1 |f||\Phi_{m,\ep}|\d y,
\end{equation*}
and Young's inequality to absorb the potential term. We omit the details.
\end{proof}

With the above lemma we are in position to provide bounds on the solution to \eqref{eq:inhom TG}.
\begin{proposition}\label{L2 bounds inhom TG solution}
Let $\Phi_{m,\ep}$ be the solution to \eqref{eq:inhom TG}. Then
\begin{itemize}
\item If $m|y-y_0|\leq 3\b$ and $\b^2\neq 1/4$, then 
\begin{equation*}
|y-y_0+i\ep|^{-\frac12+\mu} |\Phi_{m,\ep}(y,y_0)|+ |y-y_0+i\ep|^{\frac12+\mu} |\p_y \Phi_{m,\ep}(y,y_0)|\lesssim \frac{1}{m^{1+\mu}}\Vert f \Vert_{L^2_y}.
\end{equation*}
\item If $m|y-y_0|\leq 3\b$ and $\b^2 = 1/4$, then 
\begin{equation*}
|y-y_0+i\ep|^{-\frac12} |\Phi_{m,\ep}(y,y_0)|+ |y-y_0+i\ep|^{\frac12} |\p_y \Phi_{m,\ep}(y,y_0)|\lesssim \frac{1}{m} \l( 1 + \big| \log \l(m|y-y_0\pm i\ep|\r) \big| \r)  \Vert f \Vert_{L^2_y}.
\end{equation*}
\item If $m|y-y_0|\geq 3\b$ then 
\begin{equation*}
m\Vert \Phi_{m,\ep}(y,y_0)\Vert_{L^2_y(J_3^c)}+\Vert\p_y \Phi_{m,\ep}(y,y_0)\Vert_{L^2_y(J_3^c)}\lesssim \frac{1}{m}\Vert f \Vert_{L^2_y}
\end{equation*}
and
\begin{equation*}
|\p_y \Phi_{m,\ep}(y,y_0)|\lesssim \Vert f \Vert_{L^2_y}.
\end{equation*}
\end{itemize}
\end{proposition}

\begin{proof}
The first part is a straightforward application of the bounds on the Green's function from Theorem \ref{L2 bounds of G} and the Cauchy-Schwartz inequality, once we write $\Phi_{m,\ep}(y,y_0)=\int_0^1 \G_{m,\ep}^+(y,y_0,z)f(z,y_0)\d z$. The second part of the proposition follows from the first part, which gives $m \Vert \Phi_{m,\ep} \Vert_{L^2_y(J_2^c\cap J_3)}\lesssim \frac{1}{m}\Vert f \Vert_{L^2_y}$ and Lemma \ref{lemma:entangle ineq}. For the pointwise bound, assume without loss of generality that $y_0+\frac{3\b}{m}<y\leq 1$. Then, let $y_3=y_0+\frac{3\b}{m}$ and write
\begin{equation*}
\p_y\Phi_{m,\ep}(y,y_0)=\p_y \Phi_{m,\ep}(y_3,y_0) + \int_{y_3}^y \l[\l(m^2-\b^2\frac{1}{(y'-y_0+i\ep)^2}\r)\Phi_{m,\ep}(y',y_0) + f(y')\r] \d y'.
\end{equation*}
Now, $|y_3-y_0|=\frac{3\b}{m}$ so that we estimate $|\p_y \Phi_{m,\ep}(y_3,y_0)|\lesssim \frac{1}{m^{1+\mu}}\l|\frac{\b}{m}\r|^{-\frac12-\mu}\lesssim m^{-\frac12}$. Similarly, we use the second part of the proposition to estimate the remaining terms in $L^2_y(J_3^c)$ and obtain the desired conclusion.
\end{proof}

\section{Boundary terms estimates}\label{sec: boundary term estimates}
The purpose of this section is to obtain estimates on the boundary terms that appear in the expressions for $\p_{y_0}\psi_{m,\ep}^\pm(y,y_0)$ and other related derivatives. We begin by recording the following results, which will be used throughout the entire section.
\begin{proposition}\label{prop: pointwise dzG bounds}
Let $\b^2\neq 1/4$. There exists $\ep_0>0$ such that for all $y,y_0\in [0,1]$ with $m|y-y_0|\leq 3\b$ there holds
\begin{equation*}
|y-y_0\pm i \ep|^{-\frac12-\mu}|\p_z\G_{m,\ep}^\pm(y,y_0,0)| +  |y-y_0\pm i\ep|^{\frac12+\mu}|\p_y\p_z \G_{m,\ep}^\pm(y,y_0,0)|\lesssim m^{\frac12-\mu}\frac{1}{|M_-(y_0\mp i \ep))|},
\end{equation*}
and 
\begin{equation*}
|y-y_0\pm i \ep|^{-\frac12-\mu}|\p_z\G_{m,\ep}^\pm(y,y_0,1)| +  |y-y_0\pm i\ep|^{\frac12+\mu}|\p_y\p_z \G_{m,\ep}^\pm(y,y_0,1)|\lesssim m^{\frac12-\mu}\frac{1}{|M_-(1-y_0\pm i \ep))|},
\end{equation*}
for all $0\leq \ep \leq \ep_0$.
\end{proposition}
\begin{proof}
For $z=0$, note that we have the explicit expression
\begin{equation*}
\p_z\G_{m,\ep}^\pm(y,y_0,0)=\frac{M_+(1-y_0\pm i\ep)M_-(y-y_0\pm i\ep)-M_-(1-y_0\pm i\ep)M_+(y-y_0\pm i\ep)}{M_+(1-y_0\pm i\ep)M_-(-y_0\pm i\ep)-M_-(1-y_0\pm i\ep)M_+(-y_0\pm i\ep)},
\end{equation*}
so that 
\begin{equation*}
\p_y \p_z\G_{m,\ep}^\pm(y,y_0,0)=2m \frac{M_+(1-y_0\pm i\ep)M_-'(y-y_0\pm i\ep)-M_-(1-y_0\pm i\ep)M_+'(y-y_0\pm i\ep)}{M_+(1-y_0\pm i\ep)M_-(-y_0\pm i\ep)-M_-(1-y_0\pm i\ep)M_+(-y_0\pm i\ep)}.
\end{equation*}
If $m|y-y_0|\leq 3\b$, we use the bounds on the Wronskian from Proposition \ref{Bounds Small y Small z}. For $\b^2>1/4$, the conclusion is straightforward. For $\b^2<1/4$, we take a closer look to the Wronskian estimates obtained on the proof of Proposition \ref{bounds real G small y small z}. The bounds are a consequence of Lemma \ref{lower bounds M}, \ref{growth bounds real M large argument}-\ref{Comparison bounds real M order one argument}. The argument for $z=1$ is similar, we omit the details.
\end{proof}

\begin{proposition}\label{prop: special pointwise dzG bounds}
Let $\b^2= 1/4$. There exists $\ep_0>0$ such that for all $y,y_0\in [0,1]$ with $m|y-y_0|\leq 3\b$ there holds
\begin{equation*}
|y-y_0\pm i \ep|^{-\frac12}|\p_z\G_{m,\ep}^\pm(y,y_0,0)| +  |y-y_0\pm i\ep|^{\frac12}|\p_y\p_z \G_{m,\ep}^\pm(y,y_0,0)|\lesssim m^{\frac12}\frac{1+ \big| \log \l( m|y-y_0\pm i\ep|\r) \big|}{|M_0(y_0 \mp i \ep))|},
\end{equation*}
and 
\begin{equation*}
|y-y_0\pm i \ep|^{-\frac12}|\p_z\G_{m,\ep}^\pm(y,y_0,1)| +  |y-y_0\pm i\ep|^{\frac12}|\p_y\p_z \G_{m,\ep}^\pm(y,y_0,1)|\lesssim m^{\frac12}\frac{1+ \big| \log \l( m|y-y_0\pm i\ep| \r) \big|}{|M_0(1-y_0\mp i \ep))|},
\end{equation*}
for all $0\leq \ep\leq \ep_0$.
\end{proposition}

\begin{proof}
Since $m|y-y_0|\leq 3\b$, the proof follows the same ideas to show Proposition \ref{Bounds special G Small y Small z}, with the help of Lemma \ref{lower bounds M}, \ref{growth bounds special M large argument}-\ref{Comparison bounds special M order one argument}, we omit the details.
\end{proof}

\subsection{Estimates for first order boundary terms}
This subsection is devoted to obtain estimates on
$$\B_{m,\ep}^\pm(y,y_0,z)=\p_z\G_{m,\ep}^\pm(y,y_0,z)\p_z\varphi_{m,\ep}^\pm(z,y_0)$$
for $z=0$ and $z=1$ under the assumption that $m|y-y_0|\leq 3\b$. In what follows, we shall argue for $z=0$, the statements and proofs for $z=1$ are similar and we thus omit them.
We begin by providing bounds for $\p_z\varphi_{m,\ep}^\pm(0,y_0)$.

\begin{proposition}\label{prop: dyVarphi0 estimates}
Let $y_0\in[0,1]$, we have the following.
\begin{itemize}
\item If $my_0\leq 3\b$, then
$
|\p_y\varphi_{m,\ep}^\pm(0,y_0)|\lesssim m^{-\frac12} Q_{0,m}.
$
\item If $my_0\geq 3\b$, then
$
|\p_y\varphi_{m,\ep}^\pm(0,y_0)|\lesssim Q_{0,m}
$.
\end{itemize}
\end{proposition}

For the proof, we assume that $y_0<1/2$. Otherwise, the proposition follows from Proposition \ref{L2 bounds inhom TG solution}. Note that from  \eqref{eq:Fdata} and \eqref{eq:defvarphi}, there holds
\begin{align*}
\p_y\varphi_{m,\ep}^\pm(0,y_0)&= \int_0^1 \p_y\G_{m,\ep}^\pm(0,y_0,z) F_{m,\ep}^\pm(z,y_0) \d z \\
&= \int_0^1 \p_y\G_{m,\ep}^\pm(0,y_0,z)\l(F_{m}(z)+ \frac{y_0\mp i\ep}{\b^2}\D_m\o_m^0(z)\r) \d z,
\end{align*}
Further observe that, due to \eqref{eq:specassume} we have 
\begin{equation*}
\begin{aligned}
\int_0^1 \p_y\G_{m,\ep}^\pm(0,y_0,z) F_{m}^\pm(z,0) \d z &= -4(\mu +i\nu)m \int_0^1 \frac{\phi_{u,m,\ep}^\pm(z,y_0)}{\W_{m,\ep}^\pm(y_0)}F_{m}(z,0)\d z \\
&= -4(\mu + i \nu)m\int_0^1 \frac{\phi_{u,m,\ep}^\pm(z,y_0)- \phi_{u,m}(z)}{\W_{m,\ep}^\pm(y_0)}F_{m}(z,0)\d z
\end{aligned}
\end{equation*}
and we define 
$$f_{m,\ep}^\pm(z,y_0):=\phi_{u,m,\ep}^\pm(z,y_0) - \phi_{u,m}(z).$$ 

\subsubsection{Estimates on $f_{m,\ep}^\pm$ for $\b^2\neq 1/4$}
From the explicit formulas \eqref{eq:homoup} and \eqref{eq: phium}, we have
\begin{equation*}
\begin{aligned}
f_{m,\ep}^\pm(z,y_0)&=M_+(1-y_0\pm i\ep)M_-(z-y_0\pm i\ep) - M_+(1)M_-(z) \\
&\quad - M_-(1-y_0\pm i\ep)M_+(z-y_0\pm i\ep) + M_-(1)M_+(z)
\end{aligned}
\end{equation*}
and we can obtain the next result.
\begin{proposition}\label{f bounds small z}
Let  $z,y_0\in[0,1]$ such that $my_0\leq 3\b$ and $mz\leq 6\b$. Let $0\leq \ep\leq \min\l( \frac{\b}{m},\frac{1}{2m}\r)$. Then,
\begin{equation*}
|f_{m,\ep}(z,y_0)|\lesssim m^{\frac12-\mu}|y_0\pm i\ep|^{\frac12-\mu}|M_+(1-y_0\pm i\ep)|.
\end{equation*}
In particular, $\Vert f_{m,\ep}\Vert_{L^2_y(J)}\lesssim m^{-\mu}|y_0\pm i\ep|^{\frac12-\mu}|M_+(1-y_0\pm i\ep)|$.
\end{proposition}
\begin{proof}
We shall assume $\b^2<1/4$, the case $\b^2>1/4$ is analogous and easier. We write 
\begin{equation*}
\begin{aligned}
M_+(1-y_0\pm i\ep)M_-(z-y_0\pm i\ep) - M_+(1)M_-(z)& =M_+(1-y_0\pm i\ep)\Big(M_-(z-y_0\pm i\ep) - M_-(z)\Big) \\
&\quad + M_-(z)\Big( M_+(1-y_0\pm i\ep) - M_+(1)\Big)
\end{aligned}
\end{equation*}
and
\begin{equation*}
\begin{aligned}
M_-(1-y_0\pm i\ep)M_+(z-y_0\pm i\ep) - M_-(1)M_+(z)& =M_-(1-y_0\pm i\ep)\Big(M_+(z-y_0\pm i\ep) - M_+(z)\Big) \\
&\quad + M_+(z)\Big( M_-(1-y_0\pm i\ep) - M_-(1)\Big).
\end{aligned}
\end{equation*}
Firstly, we estimate
\begin{equation*}
M_+(1-y_0\pm i\ep) - M_+(1)=\int_0^1 \frac{\d}{\d s}M_+(1+s(-y_0\pm i\ep))\d s = (-y_0\pm i\ep)\int_0^1 M_+'(1+s(-y_0\pm i\ep))\d s
\end{equation*}
and we divide our argument as follows. Let $N_{\mu,0}$ be given as in Lemma \ref{growth bounds real M large argument}.

For $m\leq N_{\mu,0}$, we use Lemma \ref{asymptotic expansion M} and the fact that $y_0\leq 1/2$ to bound
\begin{equation*}
\begin{aligned}
|M_+(1-y_0\pm i\ep) - M_+(1)|&\lesssim m^{\frac12+\mu}|y_0\pm i\ep|\int_0^1\frac{\d s}{|1+s(-y_0\pm i\ep)|^{\frac12 - \mu}} \\
&\lesssim m^{\frac12+\mu}|y_0\pm i\ep| \\
&\lesssim |y_0\pm i\ep||M_+(1-y_0 \pm i\ep)|.
\end{aligned}
\end{equation*}
In the last inequality, we have used Lemma \ref{lower bounds M}, \ref{Comparison bounds real M order one argument} and \ref{growth bounds real M large argument}. Similarly,
\begin{equation*}
|M_-(1-y_0\pm i\ep) - M_-(1)|\lesssim |y_0\pm i\ep||M_-(1-y_0 \pm i\ep)|\lesssim |y_0\pm i\ep||M_+(1-y_0 \pm i\ep)|,
\end{equation*}
where we have used Lemma \ref{lower bounds M} and Lemma \ref{Comparison bounds real M order one argument} to deduce $|M_-(1-y_0 \pm i\ep)|\lesssim |M_+(1-y_0 \pm i\ep)|$.

For $m\geq N_{\mu,0}$, we claim that
\begin{equation*}
\l| \frac{M_+'(1+s(-y_0\pm i\ep))}{M_+(1-y_0\pm i\ep)}   \r|\lesssim 1.
\end{equation*}
Indeed, this follows from 
\begin{equation*}
\l| \frac{M_+'(1+s(-y_0\pm i\ep))}{M_+(1-y_0\pm i\ep)}   \r| = \l| \frac{M_+'(1+s(-y_0\pm i\ep))}{M_+(1+s(-y_0\pm i\ep))}\r| \l|\frac{M_+(1+s(-y_0\pm i\ep)) }{M_+(1-y_0\pm i\ep)}   \r|
\end{equation*}
and the corresponding bounds from Lemma \ref{growth bounds M large argument} since $2m(1-y_0)\geq m\geq N_{\mu,0}$. Hence, we have
\begin{equation*}
\begin{aligned}
|M_+(1-y_0\pm i\ep) - M_+(1)|&\leq 2m|y_0\pm i\ep|\int_0^1 \l|M_+'(1+s(-y_0\pm i\ep))\r| \d s \\
&\lesssim m|y_0\pm i\ep|\l| M_+(1-y_0\pm i\ep) \r|. \\
\end{aligned}
\end{equation*}
Similarly, we also have 
\begin{equation*}
|M_-(1-y_0\pm i\ep) - M_-(1)|\lesssim m|y_0\pm i\ep|\l| M_-(1-y_0\pm i\ep) \r| \lesssim m|y_0\pm i\ep|\l| M_+(1-y_0\pm i\ep) \r|,
\end{equation*}
where we have used Lemma \ref{growth bounds real M large argument} to deduce $ \l| M_-(1-y_0\pm i\ep) \r| \lesssim \l| M_+(1-y_0\pm i\ep) \r|$. We next turn our attention to the bounds for $M_-(z-y_0\pm i\ep) - M_-(z)$. As before, we consider two cases.

\bullpar{Case 1} 
For $2y_0\leq z$ we estimate
\begin{equation*}
M_-(z-y_0\pm i\ep) - M_-(z)=\int_0^1 \frac{\d}{\d s}M_-(z+s(-y_0\pm i\ep))\d s = (-y_0\pm i\ep)\int_0^1 M_-'(z+s(-y_0\pm i\ep))\d s.
\end{equation*}
From Lemma \ref{asymptotic expansion M}, $M_-'(\zeta)\lesssim \zeta^{-\frac12-\mu} m^{\frac12-\mu}$, and since $2y_0\leq z$, we have that $s|y_0\pm i\ep|\leq |z+s(-y_0\pm i\ep)|$, for all $s\in(0,1)$. Thus, 
\begin{equation*}
\begin{aligned}
|M_-(z-y_0\pm i\ep) - M_-(z)|&\lesssim m^{\frac12-\mu}|y_0\pm i\ep|\int_0^1\frac{\d s}{|z+s(-y_0\pm i\ep)|^{\frac12 + \mu}} \\
&\lesssim m^{\frac12-\mu}|y_0\pm i\ep|\int_0^1\frac{\d s}{|s(y_0\pm i\ep)|^{\frac12 + \mu}} \\
&\lesssim m^{\frac12-\mu}|y_0\pm i\ep|^{\frac12-\mu}.
\end{aligned}
\end{equation*}

\bullpar{Case 2} 
For $z\leq 2y_0$, we directly estimate using Lemma \ref{asymptotic expansion M}, that is,
\begin{equation*}
\begin{aligned}
|M_-(z-y_0\pm i\ep) - M_-(z)|\leq |M_-(z-y_0\pm i\ep)| + |M_-(z)|&\lesssim m^{\frac12 -\mu} \l( |z-y_0\pm i\ep|^{\frac12-\mu} + |z|^{\frac{1}{2}-\mu}\r) \\
&\lesssim m^{\frac12 -\mu}|y_0\pm i\ep|^{\frac12-\mu}.
\end{aligned}
\end{equation*}
\end{proof}

From this localised estimates, we are able to obtain bounds on $f_{m,\ep}(z,y_0)$ for $mz\geq 6\b$. For this, we first deduce useful estimates on $\phi_{u,m}^\pm(z)$.
\begin{lemma}\label{phi up bounds}
The function $\phi_{u,m}(z)=M_+(1)M_-(z) - M_-(1)M_+(z)$ satisfies 
\begin{equation*}
\D_m \phi_{u,m}^\pm(z) + \b^2 \frac{\phi_{u,m}^\pm(z)}{z^2}=0, \qquad \phi_{u,m}^\pm(1)=0.
\end{equation*}
For $J_6=\lbrace z\in[0,1]: mz\leq 6\b \rbrace$ and $J_6^c=[0,1]\setminus J_6$, it is such that
\begin{equation*}
\Vert \phi_{u,m}\Vert _{L^\infty(J_6)}\lesssim m^{\frac12-\mu}|z|^{\frac{1}{2}-\mu}|M_+(1)|, \quad \Vert \phi_{u,m}\Vert_{L^2_y(J_6)}\lesssim m^{-\frac12}|M_+(1)|
\end{equation*}
and
\begin{equation*}
\Vert \p_z\phi_{u,m}\Vert_{L^2_y(J_6^c)} + m\Vert \phi_{u,m}\Vert_{L^2_y(J_6^c)} \lesssim m^{\frac12}|M_+(1)|.
\end{equation*}
\end{lemma}

\begin{proof}
The statements for $\Vert \phi_{u,m}\Vert _{L^\infty(J_6)}$ and $\Vert \phi_{u,m}\Vert_{L^2_y(J_6)}$ follow from the asymptotic expansions for small argument given by Lemma \ref{asymptotic expansion M}. The integral estimates follow from the $\Vert \phi_{u,m}\Vert_{L^2_y(J_6)}$ bounds using Lemma \ref{lemma:entangle ineq}.
\end{proof}

The following proposition obtains $L^2$ bounds on $f_{m,\ep}^\pm(\cdot,y_0)$ from the localized bounds of Proposition \ref{f bounds small z} and the above lemma.
\begin{proposition}\label{f bounds large z}
We have that
\begin{equation*}
\Vert f_{m,\ep}^\pm(\cdot, y_0)\Vert_{L^2}\lesssim m^{-\frac12}(m|y_0-i\ep|)^{\frac12-\mu}|M_+(1-y_0+i\ep)|.
\end{equation*}
\end{proposition}
\begin{proof}
It is straightforward to see that $f_{m,\ep}^\pm(z,y_0)$ solves 
\begin{equation*}
\D_m f_{m,\ep}^\pm + \b^2 \frac{f_{m,\ep}^\pm}{(z-y_0\pm i\ep)^2}=\b^2(-y_0\pm i\ep)\l(\frac{2}{z(z-y_0\pm i\ep)^2} + \frac{-y_0\pm i\ep}{z^2(z-y_0\pm i\ep)^2}\r)\phi_{u,m} 
\end{equation*}
and $f_{m,\ep}^\pm(1,y_0)=0$. Hence, using the same strategy from Lemma \ref{Bounds G Small z Large y}, we have that
\begin{equation}\label{energy ineq f}
\begin{aligned}
\frac{1}{2}\int_{\frac{6\b}{m}}^1 |\p_z f_{m,\ep}^\pm|^2 + m^2 |f_{m,\ep}^\pm|^2\d z &\leq \frac{m^2}{\b^2}\int_{\frac{5\b}{m}}^{\frac{6\b}{m}}|f^\pm_{m,\ep}|^2\d z \\
&+\b^2|-y_0\pm i\ep|\int_{\frac{5\b}{m}}^{1}\l(\frac{2}{z} + \frac{|-y_0\pm i\ep|}{z^2}\r)\frac{|\phi_{u,m}(z)|f_{m,\ep}^\pm(z)|}{|z-y_0\pm i\ep|^2}\d z.
\end{aligned}
\end{equation}
Now, from Proposition \ref{f bounds small z}, we have
\begin{equation*}
\frac{m^2}{\b^2}\int_{\frac{5\b}{m}}^{\frac{6\b}{m}}|f^\pm_{m,\ep}|^2\d z \lesssim \frac{m}{\b} \l(m|y_0-i\ep|\r)^{1-2\mu}|M_+(1-y_0+i\ep)|^2,
\end{equation*}
while we write
\begin{equation*}
\b^2|-y_0\pm i\ep|^2\int_{\frac{5\b}{m}}^{1}\frac{|\phi_{u,m}(z)|f_{m,\ep}^\pm(z)|}{z^2|z-y_0\pm i\ep|^2}\d z = \b^2|-y_0\pm i\ep|^2\l( \int_{\frac{5\b}{m}}^{\frac{6\b}{m}} + \int_{\frac{6\b}{m}}^1 \r)\frac{|\phi_{u,m}(z)|f_{m,\ep}^\pm(z)|}{z^2|z-y_0\pm i\ep|^2}\d z.
\end{equation*}
For example, with the bounds of Proposition \ref{f bounds small z} and Lemma \ref{phi up bounds}, and the fact that $z\geq \frac{5\b}{m}$ and $y_0\leq \frac{3\b}{m}$, we have $|z-y_0\pm i\ep|^{-2}\lesssim m^2$ and
\begin{equation*}
\begin{aligned}
\b^2|-y_0\pm i\ep|^2 \int_{\frac{5\b}{m}}^{\frac{6\b}{m}}\frac{|\phi_{u,m}(z)|f_{m,\ep}^\pm(z)|}{z^2|z-y_0\pm i\ep|^2}\d z 
&\lesssim {m^2}|y_0-i\ep|^2\Vert f_{m,\ep}^\pm\Vert_{L^\infty(J)}|M_+(1)|^2 \int_{y_2}^{y_2+\frac{\b}{m}}m^{\frac12-\mu}z^{-\frac32-\mu}   \d z \\
&\lesssim m^{\frac72 -\mu}|y_0\pm i\ep|^{\frac52-\mu}|M_+(1)||M_+(1-y_0\pm i\ep)|.
\end{aligned}
\end{equation*}
On the other hand, Young's inequality and Lemma \ref{phi up bounds} gives
\begin{equation*}
\begin{aligned}
\b^2|-y_0\pm i\ep|^2\int_{\frac{6\b}{m}}^1\frac{|\phi_{u,m}(z)|f_{m,\ep}^\pm(z)|}{z^2|z-y_0\pm i\ep|^2}\d z &\leq \frac{m^2}{8}\int_{y_2+\frac{\b}{m}}^1|f_{m,\ep}^\pm(z)|^2 \d z + C m^5|y_0-i\ep|^4 {|M_+(1)|^2},
\end{aligned}
\end{equation*}
for some $C>0$ large enough. Similarly, we bound
\begin{equation*}
\begin{aligned}
\b^2|-y_0\pm i\ep|&\int_{\frac{5\b}{m}}^{\frac{6\b}{m}}\frac{|\phi_{u,m}(z)|f_{m,\ep}^\pm(z)|}{z|z-y_0\pm i\ep|^2}\d z \lesssim m^{\frac52 -\mu}|y_0\pm i\ep|^{\frac32-\mu}|M_+(1)||M_+(1-y_0\pm i\ep)|
\end{aligned}
\end{equation*}
and
\begin{equation*}
\begin{aligned}
\b^2|-y_0\pm i\ep|\int_{\frac{6\b}{m}}^1\frac{|\phi_{u,m}(z)|f_{m,\ep}^\pm(z)|}{z|z-y_0\pm i\ep|^2}\d z &\leq \frac{m^2}{8}\int_{\frac{6\b}{m}}^1|f_{m,\ep}^\pm(z)|^2 \d z + C m^3|y_0-i\ep|^2 {|M_+(1)|^2},
\end{aligned}
\end{equation*}
for some $C>0$ large enough. Hence, we absorb the potential term on the left hand side of \eqref{energy ineq f} and conclude that
\begin{equation*}
\frac{1}{4}\int_{y_2+\frac{\b}{m}}^1 |\p_z f_{m,\ep}^\pm|^2 + m^2 |f_{m,\ep}^\pm|^2\d z \lesssim m(m|y_0-i\ep|)^{1-2\mu}|M_+(1-y_0+i\ep)|^2.
\end{equation*}
and the lemma follows.
\end{proof}

\subsubsection{Estimates on $f_{m,\ep}^\pm$ for $\b^2=1/4$}
From the explicit formulas \eqref{eq:homoupspecial} and \eqref{eq: phium special}, we now have
\begin{equation*}
\begin{aligned}
f_{m,\ep}^\pm(z,y_0)&=W_0(1-y_0\pm i\ep)M_0(z-y_0\pm i\ep) - W_0(1)M_0(z) \\
&\quad - M_0(1-y_0\pm i\ep)W_0(z-y_0\pm i\ep) + M_0(1)W_0(z)
\end{aligned}
\end{equation*}
from which we obtain the following result.
\begin{proposition}\label{special f bounds small z}
Let  $z,y_0\in[0,1]$ such that $my_0\leq 3\b$ and $mz\leq 6\b$. Let $0\leq \ep\leq \min\l( \frac{\b}{m},\frac{1}{2m}\r)$. Then,
\begin{equation*}
|f_{m,\ep}(z,y_0)|\lesssim (m|y_0\pm i\ep|)^\frac12 \l( 1 + \big| \log \l( my_0 \r) \big| \r)|M_0(1-y_0\pm i\ep)|.
\end{equation*}
In particular, $\Vert f_{m,\ep}\Vert_{L^2_y(J)}\lesssim |y_0\pm i\ep|^{\frac12} \l( 1 + \big| \log \l(m|y_0\pm i\ep|\r) \big| \r)|M_0(1-y_0\pm i\ep)|$.
\end{proposition}

\begin{proof}
We write 
\begin{equation*}
\begin{aligned}
W_0(1-y_0\pm i\ep)M_0(z-y_0\pm i\ep) - W_0(1)M_0(z)& =W_0(1-y_0\pm i\ep)\Big(M_0(z-y_0\pm i\ep) - M_0(z)\Big) \\
&\quad + M_0(z)\Big( W_0(1-y_0\pm i\ep) - W_0(1)\Big)
\end{aligned}
\end{equation*}
We shall now estimate the differences involving the Whittaker function $W_0$, the estimates for the differences involving $M_0$ follow similarly as for the case $\b^2\neq 1/4$ and they are 
\begin{equation*}
\l|M_0(z-y_0\pm i\ep) - M_0(z)\r| \lesssim m^\frac12|y_0\pm i\ep|^\frac12, \quad \l|M_0(1-y_0\pm i\ep) - M_0(1)\r| \lesssim m|y_0\pm i\ep||M_0(1-y_0\pm i\ep)|.
\end{equation*}
Firstly, we estimate
\begin{equation*}
W_0(1-y_0\pm i\ep) - W_0(1)=\int_0^1 \frac{\d}{\d s}W_0(1+s(-y_0\pm i\ep))\d s = 2m(-y_0\pm i\ep)\int_0^1 W_0'(1+s(-y_0\pm i\ep))\d s
\end{equation*}
and we divide our argument as follows. Let $N_{\mu,0}$ be given as in Lemma \ref{growth bounds real M large argument}.
 
For $m\leq N_{\mu,0}$, we use Lemma \ref{asymptotic expansion W} and the fact that $y_0\leq \frac12$ to bound
\begin{equation*}
\begin{aligned}
|W_0(1-y_0\pm i\ep) - W_0(1)|&\lesssim m^{\frac12}|y_0\pm i\ep|\int_0^1\frac{1+ \big|\log m|1+s(-y_0\pm i\ep)|\big|}{|1+s(-y_0\pm i\ep)|^{\frac12 }} \d s \\
&\lesssim m^{\frac12}|y_0\pm i\ep| \\
&\lesssim |y_0\pm i\ep||M_0(1-y_0 \pm i\ep)|.
\end{aligned}
\end{equation*}
In the last inequality, we have used Lemma \ref{lower bounds M}, \ref{Comparison bounds real M order one argument} and \ref{growth bounds real M large argument}.

For $m\geq N_{\mu,0}$, we claim that
\begin{equation*}
\l| \frac{W_0'(1+s(-y_0\pm i\ep))}{M_0(1-y_0\pm i\ep)}   \r|\lesssim 1.
\end{equation*}
Indeed, this follows from 
\begin{equation*}
\l| \frac{W_0'(1+s(-y_0\pm i\ep))}{M_0(1-y_0\pm i\ep)}   \r| = \l| \frac{W_0'(1+s(-y_0\pm i\ep))}{W_0(1+s(-y_0\pm i\ep))}\r| \l|\frac{W_0(1+s(-y_0\pm i\ep)) }{W_0(1-y_0\pm i\ep)}   \r|\l|\frac{W_0(1-y_0\pm i\ep) }{M_0(1-y_0\pm i\ep)}   \r|
\end{equation*}
and the corresponding bounds from Lemma \ref{growth bounds M large argument} since $2m(1-y_0)\geq m\geq N_{\mu,0}$. Hence, we have
\begin{equation*}
\begin{aligned}
|W_0(1-y_0\pm i\ep) - W_0(1)|&\leq 2m|y_0\pm i\ep|\int_0^1 \l|W_0'(1+s(-y_0\pm i\ep))\r| \d s \\
&\lesssim m|y_0\pm i\ep|\l| M_0(1-y_0\pm i\ep) \r|. \\
\end{aligned}
\end{equation*}
We next turn our attention to the bounds for $W_0(z-y_0\pm i\ep) - W_0(z)$. As before, we consider two cases.

\bullpar{Case 1} 
For $2y_0\leq z$ we estimate
\begin{equation*}
W_0(z-y_0\pm i\ep) - W_0(z)=\int_0^1 \frac{\d}{\d s}W_0(z+s(-y_0\pm i\ep))\d s = (-y_0\pm i\ep)\int_0^1 W_0'(z+s(-y_0\pm i\ep))\d s.
\end{equation*}
From Lemma \ref{asymptotic expansion W}, $W_0'(\zeta)\lesssim m^\frac12\zeta^{-\frac12}\l( 1+ \big| \log \l(m|\zeta|\r)\big|\r)$, and since $2y_0\leq z$, we have that $s|y_0\pm i\ep|\leq |z+s(-y_0\pm i\ep)|$, for all $s\in(0,1)$. Thus, 
\begin{equation*}
\begin{aligned}
|W_0(z-y_0\pm i\ep) - W_0(z)|&\lesssim m^{\frac12}|y_0\pm i\ep|\int_0^1\frac{1+\big| \log \l(m|z+s(-y_0\pm i\ep)|\r)\big|}{|z+s(-y_0\pm i\ep)|^{\frac12}}\d s \\
&\lesssim m^{\frac12}|y_0\pm i\ep|^\frac12\int_0^1\frac{1+\big| \log \l(ms|y_0\pm i\ep|\r)\big|}{s^{\frac12}}\d s \\
&\lesssim (m|y_0\pm i\ep|)^{\frac12}\l( 1 + \big| \log \l(m|y_0\pm i\ep|\r)\big|\r).
\end{aligned}
\end{equation*}
\bullpar{Case 2} 
For $z\leq 2y_0$, we directly estimate using Lemma \ref{asymptotic expansion W}, that is,
\begin{align*}
\begin{aligned}
|M_-(z-y_0\pm i\ep) - M_-(z)|&\leq |M_-(z-y_0\pm i\ep)| + |M_-(z)| \\
&\lesssim m^\frac12|z-y_0\pm i\ep|^\frac12 \l( 1 + \big| \log \l(m|z-y_0\pm i\ep|\r) \big|\r) \\
&\quad + m^\frac12z^\frac12 \l( 1 + \big| \log \l(mz\r) \big|\r) \\
&\lesssim (m|y_0\pm i\ep|)^\frac12 \l( 1 + \big| \log \l(m|y_0\pm i\ep|\r) \big| \r).
\end{aligned}
\end{align*}
\end{proof}

From this localised estimates, we are able to obtain bounds on $f_{m,\ep}(z,y_0)$ for $mz\geq 6\b$. For this, we first deduce useful estimates on $\phi_{u,m}^\pm(z)$.
\begin{lemma}\label{special phi up bounds}
The function $\phi_{u,m}(z)=W_0(1)M_0(z) - M_0(1)W_0(z)$ satisfies 
\begin{equation*}
\D_m \phi_{u,m}^\pm(z) + \b^2 \frac{\phi_{u,m}^\pm(z)}{z^2}=0, \qquad \phi_{u,m}^\pm(1)=0.
\end{equation*}
For $J_6=\lbrace z\in[0,1]: mz\leq 6\b \rbrace$ and $J_6^c=[0,1]\setminus J_6$, it is such that
\begin{equation*}
\Vert \phi_{u,m}\Vert _{L^\infty(J)}\lesssim (mz)^\frac12 \l( 1 + \big| \log \l(mz\r) \big| \r)    |M_0(1)|, \quad \Vert \phi_{u,m}\Vert_{L^2_y(J)}\lesssim m^{-\frac12}|M_0(1)|
\end{equation*}
and
\begin{equation*}
\Vert \p_z\phi_{u,m}\Vert_{L^2_y(J^c)} + m\Vert \phi_{u,m}\Vert_{L^2_y(J^c)} \lesssim m^{\frac12}|M_+(1)|
\end{equation*}
\end{lemma}

\begin{proof}
The statement for $\Vert \phi_{u,m}\Vert _{L^\infty(J_6)}$ follows from the asymptotic expansions for small argument given by Lemma \ref{asymptotic expansion M}. For the integral estimates estimate, note that the change of variables $u=mz$ provides
\begin{align*}
\Vert \phi_{u,m}\Vert _{L^2_y(J_6)}^2 \lesssim \int_0^{\frac{6\b}{m}}(mz)\l( 1 + \big| \log \l(mz\r) \big| \r)^2  |M_0(1)|^2 \d z &= \frac{|M_0(1)|^2}{m}\int_0^{6\b} \eta\l( 1 + | \log \l(\eta\r) | \r)^2 \d \eta \\
&\lesssim \frac{|M_0(1)|^2}{m}.
\end{align*}
The result follows using Lemma \ref{lemma:entangle ineq}.
\end{proof}

The following proposition obtains $L^2(0,1)$ bounds on $f_{m,\ep}^\pm(\cdot,y_0)$ from the localized bounds of Proposition \ref{special f bounds small z} and the above Lemma. We omit its proof due to its similarity to the one for Proposition \ref{f bounds large z}.

\begin{proposition}\label{special f bounds large z}
We have that
\begin{equation*}
\Vert f_{m,\ep}^\pm(\cdot, y_0)\Vert_{L^2(0,1)}\lesssim m^{-\frac12}(m|y_0-i\ep|)^{\frac12}\l( 1 + \big| \log \l(m|y_0\pm i\ep|\r) \big|\r)   |M_0(1-y_0\pm i\ep)|
\end{equation*}
\end{proposition}

We are now able to compare $\Vert f_{m,\ep}^\pm(\cdot, y_0)\Vert_{L^2(0,1)}$ and the Wronskian $|\W_{m,\ep}^\pm(y_0)|$.

\begin{lemma}\label{lemma special f Wronskian comparison}
Let $y_0\in [0,1]$ such that $my_0\leq 3\b$. There exists $\ep_0>0$ such that
\begin{equation*}
\frac{\Vert f_{m,\ep}^\pm(\cdot, y_0)\Vert_{L^2(0,1)}}{|\W_{m,\ep}^\pm(y_0)|}\lesssim m^{-\frac32}.
\end{equation*}
\end{lemma}

\begin{proof}
Let $N_0>0$ be given by Lemma \ref{growth bounds special M large argument}, $\delta_1>0$ be given by Lemma \ref{Comparison bounds special M small argument} and $\delta_2>0$ be given by Lemma \ref{lemma: lower asymptotic bounds special W}. From Lemma \ref{lower bound special Wronskian}, there holds the following,

\bullpar{Case 1} For $m\leq N_0$ and $2m|y_0\pm i\ep|\leq \delta_1$, we have 
\begin{equation*}
|\W_{m,\ep}^\pm(y_0)|\gtrsim m|M_0(1-y_0\pm i\ep)||W_0(y_0\mp i\ep)|.
\end{equation*}
where further 
\begin{equation*}
\l| \frac{M_0(y_0\pm i\ep)}{W_0(y_0\pm i\ep)}\r| \leq \frac{1}{2\sqrt{\pi}}.
\end{equation*} 
Now, from Lemma \ref{lemma: lower asymptotic bounds special W}, if $\delta_1\leq \delta_2$, then,
\begin{equation*}
|2m(y_0\pm i\ep)|^\frac12\l( 1 + \big| \log \l(2m|y_0\pm i\ep|\r) \big| \r) \lesssim |W_0(y_0\pm i\ep)|,
\end{equation*}
and the conclusion follows. On the other hand, for $\delta_2\leq 2m|y_0\pm i\ep| \leq \delta_1$, we have that 
\begin{equation*}
m\frac{\Vert f_{m,\ep}^\pm(\cdot, y_0)\Vert_{L^2}}{|\W_{m,\ep}^\pm(y_0)|}\lesssim m^{-\frac12}\frac{(m|y_0\pm i\ep|)^\frac12}{M_0(y_0\pm i\ep)}\lesssim m^{-\frac12},
\end{equation*}
due to Lemma \ref{lemma:lower asymptotic bounds special M}.

\bullpar{Case 2} For $m\leq N_0$ and $\delta_1\leq 2my_0\leq N_0$, we have now 
\begin{equation*}
|\W_{m,\ep}^\pm(y_0)|\gtrsim m|M_0(1-y_0\pm i\ep)||M_0(y_0\mp i\ep)|,
\end{equation*}
the conclusion follows using Lemma \ref{lemma:lower asymptotic bounds special M} and the fact that $\l( 1 + \big| \log |2m(y_0\pm i\ep)| \big| \r)\lesssim 1$.

\bullpar{Case 3} For $m\geq N_0$, and $2m|y_0\pm i\ep|\leq \delta_1$, we have 
\begin{equation*}
|\W_{m,\ep}^\pm(y_0)|\gtrsim m|M_0(1-y_0\pm i\ep)||W_0(y_0\mp i\ep)|.
\end{equation*}
and also 
\begin{equation*}
\l| \frac{M_0(y_0\pm i\ep)}{W_0(y_0\pm i\ep)}\r| \leq \frac{1}{2\sqrt{\pi}},
\end{equation*} 
we proceed as in Case 1, we omit the details.

\bullpar{Case 4} For $m\geq N_0$ and $\delta_1 \leq 2my_0\leq N_0$, we have  
\begin{equation*}
|\W_{m,\ep}^\pm(y_0)|\gtrsim m|M_0(1-y_0\pm i\ep)||M_0(y_0\mp i\ep)|,
\end{equation*}
we proceed as in Case 2, we omit the details.
\end{proof}

We are now in position to prove Proposition \ref{prop: dyVarphi0 estimates}.

\begin{proof}[Proof of Proposition \ref{prop: dyVarphi0 estimates}]
For $my_0\geq 3\b$ we appeal to Proposition \ref{L2 bounds inhom TG solution} to obtain the desired bound. On the other hand, for $my_0\leq 3\b$ let us recall that we can write
\begin{equation*}
\begin{aligned}
\p_y\varphi_{m,\ep}^\pm(0,y_0)&=  \int_0^1 \p_y\G_{m,\ep}^\pm(0,y_0,z)\l(F_{m}(z)+ \frac{y_0\mp i\ep}{\b^2}\D_m\o_m^0(z)\r) \d z
\end{aligned}
\end{equation*}
For $\b^2\neq1/4$, it is straightforward to see from Proposition \ref{L2 bounds inhom TG solution} that 
\begin{equation*}
\l| \frac{y_0\mp i\ep}{\b^2}\int_0^1 \p_y\G_{m,\ep}^\pm(0,y_0,z)\D_m\o_m^0(z) \d z\r|\lesssim \frac{1}{m^{1+\mu}}|y_0\pm i\ep|^{\frac12-\mu}\Vert \o_m^0\Vert_{H^2_y},
\end{equation*}
while, thanks to Proposition \ref{f bounds large z}, the lower bounds on the Wronskian from Proposition \ref{bounds real G small y small z} and Lemma \ref{Comparison bounds real M small argument} we bound
\begin{equation*}
\begin{aligned}
\l|\int_0^1 \p_y\G_{m,\ep}^\pm(0,y_0,z) F_{m}^\pm(z,0) \d z\r| &= \l|4(\mu +i\nu)m \int_0^1 \frac{f_{m,\ep}^\pm(z,y_0)}{\W_{m,\ep}^\pm(y_0)}F_{m}^\pm(z,0)\d z\r| \\
&\lesssim  m \frac{\Vert f_{m,\ep}^\pm(\cdot,y_0)\Vert_{L^2_z}}{|\W_{m,\ep}^\pm(y_0)|}\Vert F_{m}^\pm(z,0)\Vert_{L^2_z} \\
&\lesssim m^{-\frac12}\Vert F_m \Vert_{L^2}.
\end{aligned}
\end{equation*}
Similarly, for $\b^2=1/4$, using again Proposition \ref{L2 bounds inhom TG solution},
\begin{align*}
\l| \frac{y_0\mp i\ep}{\b^2}\int_0^1 \p_y\G_{m,\ep}^\pm(0,y_0,z)\D_m\o_m^0(z) \d z\r|&\lesssim \frac{1}{m}|y_0\pm i\ep|^{\frac12}\l( 1 + \big| \log \l(m|y_0\pm i\ep|\r) \big| \r)  \Vert \o_m^0\Vert_{H^2_y}\\
&\lesssim m^{-\frac32}\Vert \o_m^0\Vert_{H^2_y},
\end{align*}
while, thanks to Lemma \ref{lemma special f Wronskian comparison} we have
\begin{equation*}
\begin{aligned}
\l|\int_0^1 \p_y\G_{m,\ep}^\pm(0,y_0,z) F_{m}^\pm(z,0) \d z\r| &= \l|4(\mu +i\nu)m \int_0^1 \frac{f_{m,\ep}^\pm(z,y_0)}{\W_{m,\ep}^\pm(y_0)}F_{m}^\pm(z,0)\d z\r| \\
&\lesssim  m \frac{\Vert f_{m,\ep}^\pm(\cdot,y_0)\Vert_{L^2_z}}{|\W_{m,\ep}^\pm(y_0)|}\Vert F_{m}^\pm(z,0)\Vert_{L^2_z} \\
&\lesssim m^{-\frac12}\  \Vert F_m \Vert_{L^2}.
\end{aligned}
\end{equation*}
With this, the proof is finished.
\end{proof}

We next provide pointwise localized bounds on $\B_{m,\ep}^\pm(y,y_0,0)$.

\begin{proposition}\label{prop:boundary term estimates}
Let $\b^2\neq1/4$ and $0\leq\ep\leq\ep_0$. Let $y,y_0\in[0,1]$ such that $m|y-y_0|\leq 3\b$. Then,
\begin{itemize}
\item If $my_0\leq 3\b$, we have
\begin{equation*}
\l|\B_{m,\ep}^\pm(y,y_0,0)\r|\lesssim m^{-\frac12}y_0^{-\frac12+\mu}|y-y_0\pm i\ep|^{\frac12-\mu}{Q_{0,m}}.
\end{equation*}
\item If $my_0\geq 3\b$, we have
\begin{equation*}
\l| \B_{m,\ep}^\pm(y,y_0,0)\r|\lesssim (m|y-y_0\pm i\ep|)^{\frac12-\mu}{Q_{0,m}}.
\end{equation*}
\end{itemize}
\end{proposition}

\begin{proposition}\label{prop: special boundary term estimates}
Let $\b^2=1/4$ and $0\leq\ep\leq\ep_0$. Let $y,y_0\in[0,1]$ such that $m|y-y_0|\leq 3\b$. Then,
\begin{itemize}
\item If $my_0\leq 3\b$, we have
\begin{equation*}
\l|\B_{m,\ep}^\pm(y,y_0,0)\r|\lesssim m^{-\frac12}y_0^{-\frac12}|y-y_0\pm i\ep|^{\frac12}\l( 1 + \big| \log \l(m|y-y_0\pm i\ep|\r) \big| \r){Q_{0,m}}.
\end{equation*}
\item If $my_0\geq 3\b$, we have
\begin{equation*}
\l| \B_{m,\ep}^\pm(y,y_0,0)\r|\lesssim (m|y-y_0\pm i\ep|)^{\frac12} \l( 1 + \big| \log \l(m|y-y_0\pm i\ep|\r) \big| \r) {Q_{0,m}}.
\end{equation*}
\end{itemize}
\end{proposition}

With the above pointwise bounds, one deduces the following integral estimates for all $\b^2>0$.
\begin{corollary}\label{cor:global boundary term estimates}
Let $y_0\in[0,1]$. Then,
\begin{itemize}
\item If $my_0\leq 3\b$, we have
\begin{equation*}
\Vert \B_{m,\ep}^\pm(\cdot,y_0,0)\Vert_{L^2_y(J_2^c\cap J_3)} \lesssim m^{-\frac32}y_0^{-\frac12}{Q_{0,m}}.
\end{equation*}
\item If $my_0\geq 3\b$, we have
\begin{equation*}
\Vert \B_{m,\ep}^\pm(\cdot,y_0,0)\Vert_{L^2_y(J_2^c\cap J_3)}\lesssim m^{-\frac12}{Q_{0,m}}.
\end{equation*}
\end{itemize}
\end{corollary}

The two propositions are a consequence of Propositions \ref{prop: pointwise dzG bounds}, \ref{prop: special pointwise dzG bounds}, the lower bounds from Lemma \ref{lemma:lower asymptotic bounds M}, \ref{lemma:lower asymptotic bounds special M}, \ref{lemma:lower asymptotic bounds real M} and the  pointwise estimates on $\p_y\varphi_{m,\ep}^\pm(0,y_0)$ from Proposition \ref{prop: dyVarphi0 estimates}.

\subsection{Boundary pointwise estimates on Green's Function's derivatives}
This subsection estimates derivatives of the Green's function $\G_{m,\ep}^\pm(y,y_0,z)$ evaluated at the boundary values $y,z\in \lbrace 0,1\rbrace$.
\begin{lemma}\label{lemma:pointwise bound dydzG00}
We have that for $my_0\geq 3\b$,
\begin{equation*}
|\p_y\p_z\G_{m,\ep}^{\pm}(0,y_0,0)|\lesssim m
\end{equation*}
while for $my_0\leq 3\b$, 
\begin{equation*}
|\p_y\p_z\G_{m,\ep}^{\pm}(0,y_0,0)|\lesssim \frac{1}{y_0}.
\end{equation*}
\end{lemma}
\begin{proof}
For $\b^2>1/4$, it follows from the proof of Proposition \ref{Bounds Small y Small z} and Lemma \ref{lemma:lower asymptotic bounds M}. For $\b^2=1/4$, it follows from Lemma \ref{lower bound special Wronskian}, Lemma \ref{lemma: lower asymptotic bounds special W} and Lemma \ref{lemma:lower asymptotic bounds special M}. For $b^2<1/4$, it follows from the proof of Proposition \ref{bounds real G small y small z} and Lemma \ref{lemma:lower asymptotic bounds real M}.
\end{proof}

\begin{lemma}\label{lemma:pointwise bound dydzG10}
For $m\geq 6\b$, we have
\begin{itemize}
\item If $my_0\leq 3\b$, then $m(1-y_0)\geq 3\b$ and 
\begin{equation*}
|\p_y\p_z\G_{m,\ep}^\pm(1,y_0,0)|\lesssim {m^{\frac12+\mu}}{y_0^{-\frac12+\mu}}.
\end{equation*}
\item If $m(1-y_0)\leq 3\b$, then $my_0\geq 3\b$ and 
\begin{equation*}
|\p_y\p_z\G_{m,\ep}^\pm(1,y_0,0)|\lesssim {m^{\frac12+\mu}}{(1-y_0)^{-\frac12+\mu}}.
\end{equation*}
\item If $my_0\geq 3\b$ and $m(1-y_0)\geq 3\b$, then
\begin{equation*}
|\p_y\p_z\G_{m,\ep}^\pm(1,y_0,0)|\lesssim m.
\end{equation*}
\end{itemize}
On the other hand, for $m\leq 6\b$, we have that
\begin{itemize}
\item If $y_0\leq \frac12$, then 
\begin{equation*}
|\p_y\p_z\G_{m,\ep}^\pm(1,y_0,0)|\lesssim {y_0^{-\frac12+\mu}}.
\end{equation*}
\item If $1-y_0\leq \frac12$, then 
\begin{equation*}
|\p_y\p_z\G_{m,\ep}^\pm(1,y_0,0)|\lesssim {(1-y_0)^{-\frac12+\mu}}.
\end{equation*}
\end{itemize}
\end{lemma}

\begin{proof}
It is straightforward 
the lower bounds on the Wronskian and Lemmas \ref{lemma:lower asymptotic bounds M}, \ref{lemma:lower asymptotic bounds real M}, \ref{lemma: lower asymptotic bounds special W}, for $\b^2>1/4$, $\b^2=1/4$ and $\b^2<1/4$, respectively.
\end{proof}

\begin{lemma}\label{lemma:pointwise bound dydzG01}
The same bounds as in Lemma \ref{lemma:pointwise bound dydzG10} hold for $|\p_y\p_z\G_{m,\ep}^\pm(0,y_0,1)|$.
\end{lemma}

\begin{lemma}\label{lemma:pointwise bound dydzG11}
We have that for $m(1-y_0)\geq 3\b$,
\begin{equation*}
|\p_y\p_z\G_{m,\ep}^{\pm}(1,y_0,1)|\lesssim m
\end{equation*}
while for $m(1-y_0)\leq 3\b$, 
\begin{equation*}
|\p_y\p_z\G_{m,\ep}^{\pm}(1,y_0,1)|\lesssim \frac{1}{1-y_0}.
\end{equation*}
\end{lemma}

\subsection{Estimates for second order boundary terms}
In what follows, we shall consider only the case $m\geq 6\b$, since the setting $m\leq 6\b$ is analogous and easier. With the pointwise derivatives bounds obtained in the four previous lemmas, we are now in position to estimate 
\begin{equation*}
\begin{aligned}
\l( \p_y + \p_{y_0}\r)^2\varphi_{m,\ep}^\pm(y,y_0)&=-\frac{2}{\b^2}\p_y\o_m^0(y) -F_{m,\ep}^\pm(y,y_0) +2\p_y\Big. \B_{m,\ep}^\pm(y,y_0,z)\Big]_{z=0}^{z=1} \\
&\quad+2 \int_0^1\p_y\G_{m,\ep}^\pm(y,y_0,z)\l(\p_zF_{m,\ep}^\pm(z,y_0)+\p_{y_0}F_{m,\ep}^\pm(z,y_0)\r)\d z
\end{aligned}
\end{equation*} 
for both $y=0$ and $y=1$. For simplicity we only discuss the case $y=0$; the results and proofs are the same for the case $y=1$.

\begin{proposition}\label{prop: tildetildeB bounds}
Let $m\geq6\b$ and $y_0\in[0,1]$. Then, we have that
\begin{itemize}
\item For $my_0\leq 3\b$ and $\b^2\neq1/4$,
\begin{equation*}
|\l( \p_y + \p_{y_0}\r)^2\varphi_{m,\ep}^\pm(y,y_0)|\lesssim \l(1+\frac{1}{m^{1+\mu}}y_0^{-\frac12-\mu} +m^{-\frac12}y_0^{-1} + m^\frac12y_0^{-\frac12}\r)Q_{1,m}
\end{equation*}
\item For $my_0\leq 3\b$ and $\b^2=1/4$,
\begin{equation*}
|\l( \p_y + \p_{y_0}\r)^2\varphi_{m,\ep}^\pm(y,y_0)|\lesssim \l(1+\frac{1}{m}y_0^{-\frac12}\l( 1+ \big| \log \l(my_0\r) \big| \r) +m^{-\frac12}y_0^{-1} + m^\frac12y_0^{-\frac12}\r)Q_{1,m}
\end{equation*}
\item For $my_0\geq 3\b$ and $m(1-y_0)\leq 3\b$,
\begin{equation*}
|\l( \p_y + \p_{y_0}\r)^2\varphi_{m,\ep}^\pm(y,y_0)|\lesssim \l( m + (1-y_0)^{-\frac12}\r)Q_{1,m}, 
\end{equation*}
\item For $my_0\geq 3\b$ and $m(1-y_0)\geq 3\b$,
\begin{equation*}
|\l( \p_y + \p_{y_0}\r)^2\varphi_{m,\ep}^\pm(y,y_0)|\lesssim mQ_{1,m}.
\end{equation*}
\end{itemize}
\end{proposition}

\begin{proof}
For $y=0$, we can estimate 
\begin{equation*}
|\p_y\o_m^0(0)|+ |F_{m,\ep}^\pm(0,y_0)|\lesssim Q_{1,m}
\end{equation*}
thanks to the Sobolev embedding. On the other hand, from Proposition \ref{L2 bounds inhom TG solution}, for $my_0\leq 3\b$ and $\b^2\neq1/4$ we have 
\begin{equation*}
\l|\int_0^1\p_y\G_{m,\ep}^\pm(0,y_0,z)\l(\p_zF_{m,\ep}^\pm(z,y_0)+\p_{y_0}F_{m,\ep}^\pm(z,y_0)\r)\d z\r|\lesssim \frac{1}{m^{1+\mu}}|y_0\pm i\ep|^{-\frac12-\mu}Q_{1,m},
\end{equation*}
while for $my_0\leq 3\b$ and $\b^2=1/4$ we have 
\begin{align*}
\l|\int_0^1\p_y\G_{m,\ep}^\pm(0,y_0,z) \l(\p_z+\p_{y_0}\r)F_{m,\ep}^\pm(z,y_0)\d z\r| \lesssim \frac{1}{m}|y_0\pm i\ep|^{-\frac12} \l( 1 + \big| \log \l(m|y_0\pm i\ep|\r) \big| \r) Q_{1,m},
\end{align*}
whereas for $my_0\geq 3\b$, we have 
\begin{equation*}
\l|\int_0^1\p_y\G_{m,\ep}^\pm(0,y_0,z)\l(\p_z+\p_{y_0}\r)F_{m,\ep}^\pm(z,y_0)\d z\r|\lesssim Q_{1,m}.
\end{equation*}
Now, for the solid boundary terms $\p_y\Big. \B_{m,\ep}^\pm(y,y_0,z)\Big]_{z=0}^{z=1}$, we shall use Proposition \ref{prop: dyVarphi0 estimates} as well as Lemmas \ref{lemma:pointwise bound dydzG00}-\ref{lemma:pointwise bound dydzG11}. Indeed, for $\p_y\B_{m,\ep}^\pm(0,y_0,0)=\p_y\p_z\G_{m,\ep}^\pm(0,y_0,0)\p_y\varphi_{m,\ep}^\pm(0,y_0)$, Lemma \ref{lemma:pointwise bound dydzG00} provides
\begin{itemize}
\item For $my_0\leq 3\b$, we have that $|\p_y\B_{m,\ep}^\pm(0,y_0,0)|\lesssim m^{-\frac12}y_0^{-1}Q_{0,m}$.
\item For $my_0\geq 3\b$, we have that $|\p_y\B_{m,\ep}^\pm(0,y_0,0)|\lesssim mQ_{0,m}$.
\end{itemize}
Similarly, for $\p_y\B_{m,\ep}^\pm(0,y_0,1)=\p_y\p_z\G_{m,\ep}^\pm(0,y_0,1)\p_y\varphi_{m,\ep}^\pm(1,y_0)$, we have from Lemma \ref{lemma:pointwise bound dydzG01} that
\begin{itemize}
\item For $my_0\leq 3\b$, then $|\p_y\B_{m,\ep}^\pm(0,y_0,1)|\lesssim m^{\frac12+\mu}y_0^{-\frac12+\mu}Q_{0,m}$.
\item For $my_0\geq 3\b$, we further distinguish 
\begin{itemize}
\item For $m(1-y_0)\leq 3\b$, then $|\p_y\B_{m,\ep}^\pm(0,y_0,1)|\lesssim m^\mu(1-y_0)^{-\frac12+\mu}Q_{0,m}$.
\item For $m(1-y_0)\geq 3\b$, then $|\p_y\B_{m,\ep}^\pm(0,y_0,1)|\lesssim mQ_{0,m}$.
\end{itemize}
\end{itemize}
As a result, for $my_0\leq 3\b$ and $\b^2\neq1/4$ we have that
\begin{equation*}
\begin{aligned}
|\l( \p_y + \p_{y_0}\r)^2\varphi_{m,\ep}^\pm(y,y_0)|&\lesssim Q_{1,m}+\frac{1}{m^{1+\mu}}y_0^{-\frac12-\mu}Q_{1,m} + m^{-\frac12}y_0^{-1}Q_{0,m} + m^{\frac12+\mu}y_0^{-\frac12+\mu}Q_{0,m} \\
&\lesssim \l(1+\frac{1}{m^{1+\mu}}y_0^{-\frac12-\mu} +m^{-\frac12}y_0^{-1} + m^\frac12y_0^{-\frac12}\r)Q_{1,m}, 
\end{aligned}
\end{equation*}
while for $my_0\leq 3\b$ and $\b^2=1/4$ we have that
\begin{equation*}
\begin{aligned}
|\l( \p_y + \p_{y_0}\r)^2\varphi_{m,\ep}^\pm(y,y_0)|&\lesssim Q_{1,m}+\frac{1}{m}y_0^{-\frac12} \l( 1 + \big| \log \l(my_0\r) \big| \r) Q_{1,m} + m^{-\frac12}y_0^{-1}Q_{0,m} + m^{\frac12}y_0^{-\frac12}Q_{0,m} \\
&\lesssim \l(1+\frac{1}{m}y_0^{-\frac12} \l( 1 + \big| \log \l(my_0\r) \big| \r) +m^{-\frac12}y_0^{-1} + m^\frac12y_0^{-\frac12}\r)Q_{1,m}.
\end{aligned}
\end{equation*}
Similarly, for $my_0\geq 3\b$ and $m(1-y_0)\leq 3\b$ we conclude
\begin{equation*}
\begin{aligned}
|\l( \p_y + \p_{y_0}\r)^2\varphi_{m,\ep}^\pm(y,y_0)|&\lesssim Q_{1,m} + m Q_{0,m} + (1-y_0)^{-\frac12}Q_{0,m} \lesssim \l(m+(1-y_0)^{-\frac12}\r)Q_{1,m}, 
\end{aligned}
\end{equation*}
whereas for $my_0\geq 3\b$ and $m(1-y_0)\geq 3\b$ we obtain
\begin{equation*}
\begin{aligned}
|\l( \p_y + \p_{y_0}\r)^2\varphi_{m,\ep}^\pm(y,y_0)|&\lesssim Q_{1,m} + m Q_{0,m}  \lesssim mQ_{1,m}
\end{aligned}
\end{equation*}
and the proof is finished.
\end{proof}

We next present estimates for 
$$\widetilde{{\B_{m,\ep}^\pm}}(y,y_0,z)=\p_z\G_{m,\ep}^\pm(y,y_0,z)\l( \p_y + \p_{y_0}\r)^2\varphi_{m,\ep}^\pm(y,y_0),$$
at $z=0$. As before, we only obtain these bounds under the assumption that $m|y-y_0|\leq 3\b$. We state them for $z=0$; the result for $z=1$ is the same and thus we omit the details. The next two Propositions are a direct consequence of Propositions \ref{prop: pointwise dzG bounds}, \ref{prop: special pointwise dzG bounds}, \ref{prop: tildetildeB bounds}, as well as Lemma \ref{lemma:lower asymptotic bounds M} \ref{lemma: lower asymptotic bounds special W}, \ref{lemma:lower asymptotic bounds special M}, \ref{lemma:lower asymptotic bounds real M}, depending on $\b^2$.

\begin{proposition}\label{prop: tildeB0 estimates}
Let $\b^2\neq1/4$ and  $y,y_0\in[0,1]$ such that $m|y-y_0|\leq 3\b$. Then,
\begin{itemize}
\item For $my_0\leq 3\b$,
\begin{equation*}
|\widetilde{{\B_{m,\ep}^\pm}}(y,y_0,0)|\lesssim m^{-\frac12}y_0^{-\frac32}(m|y-y_0\pm i\ep|)^{\frac12-\mu}Q_{1,m}
\end{equation*}
\item For $my_0\geq 3\b$ and $m(1-y_0)\leq 3\b$,
\begin{equation*}
|\widetilde{{\B_{m,\ep}^\pm}}(y,y_0,0)|\lesssim \l( m + (1-y_0)^{-\frac12}\r)(m|y-y_0\pm i\ep|)^{\frac12-\mu}Q_{1,m}, 
\end{equation*}
\item For $my_0\geq 3\b$ and $m(1-y_0)\geq 3\b$,
\begin{equation*}
|\widetilde{{\B_{m,\ep}^\pm}}(y,y_0,0)|\lesssim m(m|y-y_0\pm i\ep|)^{\frac12-\mu}Q_{1,m}.
\end{equation*}
\end{itemize}
\end{proposition}

\begin{proposition}\label{prop: special tildeB0 estimates}
Let $\b^2=1/4$ and $y,y_0\in[0,1]$ such that $m|y-y_0|\leq 3\b$. Then,
\begin{itemize}
\item For $my_0\leq 3\b$,
\begin{equation*}
|\widetilde{{\B_{m,\ep}^\pm}}(y,y_0,0)|\lesssim m^{-\frac12}y_0^{-\frac32}(m|y-y_0\pm i\ep|)^{\frac12}\l( 1 + \big| \log \l(m|y-y_0\pm i\ep|\r)\big| \r)  Q_{1,m}
\end{equation*}
\item For $my_0\geq 3\b$ and $m(1-y_0)\leq 3\b$,
\begin{equation*}
|\widetilde{{\B_{m,\ep}^\pm}}(y,y_0,0)|\lesssim \l( m + (1-y_0)^{-\frac12}\r)(m|y-y_0\pm i\ep|)^{\frac12} \l( 1 + \big| \log \l(m|y-y_0\pm i\ep|\r)\big| \r) Q_{1,m}, 
\end{equation*}
\item For $my_0\geq 3\b$ and $m(1-y_0)\geq 3\b$,
\begin{equation*}
|\widetilde{{\B_{m,\ep}^\pm}}(y,y_0,0)|\lesssim m(m|y-y_0\pm i\ep|)^{\frac12} \l( 1 + \big| \log \l(m|y-y_0\pm i\ep|\r)\big| \r) Q_{1,m}.
\end{equation*}
\end{itemize}
\end{proposition}

We upgrade the above pointwise estimates to integral bounds for $y\in[0,1]$ such that $2\b\leq m|y-y_0|\leq 3\b$, which will be useful later on. 

\begin{corollary}\label{cor: integral tildeB0 estimates}
Let $0\leq \ep\leq \ep_0\leq \frac{\b}{m}$. Then,
\begin{itemize}
\item For $my_0\leq 3\b$,
\begin{equation*}
\Vert \widetilde{{\B_{m,\ep}^\pm}}(\cdot,y_0,0)\Vert_{L^2_y(J_2^c\cap J_3)}\lesssim m^{-1}y_0^{-\frac32}Q_{1,m}
\end{equation*}
\item For $my_0\geq 3\b$ and $m(1-y_0)\leq 3\b$,
\begin{equation*}
\Vert \widetilde{{\B_{m,\ep}^\pm}}(\cdot,y_0,0)\Vert_{L^2_y(J_2^c\cap J_3)}\lesssim m^{-\frac12}\l( m + (1-y_0)^{-\frac12}\r)Q_{1,m}, 
\end{equation*}
\item For $my_0\geq 3\b$ and $m(1-y_0)\geq 3\b$,
\begin{equation*}
\Vert \widetilde{{\B_{m,\ep}^\pm}}(\cdot,y_0,0)\Vert_{L^2_y(J_2^c\cap J_3)}\lesssim m^{\frac12}Q_{1,m}.
\end{equation*}
\end{itemize}
\end{corollary}

We finish the section with estimates for 
\begin{equation*}
\p_y\B_{m,\ep}^\pm(y,y_0,z)=\p_y\p_z\G_{m,\ep}^\pm(y,y_0,z)\p_z\varphi_{m,\ep}^\pm(z,y_0)
\end{equation*} 
for $z=0$ and $z=1$ under the localizing assumption that $m|y-y_0|\leq 3\b$. The next two results follows directly from Proposition \ref{prop: pointwise dzG bounds}, Proposition \ref{prop: special pointwise dzG bounds} and Proposition \ref{prop: dyVarphi0 estimates}.

\begin{proposition}\label{prop: dyB0 estimates}
Let $\b^2\neq1/4$. Let $y,y_0\in[0,1]$ such that $m|y-y_0|\leq 3\b$. Then,
\begin{itemize}
\item For $my_0\leq 3\b$, we have that
\begin{equation*}
|\p_y\B_{m,\ep}^\pm(y,y_0,0)|\lesssim m^{-\frac12}y_0^{-\frac12+\mu}|y-y_0\pm i\ep|^{-\frac12-\mu}Q_{0,m}.
\end{equation*}
\item For $my_0\geq 3\b$, we have that
\begin{equation*}
|\p_y\B_{m,\ep}^\pm(y,y_0,0)|\lesssim m^{\frac12-\mu}|y-y_0\pm i\ep|^{-\frac12-\mu}Q_{0,m}.
\end{equation*}
\end{itemize}
\end{proposition}

\begin{proposition}\label{prop: special dyB0 estimates}
Let $b^2=1/4$ and  $y,y_0\in[0,1]$ such that $m|y-y_0|\leq 3\b$. Then,
\begin{itemize}
\item For $my_0\leq 3\b$, we have that
\begin{equation*}
|\p_y\B_{m,\ep}^\pm(y,y_0,0)|\lesssim m^{-\frac12}y_0^{-\frac12}|y-y_0\pm i\ep|^{-\frac12} \l( 1 + \big| \log \l(m|y-y_0\pm i\ep|\r)\big| \r) Q_{0,m}.
\end{equation*}
\item For $my_0\geq 3\b$, we have that
\begin{equation*}
|\p_y\B_{m,\ep}^\pm(y,y_0,0)|\lesssim m^{\frac12}|y-y_0\pm i\ep|^{-\frac12} \l( 1 + \big| \log \l(m|y-y_0\pm i\ep|\r)\big| \r) Q_{0,m}.
\end{equation*}
\end{itemize}
\end{proposition}

Finally, we state the integral bounds that are deduced from the above estimates.

\begin{corollary}\label{cor: integral dyB0 estimates}
Let   $0\leq \ep\leq \ep_0$. Then, 
\begin{itemize}
\item For $my_0\leq 3\b$, we have that
\begin{equation*}
\Vert \p_y\B_{m,\ep}^\pm(\cdot,y_0,0) \Vert_{L^2_y(J_2^c\cap J_3)}\lesssim (my_0)^{-\frac12}Q_{0,m}.
\end{equation*}
\item For $my_0\geq 3\b$, we have that
\begin{equation*}
\Vert \p_y\B_{m,\ep}^\pm(\cdot,y_0,0) \Vert_{L^2_y(J_2^c\cap J_3)}\lesssim m^{\frac12}Q_{0,m}.
\end{equation*}
\end{itemize}
\end{corollary}

\section{Estimates for the Generalized Stream-functions}\label{sec: estimates generalized stream}
This section is devoted to obtaining estimates for the generalized stream-functions $\psi_{m,\ep}^\pm(y,y_0)$ and densities $\rho_{m,\ep}^{\pm}(y,y_0,z)$, as well as for some of their derivatives. Moreover, we define
\begin{equation*}
\widetilde{\psi_m}(y,y_0):=\lim_{\ep\rightarrow 0}\psi_{m,\ep}^-(y,y_0) - \psi_{m,\ep}^+(y,y_0),
\end{equation*}
and similarly
\begin{equation*}
\widetilde{\rho_m}(y,y_0):=\lim_{\ep\rightarrow 0}\rho_{m,\ep}^-(y,y_0) - \rho_{m,\ep}^+(y,y_0).
\end{equation*} 

We state the following proposition regarding estimates for $\p_{y_0}\varphi_{m,\ep}^\pm(y,y_0)$ and $\p_{y,y_0}^2\varphi_{m,\ep}^\pm(y,y_0)$, from which one obtains the corresponding estimates for $\p_{y_0}\widetilde{\psi_m}(y,y_0)$ and $\p_{y,y_0}^2\widetilde{\psi_m}(y,y_0)$, respectively.
\begin{proposition}\label{prop: bounds dy0 varphi}
The following holds true.
\begin{itemize}
\item[$(i)$] For $m|y-y_0|\leq 3\b$ and $\b^2\neq1/4$, we have that
\begin{equation*}
|\p_{y_0}\varphi_{m,\ep}^\pm(y,y_0)|\lesssim \frac{1}{m^{1+\mu}}|y-y_0|^{-\frac12-\mu}Q_{1,m} + \sum_{\sigma=0,1}|\B_{m,\ep}^\pm(y,y_0,\sigma)|,
\end{equation*}
and
\begin{equation*}
|\p_{y,y_0}^2\varphi_{m,\ep}^\pm(y,y_0)|\lesssim \frac{1}{m^{1+\mu}}|y-y_0|^{-\frac32-\mu}Q_{1,m} + \sum_{\sigma=0,1}|\p_y\B_{m,\ep}^\pm(y,y_0,\sigma)|,
\end{equation*}
where the bounds for $|\B_{m,\ep}^\pm(y,y_0,\sigma)|$ and $|\p_y\B_{m,\ep}^\pm(y,y_0,\sigma)|$ for $\sigma=0,1$ are given in Propositions \ref{prop:boundary term estimates} and \ref{prop: dyB0 estimates}, respectively.
\item[$(ii)$] For $m|y-y_0|\leq 3\b$ and $\b^2=1/4$, we have that
\begin{equation*}
|\p_{y_0}\varphi_{m,\ep}^\pm(y,y_0)|\lesssim \frac{1}{m} |y-y_0|^{-\frac12} \l( 1 + \big| \log \l(m|y-y_0 \pm i\ep|\r) \big| \r) Q_{1,m} + \sum_{\sigma=0,1}|\B_{m,\ep}^\pm(y,y_0,\sigma)|,
\end{equation*}
and
\begin{equation*}
|\p_{y,y_0}^2\varphi_{m,\ep}^\pm(y,y_0)|\lesssim \frac{1}{m}|y-y_0|^{-\frac32}  \l( 1 + \big| \log \l(m|y-y_0\pm i\ep|\r)\big| \r) Q_{1,m} + \sum_{\sigma=0,1}|\p_y\B_{m,\ep}^\pm(y,y_0,\sigma)|,
\end{equation*}
where the bounds for $|\B_{m,\ep}^\pm(y,y_0,\sigma)|$ and $|\p_y\B_{m,\ep}^\pm(y,y_0,\sigma)|$ for $\sigma=0,1$ are given in Propositions \ref{prop: special boundary term estimates} and \ref{prop: special dyB0 estimates}, respectively.
\item[$(iii)$] For $m|y-y_0|\geq 3\b$, we have that
\begin{equation*}
\Vert \p_{y}\p_{y_0} \varphi_{m,\ep}^\pm \Vert_{L^2_y(J_3^c)}^2 + m^2 \Vert \p_{y_0} \varphi_{m,\ep}^\pm \Vert_{L^2_y(J_3^c)}^2 \lesssim Q_{1,m}^2 + m^2 \sum_{\sigma=0,1}\Vert \B_{m,\ep}^\pm(\cdot,y_0,\sigma)\Vert_{L^2_y(J_2^c\cap J_3)}^2,
\end{equation*}
where the bounds for $\Vert \B_{m,\ep}^\pm(\cdot,y_0,\sigma)\Vert_{L^2_y(J_2^c\cap J_3)}$ are given in Corollary \ref{cor:global boundary term estimates}.
\end{itemize}
In particular, these bounds also apply to $\p_{y_0}\widetilde{\psi_{m}}(y,y_0)$ and $\p_{y,y_0}^2\widetilde{\psi_{m}}(y,y_0)$.
\end{proposition}

\begin{proof}
Both $(i)$ and $(ii)$ follows from Proposition \ref{stream derivatives formulae} and Proposition \ref{L2 bounds inhom TG solution}. As for $(iii)$, we argue assuming that $\b^2\neq1/4$. Taking a $\p_{y_0}$ derivative in \eqref{eq:varphieq}, we see that
\begin{equation*}
\textsc{TG}_{m,\ep}^\pm\p_{y_0}\varphi_{m,\ep}^\pm = \p_{y_0}F_{m,\ep}^\pm-2\b^2\frac{1}{(y-y_0\pm i\ep)^3}\varphi_{m,\ep}^\pm.
\end{equation*}
In order to use Lemma \ref{lemma:entangle ineq}, we need to control $\l\Vert\p_{y_0}\varphi_{m,\ep}^\pm\r\Vert_{L^2_y(J_2^c\cap J_3)}$ and $\l\Vert\frac{1}{(y-y_0\pm i\ep)^3}\varphi_{m,\ep}^\pm\r\Vert_{L^2_y(J_2^c)}$. We begin by estimating
\begin{equation*}
\int_{y_0+\frac{2\b}{m}}^{y_0+\frac{3\b}{m}}|\p_{y_0}\varphi_{m,\ep}^\pm(y,y_0)|^2\d y \lesssim \sum_{\sigma=0,1}\int_{y_0+\frac{2\b}{m}}^{y_0+\frac{3\b}{m}}|\B_{m,\ep}^\pm(y,y_0,\sigma)|^2\d y + Q_{1,m}^2 \int_{y_0+\frac{2\b}{m}}^{y_0+\frac{3\b}{m}} \frac{1}{m^{2+2\mu}}|y-y_0|^{-1-2\mu}\d y
\end{equation*}
Now, for $\b^2< 1/4$ we have $\mu\neq 0$ and 
\begin{equation*}
Q_{1,m}^2 \int_{y_0+2\frac{\b}{m}}^{y_0+3\frac{\b}{m}} \frac{1}{m^{2+2\mu}}|y-y_0|^{-1-2\mu}\d y \lesssim \frac{1}{m^2}Q_{1,m}^2,
\end{equation*}
while for $\b^2>1/4$, we have $\mu=0$ and therefore the bound still becomes
\begin{equation*}
Q_{1,m}^2 \int_{y_0+2\frac{\b}{m}}^{y_0+3\frac{\b}{m}} \frac{1}{m^{2}}|y-y_0|^{-1}\d y \lesssim \frac{1}{m^2}Q_{1,m}^2\l(\log\l(\frac{3\b}{m}\r) - \log\l(\frac{2\b}{m}\r)\r)\lesssim \frac{1}{m^2}Q_{1,m}^2.
\end{equation*}
Therefore, we conclude that 
\begin{equation*}
\l\Vert\p_{y_0}\varphi_{m,\ep}^\pm\r\Vert_{L^2_y(J_2^c\cap J_3)}\lesssim \frac{1}{m}Q_{1,m} + \sum_{\sigma=0,1}\Vert \B_{m,\ep}^\pm(\cdot,y_0,\sigma)\Vert_{L^2_y(J_2^c\cap J_3)}
\end{equation*}

On the other hand, we use Proposition \ref{L2 bounds inhom TG solution} applied to $\varphi_{m,\ep}^\pm(y,y_0)$ to estimate
\begin{equation*}
\l\Vert\frac{1}{(y-y_0\pm i\ep)^3}\varphi_{m,\ep}^\pm\r\Vert_{L^2_y(J_2^c\cap J_3)}^2\lesssim \int_{2\b\leq m|y-y_0|\leq 3\b} \frac{1}{m^{2+2\mu}}|y-y_0\pm i\ep|^{-5-2\mu}Q_{0,m} \d y \lesssim m^2Q_{0,m}
\end{equation*}
and
\begin{equation*}
\l\Vert\frac{1}{(y-y_0\pm i\ep)^3}\varphi_{m,\ep}^\pm\r\Vert_{L^2_y(J_3^c)}^2\lesssim m^6\Vert \varphi_{m,\ep}^\pm\Vert_{L^2_y(J_3^c)}\lesssim m^2Q_{0,m}.
\end{equation*}
The result follows from applying Lemma \ref{lemma:entangle ineq}.
\end{proof}

The next proposition gives bounds on $\p_{y_0}^2\varphi_{m,\ep}^\pm$ and therefore also on $\p_{y_0}^2\widetilde{\psi_{m}}(y,y_0)$.

\begin{proposition}\label{prop: bounds dy0y0 varphi}
The following holds true.
\begin{itemize}
\item For $m|y-y_0|\leq 3\b$ and $\b^2\neq1/4$, we have that
\begin{equation*}
|\p_{y_0}^2\varphi_{m,\ep}^\pm(y,y_0)|\lesssim \frac{1}{m^{1+\mu}}|y-y_0|^{-\frac32-\mu}Q_{2,m} + \sum_{\sigma=0,1}\Big( |\p_y\B_{m,\ep}^\pm(y,y_0,\sigma)| + |\widetilde{\B_{m,\ep}^\pm}(y,y_0,\sigma)|\Big),
\end{equation*}
where the bounds for $|\p_y\B_{m,\ep}^\pm(y,y_0,\sigma)|$ and $|\widetilde{\B_{m,\ep}^\pm}(y,y_0,\sigma)|$ are given in Propositions \ref{prop: dyB0 estimates} and \ref{prop: tildeB0 estimates}, respectively.
\item For $m|y-y_0|\leq 3\b$ and $\b^2=1/4$, we have that
\begin{align*}
|\p_{y_0}^2\varphi_{m,\ep}^\pm(y,y_0)|&\lesssim \frac{1}{m}|y-y_0|^{-\frac32} \l( 1 + \big| \log \l(m|y-y_0\pm i\ep|\r)\big| \r) Q_{2,m} \\
&\quad + \sum_{\sigma=0,1}\Big( |\p_y\B_{m,\ep}^\pm(y,y_0,\sigma)| + |\widetilde{\B_{m,\ep}^\pm}(y,y_0,\sigma)|\Big),
\end{align*}
where the bounds for $|\p_y\B_{m,\ep}^\pm(y,y_0,\sigma)|$ and $|\widetilde{\B_{m,\ep}^\pm}(y,y_0,\sigma)|$ are given in Propositions \ref{prop: special dyB0 estimates} and \ref{prop:  special tildeB0 estimates}, respectively.
\item For $m|y-y_0|\geq 3\b$, we have that
\begin{align*}
\Vert  \p_{y_0}^2\varphi_{m,\ep}^\pm\Vert_{L^2_y(J_3^c)}& \lesssim Q_{2,m} + \sum_{\sigma=0,1}\l( \Vert \p_y\B_{m,\ep}^\pm(\cdot,y_0,\sigma)\Vert_{L^2_y(J_2^c\cap J_3)} + \Vert \widetilde{\B_{m,\ep}^\pm}(\cdot,y_0,\sigma)\Vert_{L^2_y(J_2^c\cap J_3)}\r) \\
&\quad +m \sum_{\sigma=0,1} \Vert B_{m,\ep}^\pm(\cdot,y_0,\sigma)\Vert_{L^2_y(J_2^c\cap J_3)},
\end{align*}
where the estimates for $\Vert \B_{m,\ep}^\pm(\cdot,y_0,\sigma)\Vert_{L^2_y(J_2^c\cap J_3)}$,  $\Vert \p_y\B_{m,\ep}^\pm(\cdot,y_0,\sigma)\Vert_{L^2_y(J_2^c\cap J_3)}$ and $\Vert \widetilde{\B_{m,\ep}^\pm}(\cdot,y_0,\sigma)\Vert_{L^2_y(J_2^c\cap J_3)}$ are given in Corollaries \ref{cor:global boundary term estimates}, \ref{cor: integral tildeB0 estimates} and \ref{cor: integral dyB0 estimates}, respectively.
\end{itemize}
In particular, these bounds also apply to $\p_{y_0}^2\widetilde{\psi_{m}}(y,y_0)$.
\end{proposition}

\begin{proof}
The first two statements of the proposition follow from Proposition \ref{stream derivatives formulae} and Proposition \ref{L2 bounds inhom TG solution}. For the third part of the proposition, we argue for $\b^2\neq1/4$. Taking $\p_{y_0}^2$ derivatives to \eqref{eq:varphieq}, we see that $\p_{y_0}^2\varphi_{m,\ep}^\pm(y,y_0)$ solves
\begin{equation*}
\textsc{TG}_{m,\ep}^\pm \p_{y_0}^2\varphi_{m,\ep}^\pm = \p_{y_0}^2F_{m,\ep}^\pm -4\b^2\frac{\p_{y_0}\varphi_{m,\ep}^\pm}{(y-y_0\pm i\ep)^3} - 6\b^2 \frac{\varphi_{m,\ep}^\pm}{(y-y_0\pm i\ep)^4}.
\end{equation*}
In order to use Lemma \ref{lemma:entangle ineq}, we need to bound $\Vert \p_{y_0}^2\varphi_{m,\ep}^\pm\Vert_{L^2_y(J_2^c\cap J_3)}$, as well as $\Vert \frac{\p_{y_0}\varphi_{m,\ep}^\pm}{(y-y_0\pm i\ep)^3}\Vert_{L^2_y(J_2^c\cap J_3)}$ and $\Vert \frac{\varphi_{m,\ep}^\pm}{(y-y_0\pm i\ep)^4}\Vert_{L^2_y(J_2^c\cap J_3)}$. We estimate
\begin{align*}
\int_{2\b\leq m|y-y_0|\leq 3\b} &|\p_{y_0}^2\varphi_{m,\ep}^\pm(y,y_0)|^2\d y \\
&\lesssim Q_{2,m}\int_{2\b\leq |y-y_0|\leq 3\b} \frac{1}{m^{2+2\mu}}|y-y_0|^{-3-2\mu}\d y \\
&\qquad+ \sum_{\sigma=0,1}\Big( \Vert \p_y\B_{m,\ep}^\pm(\cdot,y_0,\sigma)\Vert_{L^2_y(J_2^c\cap J_3)} + \Vert \widetilde{\B_{m,\ep}^\pm}(\cdot,y_0,\sigma)\Vert_{L^2_y(J_2^c\cap J_3)}\Big) \\
&\lesssim Q_{2,m} + \sum_{\sigma=0,1}\Big( \Vert \p_y\B_{m,\ep}^\pm(\cdot,y_0,\sigma)\Vert_{L^2_y(J_2^c\cap J_3)} + \Vert \widetilde{\B_{m,\ep}^\pm}(\cdot,y_0,\sigma)\Vert_{L^2_y(J_2^c\cap J_3)}\Big).
\end{align*}
Similarly, from Proposition \ref{prop: bounds dy0 varphi} we have that
\begin{equation*}
\l\Vert \frac{\p_{y_0}\varphi_{m,\ep}^\pm}{(y-y_0\pm i\ep)^3}\r\Vert_{L^2_y(J_2^c)}\lesssim m^3\Vert {\p_{y_0}\varphi_{m,\ep}^\pm}\Vert_{L^2_y(J_2^c\cap J_3)}\lesssim m^2Q_{1,m} + m^3\sum_{\sigma=0,1}\Vert \B_{m,\ep}^\pm(\cdot,y_0,\sigma)\Vert_{L^2_y(J_2^c\cap J_3)},
\end{equation*}
while using Proposition \ref{L2 bounds inhom TG solution} we obtain
\begin{equation*}
\l\Vert \frac{\varphi_{m,\ep}^\pm}{(y-y_0\pm i\ep)^4}\r\Vert_{L^2_y(J_2^c)}\lesssim m^4\Vert {\varphi_{m,\ep}^\pm}\Vert_{L^2_y(J_2^c\cap J_3)}\lesssim m^2Q_{0,m}.
\end{equation*}
With this, the proof is complete.
\end{proof}

We finish the subsection by providing the estimates for $\widetilde{\rho_m}$ and $\p_{y_0}\widetilde{\rho_m}$.
\begin{proposition}\label{prop: bounds rho and dy0 rho}
The following holds true.
\begin{itemize}
\item For $m|y-y_0|\leq 3\b$ and $\b^2\neq1/4$, we have that
\begin{equation*}
|\widetilde{\rho_{m}}(y,y_0)|\lesssim \frac{1}{m^{1+\mu}}|y-y_0|^{-\frac12-\mu}Q_{0,m}.
\end{equation*}
and
\begin{equation*}
|\p_{y_0}\widetilde{\rho_{m}}(y,y_0)|\lesssim \frac{1}{m^{1+\mu}}|y-y_0|^{-\frac32-\mu}Q_{1,m} + \sup_{0\leq \ep \leq \ep_0} \sum_{\sigma=0,1}\sum_{\kappa\in\lbrace +,-\rbrace} |y-y_0 +\kappa i\ep|^{-1} |\B_{m,\ep}^\kappa(y,y_0,\sigma)|,
\end{equation*}
where the bounds for $|\B_{m,\ep}^\pm(y,y_0,\sigma)|$ for $\sigma=0,1$ are given in Proposition \ref{prop:boundary term estimates}.
\item For $m|y-y_0|\leq 3\b$ and $\b^2=1/4$, we have that
\begin{equation*}
|\widetilde{\rho_{m}}(y,y_0)|\lesssim \frac{1}{m}|y-y_0|^{-\frac12} \l( 1 + \big| \log \l(m|y-y_0\pm i\ep|\r)\big| \r) Q_{0,m}.
\end{equation*}
and
\begin{align*}
|\p_{y_0}\widetilde{\rho_{m}}(y,y_0)|&\lesssim \frac{1}{m}|y-y_0|^{-\frac32} \l( 1 + \big| \log \l(m|y-y_0\pm i\ep|\r)\big| \r) Q_{1,m} \\
&\quad+ \sup_{0\leq \ep \leq \ep_0} \sum_{\sigma=0,1}\sum_{\kappa\in\lbrace +,-\rbrace} |y-y_0 +\kappa i\ep|^{-1} |\B_{m,\ep}^\kappa(y,y_0,\sigma)|,
\end{align*}
where the bounds for $|\B_{m,\ep}^\pm(y,y_0,\sigma)|$ for $\sigma=0,1$ are given in Proposition \ref{prop: special boundary term estimates}.
\item For $m|y-y_0|\geq 3\b$, we have that
\begin{equation*}
\Vert \widetilde{\rho_m}\Vert_{L^2_y(J_3^c)}\lesssim \frac{1}{m}Q_{0,m}
\end{equation*}
and
\begin{equation*}
\Vert \p_{y_0}\widetilde{\rho_m}\Vert_{L^2_y(J_3^c)}\lesssim Q_{1,m} + m\sum_{\sigma=0,1}\sum_{\kappa\in\lbrace +,-\rbrace} \Vert \B_{m,\ep}^\kappa(y,y_0,\sigma)\Vert_{L^2_y(J_2^c\cap J_3)}
\end{equation*}
\end{itemize}
\end{proposition}

\begin{proof}
The bounds follow directly from Proposition \ref{density derivative formulae}, Proposition \ref{L2 bounds inhom TG solution} and Proposition \ref{prop: bounds dy0 varphi}.
\end{proof}

\section{Time-decay estimates}\label{sec: time decay estimates}
This section is devoted to the proof of the time decay estimate rates of the stream function $\psi_{m}(t,y)$, its $\p_y\psi_m(t,y)$ derivative and the density $\rho_m(t,y)$.  Let us recall that we can write
\begin{equation*}
\psi_m(t,y)=\frac{1}{2\pi i}\lim_{\ep\rightarrow 0} \int_0^1 \e^{-imy_0t}\l( \psi^-_m(y,y_0) - \psi_m^+(y,y_0)\r) \d y_0,
\end{equation*}
and
\begin{equation*}
\rho_m(t,y)=\frac{1}{2\pi i}\lim_{\ep\rightarrow 0} \int_0^1 \e^{-imy_0t}\l( \rho^-_m(y,y_0) - \rho_m^+(y,y_0)\r) \d y_0.
\end{equation*}
A simple integration by parts provides
\begin{equation*}
\begin{aligned}
\psi_m(t,y)&=-\frac{1}{2\pi i}\frac{1}{imt}\lim_{\ep\rightarrow 0}\big[ \e^{-imy_0t}\l(\psi_{m,\ep}^-(y,y_0)-\psi_{m,\ep}^+(y,y_0)\r)\big]_{y_0=0}^{y_0=1} \\
&\quad +\frac{1}{2\pi i}\frac{1}{imt}\lim_{\ep\rightarrow 0}\int_0^1 \e^{-imy_0t}\l( \p_{y_0}\psi_{m,\ep}^-(y,y_0)- \p_{y_0}\psi_{m,\ep}^+(y,y_0)\r) \d y_0 \\ 
&=\frac{1}{2\pi i}\frac{1}{imt}\lim_{\ep\rightarrow 0}\int_0^1 \e^{-imy_0t}\l( \p_{y_0}\psi_{m,\ep}^-(y,y_0)- \p_{y_0}\psi_{m,\ep}^+(y,y_0)\r) \d y_0
\end{aligned}
\end{equation*}
where we use Theorem \ref{thm: boundary limiting absorption principle} to show that the solid terms associated to the spectral boundary vanish. Throughout the entire section, let us consider $\b^2\neq1/4$, unless we state otherwise. 

We begin proving the following result.
\begin{proposition}\label{prop: invdamp psi}
Let $t\geq 1$. Then,
\begin{equation*}
\begin{aligned}
\Vert \psi_m(t) \Vert_{L^2_y}&\lesssim m^{-\frac32}t^{-\frac32+\mu} Q_{2,m} \\
\end{aligned}
\end{equation*}
\end{proposition}

\begin{proof}
We write
\begin{equation*}
\begin{aligned}
\psi_m(t,y)&=\frac{1}{2\pi i}\frac{1}{imt}\lim_{\ep\rightarrow 0}\int_0^1 \e^{-imy_0t}\l( \p_{y_0}\psi_{m,\ep}^-(y,y_0)- \p_{y_0}\psi_{m,\ep}^+(y,y_0)\r) \d y_0.
\end{aligned}
\end{equation*}
Let us denote $\delta_0:=\min\l( \frac{3\b}{m},\frac{1}{2}\r)$ and let $\delta\in \l(0,\frac{\delta_0}{2}\r)$. In particular, we note that $m\delta\leq3\b$, it is bounded. We shall first show the decay rates for $\Vert \psi_m(t) \Vert_{L^2_y(\delta,1-\delta)}$ and then for $\Vert \psi_m(t) \Vert_{L^2_y(0,\delta)}$ and $\Vert \psi_m(t) \Vert_{L^2_y(1-\delta,1)}$. 
\bullpar{Step 1} For $y\in (\delta, 1-\delta)$, we write
\begin{equation*}
\begin{aligned}
\psi_m(t,y)&=\frac{1}{2\pi i}\frac{1}{imt}\lim_{\ep\rightarrow 0}\l(\int_0^{y-\frac{\delta}{2}} +\int_{y-\frac{\delta}{2}}^{y+\frac{\delta}{2}} + \int_{y+\frac{\delta}{2}}^1\r) \e^{-imy_0t}\l( \p_{y_0}\psi_{m,\ep}^-(y,y_0)- \p_{y_0}\psi_{m,\ep}^+(y,y_0)\r) \d y_0 \\
&=\mathcal{T}_1 + \mathcal{T}_2 + \mathcal{T}_3.
\end{aligned}
\end{equation*}
and we begin with estimating $\mathcal{T}_2$. There, we have that $|y-y_0|\leq \frac{\delta}{2}\leq \frac{\delta_0}{4}$ and we can use the bounds from Proposition \ref{prop: bounds dy0 varphi} to bound
\begin{align*}
|\mathcal{T}_2|&\lesssim \frac{1}{mt}\int_{y-\frac{\delta}{2}}^{y+\frac{\delta}{2}}\frac{1}{m^{1+\mu}}|y-y_0|^{-\frac12+\mu}Q_{1,m}\d y_0 + \frac{1}{mt}\int_{y-\frac{\delta}{2}}^{y+\frac{\delta}{2}} \sum_{\sigma=0,1}\sum_{\kappa\in\lbrace +,-\rbrace}|\B_{m,\ep}^\kappa(y,y_0,\sigma)| \d y_0  \\
&=\mathcal{T}_{2,1} + \mathcal{T}_{2,2}.
\end{align*}
We can integrate directly to obtain
\begin{equation*}
|\mathcal{T}_{2,1}|\lesssim  \frac{1}{m^{2+\mu}t}\delta^{\frac12-\mu}Q_{1,m}.
\end{equation*}
For $\mathcal{T}_{2,2}$, for $\sigma=0$ we decompose
\begin{align*}
\frac{1}{mt}\int_{y-\frac{\delta}{2}}^{y+\frac{\delta}{2}} |\B_{m,\ep}^\pm(y,y_0,\sigma)| \d y_0 = \frac{1}{mt}\int_{y-\frac{\delta}{2}}^{y+\frac{\delta}{2}} \l(\chi_{y_0\leq \frac{3\b}{m}} + \chi_{y_0> \frac{3\b}{m}}\r) |\B_{m,\ep}^\pm(y,y_0,\sigma)| \d y_0 = \mathcal{T}_{2,2,1} + \mathcal{T}_{2,2,2}.
\end{align*}
We use Proposition \ref{prop:boundary term estimates} to compute
\begin{equation*}
\mathcal{T}_{2,2,1}\lesssim \frac{1}{mt}\int_{y-\frac{\delta}{2}}^{y+\frac{\delta}{2}}m^{-1+\mu}y_0^{-\frac12+\mu}Q_{0,m}\d y_0\lesssim \frac{1}{m^{2-\mu}t}\delta^{\frac12+\mu}Q_{0,m}
\end{equation*}
and
\begin{equation*}
\mathcal{T}_{2,2,2}\lesssim \frac{1}{mt}\int_{y-\frac{\delta}{2}}^{y+\frac{\delta}{2}}Q_{0,m}\d y_0\lesssim \frac{1}{mt}\delta Q_{0,m}.
\end{equation*}
The bounds for the terms of $T_{2,2}$ for $\sigma=1$ are the same, we omit the details. We summarize these estimates into
\begin{equation*}
\Vert \cT_2 \Vert_{L^2_y(\delta,1-\delta)}\lesssim \frac{1}{mt}\l( m^{-\frac32}(m\delta)^{\frac12-\mu} + m^{-\frac32}(m\delta)^{\frac{1}{2}+\mu} + \delta\r)Q_{1,m}
\end{equation*}

We shall next estimate $\mathcal{T}_1$, the bounds of $\mathcal{T}_3$ are the same and the arguments to prove them are identical. For $\mathcal{T}_1$, note that we can further integrate by parts, 
\begin{equation*}
\begin{aligned}
\mathcal{T}_1&=-\frac{1}{2\pi i}\frac{1}{m^2t^2}\lim_{\ep\rightarrow 0}\big[ \e^{-imy_0t}\p_{y_0}\l(\psi_{m,\ep}^-(y,y_0)-\psi_{m,\ep}^+(y,y_0)\r)\big]_{y_0=\frac{\delta}{2}}^{y_0=y-\frac{\delta}{2}} \\
&\quad+ \frac{1}{2\pi i}\frac{1}{imt}\lim_{\ep\rightarrow 0}\int_0^{\frac{\delta}{2}} \e^{-imy_0t}\p_{y_0}\l(\psi_{m,\ep}^-(y,y_0)- \psi_{m,\ep}^+(y,y_0)\r) \d y_0 \\
&\quad+ \frac{1}{2\pi i}\frac{1}{m^2t^2}\lim_{\ep\rightarrow 0}\int_{\frac{\delta}{2}}^{y-\frac{\delta}{2}} \e^{-imy_0t}\p_{y_0}^2\l(\psi_{m,\ep}^-(y,y_0)- \psi_{m,\ep}^+(y,y_0)\r) \d y_0 \\
&=\mathcal{T}_{1,1}+ \mathcal{T}_{1,2}+\mathcal{T}_{1,3}.
\end{aligned}
\end{equation*}
We shall treat each $\mathcal{T}_{1,i}$, for $i=1,2,3$ separately.

\diampar{Estimates for $\mathcal{T}_{1,1}$}
For the boundary terms of $\mathcal{T}_{1,1}$, consider first $y_0=y-\frac{\delta}{2}$. Then, $|y-y_0|=\frac{\delta}{2}\leq\frac{\delta_0}{4}$, so that from Proposition \ref{prop: bounds dy0 varphi} we have
\begin{equation*}
\frac{1}{m^2t^2}\l|\p_{y_0}\widetilde{\psi_m}\l( y,y-\tfrac{\delta}{2}\r)\r|\lesssim \frac{1}{m^2t^2}\l(\frac{1}{m^{1+\mu}}\delta^{-\frac12-\mu}Q_{1,m} + \sum_{\sigma=0,1}\sum_{\kappa \in \lbrace +,- \rbrace} \l|\B_{m,\ep}^\kappa \l(y,y-\tfrac{\delta}{2},\sigma \r)\r|\r).
\end{equation*}
Now, from Proposition \ref{prop:boundary term estimates} we have
\begin{equation*}
\sum_{\sigma=0,1}\sum_{\kappa \in \lbrace +,- \rbrace}\l|\B_{m,\ep}^\kappa \l(y,y-\tfrac{\delta}{2},\sigma \r)\r| \lesssim \l( 1 + m^{-\frac12}(m\delta)^{-\frac12+\mu}\r)Q_{0,m},
\end{equation*}
since $y\in(\delta,1-\delta)$ ensures $y-\frac{\delta}{2}>\frac{\delta}{2}$.

For the boundary term $\mathcal{T}_{1,1}$ associated to $y_0=\frac{\delta}{2}$, since $y\in (\delta, 1-\delta)$, we have that $1-\frac{\delta}{2}\geq y-y_0\geq \frac{\delta}{2}$. Hence, for those $y\in(\delta,1-\delta)$ such that  $m|y-y_0|\leq 3\b$, we use Proposition \ref{prop: bounds dy0 varphi} to pointwise estimate
\begin{equation*}
\frac{1}{m^2t^2}\l|\p_{y_0}\widetilde{\psi_m}\l( y,\tfrac{\delta}{2}\r)\r| \lesssim \frac{1}{m^2t^2}\l(\frac{1}{m^{1+\mu}}\delta^{-\frac12-\mu}Q_{1,m} + \sum_{\sigma=0,1}\sum_{\kappa \in \lbrace +,- \rbrace} \l|\B_{m,\ep}^\kappa \l( y,\tfrac{\delta}{2},\sigma \r)\r|\r),
\end{equation*}
where we further have from Proposition \ref{prop:boundary term estimates} that
\begin{equation*}
\sum_{\sigma=0,1}\sum_{\kappa \in \lbrace +,- \rbrace}\Vert \B_{m,\ep}^\pm \l( y,\tfrac{\delta}{2},0 \r)\Vert_{L^2_y(J_3)}\lesssim \l( m^{-\frac32}\delta^{-\frac12} + m^{-\frac12}\r)  Q_{0,m}.
\end{equation*}
Next, for those $y\in(\delta,1-\delta)$ such that $m|y-y_0|\geq 3\b$ we can directly estimate in $L^2_y$ using Proposition \ref{prop: bounds dy0 varphi} to deduce that 
\begin{equation*}
\frac{1}{m^2t^2}\l\Vert \p_{y_0}\widetilde{\psi_m}\l( y,\tfrac{\delta}{2}\r)\r\Vert_{L^2_y(J_3^c)}\lesssim \frac{1}{m^2t^2}\l(\frac{1}{m}Q_{1,m} + \sum_{\sigma=0,1}\sum_{\kappa\in \lbrace +,-\rbrace}\Vert \B_{m,\ep}^\kappa \l(y, \tfrac{\delta}{2}, \sigma \r) \Vert_{L^2_y(J_2^c\cap J_3)}\r),
\end{equation*}
while from Corollary \ref{cor:global boundary term estimates} we are able to bound
\begin{equation*}
\Vert \B_{m,\ep}^\pm \l(y, \tfrac{\delta}{2}, 0 \r) \Vert_{L^2_y(J_2^c\cap J_3)}\lesssim m^{-\frac32}\delta^{-\frac12}Q_{0,m}, \quad \Vert \B_{m,\ep}^\pm \l(y, \tfrac{\delta}{2}, 1 \r) \Vert_{L^2_y(J_2^c\cap J_3)}\lesssim m^{-\frac12}Q_{0,m}.
\end{equation*}
Therefore, we have
\begin{equation*}
\Vert \cT_{1,1}\Vert_{L^2_y(\delta,1-\delta)}\lesssim \frac{1}{m^2t^2}\l(1 +  m^{-\frac12}(m\delta)^{-\frac12-\mu} + m^{-\frac12}(m\delta)^{-\frac12+\mu} + m^{-\frac32}\delta^{-\frac12} + m^{-\frac12} \r)Q_{0,m}
\end{equation*}
This concludes the analysis of $\mathcal{T}_{1,1}$.

\diampar{Estimates for $\mathcal{T}_{1,2}$}
We begin by splitting
\begin{equation*}
|\mathcal{T}_{1,2}|\lesssim \frac{1}{mt}\int_0^{\frac{\delta}{2}}\l( \chi_{m|y-y_0|\leq 3\b} + \chi_{m|y-y_0|>3\b} \r)|\p_{y_0}\widetilde{\psi_m}(y,y_0)|\d y_0 = \mathcal{T}_{1,2,1} + \mathcal{T}_{1,2,2}
\end{equation*}
We use Proposition \ref{prop: bounds dy0 varphi} to estimate
\begin{equation*}
|\cT_{1,2,1}|\lesssim \frac{1}{mt}\int_0^{\frac{\delta}{2}} \l(\frac{1}{m^{1+\mu}}|y-y_0|^{-\frac12-\mu}Q_{1,m} + \chi_{m|y-y_0|\leq 3\b}\sum_{\sigma=0,1}\sum_{\kappa\in \lbrace +,- \rbrace}|\B_{m,\ep}^\pm(y,y_0,\sigma)|\r) \d y_0
\end{equation*}
Now, since $y\in (\delta,1-\delta)$ and $y_0\leq \frac{\delta}{2}$, we have $|y-y_0|\geq \frac{\delta}{2}$. Hence,
\begin{equation*}
\int_0^{\frac{\delta}{2}} \frac{1}{m^{1+\mu}}|y-y_0|^{-\frac12-\mu}Q_{1,m}\d y_0 \lesssim \frac{1}{m^{1+\mu}}\delta^{\frac{1}{2}-\mu}Q_{1,m}.
\end{equation*}
Moreover, Proposition \ref{prop:boundary term estimates} provides
\begin{align*}
\sum_{\sigma=0,1}\sum_{\kappa\in \lbrace +,- \rbrace} &\l\Vert \int_0^{\frac{\delta}{2}} \chi_{m|y-y_0|\leq 3\b} |\B_{m,\ep}^\pm(y,y_0,\sigma)| \d y_0 \r\Vert_{L^2_y(\delta,1-\delta)} \\
&\lesssim \sum_{\sigma=0,1}\sum_{\kappa\in \lbrace +,- \rbrace}  \int_0^{\frac{\delta}{2}} \l\Vert \B_{m,\ep}^\pm(y,y_0,\sigma)\r\Vert_{L^2_y(J_3)} \d y_0 \\
&\lesssim  m^{-\frac12}\delta^{\frac12}\l( m^{-1} + \delta^\frac12\r)Q_{0,m}.
\end{align*}
As a result, we are able to bound
\begin{equation*}
\Vert \cT_{1,2,1}\Vert_{L^2_y(\delta,1-\delta)}\lesssim \frac{1}{mt}\l(m^{-1-\mu}\delta^{\frac12-\mu} + m^{-\frac32}\delta^\frac12 + m^{-\frac12}\delta\r) Q_{1,m}.
\end{equation*}
We again use Proposition \ref{prop: bounds dy0 varphi} and Corollary \ref{cor:global boundary term estimates} to estimate
\begin{align*}
\Vert \cT_{1,2,2}\Vert_{L^2_y(\delta,1-\delta)}&\lesssim \frac{1}{mt}\int_0^{\frac{\delta}{2}}\Vert \p_{y_0}\widetilde{\psi_m}(y,y_0)\Vert_{L^2_y(J_3^c)}\d y_0 \\
&\lesssim \frac{1}{mt}\int_0^{\frac{\delta}{2}} \l( \frac{1}{m}Q_{1,m} + \sum_{\sigma=0,1}\sum_{\kappa\in \lbrace +,-\rbrace}\Vert \B_{m,\ep}^\kappa \l(y, \tfrac{\delta}{2}, \sigma \r) \Vert_{L^2_y(J_2^c\cap J_3)}\r)\d y_0 \\
&\lesssim \frac{1}{mt}\l( m^{-1}\delta + m^{-\frac32}\delta^{\frac12} + m^{-\frac12}\delta\r) Q_{1,m}
\end{align*}
so that we can conclude
\begin{equation*}
\Vert \cT_{1,2} \Vert_{L^2_y(\delta,1-\delta)}\lesssim \frac{1}{mt}\l( m^{-1-\mu}\delta^{\frac12-\mu} + m^{-1}\delta + m^{-\frac32}\delta^{\frac12} + m^{-\frac12}\delta\r) Q_{1,m}.
\end{equation*}

\diampar{Estimates for $\mathcal{T}_{1,3}$}
We shall split again
\begin{equation*}
|\cT_{1,3}|\lesssim \frac{1}{m^2t^2}\int_{\frac{\delta}{2}}^{y-\frac{\delta}{2}}\l( \chi_{m|y-y_0|\leq 3\b} + \chi_{m|y-y_0|>3\b} \r)|\p_{y_0}^2\widetilde{\psi_m}(y,y_0)|\d y_0 = \mathcal{T}_{1,3,1} + \mathcal{T}_{1,3,2}.
\end{equation*}
Now, we use Proposition \ref{prop: bounds dy0y0 varphi} to estimate
\begin{align*}
|\cT_{1,3,1}|&\lesssim \frac{1}{m^2t^2}\int_{\frac{\delta}{2}}^{y-\frac{\delta}{2}} \frac{1}{m^{1+\mu}}|y-y_0|^{-\frac32-\mu}Q_{2,m}\d y_0 \\
&\quad+\frac{1}{m^2t^2} \sum_{\sigma=0,1} \sum_{\kappa\in \lbrace +,- \rbrace} \int_{\frac{\delta}{2}}^{y-\frac{\delta}{2}} \chi_{m|y-y_0|\leq 3\b}\Big( |\p_y\B_{m,\ep}^\kappa(y,y_0,\sigma)| + |\widetilde{\B_{m,\ep}^\kappa}(y,y_0,\sigma)|\Big)  \d y_0.
\end{align*}
Clearly, since $y\in(\delta,1-\delta)$, we have that
\begin{equation*}
\int_{\frac{\delta}{2}}^{y-\frac{\delta}{2}} \frac{1}{m^{1+\mu}}|y-y_0|^{-\frac32-\mu}Q_{2,m}\d y_0 \lesssim \frac{1}{m^{1+\mu}}\l[(y-y_0)^{-\frac12-\mu}\r]_{y_0=\frac{\delta}{2}}^{y_0=y-\frac{\delta}{2}}Q_{2,m}\lesssim \frac{1}{m^{1+\mu}} \delta^{-\frac12-\mu}Q_{2,m}.
\end{equation*}
Similarly, Proposition \ref{prop: tildeB0 estimates} provides
\begin{align*}
\sum_{\sigma=0,1}\sum_{\kappa\in \lbrace +,- \rbrace} &\l\Vert \int_{\frac{\delta}{2}}^{y-\frac{\delta}{2}} \chi_{m|y-y_0|\leq 3\b} |\widetilde{\B_{m,\ep}^\kappa}(y,y_0,\sigma)| \d y_0 \r\Vert_{L^2_y(\delta,1-\delta)} \\
&\lesssim \sum_{\sigma=0,1}\sum_{\kappa\in \lbrace +,- \rbrace}  \int_{\frac{\delta}{2}}^{1-\frac{\delta}{2}} \l\Vert \widetilde{\B_{m,\ep}^\kappa}(y,y_0,\sigma)\r\Vert_{L^2_y(J_3)} \d y_0 \\
&\lesssim  \l(m^{-1}\delta^{-\frac12} +m^{-\frac12}\delta^\frac12+  m^\frac12\r)    Q_{1,m},
\end{align*}
while Proposition \ref{prop: dyB0 estimates} gives
\begin{equation*}
\sum_{\sigma=0,1}\sum_{\kappa\in \lbrace +,- \rbrace}  \int_{\frac{\delta}{2}}^{y-\frac{\delta}{2}} \l| \p_y{\B_{m,\ep}^\kappa}(y,y_0,\sigma)\r|\d y_0 \lesssim \l( m^{\frac12-\mu} + (m\delta)^{\frac12-\mu} + m^{-\frac12}\delta^{-\frac{1}{2}+ \mu}\r)Q_{0,m}.
\end{equation*}
Therefore, we have
\begin{equation*}
\Vert \cT_{1,3,1}\Vert_{L^2_y(\delta,1-\delta)}\lesssim \frac{1}{m^2t^2}\l( m^\frac12 + m^{-\frac12}(m\delta)^{-\frac12-\mu}  + m(m\delta)^{\frac12-\mu} + (m\delta)^{-\frac12}\r) Q_{2,m}.
\end{equation*}
For $\cT_{1,3,2}$, we use Minkowski inequality, Proposition \ref{prop: bounds dy0y0 varphi} and Corollaries \ref{cor:global boundary term estimates}, \ref{cor: integral tildeB0 estimates} and \ref{cor: integral dyB0 estimates} to estimate
\begin{align*}
\Vert \cT_{1,3,2}\Vert_{L^2_y(\delta,1-\delta)}&\lesssim \frac{1}{m^2t^2}\int_{\frac{\delta}{2}}^{1-\frac{\delta}{2}} \l( Q_{2,m} + m \sum_{\sigma=0,1}\sum_{\kappa\in \lbrace +, - \rbrace} \Vert \B_{m,\ep}^\kappa(y,y_0,\sigma)\Vert_{L^2_y(J_2^c\cap J_3)} \r) \d y_0 \\
&\quad + \frac{1}{m^2t^2} \sum_{\sigma=0,1}\sum_{\kappa\in \lbrace +, - \rbrace} \int_{\frac{\delta}{2}}^{1-\frac{\delta}{2}}  \Vert \p_y\B_{m,\ep}^\kappa(y,y_0,\sigma)\Vert_{L^2_y(J_2^c\cap J_3)} \d y_0 \\
&\quad + \frac{1}{m^2t^2} \sum_{\sigma=0,1}\sum_{\kappa\in \lbrace +, - \rbrace} \int_{\frac{\delta}{2}}^{1-\frac{\delta}{2}} \Vert \widetilde{B_{m,\ep}^\kappa}(y,y_0,\sigma)\Vert_{L^2_y(J_2^c\cap J_3)} \d y_0 \\
&\lesssim \frac{1}{m^2t^2}\l( m^\frac12 + m^{-\frac12}\delta^\frac12 + m^{-1}\delta^{-\frac12}\r)Q_{2,m}.
\end{align*}
Hence, we conclude that
\begin{equation*}
\Vert \cT_{1,3}\Vert_{L^2_y(\delta,1-\delta)} \lesssim \frac{1}{m^2 t^2}\l( m^\frac12 + m^{-\frac12}(m\delta)^{-\frac12-\mu} + m(m\delta)^{\frac12-\mu} + m^{-\frac12}\delta^\frac12\r)Q_{2,m}
\end{equation*}
and thus
\begin{equation*}
\Vert \cT_1 \Vert_{L^2_y(\delta,1-\delta)}\lesssim \frac{1}{mt}m^{-\frac32}(m\delta)^{\frac12}Q_{2,m} + \frac{1}{m^2 t^2}\l(m^\frac12 + m^{-\frac12}(m\delta)^{-\frac12-\mu}\r)Q_{2,m}. 
\end{equation*}
In particular, gathering the estimates for $\cT_{2}$ and $\cT_{1,i}$, for $i=1,2,3$, we obtain
\begin{equation*}
\Vert \psi_m(t,y)\Vert_{L^2_y(\delta,1-\delta)}\lesssim \frac{1}{m^2 t}(m\delta)^{\frac12-\mu}Q_{2,m} + \frac{1}{m^2t^2}\l( m^\frac12 + m^{-\frac12}(m\delta)^{-\frac12-\mu}\r) Q_{2,m}
\end{equation*}

\bullpar{Step 2} For $y\in(0,\delta)$, we have that
\begin{equation*}
\psi_m(t,y)=\frac{1}{2\pi i}\frac{1}{imt}\lim_{\ep\rightarrow0}\l( \int_0^{y+\frac{\delta}{2}} + \int_{y+\frac{\delta}{2}}^1 \r) e^{imy_0t}\l( \p_{y_0}\psi_m^-(y,y_0) - \psi_m^+(y,y_0) \r) \d y_0 = \widetilde{\cT}_1 + \widetilde{\cT}_2.
\end{equation*}
One can see that the bounds for $\widetilde{\cT}_2$ here are the same as the ones for $\cT_3$, the procedure to obtain them is the same, we thus omit the details. On the other hand, for $\widetilde{\cT}_1$ we argue as follows. Note that for $0\leq y_0\leq y+\frac{\delta}{2}$, we have that $|y-y_0|\leq \delta\leq \frac{3\b}{m}$ and therefore we have from Proposition \ref{prop: bounds dy0 varphi},
\begin{equation*}
|\widetilde{\cT}_1|\lesssim \frac{1}{mt}\int_0^{y+\frac{\delta}{2}}\l( \frac{1}{m^{1+\mu}}|y-y_0|^{-\frac12-\mu}Q_{1,m} + \sum_{\sigma=0,1}\sum_{\kappa\in \lbrace +,- \rbrace}\l| \B_{m,\ep}^\kappa(y,y_0,\sigma)\r| \r) \d y_0.
\end{equation*}
Since $y\in (0,\delta)$, we trivially have that
\begin{equation*}
\frac{1}{mt}\int_0^{y+\frac{\delta}{2}} \frac{1}{m^{1+\mu}}|y-y_0|^{-\frac12-\mu}Q_{1,m}\d y_0\lesssim \frac{1}{mt}m^{-\frac32}(m\delta)^{\frac12-\mu}Q_{1,m}.
\end{equation*}
Similarly, using the bounds from Proposition \ref{prop:boundary term estimates},
\begin{align*}
\sum_{\sigma=0,1}\sum_{\kappa\in \lbrace +,- \rbrace} \int_0^{y+\frac{\delta}{2}} \l| \B_{m,\ep}^\kappa(y,y_0,\sigma)\r| \d y_0 &\lesssim \int_0^{y+\frac{\delta}{2}}\l( 1+ m^{-1+\mu}y_0^{-\frac12 + \mu}\r)\d y_0  \\
&\lesssim y+\tfrac{\delta}{2} + m^{-1+\mu}\l( y+\tfrac{\delta}{2}\r)^{\frac12+\mu}.
\end{align*}
As a result, we compute
\begin{equation*}
\Vert \widetilde{\cT}_1\Vert_{L^2_y(0,\delta)}\lesssim \frac{1}{mt}\l( m^{-2}(m\delta)^{1-\mu} + \delta^\frac32 + m^{-2}(m\delta)^{1+\mu}\r)Q_{1,m} \lesssim \frac{1}{mt}m^{-\frac32}(m\delta)^{1-\mu}Q_{1,m},
\end{equation*}
and thus we obtain
\begin{equation*}
\Vert \psi_m(t)\Vert_{L^2_y(0,\delta)}\lesssim \frac{1}{mt}m^{-\frac32}(m\delta)^\frac12Q_{2,m} + \frac{1}{m^2t^2}\l(m^\frac12 + m^{-\frac12}(m\delta)^{-\frac12-\mu}\r) Q_{2,m}.
\end{equation*}
The same bounds can be achieved for $\Vert \psi_m(t)\Vert_{L^2_y(1-\delta,1)}$, the proof of which follows along the same ideas, we omit the details. With all these bounds, we are now able to estimate
\begin{align*}
\Vert \psi_m(t)\Vert_{L^2_y}&\lesssim \frac{1}{mt}m^{-1}(m\delta)^{\frac12-\mu}Q_{2,m} + \frac{1}{m^2t^2}\l(m^\frac12 + m^{-\frac12}(m\delta)^{-\frac12-\mu}\r)Q_{2,m} \\
&\lesssim m^{-\frac32}t^{-\frac32+\mu}Q_{2,m}
\end{align*}
once we choose $\delta=\frac{c_0}{mt}$, for $c_0=\frac{1}{1000}\min(\b,1)$. The proof is complete.
\end{proof}

\begin{proposition}\label{prop: invdamp dy psi}
Let  $t\geq 1$. Then, 
\begin{equation*}
\Vert \p_y\psi_m(t)\Vert_{L^2_y}\lesssim m^{-\frac12}t^{-\frac12 + \mu} Q_{1,m}.
\end{equation*}
\end{proposition}

\begin{proof}
The argument follows the same lines as the proof for $\Vert \psi_m(t,y)\Vert_{L^2_y}$. Hence, we shall only give the bounds for the terms involved in the computation. Their proof have already been carried out in the proof of Proposition \ref{prop: invdamp psi}.

\bullpar{Step 1} For $y\in(\delta,1-\delta)$ we shall write
\begin{align*}
\p_y\psi_m(t,y)&=\frac{1}{2\pi i}\lim_{\ep\rightarrow 0} \l( \int_0^{y-\frac{\delta}{2}} + \int_{y-\frac{\delta}{2}}^{y+\frac{\delta}{2}} + \int_{y+\frac{\delta}{2}}^1 \r) e^{-imy_0t}\l( \p_y\psi_{m,\ep}^-(y,y_0) - \p_y\psi_{m,\ep}^+(y,y_0)\r) \d y_0 \\
&=\cI_1 + \cI_2 + \cI_3
\end{align*}
We begin by using Proposition \ref{L2 bounds inhom TG solution} to bound
\begin{equation*}
\l\Vert \cI_2 \r\Vert_{L^2_y(\delta,1-\delta)}\lesssim  \l\Vert \int_{y-\frac{\delta}{2}}^{y+\frac{\delta}{2}} \frac{1}{m^{1+\mu}}|y-y_0|^{-\frac12-\mu} \d y_0 \r\Vert_{L^2_y(\delta,1-\delta)}\lesssim  m^{-\frac32}(m\delta)^{\frac12-\mu}Q_{0,m}.
\end{equation*}
As before, for $\cI_1$ we split it into 
\begin{align*}
\cI_1 &= -\frac{1}{2\pi i} \frac{1}{imt}\lim_{\ep\rightarrow0} \l[ e^{-imy_0t}\l(\p_y\psi_{m,\ep}^-(y,y_0) - \p_y\psi_{m,\ep}^+(y,y_0)\r) \r]_{y_0=\frac{\delta}{2}}^{y_0=y-\frac{\delta}{2}} \\
&\quad+ \frac{1}{2\pi i}\lim_{\ep\rightarrow 0}\int_0^{\frac{\delta}{2}} \e^{-imy_0t}\p_{y}\l(\psi_{m,\ep}^-(y,y_0)- \psi_{m,\ep}^+(y,y_0)\r) \d y_0 \\
&\quad+ \frac{1}{2\pi i}\frac{1}{imt}\lim_{\ep\rightarrow 0}\int_{\frac{\delta}{2}}^{y-\frac{\delta}{2}} \e^{-imy_0t}\p_{y,y_0}^2\l(\psi_{m,\ep}^-(y,y_0)- \psi_{m,\ep}^+(y,y_0)\r) \d y_0 \\
&=\mathcal{I}_{1,1}+ \mathcal{I}_{1,2}+\mathcal{I}_{1,3}.
\end{align*}
From Proposition \ref{L2 bounds inhom TG solution} we see that
\begin{equation*}
\Vert \cI_{1,1}\Vert_{L^2_y(\delta,1-\delta)}\lesssim \frac{1}{mt} m^{-\frac12}(m\delta)^{-\frac12-\mu}Q_{0,m}, \quad \Vert \cI_{1,2}\Vert_{L^2_y(\delta,1-\delta)}\lesssim m^{-\frac32}(m\delta)^{\frac12-\mu}Q_{0,m}.
\end{equation*}
Similarly, from Proposition \ref{prop: bounds dy0 varphi}, we obtain
\begin{equation*}
\Vert \cI_{1,3}\Vert_{L^2_y(\delta,1-\delta)}\lesssim \frac{1}{mt}\l(m^\frac12 + (m\delta)^{-\frac12-\mu}\r)Q_{1,m}.
\end{equation*}
The bounds for $\cI_3$ are the same as the ones for $\cI_1$, we omit the details. Recovering all terms, we conclude that
\begin{equation*}
\Vert \p_y\psi_m(t)  \Vert_{L^2_y(\delta,1-\delta)} \lesssim m^{-\frac32}(m\delta)^{\frac12-\mu}Q_{0,m} + \frac{1}{mt}\l(m^\frac12 + (m\delta)^{-\frac12-\mu}\r)Q_{1,m}.
\end{equation*}

\bullpar{Step 2} For $y\in(0,\delta)$, we shall split now
\begin{equation*}
\p_y\psi_m(t,y)=\frac{1}{2\pi i}\lim_{\ep\rightarrow0}\l( \int_0^{y+\frac{\delta}{2}} + \int_{y+\frac{\delta}{2}}^1 \r) e^{-imy_0t}\l( \p_{y}\psi_m^-(y,y_0) - \p_y\psi_m^+(y,y_0) \r) \d y_0 = \widetilde{\cI}_1 + \widetilde{\cI}_2.
\end{equation*}
As before, the bound for $\widetilde{\cI}_2$ is the same as the bound for $\cI_3$.  For $\widetilde{\cI}_1$, note that $|y-y_0|\leq \delta\leq \frac{3\b}{m}$ so that we shall use Proposition \ref{L2 bounds inhom TG solution} to find that
\begin{equation*}
\Vert \widetilde{\cI}_1 \Vert_{L^2_y(0,\delta)}\lesssim m^{-2}(m\delta)^{1-\mu}Q_{0,m}.
\end{equation*}
Gathering the previous bound, we obtain
\begin{equation*}
\Vert \p_y\psi_m(t)  \Vert_{L^2_y}\lesssim m^{-\frac32}(m\delta)^{\frac12-\mu}Q_{0,m} + \frac{1}{mt}\l(m^\frac12 + (m\delta)^{-\frac12-\mu}\r)Q_{1,m}.
\end{equation*}
As before, the conclusion follows for $\delta=\frac{c_0}{mt}$, with $c_0=\frac{1}{1000}\min\l(\b,1\r)$.
\end{proof}

We next obtain the decay rates for the perturbed density.
\begin{proposition}\label{prop: invdamp rho}
Let   $t\geq 1$. Then,
\begin{equation*}
\Vert \rho_m(t)\Vert_{L^2_y}\lesssim m^{-\frac12}t^{-\frac12 + \mu} Q_{1,m}.
\end{equation*}
\end{proposition}

\begin{proof}
The demonstration also follows the same strategy as the proof for Proposition \ref{prop: invdamp psi}, we just present the main ideas and bounds.

\bullpar{Step 1} For $y\in(\delta,1-\delta)$ we write 
\begin{align*}
\rho_m(t,y)&=\frac{1}{2\pi i}\lim_{\ep\rightarrow 0} \l( \int_0^{y-\frac{\delta}{2}} + \int_{y-\frac{\delta}{2}}^{y+\frac{\delta}{2}} + \int_{y+\frac{\delta}{2}}^1 \r) e^{-imy_0t}\l( \rho_{m,\ep}^-(y,y_0) - \rho_{m,\ep}^+(y,y_0)\r) \d y_0 \\
&=\cS_1 + \cS_2 + \cS_3
\end{align*}
We use Proposition \ref{prop: bounds rho and dy0 rho} to bound
\begin{equation*}
\Vert \cS_2\Vert_{L^2_y(\delta,1-\delta)}\lesssim m^{-\frac32}(m\delta)^{\frac12-\mu}.
\end{equation*}
As before, both the bounds for $\cS_3$ and $\cS_1$ and the manner to show them are the same, we just comment on $\cS_1$, which we split as follows. 
\begin{align*}
\cS_1 &= -\frac{1}{2\pi i} \frac{1}{imt}\lim_{\ep\rightarrow0} \l[ \e^{-imy_0t}\l(\rho_{m,\ep}^-(y,y_0) - \rho_{m,\ep}^+(y,y_0)\r) \r]_{y_0=\frac{\delta}{2}}^{y_0=y-\frac{\delta}{2}} \\
&\quad+ \frac{1}{2\pi i}\lim_{\ep\rightarrow 0}\int_0^{\frac{\delta}{2}} \e^{-imy_0t}\l(\rho_{m,\ep}^-(y,y_0)- \rho_{m,\ep}^+(y,y_0)\r) \d y_0 \\
&\quad+ \frac{1}{2\pi i}\frac{1}{imt}\lim_{\ep\rightarrow 0}\int_{\frac{\delta}{2}}^{y-\frac{\delta}{2}} \e^{-imy_0t}\p_{y_0}\l(\rho_{m,\ep}^-(y,y_0)- \rho_{m,\ep}^+(y,y_0)\r) \d y_0 \\
&=\mathcal{S}_{1,1}+ \mathcal{S}_{1,2}+\mathcal{S}_{1,3}.
\end{align*}
From Proposition \ref{prop: bounds rho and dy0 rho} we easily deduce
\begin{equation*}
\Vert \cS_{1,1} \Vert_{L^2_y(\delta,1-\delta)}\lesssim \frac{1}{mt} m^{-\frac12}(m\delta)^{-\frac{1}{2}-\mu}Q_{0,m}, \quad \Vert \cS_{1,2} \Vert_{L^2_y(\delta,1-\delta)}\lesssim m^{-\frac32}(m\delta)^{\frac12-\mu}Q_{0,m}.
\end{equation*}
On the other hand, Proposition \ref{prop: bounds rho and dy0 rho} also yields
\begin{equation*}
\Vert \cS_{1,3} \Vert_{L^2_y(\delta,1-\delta)}\lesssim \frac{1}{mt}\l(  m^{\frac12} + (m\delta)^{-\frac12-\mu}\r)Q_{1,m}.
\end{equation*}
Gathering the bounds we get
\begin{equation*}
\Vert \rho_m(t) \Vert_{L^2_y(\delta,1-\delta)}\lesssim m^{-\frac32}(m\delta)^{\frac12-\mu}Q_{0,m} + \frac{1}{mt}\l(  m^{\frac12} + (m\delta)^{-\frac12-\mu}\r)Q_{1,m}.
\end{equation*}

\bullpar{Step 2} For $y\in(0,\delta)$ we shall now consider 
\begin{equation*}
\rho_m(t,y)=\frac{1}{2\pi i}\lim_{\ep\rightarrow0}\l( \int_0^{y+\frac{\delta}{2}} + \int_{y+\frac{\delta}{2}}^1 \r) e^{-imy_0t}\l( \rho_m^-(y,y_0) - \rho_m^+(y,y_0) \r) \d y_0 = \widetilde{\cS}_1 + \widetilde{\cS}_2.
\end{equation*}
The bounds for $\widetilde{\cS_2}$ are the same as the ones for $\cS_3$, we just focus on $\widetilde{\cS_1}$. From Proposition \ref{prop: bounds rho and dy0 rho}, we see that
\begin{equation*}
\Vert \widetilde{\cS_1} \Vert_{L^2_y(0,\delta)}\lesssim m^{-2}(m\delta)^{1-\mu}Q_{0,m}.
\end{equation*}
With this, it follows that
\begin{equation*}
\Vert \rho_m(t)\Vert_{L^2_y(0,1)}\lesssim m^{-\frac32}(m\delta)^{\frac12-\mu}Q_{0,m} + \frac{1}{mt}\l(  m^{\frac12} + (m\delta)^{-\frac12-\mu}\r)Q_{1,m}
\end{equation*}
and thus the Proposition is proved once we choose $\delta=\frac{c_0}{mt}$, with $c_0=\frac{1}{1000}\min\l(\b,1\r)$.
\end{proof}

We next prove the inviscid damping decay estimates for the case $\b^2=1/4$. The precise bounds are recorded in the following proposition.

\begin{proposition}\label{prop: invdamp special}
Let   $t\geq 1$. Then,
\begin{equation*}
\begin{aligned}
\Vert \psi_m(t) \Vert_{L^2_y}&\lesssim m^{-\frac32}t^{-\frac32}(1 + \log(t) )\l(\Vert \rho_m^0 \Vert_{H^{4}_y} + \Vert \o_m^0 \Vert_{H^{4}_y}\r) \\
\Vert \p_y\psi_m(t)\Vert_{L^2_y}&\lesssim m^{-\frac12}t^{-\frac12}(1 + \log(t) )\l(\Vert \rho_m^0 \Vert_{H^{3}_y} + \Vert \o_m^0 \Vert_{H^{3}_y}\r) \\
\Vert \rho_m(t)\Vert_{L^2_y}&\lesssim m^{-\frac12}t^{-\frac12}(1 + \log(t) )\l(\Vert \rho_m^0 \Vert_{H^{3}_y} + \Vert \o_m^0 \Vert_{H^{3}_y}\r).
\end{aligned}
\end{equation*}
\end{proposition}

\begin{proof}
The proof follows along the same lines for the case $\b^2\neq 1/4$, the only difference is the logarithmic singularity present in the bounds of several quantities. For this, we note that for $m\delta< 1$, 
\begin{align*}
\int_{y-\frac{\delta}{2}}^{y+\frac{\delta}{2}}\frac{1}{m}|y-y_0|^{-\frac12} &\l( 1+ \big| \log \l(m|y-y_0|\r)\big| \r)\d y_0 \\
&\lesssim m^{-\frac12}\int_{y-\frac{\delta}{2}}^{y+\frac{\delta}{2}}(m|y-y_0|)^{-\frac12}\l( 1+ \big| \log \l(m|y-y_0|\r)\big| \r)\d y_0 \\
&\lesssim m^{-\frac32} \int_{-m\frac{\delta}{2}}^{m\frac{\delta}{2}} |\eta|^{-\frac12} (1 - \log |\eta|)\d \eta \\
&\lesssim m^{-\frac32}(m\delta)^{\frac12}\l( 1 + \big| \log \l(m\delta\r) \big| \r).
\end{align*}
Here, we have used that, for $0<m\delta\leq 1$, 
\begin{align*}
-\int_0^{m\delta} \eta^{-\frac12}\log(\eta) \d \eta= -\int_0^{m\delta} 2\p_\eta(\eta^\frac12)\log(\eta)\d \eta &= \l[-2\eta^\frac12\log(\eta)\r]_{\eta=0}^{\eta=m\delta}+ 2\int_0^{m\delta} \eta^{-\frac12}\d \eta \\
&\lesssim (m\delta)^\frac12\l( 1 + \big| \log \l(m\delta\r) \big|\r).
\end{align*}
The same argument also yields
\begin{equation*}
\int_{y-\frac{\delta}{2}}^{y+\frac{\delta}{2}}\frac{1}{m}|y-y_0|^{-\frac32} \l( 1+ \big| \log \l(m|y-y_0|\r)\big| \r)\d y_0\lesssim m^{-\frac12}(m\delta)^{-\frac12}\l( 1+ \big| \log \l(m\delta\r)\big| \r).
\end{equation*}
With this, the result follows thanks to the estimates obtained in Propositions \ref{prop: bounds dy0 varphi}-\ref{prop: bounds rho and dy0 rho}, we omit the details.
\end{proof}

Finally, Theorem \ref{thm:mainchan} is a direct consequence of Propositions \ref{prop: invdamp psi}-\ref{prop: invdamp special} together with Parseval identity.

\appendix

\section{Properties of the Whittaker functions}\label{app:Whittaker}
Here we state and prove some properties of the Whittaker functions that are used throughout the paper, we refer to \cite{NIST} for a complete description of the Whittaker functions. 

\subsection{Basic definitions and asymptotic expansions}
For $\gamma, \zeta\in \C$, the Whittaker function $M_{0,\gamma}(\zeta)$ is given by
\begin{equation*}
M_{0,\gamma}(\zeta)=e^{-\frac12\zeta}\zeta^{\frac12 + \gamma}M\l( \tfrac12 + \gamma, 1+2\gamma, \zeta\r), \quad M(a,b,\zeta) = \sum_{s=0}^\infty \frac{(a)_s}{(b)_s s!}\zeta^s,
\end{equation*}
where $(a)_s=a(a+1)(a+2)\dots (a+s-1)$ denotes the Pochhammer symbol. For $\gamma=0$, we also introduce 
\begin{equation*}
W_{0,0}(\zeta)=e^{-\frac12\zeta}\sqrt{\frac{\zeta}{\pi}}\sum_{s=0}^\infty \frac{\l( \tfrac12 \r)_s}{(s!)^2}\zeta^s \l( 2\frac{\Gamma'(1+s)}{\Gamma(1+s)} - \frac{\Gamma'(\tfrac12+s)}{\Gamma(\tfrac12+s)} - \log(\zeta)\r),
\end{equation*}
where $\Gamma(x)$ denotes the Gamma function.

We recall that $\mu=\Re\l(\sqrt{1/4-\b^2}\r)$ and $\nu=\Im\l(\sqrt{1/4-\b^2}\r)$, and set $\gamma=\mu+i\nu$. We begin by recording some basic properties regarding complex conjugation for $M_{0,\gamma}(\zeta)$, which can be deduce from the series definition of $M_{0,\gamma}$ and $W_{0,0}$.
\begin{lemma}
We have the following
\begin{itemize}
\item For $\b^2>1/4$, then $M_{0,i\nu}(\zeta)=\overline{M_{0,-i\nu}\l(\overline{\zeta}\r)}$.
\item For $\b^2\leq 1/4$, then $M_{0,\mu}(\zeta) = \overline{M_{0,\mu}\l(\overline{\zeta}\r)}$. Additionally, for $x\in\R$ then $M_{0,\mu}(x),W_{0,0}(x)\in \R$. 
\end{itemize}
\end{lemma}

We next state an analytic continuation property, which is key in studying the Wronskian of the Green's function and is directly determined by the analytic continuation of the nonentire term of $M_{0,\gamma}(\zeta)$, which is $\zeta^{\frac12+\gamma}$, for $\zeta\in\C$. 

\begin{lemma}[\cite{NIST}]\label{analytic continuation}
Let $\b^2>0$. Then
\begin{equation*}
M_{0,\gamma}(\zeta \e^{\pm \pi i}) = \pm i \e^{\pm \gamma \pi i}M_{0,\gamma}(\zeta), \quad W_{0,0}(\zeta e^{\pm i\pi})=\sqrt{\pi}M_{0,0}(\zeta) \pm iW_{0,0}(\zeta)
\end{equation*}
for all $\zeta\in\C$.
\end{lemma}

The next result gives a precise description of the asymptotic expansion of $M_\pm(\zeta)$ and its derivatives, for $\zeta$ in a bounded domain. 
\begin{lemma}\label{asymptotic expansion M}
Let $\zeta\in\C$. Let $B_R\subset\C$ denote the closed unit ball of radius $R>0$ centered in the origin. Then,
\begin{equation*}
\begin{aligned}
M_{0,\pm\gamma}(\zeta)&=\zeta^{\frac12 \pm\gamma}\mathcal{E}_{0,\pm\gamma}(\zeta),\quad M_{0,\pm\gamma}'(\zeta)=\zeta^{-\frac12 \pm\gamma}\mathcal{E}_{1,\pm\gamma}(\zeta),
\end{aligned}
\end{equation*}
where $\mathcal{E}_{j,\pm\gamma}\in L^\infty(B_R)$ and $\Vert \mathcal{E}_{j,\pm\gamma}\Vert_{L^\infty(B_R)} \lesssim_{\gamma,R} 1$, for all $j\in \l\lbrace 0,1,2 \r\rbrace$.

Moreover, for $R_m:=\frac{R}{2m}$ and $M_{\pm}(\zeta)=M_{0,\pm\gamma}(2m\zeta)$, let $B_{R_m}\subset\C$ denote the closed ball centered in the origin of radius $R_m$. We have that
\begin{equation*}
\begin{aligned}
M_\pm(\zeta)&=\zeta^{\frac12 \pm\gamma}\mathcal{E}_{m,0,\pm\gamma}(\zeta),\quad M_\pm'(\zeta)=\zeta^{-\frac12 \pm\gamma}\mathcal{E}_{m,1,\pm\gamma}(\zeta),
\end{aligned}
\end{equation*}
where $\mathcal{E}_{m,j,\pm\gamma}\in L^\infty(B_{R_m})$ and $\Vert \mathcal{E}_{m,j,\pm\gamma}\Vert_{L^\infty(B_{R_m})} \lesssim_\gamma (2m)^{\frac12\pm\mu}$, for all $j\in \l\lbrace 0,1,2 \r\rbrace$.
\end{lemma}
\begin{proof}
Firstly, from \cite{NIST} we know that
\begin{equation*}
\begin{aligned}
M_{0,\pm\gamma}(\zeta)&=\e^{-\frac12\zeta} \zeta^{\frac12 \pm \gamma} M\l(\frac12 \pm \gamma, 1\pm2\gamma,\zeta\r) =\zeta^{\frac12 \pm\gamma}\mathcal{E}_{0,\pm\gamma}(\zeta),
\end{aligned}
\end{equation*}
where $\mathcal{E}_{0,\pm\gamma}(\zeta)$ is entire and $\Vert \mathcal{E}_{0,\pm\gamma}\Vert_{L^\infty(B_R)} \lesssim_{\gamma,R} 1$. On the other hand, note that
\begin{equation}\label{eq W'}
M_{0,\pm\gamma}'(\zeta)=-\frac{1}{2}M_{0,\pm\gamma}(\zeta)+\l(\frac{1}{2}\pm\gamma\r)\frac{M_{0,\pm\gamma}(\zeta)}{\zeta}+\frac12\zeta^{-\frac12}M_{-\frac{1}{2},\frac{1}{2}\pm\gamma}(\zeta),
\end{equation}
where further
\begin{equation*}
M_{-\frac12,\frac12\pm\gamma}(\zeta)=\e^{-\frac12\zeta}\zeta^{1\pm\gamma}M\l(\frac32\pm\gamma,2\pm2\gamma,\zeta\r)=\zeta^{1\pm\gamma}\mathcal{H}_{\pm\gamma}(\zeta),
\end{equation*}
with $\mathcal{H}_{\pm\gamma}(\zeta)$ entire and thus uniformly bounded in $B_R$. Hence,  
\begin{equation*}
M_{0,\pm\gamma}'(\zeta)=\zeta^{-\frac12-\gamma}\l( \l(\frac12\pm\gamma\r) \mathcal{E}_{0,\pm\gamma}(\zeta)+ \frac12 \zeta \l(\mathcal{H}_{\pm\gamma}(\zeta) - \mathcal{E}_{0,\pm\gamma}(\zeta)\r)\r) = \zeta^{-\frac12-\gamma}\mathcal{E}_{1,\pm\gamma}(\zeta),
\end{equation*}
with $\Vert \mathcal{E}_{1,\pm\gamma}(\zeta) \Vert_{L^\infty(B_R)}\lesssim_{\gamma,R} 1$. The formulas and bounds for $M_\pm(\zeta)=M_{0,\pm\gamma}(2m\zeta)$ and its derivatives follow from those for $M_{0,\pm\gamma}$, the chain rule and the observation that $2m\zeta\in B_R$ provided that $\zeta\in B_{R_m}$.
\end{proof}

\begin{lemma}\label{asymptotic expansion W}
Let $\b^2=1/4$ and $\zeta\in\C$. Let $B_R\subset\C$ denote the closed ball of radius $R>0$ centered at the origin. Then,
\begin{equation*}
W_{0,0}(\zeta) = \zeta^\frac12 \big( \E_{0,1}(\zeta) - \log (\zeta) \E_{0,2}(\zeta)\big), \quad W_{0,0}'(\zeta) = \zeta^{-\frac12} \big( \E_{1,1}(\zeta) - \log (\zeta) \E_{1,2}(\zeta)\big),
\end{equation*}
where $\E_{j,k}(\zeta)$ are entire functions in $\C$ and $\Vert \E_{j,k}\Vert_{L^\infty(B_R)}\lesssim 1$, for $j=0,1$ and $k=1,2$.
\end{lemma}
\begin{proof}
We begin with noting that $W_{0,0}(2\zeta)=\sqrt{\frac{2\zeta}{\pi}}K_0(\zeta)$, where $K_0(\cdot)$ is the modified Bessel function of second kind of order 0. Moreover, we have that
\begin{equation*}
K_0(\zeta)=-\l( \ln \l(\frac{\zeta}{2}\r) + \varsigma\r) I_0(\zeta) +2\sum_{k=1}^\infty \frac{I_{2k}(\zeta)}{k}, 
\end{equation*} 
where
\begin{equation*}
I_{2k}(\zeta)=\l(\frac{\zeta}{2}\r)^{2k}\sum_{j\geq 0}\frac{\l( \frac{\zeta^2}{4}\r)^j}{j!(2k+j)!}.
\end{equation*}
Here, $I_j(\zeta)$ denotes the modified Bessel function of first kind of order $j\in \N$. In particular, one observes that $|I_{2k}(\zeta)|\leq I_{2k}(|\zeta|)$. Additionally, since $\cosh (|\zeta|) = I_0(|\zeta|) + 2\sum_{k=1}^\infty I_{2k}(|\zeta|)$, see \cite{NIST}, we can bound
\begin{equation*}
\l| 2\sum_{k=1}^\infty \frac{I_{2k}(z)}{k} \r| \leq 2 \sum_{k=1}^\infty I_{2k}(|\zeta|)= \cosh (|\zeta|) -I_0(|\zeta|)< \cosh(|\zeta|).
\end{equation*}
Therefore, since $I_j(\zeta)$ is analytic in $\C$ for all $j\in\N$ and $\frac12\zeta\in B_R$, we can write
\begin{equation*}
W_{0,0}(\zeta) = \zeta^\frac12 \big( \E_{0,1}(\zeta) - \log (\zeta) \E_{0,2}(\zeta)\big),
\end{equation*}
where 
\begin{equation*}
\E_{0,1}(\zeta)=\l( \log (2) - \varsigma \r) I_0(\zeta) + 2\sum_{k=1}^{+\infty}\frac{I_{2k}(\zeta)}{k}, \quad \E_{0,2}(\zeta) = I_0(\zeta) 
\end{equation*}
and they are such that $\Vert \E_{0,j}(\zeta)\Vert_{L^\infty(B_R)}\lesssim 1$, for $j=1,2$.

For $W_{0,0}'(\zeta)$, note that $W_{0,0}'(\zeta)=\frac{1}{2\sqrt{\pi\zeta}}\l( K_0(\zeta/2)+\zeta K_0'(\zeta/2)\r)$. As before, we can write
\begin{equation*}
K_0'(\zeta)=-K_1(\zeta)= -\l[ \frac{1}{\zeta}I_0(\zeta) + \l( \log \l(\frac12\zeta\r) + \varsigma -1 \r) I_1(\zeta) - \sum_{k\geq 1}\frac{(1+2k)}{k(1+k)}I_{1+2k}(\zeta)\r].
\end{equation*}
Since $\sinh (\zeta) = I_1(\zeta) + 2\sum_{k\geq 1}I_{1+2k}(\zeta)$, confer \cite{NIST}, we bound
\begin{equation*}
\l|\sum_{k\geq 1}\frac{(1+2k)}{k(1+k)}I_{1+2k}(\zeta)\r| \leq \sinh (|\zeta|)-I_1(|\zeta|)\leq \sinh(|\zeta|).
\end{equation*}
and we conclude the existence of two entire functions $\E_{1,1}(\zeta)$ and $\E_{1,2}(\zeta)$ such that $\Vert \E_{1,j}(\zeta)\Vert_{L^\infty(B_R)}\lesssim 1$, for $j=1,2$ and for which 
\begin{equation*}
W_{0,0}'(\zeta)=\zeta^{-\frac12}\big( \E_{1,1}(\zeta) -\log(\zeta)\E_{0,2}(\zeta)\big).
\end{equation*}
\end{proof}

\subsection{Lower bounds for Whittaker functions}
The next lemma shows that there are no zeroes of $M_+(x)$, for any $x\in(0,\infty)$.
\begin{lemma}\label{lower bounds M}
Let $x>0$. We have the following.
\begin{itemize}
\item For $\b^2\leq1/4$, then $M_{0,\mu}(x)$ is monotone increasing and 
\begin{equation*}
M_{0,\mu}(x)>x^{\frac12+\mu}, \quad M\l(\tfrac12 + \mu,1+2\mu,x \r)\geq e^{\frac12 x}.
\end{equation*}
\item For $\b^2>1/4$, then $|M_{0,i\nu}(x)|$ is monotone increasing and 
\begin{equation*}
x|\Gamma(1+i\nu)|^2\frac{\sinh(\nu\pi)}{\nu\pi}\leq |M_{0,i\nu}(x)|^2 \leq x\cosh(x)|\Gamma(1+i\nu)|^2\frac{\sinh(\nu\pi)}{\nu\pi}, 
\end{equation*}
with also
\begin{equation*}
\l| M\l( \tfrac12 +i\nu, 1+ 2i\nu, x \r) \r| \geq \e^{\frac12x}|\Gamma(1+i\nu)|\sqrt{\frac{\sinh(\nu\pi)}{\nu\pi}}.
\end{equation*}
\end{itemize}
\end{lemma}

\begin{proof}
From \cite{NIST}, we have
\begin{equation*}
M_{0,\gamma}\left(2x\right)=2^{2\gamma+\frac{1}{2}}\Gamma\left(1+\gamma\right)\sqrt{x}I_{\gamma}\left(x\right).
\end{equation*}
For $\b^2\leq1/4$, we have $\gamma=\mu$ and the conclusion is straightforward, since we can use the power series representation of $I_\mu(x)$ to obtain
\begin{equation*}
M_{0,\mu}(2x)>(2x)^{\frac12+\mu}.
\end{equation*}
On the other hand, for $\b^2> 1/4$, we have $\gamma=i\nu$ and 
\begin{equation*}
M_{0,i\nu}\left(2x\right)=2^{2i\nu+\frac{1}{2}}\Gamma\left(1+i\nu\right)\sqrt{x}I_{i\nu}\left(x\right).
\end{equation*}
Therefore,
\begin{equation*}
\begin{aligned}
|M_{0,i\nu}(2x)|^2 &= 2x|\Gamma(1+i\nu)|^2I_{i\nu}(x)I_{-i\nu}(x)=2x|\Gamma(1+i\nu)|^2 \frac{2}{\pi}\int_0^{\frac{\pi}{2}} I_0(2x\cos\theta)\cosh(2\nu\theta)\d\theta
\end{aligned}
\end{equation*}
The upper and lower bound follow from the fact that $1 \leq I_0(x) \leq \cosh(x)$, for all $x\geq 0$. See \cite{NIST} for the product formula for $I_{i\nu}(x)I_{-i\nu}(x)$.
\end{proof}

\subsection{Growth bounds and comparison estimates for $\b^2>1/4$}
In this subsection we treat the case $\b^2>1/4$, so that $\mu=0$ and $\nu=\sqrt{\b^2-1/4}$.
\begin{lemma}\label{growth bounds M large argument}
Denote $a:=\frac12 + i\nu$ and $b:=2a$. Then, there exists $C>0$ and $N_{\nu,0}>0$ such that 
\begin{equation*}
\e^{-\frac18\nu\pi}\e^{\frac12\Re\zeta}\leq\l|\frac{\Gamma(a)}{\Gamma(b)}M_+(\zeta)\r|\leq \e^{\frac18\nu\pi}\e^{\frac12\Re\zeta},
\end{equation*} 
and 
\begin{equation*}
\l| \frac{M_\pm'(\zeta)}{M_\pm(\zeta)}-\frac12\r| \leq \frac14, 
\end{equation*}
for all $\Re\zeta \geq N_{\nu,0}$.
\end{lemma}
\begin{proof}
Let $\zeta\in\C$. We recall that
\begin{equation*}
\begin{aligned}
M_+(\zeta)&=\e^{-\frac12\zeta}\zeta^aM(a,b,\zeta) =\e^{-\frac12\zeta}\zeta^a\frac{\Gamma(b)}{\Gamma(a)}\l( \e^{-i\pi a}U(a,b,\zeta) + \e^{i\pi a}\e^\zeta U(a,b,\e^{i\pi} \zeta)\r).
\end{aligned}
\end{equation*}
Moreover, we have that $U(a,b,\zeta) = \zeta^{-a}+ \mathcal{E}_1(\zeta)$, where further
\begin{equation*}
|\zeta^a\E_1(\zeta)|\leq \frac{2\b^2}{|\zeta|}\e^{\frac{2\b^2}{|\zeta|}}.
\end{equation*}
In the sequel, we write $x:=\frac{2\b^2}{|\zeta|}$. Therefore, we can write
\begin{equation}\label{M large argument}
M_+(\zeta)=\frac{\Gamma(b)}{\Gamma(a)}\e^{\frac12\zeta}\l(  \l[1 + (\e^{i\pi}\zeta)^a\E_1(\e^{i\pi}\zeta)\r] + \e^{-\zeta}\e^{-i\pi a}\l[1+ \zeta^a\E_1(\zeta)\r]\r)
\end{equation}
We shall focus on obtaining upper and lower bound estimates for 
\begin{equation*}
\l|  \l[1 + (\e^{i\pi}\zeta)^a\E_1(\e^{i\pi}\zeta)\r] + \e^{-\zeta}\e^{-i\pi a}\l[1+ \zeta^a\E_1(\zeta)\r]\r|
\end{equation*}
when $\Re\zeta$ is large. To this end, we note that $|\zeta^a\E_1(\zeta)|\leq 2x$, for $x\leq \frac12$. Moreover,
\begin{equation*}
|1+\zeta^a\E_1(\zeta)| \leq \e^{\nu\frac{\pi}{16}},
\end{equation*}
provided that $x\leq \frac12\min\l\lbrace 1, \e^{\nu\frac{\pi}{16}}-1\r\rbrace$. Similarly, we also have that 
\begin{equation*}
1+\e^{-\zeta} \e^{\nu\pi}\leq \e^{\nu\frac{\pi}{16}},
\end{equation*}
for all $\Re \zeta > \nu\pi-\log (\e^{\frac{\nu\pi}{16}}-1)$. Hence,
\begin{equation*}
\l|  \l[1 + (\e^{i\pi}\zeta)^a\E_1(\e^{i\pi}\zeta)\r] + \e^{-\zeta}\e^{-i\pi a}\l[1+ \zeta^a\E_1(\zeta)\r]\r| \leq \e^{\frac{\nu\pi}{8}}.
\end{equation*}
On the other hand, for $x\leq \min\l\lbrace \frac{1}{4}\l( 1- \e^{-\frac1{8}\nu\pi}\r), \frac12\r\rbrace$, we have that 
\begin{equation*}
|\zeta^a\E_1(\zeta)| \leq \frac12\l(1-\e^{-\nu\frac{\pi}{8}}\r),
\end{equation*}
and also
\begin{equation*}
\l| \e^{-\zeta}\e^{\nu\pi}\l(1+ \zeta^a\E_1(\zeta)\r)\r|\leq \frac12\l(1-\e^{-\nu\frac{\pi}{8}}\r),
\end{equation*}
provided that $\Re\zeta>\nu\pi +\log\l(\frac{4}{1-\e^{-\frac18\nu\pi}}\r)$. Therefore, we can lower bound
\begin{equation*}
\l|  \l[1 + (\e^{i\pi}\zeta)^a\E_1(\e^{i\pi}\zeta)\r] + \e^{-\zeta}\e^{-i\pi a}\l[1+ \zeta^a\E_1(\zeta)\r]\r| \geq \e^{-\frac{\nu\pi}{8}}.
\end{equation*}
We choose $N_\nu>0$ so that all the above conditions are satisfied when $\Re\zeta \geq N_\nu$. For the second part of the lemma, we take a $\frac{\d}{\d \zeta}$ derivative in \eqref{M large argument} to obtain
\begin{equation*}
\begin{aligned}
M_+'(\zeta)&=\frac12M_+(\zeta) + \frac{\Gamma(b)}{\Gamma(a)}\e^{\frac12\zeta}\l(\frac{a}{\zeta}(\e^{i\pi}\zeta)^a\E_1(\e^{i\pi}\zeta) + \e^{i\pi}(\e^{i\pi}\zeta)^a\E_1'(\e^{i\pi}\zeta)\r) \\ 
&\quad + \frac{\Gamma(b)}{\Gamma(a)}\e^{-\frac12\zeta}\e^{i\pi a}\l( \frac{a}{\zeta}\zeta^a\E_1(\zeta) + \zeta^a\E_1'\zeta) -1 -\zeta^a\E_1(\zeta)\r).
\end{aligned}
\end{equation*}
Since $|\zeta^a\E_1(\zeta)|\leq \frac{2\b^2}{|\zeta|}\e^{\frac{2\b^2}{|\zeta|}}$ and $|\zeta^a\E_1'(\zeta)|\leq \frac{4\b^2}{|\zeta|}\e^{\frac{2\b^2}{|\zeta|}}$, confer \cite{NIST},  we find that
\begin{equation*}
\l|\frac{a}{\zeta}(\e^{i\pi}\zeta)^a\E_1(\e^{i\pi}\zeta) + \e^{i\pi}(\e^{i\pi}\zeta)^a\E_1'(\e^{i\pi}\zeta) \r| \leq \l(\frac{|a|}{|\zeta|}+2\r)\frac{2\b^2}{|\zeta|}\e^{\frac{2\b^2}{|\zeta|}}\leq 6x,
\end{equation*}
and
\begin{equation*}
\l|\e^{i\pi a}\l( \frac{a}{\zeta}\zeta^a\E_1(\zeta) + \zeta^a\E_1'\zeta) -1 -\zeta^a\E_1(\zeta)\r)\r|\leq \l(1 + \l(\frac{|a|}{|\zeta|}+3\r)\frac{2\b^2}{|\zeta|}\e^{\frac{2\b^2}{|\zeta|}} \r)\leq 5,
\end{equation*}
for $|\zeta|\geq |a|$ and $x\leq \frac12$. Therefore, 
\begin{equation*}
\begin{aligned}
\l|\frac{M_+'(\zeta)}{M_+(\zeta)}-\frac12\r| &\leq \l|\frac{\frac{a}{\zeta}(\e^{i\pi}\zeta)^a\E_1(\e^{i\pi}\zeta) + \e^{i\pi}(\e^{i\pi}\zeta)^a\E_1'(\e^{i\pi}\zeta)}{1 + (\e^{i\pi}\zeta)^a\E_1(\e^{i\pi}\zeta) + \e^{-\zeta}\e^{-i\pi a}\l(1+ \zeta^a\E_1(\zeta)\r)}\r| \\
&+ \e^{-\Re(\zeta)}\l|\frac{\e^{i\pi a}\l( \frac{a}{\zeta}\zeta^a\E_1(\zeta) + \zeta^a\E_1'\zeta) -1 -\zeta^a\E_1(\zeta)\r)}{1 + (\e^{i\pi}\zeta)^a\E_1(\e^{i\pi}\zeta) + \e^{-\zeta}\e^{-i\pi a}\l(1+ \zeta^a\E_1(\zeta)\r)}\r|,
\end{aligned}
\end{equation*}
which can be made arbitrarily small due to the previous bounds for $\Re(\zeta)$ sufficiently large.
\end{proof}

\begin{lemma}\label{Comparison bounds M small argument}
Let $y_0\in[0,1]$ such that $2my_0\leq N_{\nu,0}$. Then, there exists $\ep_0>0$ such that 
\begin{equation*}
\l| \frac{M_+(y_0-i\ep)}{M_+(y_0+i\ep)}\r| \leq \e^{\frac54\nu\pi},
\end{equation*}
for all $\ep\leq\ep_0$.
\end{lemma}
\begin{proof}
Let $\theta=\arg(y_0-i\ep)\in \l[ -\frac{\pi}{2}, 0\r]$. Recall that for $a=\frac12+i\nu$ and $b=2a$, for $\zeta\in\C$,
\begin{equation*}
\begin{aligned}
M_+(\zeta)&=\e^{-\frac12\zeta}\zeta^aM(a,b,\zeta).
\end{aligned}
\end{equation*} 
Therefore, we can estimate
\begin{equation*}
\l| \frac{M_+(y_0-i\ep)}{M_+(y_0+i\ep)}\r| = \l| \e^{-i\ep}\e^{i\theta}\e^{-2\nu\theta} \frac{M(a,b,2m(y_0-i\ep))}{M(a,b,2m(y_0+i\ep))}\r|\leq \e^{\nu\pi}\l|\frac{M(a,b,2m(y_0-i\ep))}{M(a,b,2m(y_0+i\ep))}\r|.
\end{equation*}
Now, since $\frac{\d}{\d\zeta}M(a,b,\zeta)=\frac12M(a+1,b+1,\zeta)$, which is entire in $\zeta\in\C$, we have that
\begin{equation*}
M(a,b,2m(y_0+i\ep))=M(a,b,2my_0) + \int_{0}^{\ep}imM(a+1,b+1,2m(y_0+is))\d s.
\end{equation*}
We can further bound the error term by noting that $|2m(y_0+is)|\leq N_\nu+10\b^2$, for all $|s|\leq |\ep|$. As a result, there exists $C_\nu$ such that $|M(a+1,b+1,2m(y_0+is))|\leq C_\nu$, for all $|s|\leq |\ep|$. Therefore,
\begin{equation*}
\begin{aligned}
\l|\int_{0}^{\ep}imM(a+1,b+1,2m(y_0+is))\d s\r| &\leq |M(a,b,2my_0)|\frac{C_\nu m|\ep|}{|M(a,b,2my_0)|} \\
&\leq C_\nu|M(a,b,2my_0)|\e^{-my_0}\sqrt{\frac{\nu\pi\cosh\nu\pi}{\sinh\nu\pi}}m|\ep| \\
&\leq (1-\e^{-\frac18\nu\pi})|M(a,b,2my_0)|,
\end{aligned}
\end{equation*}
for all $0\leq|\ep|\leq \ep_0=\frac{1-\e^{-\frac18\nu\pi}}{m}\sqrt{\frac{\sinh\nu\pi}{\nu\pi\cosh\nu\pi}} C_\nu$. Consequently, we have that 
\begin{equation*}
\l| \frac{M(a,b,2m(y_0-i\ep))}{M(a,b,2m(y_0+i\ep))}\r| \leq \e^{\frac14\nu\pi}.
\end{equation*}
\end{proof}
\begin{lemma}\label{Comparison bounds M order one argument}
Let $N_{\nu,0}$ be given as above and $N_{\nu,1}>0$. Let $\sigma\in\lbrace +,-\rbrace$. If $N_{\nu,1}< N_{\nu,0}$, then, there exists $\ep_0>0$ such that
\begin{equation*}
|M_\sigma(y_0)|\e^{-\frac18\nu\pi} \leq |M_\sigma(y_0+i\ep)|\leq |M_\sigma(y_0)|\e^{\frac18\nu\pi}, 
\end{equation*} 
for all $y_0\in[0,1]$ such that $N_{\nu,1}\leq 2my_0 \leq N_{\nu,0}$, and all $0\leq |\ep|\leq \ep_0$.
\end{lemma}
\begin{proof}
The result follows from the Fundamental Theorem of Calculus, the asymptotic expansions of $M_\sigma$ and $M_\sigma'$ for small arguments from Lemma tal and the lower bounds on $|M_\sigma|$ from Lemma Qual. More precisely, assume without loss of generality that $0\leq \ep$ and note that
\begin{equation*}
M_\sigma(y_0+i\ep)=M_\sigma(y_0) + \int_0^\ep \frac{\d}{\d s}M_\sigma(y_0+is)\d s.
\end{equation*}
Thanks to the asymptotic expansions for small arguments we next estimate
\begin{equation*}
\l| \int_0^\ep \frac{\d}{\d s}M_\sigma(y_0+is)\d s \r| \leq \int_0^\ep |M_\sigma'(y_0+is)| \d s \leq C_\nu\int_0^\ep \frac{(2m)^\frac12}{|y_0+is|^\frac12}\d s \leq C_\nu(2m\ep)(2my_0)^{-\frac12}.
\end{equation*}
Using the lower bound $(2my_0)^{\frac12}\leq \sqrt{\frac{\nu\pi\cosh\nu\pi}{\sinh\nu\pi}}|M_\sigma(y_0)|$ we have that
\begin{equation*}
\l| \int_0^\ep \frac{\d}{\d s}M_\sigma(y_0+is)\d s \r| \leq C_\nu|M_\sigma(y_0)|\frac{\ep}{y_0}\leq C_\nu N_{\nu,2}^{-1}2m\ep.
\end{equation*}
We now choose $\ep_0=\frac{N_{\nu,2}}{2mC_\nu}(1-\e^{-\frac18\nu\pi})$. The conclusion of the lemma follows swiftly for all $\ep\leq \ep_0$.
\end{proof}

\begin{lemma}\label{lemma:lower asymptotic bounds M}
Let $y_0\in[0,1]$ and $0\leq \ep \leq \frac{\b}{m}$. Then, 
\begin{itemize}
\item If $my_0\leq 3\b$, there exists $\ep_0>0$ such that
\begin{equation*}
\l(m|y_0+ i\ep|\r)^\frac12\lesssim |M_\pm(y_0+ i\ep)|
\end{equation*}
for all $0\leq \ep\leq \ep_0$.
\item If $my_0\geq 3\b$, 
\begin{equation*}
1\lesssim {|M_\pm(y_0+i\ep)|}.
\end{equation*}
\end{itemize}
\end{lemma}
\begin{proof}
For the first part of the Lemma, recall $M_+(\zeta)=e^{-\frac12\zeta}\zeta^{a}M\l(a, b,\zeta\r)$, the lemma follows once we obtain lower bounds on $e^{-\frac12\zeta}M\l(a,b,\zeta\r)$. For this, note that since $\frac{\d}{\d\zeta}M(a,b,\zeta)=\frac12M(a+1,b+1,\zeta)$, which is entire in $\zeta\in\C$, we have that
\begin{equation*}
M(a,b,2m(y_0+i\ep))=M(a,b,2my_0) + \int_{0}^{\ep}imM(a+1,b+1,2m(y_0+is))\d s.
\end{equation*}
We further bound the error term by noting that $|2m(y_0+is)|\leq 10\b$, for all $|s|\leq |\ep|$. As a result, there exists $C>0$ such that $|M(a+1,b+1,2m(y_0+is))|\leq C$, for all $|s|\leq |\ep|$. Therefore, using the lower bounds on $|M(a,b,2my_0)|$ from Lemma \ref{lower bounds M},
\begin{equation*}
\begin{aligned}
\l|\int_{0}^{\ep}imM(a+1,b+1,2m(y_0+is))\d s\r| &\leq |M(a,b,2my_0)|\frac{C m|\ep|}{|M(a,b,2my_0)|} \\
&\leq C|M(a,b,2my_0)||\Gamma(1+i\nu)|\sqrt{\frac{\nu\pi}{\sinh\nu\pi}}m|\ep|. \\
\end{aligned}
\end{equation*}
In particular, there exists $\ep_0>0$ such that for all $0\leq\ep\leq \ep_0$,
\begin{equation*}
\begin{aligned}
e^{-my_0}|M(a,b,2m(y_0+i\ep))|&\geq e^{-my_0}|M(a,b,2my_0)|\l( 1-Cm\frac{\ep}{|M(a,b,2my_0)|}\r) \\
&\geq \frac12\frac{1}{|\Gamma(1+i\nu)|}\sqrt{\frac{\sinh(\nu\pi)}{\nu\pi}},
\end{aligned}
\end{equation*}
and the first part of the lemma follows. As for the second statement, it is a direct consequence of Lemma \ref{growth bounds M large argument} and the fact that $|M_\pm(\cdot)|$ is bounded in compact domains (it is entire).
\end{proof}

\subsection{Growth bounds and comparison estimates for $\b^2=1/4$}
\begin{lemma}\label{growth bounds special M large argument}
Let $\b^2 = 1/4$ and let $\mu=\sqrt{1/4-\b^2}$. Denote $a:=\frac12$ and $b:=2a=1$. Then, there exists $N_0>0$ such that 
\begin{equation*}
\l|\frac{W_{0,0}(\zeta)}{M_{0,0}(\zeta)}\r|\leq 2\sqrt{\pi}\e^{-\Re\zeta}, \quad \l| \frac{W_{0,0}'(\zeta)}{W_{0,0}(\zeta)} + \frac12 \r| \leq \frac14, \quad \frac12e^{-\frac12\Re\zeta}\leq |W_{0,0}(\zeta)| \leq \frac32\e^{-\frac12\Re\zeta}
\end{equation*} 
for all $\Re\zeta \geq N_0$.
\end{lemma}
\begin{proof}
Let $\zeta\in\C$ and $\delta>0$. We recall that
\begin{equation*}
\begin{aligned}
M_{0,0}(\zeta)&=\e^{-\frac12\zeta}\zeta^{\frac12} M(1/2,1,\zeta) \\
&=\e^{-\frac12\zeta}\zeta^{1/2}\frac{\Gamma(1)}{\Gamma(1/2)}\l( -iU(1/2,1,\zeta) + i\e^\zeta U(1/2,1,\e^{i\pi} \zeta)\r),
\end{aligned}
\end{equation*}
while also
\begin{equation*}
W_{0,0}(\zeta)=\e^{-\frac12\zeta}\zeta^\frac12U(1/2,1,\zeta).
\end{equation*}
Thus, we have that
\begin{equation*}
\frac{W_{0,0}(\zeta)}{M_{0,0}(\zeta)}=-i\sqrt{\pi}\frac{U(1/2,1,\zeta)}{\e^\zeta U(1/2,1,\e^{i\pi}\zeta) - U(1/2,1,\zeta)}=-i\sqrt{\pi}\frac{1}{\e^\zeta\frac{U(1/2,1,\e^{i\pi}\zeta)}{U(1/2,1,\zeta)}-1}.
\end{equation*}
Now, we also recall that $U(1/2,1,\zeta)=\zeta^{-\frac12}\l( 1 + \zeta^\frac12\E_1(\zeta)\r)$, with $|\zeta^\frac12\E_1(\zeta)|\leq \frac{1}{2|\zeta|}\e^{\frac{1}{2|\zeta|}}$. Therefore, we have the lower bound
\begin{equation*}
\l| \frac{U(1/2,1,\e^{i\pi}\zeta)}{U(1/2,1,\zeta)} \r| \geq \frac34,
\end{equation*}
for $|\zeta|$ sufficiently large. Moreover, $\frac34 \e^{\Re\zeta}-1\geq \frac12 \e^{\Re\zeta}$, for all $\Re\zeta \geq 2\ln 2$. The desired conclusion follows.

For the second part of the Lemma, since $W_{0,0}(\zeta)=e^{-\frac12\zeta}\zeta^\frac12\l( \zeta^{-\frac12} + \E_1(\zeta)\r)$, we note that
\begin{equation*}
W_{0,0}'(\zeta)= -\frac12 W_0(\zeta) +\frac12\frac{W_0(\zeta)}{\zeta} + e^{-\frac12\zeta}\l(-\frac{1}{2\zeta} + \zeta^\frac12\E_1'(\zeta)\r),
\end{equation*}
where we recall that $\l| \zeta^\frac12\E_1'(\zeta)\r|\leq \frac{1}{4|\zeta|}e^{\frac{1}{2|\zeta|}}\leq x$, for $x=\frac{1}{2|\zeta|}\leq \frac12$. Hence,
\begin{equation*}
\l| \frac{W_{0,0}'(\zeta)}{W_{0,0}(\zeta)} + \frac12\r| \leq x + 2x\l|\frac{1}{1+\zeta^\frac12\E_1(\zeta)}\r|\leq 2x\leq \frac14,
\end{equation*}
for all $x\leq\frac18$.

For the third statement of the Lemma, note that $W_{0,0}(\zeta)=e^{-\frac12\zeta}\l(1+ \zeta^\frac12\E_1(\zeta)\r)$, the conclusion follows for $|\zeta|$ large enough so that $\l| \zeta^\frac12\E_1(\zeta) \r| \leq \frac12$.
\end{proof}

\begin{lemma}\label{Comparison bounds special M small argument}
Let $\b^2 = 1/4$. Denote $a:=\frac12 \pm \mu$ and $b:=2a=1$. Then, for all $\epsilon>0$ there exists $\delta_0>0$ such that 
\begin{equation*}
\l| \frac{M_{0,0}(\zeta)}{W_{0,0}(\zeta)}\r| \leq \epsilon,
\end{equation*}
for all $\zeta\in\C$ such that $|\zeta|\leq \delta$.
\end{lemma}

\begin{proof}
We use the functional relation between the Whittaker functions and the modified Bessel functions in order to extract the correct asymptotic behaviour of the functions near the origin and estimate the quotient precisely. In this direction, recall that
\begin{equation*}
M_{0,0}(2\zeta)=\sqrt{2\zeta}I_0(\zeta), \quad W_{0,0}(2\zeta) = \sqrt{\frac{2\zeta}{\pi}}K_0(\zeta),
\end{equation*}
where $I_0(\zeta)$ and $K_0(\zeta)$ denote the modified Bessel functions of order 0. Moreover, we have that
\begin{equation*}
K_0(\zeta)=-\l( \ln \l(\frac{\zeta}{2}\r) + \varsigma\r) I_0(\zeta) +2\sum_{k=1}^\infty \frac{I_{2k}(\zeta)}{k}, 
\end{equation*} 
where
\begin{equation*}
I_{2k}(\zeta)=\l(\frac{\zeta}{2}\r)^{2k}\sum_{j\geq 0}\frac{\l( \frac{\zeta^2}{4}\r)^j}{j!(2k+j)!}.
\end{equation*}
In particular, one observes that $|I_{2k}(\zeta)|\leq I_{2k}(|\zeta|)$. Moreover, under the observation that $(2k+j)!\geq (2k)!j!$, we can bound
\begin{equation*}
\l| 2\sum_{k=1}^\infty \frac{I_{2k}(z)}{k} \r| \leq 2 \sum_{k=1}^\infty I_0(|\zeta|)\frac{\l( \frac12 |\zeta|\r)^{2k}}{k(2k)!}\leq 2I_0(|\zeta|)\l( \cosh \frac12|\zeta| -1\r).
\end{equation*}
With this, together with the fact that $I_0(\cdot)$ is analytic in $\C$ and  $I_0(\zeta)\rightarrow 1$ when $\zeta\rightarrow 0$, we have that
\begin{equation*}
|K_0(\zeta)|\geq -\frac{1}{2}\ln \l(\frac12|\zeta|\r),
\end{equation*}
for $|\zeta|$ sufficiently small. The conclusion follows, since for $|\zeta|$ sufficiently small we have
\begin{equation*}
\l| \frac{M_{0,0}(\zeta)}{W_{0,0}(\zeta)}\r|\leq -\frac{3}{\ln \l( \frac14 |\zeta| \r)}
\end{equation*}
\end{proof}

\begin{lemma}\label{Comparison bounds special M order one argument}
Let $\b^2=1/4$ and let $y_0\in[0,1]$ such that $N_2\leq 2my_0\leq N_1$. Then, for all $\epsilon>0$ there exists $\ep_0>$ such that
\begin{equation*}
\l| \frac{W_0(y_0-i\ep)}{M_0(y_0-i\ep)}-\frac{W_0(y_0)}{M_0(y_0)}\r|\leq \epsilon, \quad \l| \frac{W_0(y_0-i\ep)}{M_0(y_0-i\ep)}\r|\leq C
\end{equation*}
for all $\ep\leq \ep_0$ and some $C>0$. In particular,
\begin{equation*}
\l| \Im \l(\frac{W_0(y_0-i\ep)}{M_0(y_0-i\ep)}\r)\r|\leq \epsilon,
\end{equation*}
\end{lemma}
\begin{proof}
It follows from the continuity of the functions involved, plus the fact that $M_{0,0}(x)$ does not vanish and $W_{0,0}(x)$ is bounded, for any $x>0$ such that $0<N_2\leq x \leq N_1< \infty$.
\end{proof}

\begin{lemma}\label{lemma: lower asymptotic bounds special W}
There exists $\delta_2>0$ such that
\begin{equation*}
|\zeta|^\frac12 \l( 1 + \big| \log |\zeta| \big| \r) \lesssim |W_{0,0}(\zeta)|,
\end{equation*}
for all $|\zeta|\leq \delta_2$.
\end{lemma}

\begin{proof}
Recall that $W_{0,0}(\zeta)= \sqrt{\frac{\zeta}{\pi}}K_0(\zeta/2)$ and the fact that $|K_0(\zeta)|\geq -\frac{1}{2}\log \l( \frac{|\zeta|}{2}\r)$ for $\zeta\rightarrow 0$. Then,
\begin{equation*}
|W_{0,0}(\zeta)|\geq \frac{1}{2\sqrt{\pi}}|\zeta|^\frac12 \big| \log |\zeta| - \log 4 \big|\geq \frac{1}{20\sqrt{\pi}}|\zeta|^\frac12\big | \log |\zeta| \big| \geq \frac{1}{40\sqrt{\pi}}|\zeta|^\frac12\l( 1 + \big| \log |\zeta| \big| \r),
\end{equation*}
for $|\zeta|$ sufficiently small.
\end{proof}

\begin{lemma}\label{lemma:lower asymptotic bounds special M}
Let $y_0\in[0,1]$ and $0\leq \ep \leq \frac{\b}{m}$. Then, 
\begin{itemize}
\item If $my_0\leq 3\b$, there exists $\ep_0>0$ such that
\begin{equation*}
\l(m|y_0+ i\ep|\r)^{\frac12}\lesssim |M_0(y_0+ i\ep)|
\end{equation*}
for all $0\leq \ep\leq \ep_0$.
\item If $my_0\geq 3\b$, 
\begin{equation*}
1\lesssim {|M_\pm(y_0+i\ep)|}.
\end{equation*}
\end{itemize}
\end{lemma}

\begin{proof}
The proof uses the ideas from Lemma \ref{lemma:lower asymptotic bounds M} together with the bounds from Lemma \ref{growth bounds real M large argument}. We omit the details.
\end{proof}

\subsection{Growth bounds and comparison estimates for $\b^2<1/4$}
In this subsection we consider the case $\b^2 < 1/4$, for which $\mu=\sqrt{1/4-\b^2}$ with $\mu\in \l( 0,\frac12\r)$ and $\nu=0$.
\begin{lemma}\label{growth bounds real M large argument}
Denote $a_\pm:=\frac12 \pm \mu$ and $b_\pm:=2a_\pm$. Then, 
\begin{equation*}
\lim_{\Re(\zeta)\rightarrow +\infty}\l|\frac{\Gamma(a)}{\Gamma(b)}\e^{-\frac12\zeta}M_+(\zeta)\r|=1.
\end{equation*} 
Moreover, let $C_\mu=2^{-4\mu}\frac{\Gamma(1-\mu)}{\Gamma(1+\mu)}$. There exists $N_{\mu,0}>0$ such that 
\begin{equation*}
\l| \frac{M_-(\zeta)}{M_+(\zeta)} - 2^{-4\mu}\frac{\Gamma(1-\mu)}{\Gamma(1+\mu)} \r| \leq \min\l( 5, \frac{\tan\mu\pi}{1+\tan\mu\pi}\r)\frac{C_\mu}{4}.
\end{equation*}
for all $\Re\zeta \geq N_{\mu,0}$.
\end{lemma}
\begin{proof}
Let $\zeta\in\C$ and $\delta>0$. We recall that
\begin{equation*}
\begin{aligned}
M_\pm(\zeta)&=\e^{-\frac12\zeta}\zeta^{a_\pm} M(a_\pm,b_\pm,\zeta) \\
&=\e^{-\frac12\zeta}\zeta^{a_\pm}\frac{\Gamma(b_\pm)}{\Gamma(a_\pm)}\l( \e^{-i\pi a_\pm}U(a_\pm,b_\pm,\zeta) + \e^{i\pi a_\pm}\e^\zeta U(a_\pm,b_\pm,\e^{i\pi} \zeta)\r).
\end{aligned}
\end{equation*}
Moreover, we have that $U(a_\pm,b_\pm,\zeta) = \zeta^{-a_\pm}+ \mathcal{E}_{\pm}(\zeta)$, where further
\begin{equation*}
|\zeta^{a_\pm}\E_\pm(\zeta)|\leq \frac{2\b^2}{|\zeta|}\e^{\frac{2\b^2}{|\zeta|}}.
\end{equation*}
In the sequel, we write $x:=\frac{2\b^2}{|\zeta|}$. Therefore, we can write
\begin{equation*}
M_\pm(\zeta)=\frac{\Gamma(b_\pm)}{\Gamma(a_\pm)}\e^{\frac12\zeta}\l(  \l[1 + (\e^{i\pi}\zeta)^{a_\pm}\E_\pm(\e^{i\pi}\zeta)\r] + \e^{-\zeta}\e^{-i\pi a_\pm}\l[1+ \zeta^{a_\pm}\E_\pm(\zeta)\r]\r)
\end{equation*}
Now, since $b_\pm=2a_\pm$, we have that
\begin{equation*}
\frac{\Gamma(b_\pm)}{\Gamma(a_\pm)}=\pi^{-\frac12}2^{2a_\pm-1}\Gamma\l( a_\pm +\frac12\r),
\end{equation*}
confer, \cite{NIST}. Therefore, 
\begin{equation*}
\frac{\displaystyle {\Gamma(b_-)}/{\Gamma(a_-)}}{\displaystyle {\Gamma(b_+)}/{\Gamma(a_+)}}=2^{-4\mu}\frac{\Gamma(1-\mu)}{\Gamma(1+\mu)}
\end{equation*}
and we note that
\begin{equation*}
\frac{M_-(\zeta)}{M_+(\zeta)}=2^{-4\mu}\frac{\Gamma(1-\mu)}{\Gamma(1+\mu)}\frac{ 1 + (\e^{i\pi}\zeta)^{a_-}\E_-(\e^{i\pi}\zeta) + \e^{-\zeta}\e^{-i\pi a_-}\l[1+ \zeta^{a_-}\E_-(\zeta)\r]}{  1 + (\e^{i\pi}\zeta)^{a_+}\E_+(\e^{i\pi}\zeta) + \e^{-\zeta}\e^{-i\pi a_+}\l[1+ \zeta^{a_+}\E_+(\zeta)\r]}.
\end{equation*}
Moreover, we observe that that $|\zeta^a\E_1(\zeta)|\leq \frac{4\b^2}{|\zeta|}$, for $|\zeta|\geq 4\b^2$. Hence, for any $\delta>0$,
\begin{equation*}
\l| (\e^{i\pi}\zeta)^{a_\pm}\E_\pm(\e^{i\pi}\zeta) + \e^{-\zeta}\e^{-i\pi a_\pm}\l( 1+\zeta^{a_\pm}\E_\pm(\zeta)\r) \r| \leq \delta
\end{equation*}
provided that $\Re\zeta> N_{\mu,0}$ for some $N_{\mu,0}>0$.
\end{proof}

\begin{lemma}\label{Comparison bounds real M small argument}
Denote $a_\pm:=\frac12 \pm \mu$ and $b_\pm:=2a_\pm$. Then, 
\begin{equation*}
\lim_{\zeta\rightarrow 0}\frac{M_+(\zeta)}{M_-(\zeta)}=0.
\end{equation*}
Therefore, there exists $\delta_{\mu,1}>0$ such that 
\begin{equation*}
\l| \frac{M_+(\zeta)}{M_-(\zeta)}\r| \leq \min\l( \frac{1}{5M(a_-,b_-,2N_{\mu,0})}, \frac{1}{3C_\mu}\r)
\end{equation*}
for all $\zeta\in\C$ such that $|\zeta|\leq \delta_{\mu,1}$.
\end{lemma}
\begin{proof}
We recall once again that
\begin{equation*}
M_\pm(\zeta)=\e^{-\frac12\zeta}\zeta^{a_\pm} M(a_\pm,b_\pm,\zeta).
\end{equation*}
Hence, we directly compute
\begin{equation*}
\frac{M_+(\zeta)}{M_-(\zeta)}=\zeta^{2\mu}\frac{M(a_+,b_+,\zeta)}{M(a_-,b_-,\zeta)}. 
\end{equation*}
Since $M(a_\pm,b_\pm,\zeta)\rightarrow 1$ for $\zeta\rightarrow 0$, and $2\mu>0$, the conclusion follows for $|\zeta|$ small enough.
\end{proof}

\begin{lemma}\label{Comparison bounds real M order one argument}
Denote $a_\pm:=\frac12 \pm \mu$ and $b_\pm:=2a_\pm$. Let $y_0\in[0,1]$ such that $N_{\mu,1}\leq 2my_0\leq N_{\mu,0}$, for some $N_{\mu,1}\in\R$. Then, there exists $\ep_0>0$ such that
\begin{equation*}
\l| \frac{M_\pm(y_0+i\ep)}{M_\pm(y_0)} -1 \r| \leq \frac{\sin\mu\pi}{5},
\end{equation*} 
for all $0<|\ep|\leq \ep_0$.
\end{lemma}
\begin{proof}
Assume without loss of generality that $\ep>0$. Then,
\begin{equation*}
M_\pm(y_0+i\ep)=M_\pm(y_0) + \int_0^\ep \frac{\d}{\d s}M_\pm(y_0+is)\d s.
\end{equation*}
Thanks to the asymptotic expansions for small arguments we next estimate
\begin{equation*}
\l| \int_0^\ep \frac{\d}{\d s}M_\pm(y_0+is)\d s \r| \leq \int_0^\ep |M_\pm'(y_0+is)| \d s \leq C_\mu\int_0^\ep \frac{(2m)^{\frac12 \pm\mu}}{|y_0+is|^{\frac12 \mp\mu}}\d s \leq C_\mu(2m\ep)(2my_0)^{-\frac12}.
\end{equation*}
Using the lower bound $(2my_0)^{\frac12\pm\mu}\leq M_\pm(y_0)$ we have that
\begin{equation*}
\l| \int_0^\ep \frac{\d}{\d s}M_\pm(y_0+is)\d s \r| \leq C_\mu M_\pm(y_0)\frac{\ep}{y_0}\leq C_\mu M_\pm(y_0)N_{2}^{-1}2m\ep.
\end{equation*}
Hence, 
\begin{equation*}
M_\pm(y_0)\l(1-C_\mu\frac{\ep}{N_2}\r) \leq |M_\pm(y_0+ i\ep)|\leq M_\pm(y_0)\l(1+C_\mu\frac{\ep}{N_2}\r)
\end{equation*}
and now choose $\ep_0>0$ sufficiently small, so that the conclusion of the lemma follows swiftly for all $\ep\leq \ep_0$.
\end{proof}

\begin{lemma}\label{lemma:lower asymptotic bounds real M}
Let $y_0\in[0,1]$ and $0\leq \ep \leq \frac{\b}{m}$. Then, 
\begin{itemize}
\item If $my_0\leq 3\b$, there exists $\ep_0>0$ such that
\begin{equation*}
\l(m|y_0+ i\ep|\r)^{\frac12\pm \mu}\lesssim |M_\pm(y_0+ i\ep)|
\end{equation*}
for all $0\leq \ep\leq \ep_0$.
\item If $my_0\geq 3\b$, 
\begin{equation*}
1\lesssim {|M_\pm(y_0+i\ep)|}.
\end{equation*}
\end{itemize}
\end{lemma}
\begin{proof}
The proof uses the ideas from Lemma \ref{lemma:lower asymptotic bounds M} together with the bounds from Lemma \ref{growth bounds real M large argument}. We omit the details.
\end{proof}

 \section*{Acknowledgments} 
The research of MCZ was partially supported by the Royal Society URF\textbackslash R1\textbackslash 191492 and EPSRC Horizon Europe Guarantee EP/X020886/1.

\bibliographystyle{abbrv}
\bibliography{CZN-InvBoussChanBiblio}
\end{document}